\documentclass[final]{siamart190516} 

\usepackage{imakeidx}

\usepackage{lipsum}
\usepackage{amsfonts}
\usepackage{graphicx}
\usepackage{epstopdf}
\usepackage{algorithmic}
  \DeclareGraphicsExtensions{.eps,.pdf,.png,.jpg}

\graphicspath{{figures/}{BrBuKaKu_SIAMReview/figures/}{BrBuKaKu_SIAMReview/}}
\numberwithin{theorem}{section}

\newcommand{\TheTitle}{Modern Koopman Theory for Dynamical Systems}
\newcommand{\TheAuthors}{S. L. Brunton, M. Budi\v{s}i\'{c}, E. Kaiser, J. N. Kutz}

\headers{\TheTitle}{\TheAuthors}

\title{{\TheTitle}\thanks{Submitted to the editors Feb. 24, 2021.
\funding{
SLB acknowledges support from the Army Research Office (W911NF-17-1-0306 and W911NF-19-1-0045).
SLB and JNK acknowledge support from the Defense Advanced Research Projects Agency (DARPA contract HR011-16-C-0016) and the UW Engineering Data Science Institute, NSF HDR award \#1934292.
}}}

\author{
  \hspace{-.1in}Steven L. Brunton\thanks{Dept. of Mechanical Engineering, University of Washington, Seattle, WA (\email{sbrunton@uw.edu}).}
  \and
    Marko Budi\v{s}i\'{c}\thanks{Dept. of Mathematics, Clarkson University, Potsdam, NY (\email{mbudisic@clarkson.edu}).}
  \and
    Eurika Kaiser\thanks{Dept. of Mechanical Engineering, University of Washington, Seattle, WA  (\email{eurika@uw.edu}).}
  \and
  J. Nathan Kutz\thanks{Dept. of Applied Mathematics, University of Washington, Seattle, WA (\email{kutz.edu}).  \hspace{-.1in}}
}

\usepackage{amsopn}

\usepackage{amsmath,amsfonts,amscd,amssymb,graphicx}
\usepackage{xfrac}
\usepackage{mathtools}
\usepackage{array}
\usepackage{booktabs}
\usepackage{multirow}
\usepackage{paralist}

\newcolumntype{L}[1]{>{\raggedright\let\newline\\\arraybackslash\hspace{0pt}}m{#1}}
\newcolumntype{C}[1]{>{\centering\let\newline\\\arraybackslash\hspace{0pt}}m{#1}}
\newcolumntype{R}[1]{>{\raggedleft\let\newline\\\arraybackslash\hspace{0pt}}m{#1}}

\newcommand{\bb}{\mathbf{b}}
\newcommand{\bc}{\mathbf{c}}
\newcommand{\bd}{\mathbf{d}}
\newcommand{\bg}{\mathbf{g}}

\newcommand{\bs}{\mathbf{s}}

\newcommand{\bu}{\mathbf{u}}
\newcommand{\bv}{\mathbf{v}}

\newcommand{\bx}{\mathbf{x}}

\newcommand{\by}{\mathbf{y}}
\newcommand{\bz}{\mathbf{z}}
\newcommand{\bA}{\mathbf{A}}
\newcommand{\bB}{\mathbf{B}}
\newcommand{\bC}{\mathbf{C}}
\newcommand{\bD}{\mathbf{D}}
\newcommand{\bF}{\mathbf{F}}
\newcommand{\bG}{\mathbf{G}}
\newcommand{\bH}{\mathbf{H}}
\newcommand{\bI}{\mathbf{I}}

\newcommand{\bL}{\mathbf{L}}
\newcommand{\bK}{\mathbf{K}}

\newcommand{\bQ}{\mathbf{Q}}
\newcommand{\bR}{\mathbf{R}}
\newcommand{\bS}{\mathbf{S}}

\newcommand{\bT}{\mathbf{T}}

\newcommand{\bW}{\mathbf{W}}
\newcommand{\bX}{\mathbf{X}}

\newcommand{\bZ}{\mathbf{Z}}

\newcommand{\bAtilde}{\mathbf{\tilde{A}}}

\newcommand{\bVtilde}{\mathbf{\tilde{V}}}
\newcommand{\bSigmatilde}{\boldsymbol{\tilde{\Sigma}}}

\newcommand{\bxi}{\boldsymbol{\xi}}

\newcommand{\bphi}{\boldsymbol{\phi}}
\newcommand{\bPhi}{\boldsymbol{\Phi}}
\newcommand{\bGamma}{\boldsymbol{\Gamma}}

\newcommand{\bTheta}{\boldsymbol{\Theta}}

\newcommand{\bLambda}{\boldsymbol{\Lambda}}
\newcommand{\bOmega}{\boldsymbol{\Omega}}

\DeclareMathOperator*{\argmin}{arg\rm{}min}

\newcommand{\flow}{\ensuremath{\mathbf{F}}}
\newcommand{\koop}{\ensuremath{\mathcal{K}}}
\newcommand{\pf}{\ensuremath{\mathcal{P}}}
\newcommand{\gen}{\ensuremath{\mathcal{L}}}

\DeclareMathOperator{\lspan}{span}
\DeclarePairedDelimiter{\avg}{\langle}{\rangle}
\DeclarePairedDelimiter{\norm}{\lVert}{\rVert}
\DeclarePairedDelimiter{\abs}{\lvert}{\rvert}
\newcommand{\corr}[1][g]{C_{#1}}
\usepackage{overpic}
\ifpdf
\hypersetup{
  pdftitle={\TheTitle},
  pdfauthor={\TheAuthors}
}
\fi

\graphicspath{{./figures/}{./}}

\usepackage{subcaption}
\usepackage{etoolbox}
\makeatletter
\patchcmd{\@addmarginpar}{\ifodd\c@page}{\ifodd\c@page\@tempcnta\m@ne}{}{}
\makeatother
\reversemarginpar

\usepackage[colorinlistoftodos]{todonotes}
\setuptodonotes{size=\scriptsize,backgroundcolor=red!15!white}

\makeatletter
\renewcommand\@tocrmarg{2.55em}
\renewcommand\@dotsep{4.5}
\renewcommand\tableofcontents{%
	\section*{\contentsname\newline
		\@mkboth{%
			\MakeUppercase\contentsname}{\MakeUppercase\contentsname}}%
	\@starttoc{toc}%
}
\newcommand*\l@part[2]{%
	\ifnum \c@tocdepth >-2\relax
	\addpenalty\@secpenalty
	\addvspace{2.25em \@plus\p@}%
	\setlength\@tempdima{3em}%
	\begingroup
	\parindent \z@ \rightskip \@pnumwidth
	\parfillskip -\@pnumwidth
	{\leavevmode
		\large \bfseries #1\hfil \hb@xt@\@pnumwidth{\hss #2}}\par
	\nobreak
	\if@compatibility
	\global\@nobreaktrue
	\everypar{\global\@nobreakfalse\everypar{}}%
	\fi
	\endgroup
	\fi}
\newcommand*\l@section[2]{%
	\ifnum \c@tocdepth >\z@
	\addpenalty\@secpenalty
	\addvspace{1.0em \@plus\p@}%
	\setlength\@tempdima{1.5em}%
	\begingroup
	\parindent \z@ \rightskip \@pnumwidth
	\parfillskip -\@pnumwidth
	\leavevmode \bfseries
	\advance\leftskip\@tempdima
	\hskip -\leftskip
	#1\nobreak\hfil \nobreak\hb@xt@\@pnumwidth{\hss #2}\par
	\endgroup
	\fi}
\newcommand*\l@subsection{\@dottedtocline{2}{1.5em}{2.3em}}
\newcommand*\l@subsubsection{\@dottedtocline{3}{3.8em}{3.2em}}
\newcommand*\l@paragraph{\@dottedtocline{4}{7.0em}{4.1em}}
\newcommand*\l@subparagraph{\@dottedtocline{5}{10em}{5em}}
\makeatother

\begin{document}
\maketitle

\begin{abstract}
The field of dynamical systems is being transformed by the mathematical tools and algorithms emerging from modern computing and data science.
First-principles derivations and asymptotic reductions are giving way to data-driven approaches that formulate models in operator theoretic or probabilistic frameworks.
Koopman spectral theory has emerged as a dominant perspective over the past decade, in which nonlinear dynamics are represented in terms of an infinite-dimensional linear operator acting on the space of all possible measurement functions of the system.
This linear representation of nonlinear dynamics has tremendous potential to enable the prediction, estimation, and control of nonlinear systems with standard textbook methods developed for linear systems.
However, obtaining finite-dimensional coordinate systems and embeddings in which the dynamics appear approximately linear remains a central open challenge.
The success of Koopman analysis is due primarily to three key factors:  1) there exists rigorous theory connecting it to classical geometric approaches for dynamical systems, 2) the approach is formulated in terms of measurements, making it ideal for leveraging big-data and machine learning techniques, and 3) simple, yet powerful numerical algorithms, such as the dynamic mode decomposition (DMD), have been developed and extended to reduce Koopman theory to practice in real-world applications.
In this review, we provide an overview of modern Koopman operator theory, describing recent theoretical and algorithmic developments and highlighting these methods with a diverse range of applications.
We also discuss key advances and challenges in the rapidly growing field of machine learning that are likely to drive future developments and significantly transform the theoretical landscape of dynamical systems.
\end{abstract}

\begin{keywords}
Dynamical systems, Koopman operator, Data-driven discovery, Control theory, Spectral theory, Operator theory, Dynamic mode decomposition, Embeddings.
\end{keywords}

\begin{AMS}
34A34, 37A30, 37C10, 37M10, 37M99, 37N35, 47A35, 47B33
\end{AMS}

\tableofcontents

\section{Introduction}\label{Sec:Introduction}

Nonlinearity is a central challenge in dynamical systems, resulting in diverse phenomena, from bifurcations to chaos, that manifest across a range of disciplines.
However, there is currently no overarching mathematical framework for the explicit and general characterization of nonlinear systems.
In contrast, linear systems are completely characterized by their spectral decomposition (i.e., eigenvalues and eigenvectors), leading to generic and computationally efficient algorithms for prediction, estimation, and control.
Importantly, {\em linear superposition} fails for nonlinear dynamical systems, leading to a variety of interesting {phenomena} including frequency shifts and the generation of harmonics.
The Koopman operator theory of dynamical systems provides a promising alternative perspective, where linear superposition may be possible even for strongly nonlinear dynamics via the infinite-dimensional, but linear, Koopman operator.
The Koopman operator is linear, advancing measurement functions of the system, and its spectral decomposition completely characterizes the behavior of the nonlinear system.
Finding tractable \emph{finite-dimensional} representations of the Koopman operator is closely related to finding effective coordinate transformations in which the nonlinear dynamics appear linear.
Koopman analysis has recently gained renewed interest {with the pioneering work of Mezi\'{c} and collaborators}~\cite{mezic2004physicad,mezic2005nd,rowley2009jfm,budisic2012chaos,budisic2012physd,mezic2013arfm} { and its strong connections to data-driven modeling~\cite{schmid2010jfm,rowley2009jfm,kutz2016book}}.
{This review provides} an overview of modern Koopman theory for dynamical systems, including an in-depth analysis of leading computational algorithms, such as the dynamic mode decomposition (DMD) {of Schmid~\cite{schmid2010jfm}}.

Koopman introduced his operator theoretic perspective of dynamical systems in 1931  to describe the evolution of measurements of Hamiltonian systems~\cite{koopman1931pnas}, and this theory was generalized by Koopman and von Neumann to systems with a continuous eigenvalue spectrum in 1932~\cite{koopman1932pnas}.
Koopman's 1931 paper was central to the celebrated proofs of the ergodic theorem by von Neumann~\cite{neumann1932pnas} and Birkhoff~\cite{birkhoff1931pnas,birkhoff1932pnas}.
The history of these developments is fraught with intrigue, as discussed by Moore~\cite{moore2015pnas}; { a comprehensive survey of the operator-theoretic developments in ergodic theory can be found in~\cite{eisner2015}}.
In his original paper~\cite{koopman1931pnas}, Koopman drew connections between the Koopman eigenvalue spectrum and conserved quantities, integrability, and ergodicity.
For Hamiltonian flows, the Koopman operator is unitary.
{Efforts in the past two decades by Mezi\'{c} and colleagues have extended this theory from Hamiltonian systems with measure-preserving dynamics to dissipative and non-smooth dynamics~\cite{mezic2004physicad,mezic2005nd,mezic2013arfm}.
Furthermore, Rowley et al.~\cite{rowley2009jfm} rigorously connected the Koopman mode decomposition, introduced by Mezi\'{c} in 2005~\cite{mezic2005nd}, with the  DMD algorithm, introduced by Schmid in the fluid mechanics community~\cite{schmid2009dynamic,schmid2010jfm}.
This serendipitous connection justified the application of DMD to nonlinear systems and provided a powerful data-driven algorithm for the approximation of the Koopman operator.
This confluence of modern Koopman theory with a simple and effective numerical realization has resulted in rapid progress in the past decade, which is the main focus of this review.
It is important to note that this is not a comprehensive history of modern developments, as this is beyond the scope of the present review.
}
{Interestingly,} DMD, the leading numerical algorithm for approximating the Koopman operator, is built on the discrete Fourier transform (DFT) and the singular value decomposition (SVD), which both provide unitary coordinate transformations~\cite{brunton2019data}.

The operator theoretic framework discussed here complements the traditional geometric and probabilistic perspectives on dynamical systems.
For example, level sets of Koopman eigenfunctions form invariant partitions of the state-space of a dynamical system~\cite{budisic2012physd}; in particular, eigenfunctions of the Koopman operator may be used to analyze the ergodic partition~\cite{mezic1999chaos,budisic2009cdc,rokhlin1966}.
Koopman analysis has also been shown to generalize the Hartman--Grobman theorem to the entire basin of attraction of a stable or unstable equilibrium point or periodic orbit~\cite{lan2013physd}.
The Koopman operator is also known as the composition operator, which is formally the pull-back operator on the space of scalar observable functions~\cite{marsdenmtaa}, and it is the dual, or left-adjoint, of the Perron--Frobenius (PF) operator, or transfer operator, which is the push-forward operator on the space of probability density functions.
When a polynomial basis is chosen to represent the Koopman operator, then it is closely related to Carleman linearization~\cite{carleman1932am,carleman1933theorie,carleman1933systemes}, which has been used extensively in nonlinear  control~\cite{steeb1980non,kowalski1991nonlinear,banks1992infinite,svoronos1994discretization}.
{This review complements several other excellent reviews published in the last ten years~\cite{budisic2012physd,mezic2013arfm,taira2017,taira2020,otto2021koopman}, and a detailed comparison is given in \Cref{sec:organization}.}

\subsection{An overview of Koopman theory}
In this review, we will consider dynamical systems of the form
\begin{align}
    \frac{d}{dt}\bx(t) = \mathbf{f}(\bx(t)),\label{Eq:ContinuousDynamicsGeneral}
\end{align}
where $\bx\in\mathcal{X}\subseteq\mathbb{R}^n$ is the state of the system{, possibly living on a submanifold $\mathcal{X}$ of a $n$-dimensional vector space $\mathbb{R}^n$,} and $\mathbf{f}$ is a vector field describing the dynamics.
In general, the dynamics may also depend on time $t$, parameters $\boldsymbol{\beta}$, and external actuation or control $\bu(t)$.
Although we omit these here for simplicity, they will be considered in later sections.

A major goal of modern dynamical systems is to find a new vector of coordinates $\bz$ such that either
\begin{equation}
\bx=\boldsymbol{\varphi}(\bz) \quad\text{ or }\quad \bz=\boldsymbol{\varphi}(\bx) \label{eq:linearizing-coordinates}
\end{equation}
where the dynamics are simplified, or ideally, linearized:
\begin{align}
    \frac{d}{dt}\bz = \bL \bz.
\end{align}
{In these new linearizing coordinates, the dynamics of $\bz$ are entirely determined by the matrix $\bL$.
The future evolution of the system in these coordinates may be fully characterized by the eigendecomposition of $\bL$. }
While in geometric dynamics, one asks for homeomorphic (continuously invertible) or even diffeomorphic coordinate maps, which trivialize the choice between the two options in~\eqref{eq:linearizing-coordinates}, there is little hope for global coordinate maps of this sort.
Rather, we contend with embeddings $\boldsymbol{\varphi}$ that lift the dynamics into a higher-dimensional space of $\bz$ variables, allowing for ``unfolding'' of nonlinearities.

In practice, we typically have access to \emph{measurement data} of our system, discretely sampled in time.
This data is governed by the discrete-time dynamical system
\begin{eqnarray}
\mathbf{x}_{k+1}=\flow(\mathbf{x}_k),\label{Eq:DiscreteDynamics}
\end{eqnarray}
where $\bx_k=\bx(t_k)=\bx(k\Delta t)$.
 Also known as a \emph{flow map}, the discrete-time dynamics are more general than the continuous-time formulation in \eqref{Eq:ContinuousDynamicsGeneral}, encompassing discontinuous and hybrid systems as well.
In this case, the goal is still to find a linearizing coordinate transform so that $\bz_{k+1}=\bK\bz_k${, where the matrix $\bK$ is the discrete-time analog of the continuous-time matrix $\bL$}.
These coordinates are given by eigenfunctions of the discrete-time Koopman operator, $\koop$, which advances a measurement function $g(\bx)$ of the state forward in time through the dynamics:
\begin{eqnarray}
\koop g(\bx_k) \coloneqq g(\flow(\bx_k)) = g(\bx_{k+1}).
\end{eqnarray}
For an eigenfunction $\varphi$ of $\koop$, corresponding to an eigenvalue $\lambda$, this becomes
\begin{eqnarray}
\koop \varphi(\bx_k) = \lambda \varphi(\bx_k) = \varphi(\bx_{k+1}).
\end{eqnarray}
Thus, a tremendous amount of effort has gone into characterizing the Koopman operator and approximating its spectral decomposition from measurement data.

The coordinates $\boldsymbol{\varphi}$ and the {matrix} $\bL$ are closely related to the continuous-time analogue $\gen{}$ of the {discrete-time} Koopman operator $\koop$, {which will both be introduced in more detail in \Cref{Sec:Koopman}}.
In particular, eigenfunctions $\varphi_j$ of $\gen{}$ provide such a linearizing coordinate system, and the matrix $\bL$ is obtained by restricting the operator $\gen{}$ to the span of these functions.
Spectral theory provides a complete description of the dynamics in terms of the eigenstructure of $\bL$.
Thus, transforming the system into coordinates where the dynamics are linear dramatically simplifies all downstream analysis and control efforts.

{In the following, we will generally use calligraphic symbols for operators (e.g., $\mathcal{L}$ or $\mathcal{K}$) and bold capital letters for matrices (e.g., $\mathbf{L}$ or $\mathbf{K}$).
It should be noted that matrices are representations of finite-dimensional linear operators in a particular basis, and so we will occasionally refer to multiplication by a matrix as a linear operator.
Much of modern Koopman theory is concerned with uncovering the intrinsic spectral properties of an operator up to conjugacy or change of coordinates.
}

\begin{figure}[t]
\begin{center}
\begin{overpic}[width=1\textwidth]{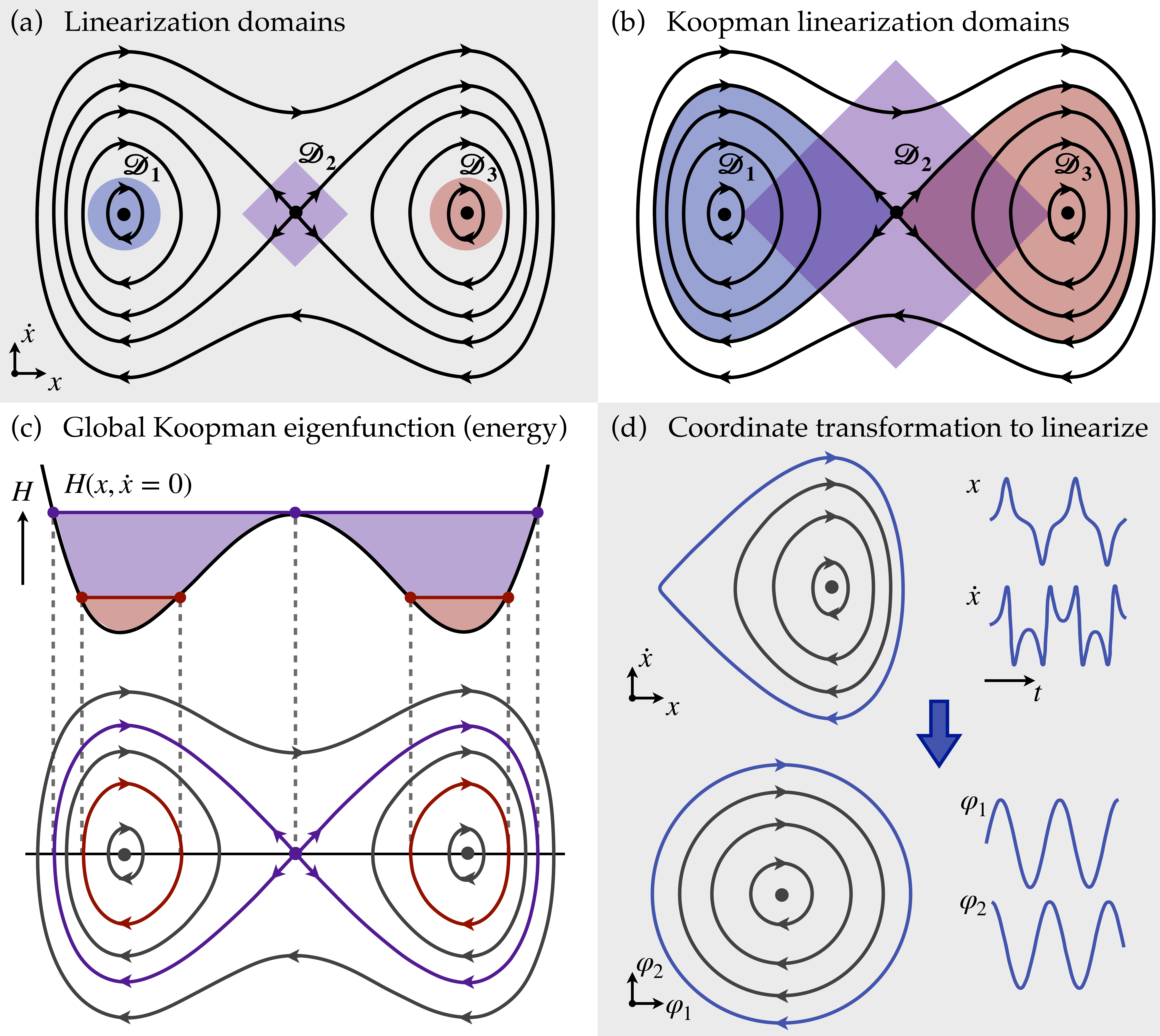}
\end{overpic}
\caption{Different Koopman perspectives for the Duffing oscillator, $\ddot{x} = x-x^3$, the equation for a particle in a double potential well.  (a) Traditional linearization near the fixed points gives small regions where the system is approximately linear.  (b) Koopman theory may extend the Hartman--Grobman theorem to enlarge the domain of linearity until the next fixed point~\cite{lan2013physd}. (c) There are also global Koopman eigenfunctions, like the Hamiltonian energy, although these lose information about which basin the solution is in.  (d) Yet a third perspective seeks a coordinate transformation to rescale space and time until dynamics live on a hypertoroid.}\label{Fig:KoopmanOverview}
\end{center}
\end{figure}

\subsection{An illustrative example: The Duffing oscillator}
Although the objective of Koopman theory is easily expressed mathematically, it is helpful to explore its application on a simple dynamical system.
Consider the nonlinear Duffing system $\ddot{x}=x-x^3$ with state space representation
\begin{subequations}
    \begin{align}
        \dot{x}_1 & = x_2\\
        \dot{x}_2 & = x_1-x_1^3,
    \end{align}
\end{subequations}
and the corresponding phase portrait in~\cref{Fig:KoopmanOverview}.
This example has three fixed points, a saddle at the origin with Jacobian eigenvalues $\lambda_{\pm}=\pm 1$ and two centers at $(x_1,x_2)=(\pm 1,0)$ with eigenvalues $\pm\sqrt{2}i$.
These results, shown in \cref{Fig:KoopmanOverview}(a), can be obtained by a local phase-plane analysis~\cite{boyce2017elementary}, and these local linearizations are valid in a small neighborhood of each fixed point, illustrated by the shaded regions.

{
The Duffing oscillator is a classic textbook example of a weakly nonlinear system, where fixed points may be identified and linearized about, and the stable and unstable manifolds emanating from these points organize the entire phase space structure.
This so-called \emph{geometric} perspective on dynamical systems has become the dominant viewpoint over the past century, and powerful techniques have been developed to analyze and visualize these systems~\cite{guckenheimer1986,Haller2002pof,Shadden2005pd,Farazmand2012chaos,HLBR_turb,tallapragada2013set,Haller2015arfm,Strogatz2018book}.
However, this example, even in its simplicity, highlights many of the challenges and subtleties of Koopman operator theory.
For example, the Duffing oscillator exhibits a continuous spectrum of frequencies, as can be seen by increasing the energy from one of the center fixed points, where the frequency is determined by the linearized eigenvalues, up to the saddle point, where the period of oscillation tends toward infinity.
However, for the Koopman operator acting on ``usual'' function spaces, the only eigenspace is the one at $\lambda=0$, containing the conserved Hamiltonian energy of the system, indicator functions associated with invariant sets of positive area, such as ``bands'' of periodic orbits, and other invariant functions.
Adding dissipation to the system regularizes the problem, at the expense of removing the continuous spectrum of frequencies.
Mezi\'{c} recently conducted an in-depth study of these subtleties that appear even in simpler systems, such as the nonlinear pendulum~\cite{mezic2019}, which we summarize in \ref{sec:koopman-spectrum}.
}

There is no homeomorphic coordinate transformation that captures the global dynamics of this system with a linear operator, since any such linear operator has either one fixed point at the origin, or a subspace of infinitely many fixed points~\cite{brunton2016plosone}, but never three isolated fixed points.
Instead, the Koopman operator can provide a system of coordinate transformations that extend the local neighborhoods where a linear model is valid to the full basin around them~\cite{lan2013physd}, as shown in \cref{Fig:KoopmanOverview}(b).

This interpretation may lead to seemingly unintuitive results, as even the most obvious observables, such as $g(\bx) = \bx$, may require different expansions in different areas of the state space, as recently explored by Page and Kerswell~\cite{page2019}.
Indeed, systems that possess multiple simple invariant solutions cannot be represented by a Koopman expansion that is globally uniformly convergent, implying that any choice of a representation of the Koopman operator $\koop$ in a system of coordinates may not hold everywhere, for any given observable of interest.

Even if there are no globally convergent linear Koopman representations, the Koopman eigenfunctions maybe be globally well defined, and even regular in the entire phase space.
For example, the Hamiltonian energy function ${H=x_2^2-x_1^2/2+x_1^4/4}$ is a Koopman eigenfunction for this conservative system, with eigenvalue $\lambda=0$, as shown in \cref{Fig:KoopmanOverview}(c).
In a sense, this global eigenfunction exploits symmetry in the dynamics to represent a global dynamic quantity, valid in all regions of phase space.
This example establishes a connection between the Koopman operator and Noether's theorem~\cite{noether1918invariante}, as a symmetry in the dynamics gives rise to a new Koopman eigenfunction with eigenvalue $\lambda = 0$.
In addition, the constant function $\varphi\equiv1$ is a trivial eigenfunction corresponding to $\lambda=0$ for every dynamical system, which is a consequence of area-preservation of the system.
{ The eigenspace at \(\lambda=0\) is rich, as the characteristic function of any invariant set of positive area, for example, an annular band around either equilibrium, acts as an eigenfunction.
  At the same time, the minimal invariant sets, corresponding to individual orbits, are of area-zero, and their characteristic functions are neither continuous nor members of any \(L^{p}\) space.
  Instead, they give rise to  Dirac-\(\delta\) distributions that can be considered as eigenfunctions only in the appropriate function space; for more on this see \Cref{sec:koopman-spectrum}.}

A final interpretation of Koopman theory is shown in \cref{Fig:KoopmanOverview}(d), related to \cref{Fig:KoopmanOverview}(b).
In this case, we see that the Koopman operator establishes a change of coordinates, in space and time, in which the original nonlinear trajectories become linear.
When obtaining these approximate coordinate systems from data, only approximately periodic dynamics will persist for long times.
Rescaling space and time to make dynamics approximately periodic is the approach taken in~\cite{lusch2017arxiv,lange2020fourier}; { Giannakis also performed a nonlinear time transformation to study mixing systems~\cite{giannakis2019}, while Bollt and coauthors studied the finite-time blowup in the context of Koopman theory~\cite{bollt2018}}.

\subsection{Dynamics in the big data era}
The recent interest in Koopman operator theory is inherently linked to a growing wealth of data.
Measurement technologies across every discipline of the engineering, biological, and physical sciences have revolutionized our ability to interrogate and extract quantities of interest from complex systems in real time with rich multi-modal and multi-fidelity time-series data.
From broadband sensors that measure at exceptionally high sampling rates to high-resolution imaging, the front lines of modern science are being transformed by unprecedented quality and quantities of data.
These advances are based on three foundational technologies: (i) improved sensors capable of measurement with higher quality and quantity, (ii) improved hardware for high-performance computing, data storage, and transfer, and (iii) improved algorithms for processing the data.
Taken together, these advances form the underpinnings of the {\em big data} era and are driving the fourth paradigm of scientific discovery~\cite{hey2009fourth}, whereby the data itself drives discovery.

Data science is now a targeted growth area across almost all academic disciplines.
However, it has been almost six decades since John Tukey, co-developer of the \emph{fast Fourier transform} with James Cooley, first advocated for data science as its own discipline~\cite{tukey1962future,donoho201750}.
Not surprisingly, the 1960s also coincided with pioneering developments of Gene Golub and co-workers on numerical algorithms for computing the {\em singular value decomposition} of a matrix~\cite{stewart1993early}, enabling one of the earliest data science exploration tools:  {\em principal component analysis} (PCA).
Thus the mathematical foundations for the {\em big data} era have been long in development.
Indeed, its maturity is reflected in the emergence of the two cultures~\cite{breiman2001statistical} of {\em statistical learning} and {\em machine learning}.
In the former, the primary focus is on the development of interpretable models of data, while in the latter, accuracy is of paramount importance.
Although accuracy and interpretability are not mutually exclusive, often the refinement of one comes at the expense of the other.
For example, modern machine learning and artificial intelligence algorithms are revolutionizing computer vision and speech processing through deep neural network (DNN) architectures.
DNNs have produced performance metrics in these fields far beyond any previous algorithms.
Although individual components of DNNs may be interpretable, the integration across many layers with nonlinear activation functions typically render them opaque and uninterpretable.
In contrast, sparsity promoting algorithms, such as the LASSO~\cite{tibshirani1996lasso}, are examples of statistical learning where interpretable variable selection is achieved.
As will be highlighted throughout this review, Koopman theory is amenable to many of the diverse algorithms developed in both {\em statistical learning} and {\em machine learning}.

Dynamical systems theory has a long history of leveraging data for improving modeling insights, promoting parsimonious and interpretable models, and generating forecasting capabilities.
In the 1960s, Kalman introduced a rigorous mathematical architecture~\cite{kalman1960general,kalman1963mathematical} whereby data and models could be combined through data assimilation techniques~\cite{evensen2009data,law2015data}, which is especially useful for forecasting and control.
Thus the integration of streaming data and dynamical models has a nearly seven decade history.
In the modern era, it is increasingly common to build the dynamical models from the data directly using machine learning~\cite{Schmidt2009science,battaglia2018relational,Raissi2019jcp,Noe2019science,raissi2020science,bar2019learning,cranmer2019learning,sanchez2020learning,lee2020model,lu2021learning,li2020fourier,li2020multipole,li2020neural,rackauckas2020universal,cranmer2020lagrangian,wang2021learning}.
This is especially important in complex systems where first principles models are not available, or where it is not even known what the correct state-space variable should be.
Biological systems, such as those that arise in neuronal recordings in the brain, are well suited for such data-driven model discovery techniques, as whole-brain imaging provides insight into how the microscale dynamics of individual neurons produces large scale patterns of spatio-temporal brain activity.
Such systems are ideally suited for leveraging the modern data-driven modeling tools of machine learning to produce dynamical models characterizing the observed data.

\subsection{Koopman objectives and applications}
The ultimate promise of Koopman spectral theory is the ability {to} analyze, predict, and control nonlinear systems with the wealth of powerful techniques from linear systems theory.
This overarching goal may be broken into several specific goals:

\paragraph{Diagnostics}
Spectral properties of the Koopman operator may be used effectively to characterize complex, high-dimensional dynamical systems.
For example, in fluid mechanics, DMD is used to approximate the Koopman mode decomposition, resulting in an expansion of the flow as a linear combination of dominant coherent structures.
Thus, Koopman mode decomposition can be used for dimensionality reduction and model reduction, generalizing the space-time separation of variables that is classically obtained via either Fourier transform or singular value decomposition~\cite{brunton2019data}.

\paragraph{Prediction}
One of the major benefits of dynamical systems over other statistical models is that they provide a mechanistic understanding of how the future evolves, based on the current state of a system.
Prediction is thus one of the central goals of any dynamical system framework.
The prediction of a nonlinear system presents numerous challenges, such as sensitive dependence on initial conditions, parameter uncertainty, and bifurcations.
Koopman operator theory provides insights into which observables are easier or harder to predict, and suggests robust numerical algorithms for such time-series prediction.

\paragraph{Estimation and control}
In many systems, we not only seek to understand the dynamics, but also to modify, or \emph{control}, their behavior for some engineering goal.
In modern applications, it is rare to have complete measurements of a high-dimensional system, necessitating the \emph{estimation} of these quantities from limited sensor measurements.
By viewing the system through the linearizing lens of Koopman eigenfunctions, it is possible to leverage decades of results from linear estimation and control theory for the analysis of nonlinear systems.
Indeed, linear systems are by far the most well studied and completely characterized class of dynamical systems for control.

\paragraph{Uncertainty quantification and management}
Uncertainty is a hallmark feature of the real world, and dynamical systems provide us with a framework to forecast, quantify, and manage uncertainty.
However, modeling uncertainty is challenging for strongly nonlinear and chaotic systems.
Integrating uncertainty quantification into {the transfer operator framework has been connected early with the modern applications of Koopman analysis~\cite{mezic2004b,mezic2008c,mezic2004b} and remains a topic of ongoing research~\cite{mukherjee2015,streif2013}.}

\paragraph{Understanding}
Beyond the practical goals value of prediction and control, dynamical systems provide a compelling framework with which to understand and model complex systems.
Normal forms, for example, distill essential structural information about the dynamics, while remaining as simple as possible.
For Koopman theory to be embraced, it will need to facilitate similar intuition and understanding.

\begin{figure}
    \centering
    \begin{overpic}[width=\textwidth]{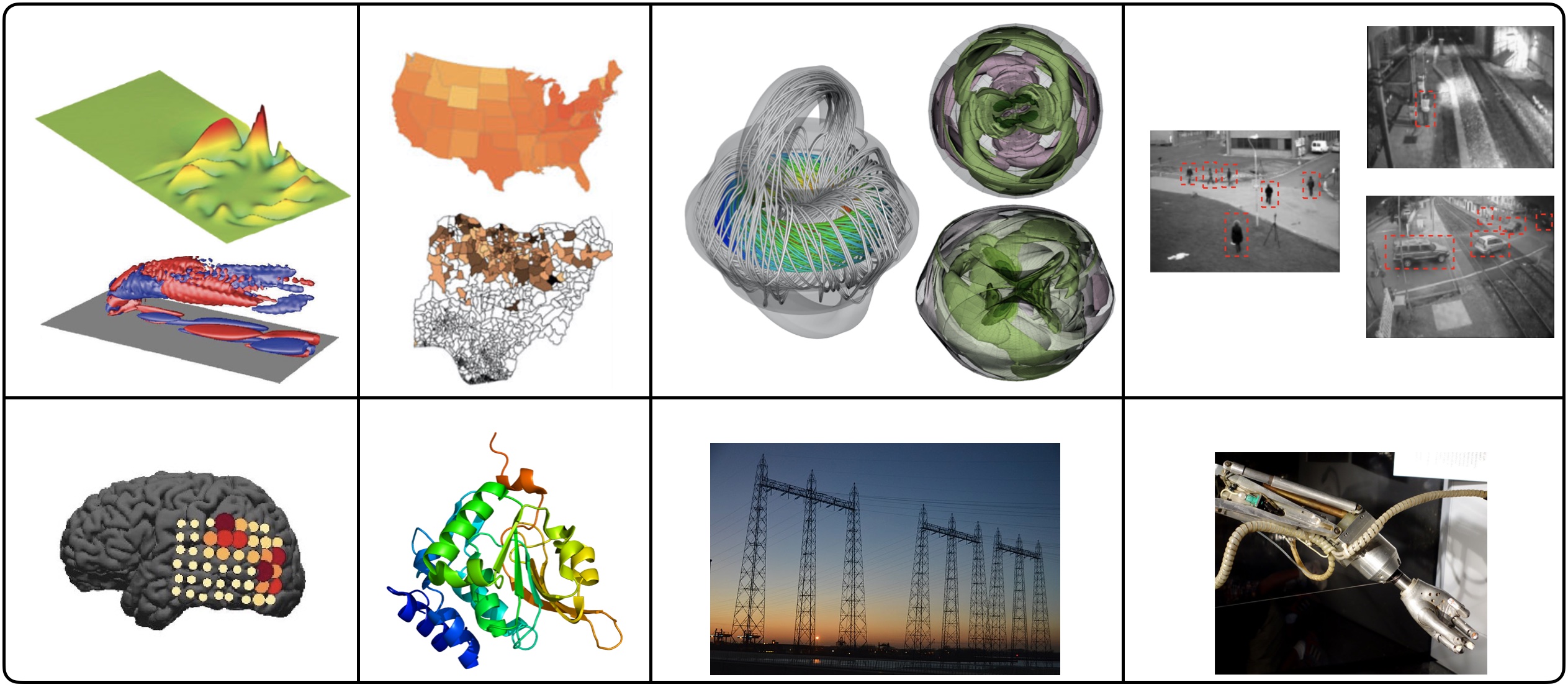}
    \footnotesize
    \put(1,41){(a) Fluids}
    \put(23.25,41){(b) Epidemiology}
    \put(42,41){(c) Plasmas}
    \put(72.5,41){(d) Video}
    \put(1,16){(e) Neuroscience}
    \put(23.25,16){(f) Chemistry}
    \put(42,16){(g) Power Grid}
    \put(72.5,16){(h) Robotics}
    \end{overpic}
    \caption{Overview of applications of data-driven Koopman analysis via DMD. Figures reproduced with permission from: (a) top~\cite{schmid2010jfm}, bottom~\cite{rowley2009jfm}; (b) \cite{proctor2015ih}; (c) \cite{kaptanoglu2020pop}; (d) \cite{erichson2016jrtp}; (e) \cite{brunton2016b}; (f) from Emw \url{https://commons.wikimedia.org/wiki/File:Protein_PCMT1_PDB_1i1n.png}; (g) from Henk Monster \url{https://commons.wikimedia.org/wiki/File:Power_grid_masts_besides_the_new_Waalbridge_Nijmegen_-_panoramio.jpg}; (h) from
Daderot \url{https://commons.wikimedia.org/wiki/File:Minsky's_robot_arm,_late_1960s,_view_2_-_MIT_Museum_-_DSC03759.JPG}.
    }
    \label{fig:Overview}
\end{figure}
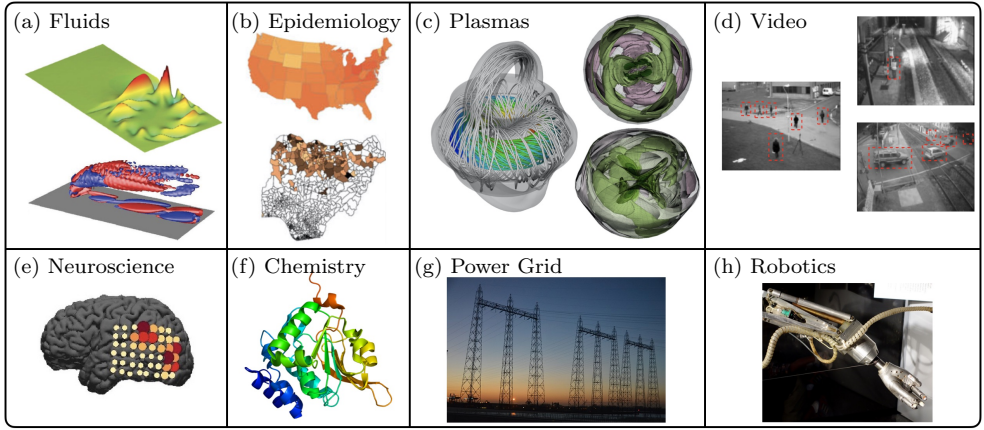

Across a wide range of applications, researchers are developing and advancing Koopman operator theory to address these challenges.
Some of these applications include fluid dynamics~\cite{schmid2010jfm,rowley2009jfm,mezic2013arfm}, epidemiology~\cite{proctor2015ih}, neuroscience~\cite{brunton2016b,alfatlawi2019arxiv}, plasma physics~\cite{taylor2017arxiv,kaptanoglu2020pop}, finance~\cite{mann2016qf}, robotics~\cite{berger2014ieee,abraham2019ieee,bruder2019proc,bruder2020arxiv}, and the power grid~\cite{susuki2011jns,susuki2011c}; a number of these are shown in \cref{fig:Overview} and will be explored more in \cref{Sec:DMD}.

\subsection{Organization and goals of review}\label{sec:organization}
A primary goal of this review is to make modern Koopman operator theory accessible to researchers in many diverse fields.
We seek to provide the reader with a big-picture overview of the state of the art, expediting the process of getting up to speed in this rapidly developing field.
Further, we explore the dual theoretical and applied perspectives for understanding Koopman operator theory.
To this end, each section will be approached from the perspective of providing a unified overview of major theoretical, methodological, numerical, and applied developments.
Furthermore, we have compiled a number of striking success stories along with several outstanding challenges to help guide readers wishing to enter this field.
Finally, it is our goal to highlight connections with existing and emerging methods from other fields, including geometric and probabilistic dynamical systems, computational science and engineering, and machine learning.

{ There have been several excellent reviews focusing on Koopman operator theory and related methodologies published during the past ten years~\cite{budisic2012physd,mezic2013arfm,taira2017,taira2020,otto2021koopman}.
This review is similar in spirit to~\cite{budisic2012physd}, but surveys the considerable body of literature published in the decade since~\cite{budisic2012physd}.
For example, the present review establishes connections between the Koopman operator theory and classical nonlinear dynamics (see \Cref{Sec:AdvancedKoopman}) and, in particular, with modern control theory (see \Cref{Sec:Control}).
The review by Mezi\'{c}~\cite{mezic2013arfm} and the two broader reviews that surveyed Koopman and DMD techniques~\cite{taira2017,taira2020} focused on developments associated with fluid dynamics; while the original applications of Koopman and DMD were to problems of fluid motion, these techniques were shown to be useful in a broader range of fields which we highlight here.
The review by Otto and Rowley~\cite{otto2021koopman} focuses mainly on connections to control theory, and although \Cref{Sec:Control} of the present review overlaps with this material, it is in the broader context of connections to nonlinear systems and data-driven numerical implementations.
Finally, we envisioned this article to serve as guide to the body of work published within the last decade, but also as a practical introduction to techniques associated with the Koopman operator  for those that are new to the field (particularly~\Cref{Sec:Koopman}).
}

This review begins with relatively standard material and builds up to advanced concepts.
\Cref{Sec:Koopman} provides a practical introduction Koopman operator theory, including simple examples, the Koopman mode decomposition, and spectral theory.
Similarly, \cref{Sec:DMD} provides an overview of dynamic mode decomposition, which is the most simple and widely used numerical algorithm to approximate the Koopman operator.
\Cref{Sec:AdvancedKoopman} then introduces many of the most important concepts in \emph{modern} Koopman theory.
Advanced numerical and data-driven algorithms for representing the Koopman operator are presented in \cref{Sec:Observables}.
Thus, \cref{Sec:AdvancedKoopman} and~\cref{Sec:Observables} provide advanced material extending \cref{Sec:Koopman} and \cref{Sec:DMD}, respectively.
\Cref{Sec:Control} investigates how Koopman theory is currently extended for advanced estimation and control.
\Cref{Sec:Discussion} provides a discussion and concluding remarks. This section also describes several ongoing challenges in the community, along with recent extensions and motivating applications.
The goal here is to provide researchers with a quick sketch of the state of the art, so that they may quickly get to the frontier of research.

\section{A practical introduction to the Koopman operator framework}\label{Sec:Koopman}
\subsection{Definitions and vocabulary}
The Koopman operator advances measurement functions of the state of a dynamical system with the flow of the dynamics.
To explain the basic properties of the Koopman operator, we begin with an autonomous ordinary differential equation \eqref{Eq:ContinuousDynamicsGeneral} on a finite-dimensional state space \(\mathcal{X}\subseteq \mathbb{R}^{n}\).
The flow map operator, or time-\(t\) map,  \(\flow^{t} \colon \mathcal{X} \to \mathcal{X}\) advances initial conditions $\bx(0)$ forward along the trajectory by a time $t$, so that trajectories evolve according to
\begin{equation}
  \bx(t) = \flow^{t}( \bx(0) ). \label{eq:flow-map}
\end{equation}
The family of Koopman operators \(\koop^{t} \colon \mathcal{G}(\mathcal{X}) \to \mathcal{G}(\mathcal{X})\), parameterized by $t$, are given by
\begin{equation}
\koop^{t} g (\bx)  = g (\flow^{t}(\bx))  \label{eq:koopman},
\end{equation}
where \(\mathcal{G}(\mathcal{X})\) is a set of \emph{measurement functions}
\(g \colon \mathcal{X} \to \mathbb{C}\).
Another name for \(g\), derived from this framework's origin in quantum mechanics, is an \emph{observable} function, although this should not be confused with the unrelated \emph{observability} from control theory.
We can interpret \eqref{eq:koopman} as defining a family of functions
\begin{equation}
g_{t} \coloneqq \koop^{t}g, \quad g_{0}\coloneqq g\label{eq:functional-trajectory}
\end{equation}
that corresponds to the trajectory \(t \mapsto g_{t}\) in the set \(\mathcal{G}(\mathcal{X})\) of measurement functions.

In most applications, the set of functions \(\mathcal{G}(\mathcal{X})\) is not defined \emph{a priori}, but is loosely specified by a set of properties it should satisfy, e.g., that it is a vector space, that it possesses an inner product, that it is complete, or that it contains certain functions of interest, such as continuous functions on \(\mathcal{X}\).
Hilbert spaces, such as \(L^{2}(\mathcal{X},d\mu)\) or reproducing kernel Hilbert spaces (RKHS), are a common choice in modern applied mathematics, although historically other Banach spaces, such as integrable functions \(L^{1}\) or continuous functions \(C(\mathcal{X})\) have also been used.
The choice of the space, whether explicit or implicit, can have consequences on the properties of the operator and its approximations.
In all cases, however, \(\mathcal{G}(\mathcal{X})\) is of significantly higher dimension than \(\mathcal{X}\), i.e., countably or uncountably infinite.

The most significant property of the Koopman operator is that it is linear when \(\mathcal{G}(\mathcal{X})\) is a linear (vector) space of functions:
 \begin{equation}\label{eq:linearity}
 \begin{aligned}
 \koop^t\left(\alpha_1g_1(\bx) + \alpha_2g_2(\bx)\right) & = \alpha_1 g_1\left(\flow^{t}(\bx)\right) +  \alpha_2 g_2\left(\flow^{t}(\bx)\right)\\
 & = \alpha_1 \koop^t g_1(\bx) + \alpha_2 \koop^t g_2(\bx).
 \end{aligned}
\end{equation}
This property holds regardless of whether \(\flow^{t} \colon \mathcal{X} \to \mathcal{X}\) is linear itself, as it is simply a consequence of definition \eqref{eq:koopman}, since the argument function \(g\) is on the ``outside'' of the composition, allowing linearity to carry over from the vector space of observables.
In this sense, the Koopman framework obtains linearity of \(\koop^{t}\) despite the nonlinearity of \(\flow^{t}\) by trading the finite-dimensional state space \(\mathcal{X}\) for an infinite-dimensional function space \(\mathcal{G}(\mathcal{X})\).

When time \(t\) is discrete, \(t\in \mathbb{N}\), and the dynamics are autonomous, then \(\flow^{t}\) is a repeated \(t\)-fold composition of $\flow\equiv\flow^1$ given by \(\flow^{t}(\bx) = \flow(\flow(\cdots (\flow(\bx))))\), so that \(\koop^{t} g\) is likewise generated by repeated application of \(\koop \equiv \koop^{1}\).
The generator \(\koop\) of the (countable) composition semigroup is then called \emph{the} Koopman operator, which results in a dynamical system
\begin{equation}
  \label{eq:koopman-dynamics-discrete}
  g_{k+1} = \koop g_{k},
\end{equation}
analogous to \(\bx_{k+1} = \flow (\bx_{k})\), except that \eqref{eq:koopman-dynamics-discrete} is linear and infinite-dimensional.

When time \(t\) is continuous, the flow map family satisfies the \emph{semigroup} property
\begin{equation}
  \label{eq:semigroup-property}
  \flow^{t+s}(\bx) = \flow^{t}( \flow^{s}(\bx) ), \forall \bx,\ t,s \geq 0,
\end{equation}
which can be strengthened to a group property \(t,s \in \mathbb{R}\) if the flow map is invertible.
The Koopman operator family \(\koop^{t}\) inherits these properties as well.
Given a continuous and sufficiently smooth dynamics, it is also possible to define the continuous-time infinitesimal generator {$\mathcal{L}$} of the Koopman operator family as
\begin{equation}
\gen{} g \coloneqq  \lim_{t\to 0} \frac{\koop^{{t}}g-g}{t} =  \lim_{t\to 0} \frac{g\circ \flow^{t}-g}{t}.
\label{Eq:Koopman:InfinitesimalGenerator}
\end{equation}
As our goal here is to provide an introduction, we omit the mathematical scaffolding that accompanies a careful definition of an operator derivative; all details in the context of Koopman operator theory are available in \cite[\S 7.5]{lasota1994}.

The generator \(\gen{}\) has been called the Lie operator~\cite{koopman1931pnas}, as it is the Lie derivative of \(g\) along the vector field \(\mathbf{f}(\bx)\) when the dynamics is given by \eqref{Eq:ContinuousDynamicsGeneral}~\cite{marsdenmtaa,chicone1999evolution}.
This follows from applying the chain rule to the time derivative of $g(\bx)$:
\begin{align}
\frac{d}{dt} g(\bx(t)) &= \nabla g \cdot \dot{\bx} (t) =  \nabla g \cdot \mathbf{f}( \bx(t) ) \label{eq:time-derivative}
\end{align}
and equating with
\begin{align}
\frac{d}{dt} g(\bx(t)) &= \lim_{\tau\to 0} \frac{g( \bx(t + \tau) ) -g( \bx(t) )}{\tau} = \gen{}(g(\bx(t))), \label{eq:Lie-derivative}
\end{align}
resulting in
\begin{align}
  \label{eq:Lie-operator}
  \gen{} g = \nabla g \cdot \mathbf{f}.
\end{align}

The adjoint of the Lie operator is called the Liouville operator, especially in Hamiltonian dynamics~\cite{gaspard1998,gaspard1995}, while the adjoint of the Koopman operator is the Perron--Frobenius operator~\cite{dellnitz1999,dellnitz2001,froyland2003,froyland2014}.
In many ways, the operator-theoretic framework for applied dynamical systems has two dual perspectives, corresponding either to the Koopman operator or the Perron--Frobenius operator.
In \cref{sec:perron-frobenius} we discuss these parallels in more detail.

Similar to \eqref{eq:koopman-dynamics-discrete}, \(\gen{}\) induces  a linear dynamical system in continuous-time:
\begin{eqnarray}
\frac{d}{dt}g = \gen{} g.\label{eq:koopman-dynamics-continuous}
\end{eqnarray}

The linear dynamical systems in~\eqref{eq:koopman-dynamics-discrete} and~\eqref{eq:koopman-dynamics-continuous} are analogous to the dynamical systems in~\eqref{Eq:ContinuousDynamicsGeneral} and~\eqref{Eq:DiscreteDynamics}, respectively.

An important special case of an observable function is the projection onto a component of the state \(g(\bx) = x_{i}\) or, with a slight abuse of notation, \(g(\bx) = \bx\).
In this case, the left-hand side of \eqref{eq:koopman-dynamics-continuous} is plainly \(\dot{\bx}\), but the right hand side \(\gen{} \bx\) may not be simple to represent in a chosen basis for the space \(\mathcal{G}(\mathcal{X})\).
It is typical in real problems that representing \(\gen{}\bx\) will involve an infinite number of terms.
For certain special structures, this may not be the case, as will be demonstrated in \cref{Sec:SlowManifold}.

In summary, the Koopman operator is linear, which is appealing, but is infinite dimensional, posing issues for representation and computation.
Instead of capturing the evolution of all measurement functions in a function space, applied Koopman analysis attempts to identify key measurement functions that evolve linearly with the flow of the dynamics.
Eigenfunctions of the Koopman operator provide just such a set of special measurements that behave linearly in time.
In fact, a primary motivation to adopt the Koopman framework is the ability to simplify the dynamics through the eigendecomposition of the operator.

\subsection{Eigenfunctions and the spectrum of eigenvalues}\label{sec:koopman-eigenfunctions}

A Koopman eigenfunction \(\varphi(\bx)\) corresponding to an eigenvalue \(\lambda\) satisfies
\begin{equation}
\varphi(\bx_{k+1}) = \koop \varphi(\bx_k) = \lambda\varphi(\bx_k).\label{Eq:KoopmanEfun:Discrete}
\end{equation}
In continuous-time, a Lie operator eigenfunction \(\varphi(\bx)\) satisfies
\begin{equation}
\frac{d}{dt}\varphi(\bx) = \gen\varphi(\bx) = \mu \varphi(\bx),\label{Eq:KoopmanEfun}
\end{equation}
where $\mu$ is a continuous-time eigenvalue.
In general, eigenvalues and eigenvectors are complex-valued scalars and functions, respectively, even when the state space $\mathcal{X}$ and dynamics \(\bF(\bx)\) are real-valued.

It is simple to show that Koopman eigenfunctions \(\varphi(\bx)\) that satisfy \eqref{Eq:KoopmanEfun:Discrete} for \(\lambda \not = 0\) are also eigenfunctions of the Lie operator, although with a different eigenvalue.
Applying the Lie operator \eqref{Eq:Koopman:InfinitesimalGenerator} to such a \(\varphi\) leads to
\begin{equation}
  \label{eq:koopman-to-lie-eigenfunctions}
  \koop^t \varphi = \lambda^t \varphi \quad \Longrightarrow \quad \gen{} \varphi = \lim_{t \to 0}\frac{\koop^{t} \varphi - \varphi}{t} = \lim_{t \to 0}\frac{\lambda^{t} - 1}{t} \varphi = \log(\lambda) \varphi.
\end{equation}
Conversely, the induced dynamics \eqref{eq:koopman-dynamics-continuous} applied to an eigenfunction of \(\gen\) leads to
\begin{equation}
  \label{eq:lie-to-koopman-eigenfunctions}
  \gen{} \varphi = \mu \varphi \quad \Longrightarrow \quad \frac{d}{dt}\varphi = \gen \varphi = \mu \varphi.
\end{equation}
An eigenfunction \(\varphi\)  of \(\gen{}\) with eigenvalue \(\mu\) is then an eigenfunction of \(\koop^{t}\) with eigenvalue \(\lambda^{t} = \exp(\mu t)\).
Thus, we will typically not make a distinction between Lie and Koopman eigenfunctions in the context of autonomous dynamical systems.

Eigenfunctions of \(\koop, \gen\) that are induced by already-linear dynamics further illustrate the connection between linear discrete dynamics \(\bx_{n+1} = \bA \bx_{n}\) with analogous concepts for \(g_{n+1} = \koop g_{n}\).
Given a \emph{left}-eigenvector \(\bxi^{\top}\bA  = \lambda \bxi^{\top}\) {of the matrix $\bA$}, we form a corresponding Koopman eigenfunction as
\begin{align}
  \label{eq:eigenvectors-to-eigenfunctions}
  \varphi(\bx) &\coloneqq  \bxi^{\top}\bx  \\ \shortintertext{since}
  \koop \varphi(\bx) &= \varphi( \bA\bx ) = \bxi^{\top} \bA \bx = \lambda \bxi^{\top}\bx = \lambda \varphi(\bx).
\end{align}
In other words, while right-eigenvectors of \(\bA\) give rise to time-invariant directions in the state space \(\mathcal{X}\), {which will be known as {Koopman modes} or {dynamic modes}}, the left-eigenvectors giver rise to Koopman eigenfunctions, which are similarly time-invariant directions in the space of observables \(\mathcal{G}(\mathcal{X})\).

In general systems, a set of Koopman eigenfunctions may be used to generate more eigenfunctions.
In discrete time, we find that the product of two eigenfunctions $\varphi_1(\bx)$ and $\varphi_2(\bx)$ is also an eigenfunction:
\begin{equation}
\begin{aligned}
\koop \left(\varphi_1(\bx)\varphi_2(\bx)\right) &= \varphi_1(\bF(\bx))\varphi_2(\bF(\bx))\\
& = \lambda_1 \lambda_2 \varphi_1(\bx)\varphi_2(\bx)
\end{aligned}\label{Eq:KoopmanLattice:DT}
\end{equation}
with a new eigenvalue \(\lambda_1\lambda_2\) given by the product of the two eigenvalues of \(\varphi_1(\bx)\) and \(\varphi_2(\bx)\).
{ This argument assumes implicitly that the product of the two eigenfunction is again an eigenfunction or, more strongly, that the space of observables is closed under multiplication.}
The corresponding relationship for \(\gen\) can be found by applying \eqref{eq:koopman-to-lie-eigenfunctions},
\begin{equation}
  \label{eq:mult-to-add-lattice}
  \lambda_{1} \lambda_{2} = e^{\mu_{1}} e^{\mu_{2}} = e^{\mu_{1} + \mu_{2}},
\end{equation}
resulting in
\begin{equation}
  \label{eq:mult-to-add-latticeL}
  \gen{}\left(\varphi_1(\bx)\varphi_2(\bx)\right)= \left(\mu_1+\mu_2\right)\varphi_1(\bx)\varphi_2(\bx).
\end{equation}
A simple consequence is that a complex conjugate pair of eigenfunctions  of \(\gen{}\), \( (\mu, \varphi)\), \((\bar \mu\) and \(\bar \varphi)\), additionally implies the existence of a real-valued eigenfunction \(\abs{\varphi} = \sqrt{ \varphi \bar \varphi}\) with an associated eigenvalue \((\mu + \bar \mu)/2 = \operatorname{Re} \mu\), thus leading to a non-oscillatory growth/decay of the eigenvalue.

Algebraically, { when the space of observables is closed under multiplication}, the set of Koopman eigenfunctions establishes a commutative monoid\footnote{Monoids are groups without guaranteed inverses, or semigroups with an identity element.} under point-wise multiplication.
Thus, depending on the dynamical system, there may be a finite set of \emph{generator} eigenfunction elements that may be used to construct all other eigenfunctions.
{ Cardinality of the set of eigenfunctions and the relationships between them are further explored in~\cite{bollt2021}.}
The corresponding Koopman eigenvalues form a multiplicative lattice, or an additive lattice for Lie eigenvalues due to \eqref{eq:koopman-to-lie-eigenfunctions}.

Observables that can be formed as linear combinations of eigenfunctions, i.e., \(g \in \lspan\{\varphi_{k}\}_{k=1}^{K}\)  have a particularly simple evolution under the Koopman operator
\begin{equation}
  \label{eq:linear-combinations-of-eigenfunctions}
  g(\bx) = \sum_{k} v_{k} \varphi_{k} \quad \Longrightarrow \quad
  \koop^{t} g(\bx) = \sum_{k} v_{k} \lambda_{k}^{t} \varphi_{k}.
\end{equation}
This implies that the subspace \(\lspan\{\varphi_{k}\}_{k=1}^{K}\) is invariant under the action of \(\koop\).

These simple relationships enable the analysis of the phase portrait of a dynamical system in terms of level sets.
They also enable the analysis of the evolution of general observables through the lens of their decomposition into eigenfunctions, {which will provide the foundation for the {Koopman mode decomposition} of \Cref{Sec:KMD} and the {dynamic mode decomposition} of \Cref{Sec:DMD}.}

\subsubsection{Level sets of eigenfunctions}
Level sets of eigenfunctions can clarify the relationships between sets in the domain \(\mathcal{X}\).
In particular, {(sub-)}level sets associated with eigenfunctions map into each other.
To see this, we write the eigenvalue and the eigenfunction in polar form
\begin{align}
  \varphi(\bx) &= R(\bx) e^{i \Theta(\bx)}  \label{eq:eigenfunction-polar}\\
  \lambda &= r e^{i \theta}   \label{eq:eigenvalue-polar}
\end{align}
and define the associated sublevel sets to be
\begin{align}
  M_{\varphi}(C) &\coloneqq  \{\bx \in \mathcal{X} \colon R(\bx) \leq C\}, \label{eq:sublevel-magnitude}\\
  A_{\varphi}(\alpha) &\coloneqq  \{\bx \in \mathcal{X} \colon \Theta(\bx) \leq \alpha \} \label{eq:sublevel-angle}.
\end{align}
Given \(\bx \in \mathcal{X}\), we can evaluate the eigenfunction on its successor \(\bx^{+} = \bF(\bx)\), to find
\begin{align}
  \varphi(\bx^{+}) = \koop \varphi( \bx ) = \lambda \varphi(\bx) = (r R(\bx)) e^{i (\theta + \Theta(\bx))},
\end{align}
so that if \(\bx \in M_{\varphi}(C)\) then \(\bx^{+} \in M_{\varphi}(rC)\) and if \(\bx \in A_{\varphi}(\alpha)\) then \({\bx^{+} \in A_{\varphi}(\alpha + \theta)}\), or more succinctly
\begin{equation}
  \label{eq:level-sets}
  \bF( M_{\varphi}(C) ) = M_{\varphi}(rC), \quad \bF( A_{\varphi}(\alpha) ) = A_{\varphi}(\alpha + \theta).
\end{equation}

For the particular case of \(\lambda = 1\) these relationships simplify to
\begin{equation}
  \label{eq:invariant-level-sets}
  \bF( M_{\varphi}(C) ) = M_{\varphi}(C), \quad \bF( A_{\varphi}(\alpha) ) = A_{\varphi}(\alpha),
\end{equation}
implying that level sets of invariant eigenfunctions are invariant sets.

When the eigenvalue is on the unit circle, with a phase \(\theta/2\pi  = {p}/{q} \in \mathbb{Q}\), then the level sets of phase satisfy
\begin{equation}
  \label{eq:periodic-level-sets}
  \bF( A_{\varphi}(\alpha) ) = A_{\varphi}\left(\alpha + \frac{2\pi p}{q}\right)  \quad \bF^{q}( A_{\varphi}(\alpha) ) = A_{\varphi}(\alpha + 2\pi p) = A_{\varphi}(\alpha),
\end{equation}
and, therefore, identify chains of \(q\)-periodic sets where the value of the phase indicates the order in which they are visited by a trajectory \(\bx_{k}\).

\begin{figure}[t]
  \centering
\includegraphics[width=0.4\linewidth]{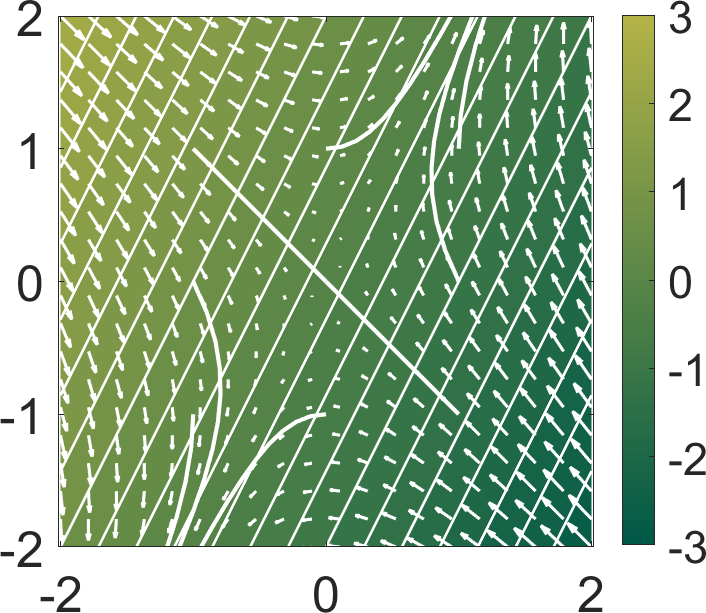}\hspace{2em}
  \includegraphics[width=0.4\linewidth]{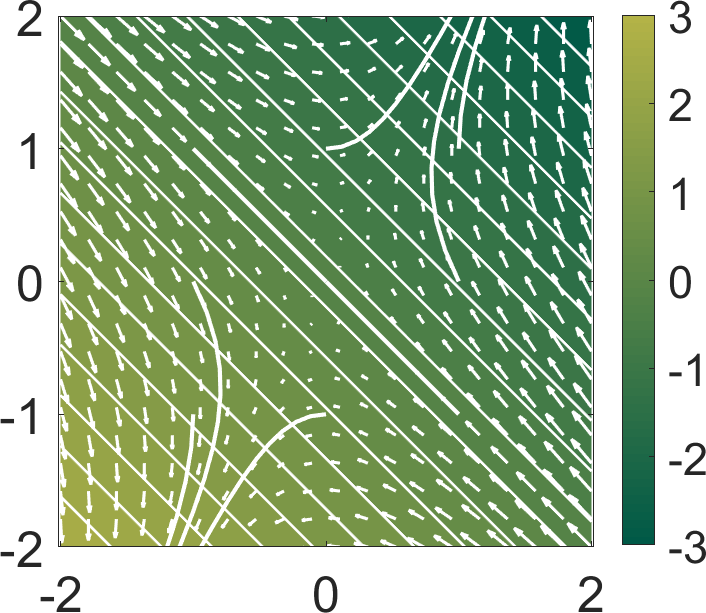}
  \caption{Koopman eigenfunctions constructed as \eqref{eq:eigenvectors-to-eigenfunctions} for a linear matrix ODE \(\dot\bx = \bA \bx\) with a saddle-type fixed point at the origin. Velocity field and sample orbits are overlaid in white.}
  \label{fig:koopman-efun-saddle}
\end{figure}

To illustrate these concepts, \cref{fig:koopman-efun-saddle} shows two eigenfunctions constructed from two left eigenvectors of a planar linear system \(\dot\bx = \bA\bx\) with saddle-type dynamics, as in \eqref{eq:eigenvectors-to-eigenfunctions}.
In contrast to right eigenvectors of \(\bA\), which act as invariant manifolds attracting/repelling the trajectories, Koopman eigenfunctions foliate the space into leaves (level sets), providing a time-ordering on the domain associated with each eigenvalue.

\begin{figure}[t]
  \centering
  \hspace{.2in}
  \begin{subfigure}[t]{0.45\linewidth}
\includegraphics[height=1.75in]{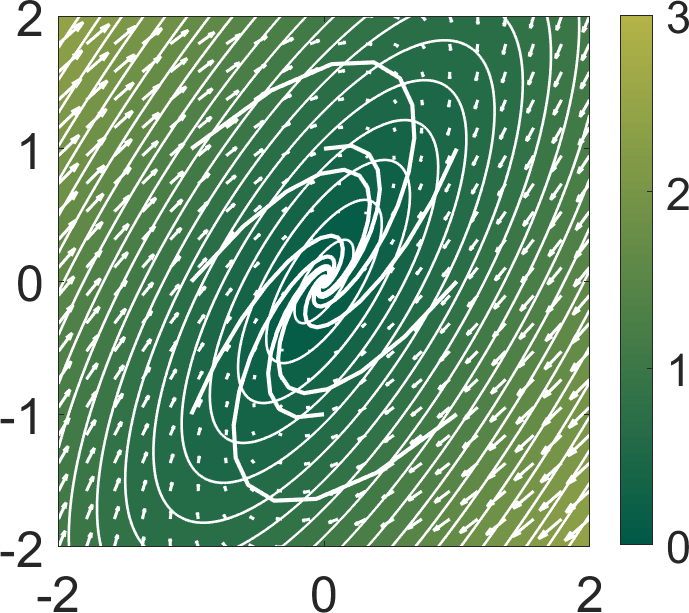}
\end{subfigure}
  \begin{subfigure}[t]{0.45\linewidth}
  \includegraphics[height=1.75in]{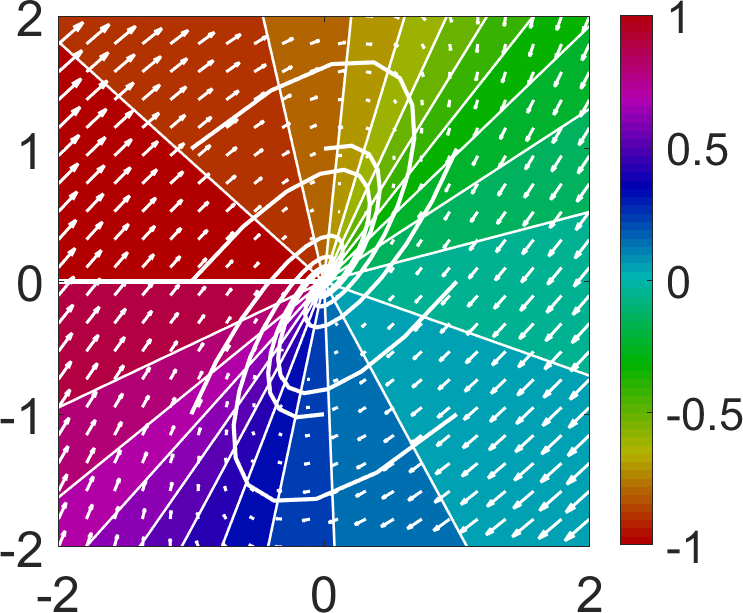}
  \end{subfigure}
  \caption{Modulus \(R(\bx)\) and argument \(\Theta(\bx)\)  of the Koopman eigenfunction constructed from a complex eigenvector of a linear matrix ODE \(\dot\bx = \bA \bx\) with a focus-type fixed point at the origin.}
  \label{fig:koopman-efun-focus}
\end{figure}

When {the} system matrix \(\bA\) has complex conjugate eigenvalues \(\lambda,\bar{\lambda}\), the associated eigenfunctions \(\bxi^{\top}\bx\) and \(\bar{\bxi}^{\top}\bx\) are themselves complex conjugates.
Visualizing the associated modulus \(R(\bx)\) and phase \(\Theta(\bx)\) \eqref{eq:eigenfunction-polar} for planar dynamics with a focus equilibrium, as in \cref{fig:koopman-efun-focus}, demonstrates that the modulus acts as a Lyapunov function for the stable focus, as sub-level sets provide a time-ordering on the plane implying that trajectories converge to the origin.
Indeed, \(\abs{\bxi^{\top} \bx}\) is itself an eigenfunction of the Koopman operator \(\koop^{t}\) associated with eigenvalues \(e^{t \Re \lambda}\), therefore capturing the exponential envelope of oscillating trajectories.
Level sets of \(\abs{\bxi^{\top} \bx}\) are called \emph{isostables}~\cite{mauroy2013}.
An analogous definition of an isostable for real-valued eigenvalues gives points that converge to the same trajectory in the stable/unstable manifold.

Level sets of the argument (angle) provide a cyclic foliation of the domain, acting as \emph{isochrons}~\cite{mauroy2013}, i.e., collections of initial conditions that converge to the origin with a common phase.
We will revisit the concepts of isochrons and isostables as tools for nonlinear reduction of order, in particular as a foundation for understanding the synchronization and control of oscillators, in \cref{sec:eigenfunctions-as-coordinates}.

\subsubsection{Computing eigenfunctions}
There are several approaches to approximate Koopman eigenfunctions computationally.
Given the evolution of an observable \(g(\bx(t))=g\left( \flow^t( \bx ) \right)\) it is possible to compute an eigenfunction associated with an eigenvalue \(e^{i\omega}\) by forming the following long-term \emph{harmonic} or \emph{Fourier average}
\begin{equation}
  \label{eq:harmonic-average}
  \tilde{g}_{\omega}(\bx) \coloneqq  \lim_{T\to\infty} \frac{1}{T}\int_{0}^{T} g\left( \flow^t( \bx ) \right) e^{-i \omega t} dt.
\end{equation}
If \(\omega = 0\), harmonic averages reduce to the trajectory (\emph{ergodic}) average.
When the space of observables \(\mathcal{G}(\mathcal{X})\) is taken as the \(L^{1}\) space of integrable functions with respect to a dynamically-preserved measure \(\mu\), then pointwise convergence \(\mu\)-almost everywhere is guaranteed by the Birkhoff ergodic theorem~\cite{birkhoff1931pnas}.

More generally, Yosida's theorem~\cite{budisic2012chaos} implies that in Banach spaces the limit not only exists, but that \(g \mapsto \tilde{g}_{\omega}\) is a projection onto the eigenspace associated with the eigenvalue \(e^{i\omega}\).
Said another way, eigenfunctions associated with eigenvalues on the unit circle \(\abs{\lambda} = 1\) can be computationally approximated through~\eqref{eq:harmonic-average}, as employed in~\cite{mezic2004physicad,mezic2005nd,budisic2012physd,levnajic2015}.
Since such eigenvalues correspond to neutrally-stable steady-state behavior, they are of practical importance to analyze dynamics.
This approach can, in principle, be extended to compute eigenfunctions associated with other eigenvalues as well, with conditions and restrictions given in~\cite{mohr2014arxiv,mohr2016a}.
{Recently, \cite{das2020a} also provided an extension and reframing of such results to RKHS, including both theoretical and computational results.
The existence and uniqueness of \(C^{k}\)-regular eigenfunctions was established in~\cite{kvalheim2021}, which has further implications on the use of such eigenfunctions in conjugacy arguments (see later \Cref{sec:eigenfunctions-as-coordinates}).}

Instead of the iteration~\eqref{eq:harmonic-average} that requires simulating the evolution \(g_{t}\), in certain cases it is possible to solve for eigenfunctions directly.
The eigenfunction relation \(\koop \varphi(\bx) = \varphi( \bF(\bx) ) = \lambda \varphi(\bx)\) is a compositional algebraic equation that is challenging to solve directly.
The relation~\eqref{eq:Lie-operator} implies that the analoguous relation for eigenfunctions of \(\gen\) is the PDE
\begin{equation}
  \nabla\varphi(\bx)\cdot\mathbf{f}(\bx)= \lambda\varphi(\bx).  \label{eq:lie-eigenfunction}
\end{equation}
With this PDE, it is possible to approximate eigenfunctions, either by solving for the Laurent series as in \cref{Sec:Koopman:Expansions} or with data via regression as in \cref{Sec:Observables}.
This formulation assumes that the dynamics are both continuous and differentiable.

\subsection{Koopman mode decomposition and finite representations}\label{Sec:KMD}
Until now, we have considered scalar measurements of a system, and we explored special \emph{eigen}-measurements (i.e., Koopman eigenfunctions) that evolve linearly in time.
However, we often take multiple measurements of a system, which we will arrange in a vector $\bg$:
\begin{align}
\bg(\bx) = \begin{bmatrix} g_1(\bx)\\ g_2(\bx)\\ \vdots \\ g_p(\bx)\end{bmatrix}.\label{Eq:Koopman:Measurement}
\end{align}
Each of the individual measurements may be expanded in terms of a basis of eigenfunctions $\varphi_j(\bx)$:
\begin{align}
g_i(\bx) = \sum_{j=1}^{\infty} v_{ij}\varphi_j(\bx).
\end{align}
Thus, the vector of observables, $\bg$, may be similarly expanded:
\begin{align}
\bg(\bx) = \begin{bmatrix} g_1(\bx)\\ g_2(\bx)\\ \vdots \\ g_p(\bx)\end{bmatrix} = \sum_{j=1}^{\infty}\varphi_j(\bx) \bv_j,\label{Eq:KoopmanMode}
\end{align}
where $\bv_j$ is {known as} the $j$-th \emph{Koopman mode} associated with the eigenfunction $\varphi_j$.

For conservative dynamical systems, such as those governed by Hamiltonian dynamics, the Koopman operator is unitary on the space \(\mathcal{G}(\mathcal{X})\) of square-integrable functions with respect to the conserved measure.
Thus, the Koopman eigenfunctions form an orthonormal basis for conservative systems, and it is possible to compute the Koopman modes $\bv_j$ directly by projection:
\begin{align}
\bv_j = \begin{bmatrix} \avg{ \varphi_j, g_1}\\ \avg{ \varphi_j, g_2 } \\ \vdots \\ \avg{ \varphi_j, g_p }\end{bmatrix},
\end{align}
where $\avg{ \cdot,\cdot }$ is the standard inner product of functions in \(\mathcal{G}(\mathcal{X})\).
{Thus the expansion of the observable function in \eqref{Eq:KoopmanMode} may be thought of as a change of basis into eigenfunction coordinates. }
These {Koopman} modes have a physical interpretation in the case of direct spatial measurements of a system, ${\bg(\bx) = \bx}$, in which case they are coherent \emph{spatial} modes that behave linearly with the same temporal dynamics (i.e., oscillations, possibly with linear growth or decay).
{The Koopman modes $\bv$ are also known as dynamic modes, and their computation is discussed in \Cref{Sec:DMD}.}

Given the decomposition in \eqref{Eq:KoopmanMode}, it is possible to represent the dynamics of the measurements $\bg$ as follows:
\begin{subequations}\label{Eq:KoopmanModeDecomposition}
\begin{align}
\bg(\bx(t)) = \mathcal{K}^t \bg(\bx_0)
& = \mathcal{K}^t\sum_{j=1}^{\infty} \varphi_j(\bx_0) \bv_j\\
& = \sum_{j=1}^{\infty} \mathcal{K}^t\varphi_j(\bx_0) \bv_j\\
& = \sum_{j=1}^{\infty} \lambda_j^t\varphi_j(\bx_0) \bv_j.
\end{align}
\end{subequations}
This sequence of triples $\{(\lambda_j, \varphi_j, \bv_j)\}_{j=1}^{\infty}$ is the \emph{Koopman mode decomposition} and was introduced by Mezi\'{c} in 2005~\cite{mezic2005nd}.
{Often, it is possible to approximate this expansion as a truncated sum of only a few dominant terms.}
The Koopman mode decomposition was later connected to data-driven regression via the dynamic mode decomposition (DMD)~\cite{rowley2009jfm}, which is explored in~\cref{Sec:DMD}.
{The DMD eigenvalues will approximate the Koopman eigenvalues $\lambda_j$, the DMD modes will approximate the Koopman modes $\bv_j$, and the DMD mode amplitudes will approximate the corresponding Koopman eigenfunctions evaluated at the initial condition $\varphi_j(\bx_0)$.
It is important to note that the Koopman modes and eigenfunctions are distinct mathematical objects, requiring different approaches for approximation.
Koopman eigenfunctions are often more challenging to compute than Koopman modes, motivating advanced techniques, such as the extended DMD algorithm~\cite{williams2015jnls} in~\Cref{Sec:EDMD}.}

\subsubsection{Invariant eigenspaces and finite-dimensional models}
Instead of capturing the evolution of all measurement functions in a Hilbert space, applied Koopman analysis approximates the evolution on an invariant subspace spanned by a finite set of measurement functions.
A \emph{Koopman-invariant subspace} is defined as the span of a set of functions $\{g_1, g_2, \cdots , g_p\}$ if all functions $g$ in this subspace
\begin{eqnarray}
g=\alpha_1g_{1}+\alpha_2g_{2}+\cdots+\alpha_pg_{p}
\end{eqnarray}
remain in this subspace after being acted on by the Koopman operator $\mathcal{K}$:
\begin{eqnarray}
\mathcal{K}g =\beta_1 g_{1} + \beta_2 g_{2} + \cdots + \beta_p g_{p}.
\end{eqnarray}
It is possible to obtain a finite-dimensional matrix representation of the Koopman operator by restricting it to an invariant subspace spanned by a finite number of functions $\{g_j\}_{j=1}^{p}$.
The matrix representation $\mathbf{K}$ acts on a vector space $\mathbb{R}^p$, with the coordinates given by the values of $g_{j}(\mathbf{x})$.
This induces a finite-dimensional linear system.

Any finite set of eigenfunctions of the Koopman operator will span an invariant subspace.
 Discovering these eigenfunction coordinates is, therefore, a central challenge, as they provide intrinsic coordinates along which the dynamics behave linearly.
In practice, it is more likely that we will identify an \emph{approximately} invariant subspace, given by a set of functions $\{g_j\}_{j=1}^{p}$, where each of the functions $g_j$ is well approximated by a finite sum of eigenfunctions: $g_j\approx \sum_{k=1}^{p}v_{jk}\varphi_k$.

\subsection{Example of a simple Koopman embedding}\label{Sec:SlowManifold}
Here, we consider an example system with a single fixed point from Tu et al.~\cite{tu2014jcd} that is explored in more detail in Brunton et al.~\cite{brunton2016plosone}, given by:
\begin{subequations}
\begin{align}
\dot x_1 &= \mu x_1\\
\dot x_2 &= \lambda(x_2-x_1^2).
\end{align}
\end{subequations}
For $\lambda<\mu<0$, the system exhibits a slow attracting manifold given by $x_2=x_1^2$.
It is possible to augment the state $\bx$ with the nonlinear measurement $g=x_1^2$, to define a three-dimensional Koopman invariant subspace.
In these coordinates, the dynamics become linear:
\begin{subequations}
\begin{align}
\frac{d}{dt}\begin{bmatrix} y_1\\ y_2\\ y_3\end{bmatrix} =
\begin{bmatrix} \mu & 0 & 0 \\ 0& \lambda& - \lambda\\ 0 & 0 & 2\mu\end{bmatrix}
\begin{bmatrix}y_1\\ y_2\\ y_3\end{bmatrix}
\quad \text{for}\quad
\begin{bmatrix}y_1\\ y_2\\ y_3\end{bmatrix} =
\begin{bmatrix}x_1\\ x_2\\x_1^2\end{bmatrix}.\label{Eq:QuadraticAttractor}
\end{align}
\end{subequations}

\begin{figure}
\vspace{-.1in}
\begin{center}
\begin{overpic}[width=.95\textwidth]{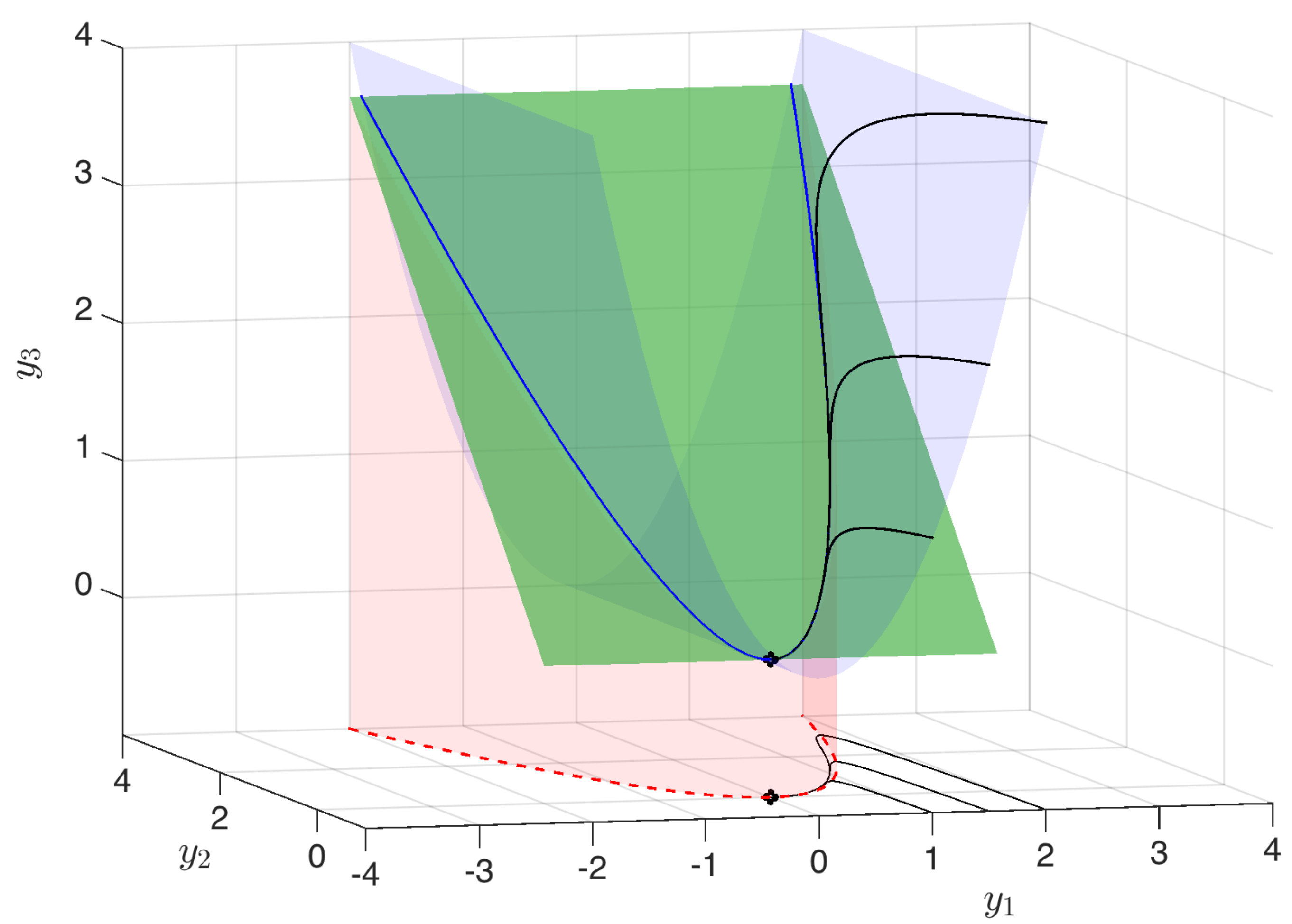}
\end{overpic}
\vspace{-.1in}
\caption{Visualization of three-dimensional linear Koopman system from \eqref{Eq:QuadraticAttractor} along with projection of dynamics onto the $x_1$-$x_2$ plane.  The attracting slow manifold is shown in red, the constraint $y_3=y_1^2$ is shown in blue, and the slow unstable subspace of \eqref{Eq:QuadraticAttractor} is shown in green.  Black trajectories of the linear Koopman system in $\mathbf{y}$ project onto trajectories of the full nonlinear system in $\mathbf{x}$ in the $y_1$-$y_2$ plane.  Here, $\mu=-0.05$ and $\lambda=1$.  \textit{Reproduced from Brunton et al.~\cite{brunton2016plosone}.}}\label{Fig:quad_manifold}
\end{center}
\vspace{-.1in}
\end{figure}

The full three-dimensional Koopman observable vector space is visualized in \cref{Fig:quad_manifold}.
Trajectories that start on the invariant manifold $y_3=y_1^2$, visualized by the blue surface, are constrained to stay on this manifold.
There is a \emph{slow} subspace, spanned by the eigenvectors corresponding to the slow eigenvalues $\mu$ and $2\mu$; this subspace is visualized by the green surface.
Finally, there is the original asymptotically attracting manifold of the original system, $y_2=y_1^2$, which is visualized as the red surface.
The blue and red parabolic surfaces always intersect in a parabola that is inclined at a $45^{\circ}$ angle in the $y_2$-$y_3$ direction.
The green surface approaches this $45^{\circ}$ inclination as the ratio of fast to slow dynamics become increasingly large.
In the full three-dimensional Koopman observable space, the dynamics produce a single stable node, with trajectories rapidly attracting onto the green subspace and then slowly approaching the fixed point.

The left eigenvectors of the Koopman operator yield Koopman eigenfunctions.
The Koopman eigenfunctions of \eqref{Eq:QuadraticAttractor} corresponding to eigenvalues $\mu$ and $\lambda$ are:
\begin{eqnarray}
\varphi_{\mu}=x_1, \quad\text{and}\quad\varphi_{\lambda} = x_2-bx_1^2 \quad\text{with}\quad b =\frac{\lambda}{\lambda-2\mu}.
\end{eqnarray}
The constant $b$ in $\varphi_\lambda$ captures the fact that for a finite ratio $\lambda/\mu$, the dynamics only shadow the asymptotically attracting slow manifold $x_2=x_1^2$, but in fact follow neighboring parabolic trajectories.
This is illustrated more clearly by the various surfaces in \cref{Fig:quad_manifold} for different ratios $\lambda/\mu$.

In this way, a set of intrinsic coordinates may be determined from the observable functions defined by the left eigenvectors of the {matrix representation of the} Koopman operator on an invariant subspace.
Explicitly,
\begin{eqnarray}
\varphi_{\alpha}(\bx) = \bxi^{\top}_{\alpha}\bz(\bx), \quad\text{where}\quad \bxi^{\top}_{\alpha}\mathbf{K} = \alpha\bxi^{\top}_{\alpha}.
\end{eqnarray}
These eigen-observables define observable subspaces that remain invariant under the Koopman operator, even after coordinate transformations.
As such, they may be regarded as intrinsic coordinates~\cite{williams2015jnls} on the Koopman-invariant subspace.

\subsection{Analytic series expansions for eigenfunctions}\label{Sec:Koopman:Expansions}
Given the dynamics in \eqref{Eq:ContinuousDynamicsGeneral}, it is possible to solve the PDE in \eqref{eq:lie-eigenfunction} using standard techniques, such as recursively solving for the terms in a Taylor or Laurent series.
A number of simple examples are explored below.

\subsubsection{Linear dynamics}
Consider the simple linear dynamics
\begin{equation}
\frac{d}{dt}x =  x.
\end{equation}
Assuming a Taylor series expansion for $\varphi(x)$:
\begin{align*}
\varphi(x) &= c_0 + c_1x + c_2x^2 + c_3x^3 + \cdots
\end{align*}
then the gradient and directional derivatives are given by:
\begin{align*}
\nabla \varphi & = c_1 + 2c_2x + 3 c_3x^2 + 4c_4x^3+\cdots\\
\nabla\varphi\cdot f & = c_1 x + 2c_2 x^2 +3c_3 x^3 + 4c_4 x^4+\cdots
\end{align*}
Solving for terms in the Koopman eigenfunction PDE \eqref{eq:lie-eigenfunction}, we see that $c_0=0$ must hold.
For any positive integer $\lambda$ in \eqref{eq:lie-eigenfunction}, only one of the coefficients may be nonzero.
Specifically, for $\lambda=k\in\mathbb{Z}^+$, then $\varphi(x) = cx^k$ is an eigenfunction for any constant $c$.
For instance, if $\lambda=1$, then $\varphi(x) =  x$.

\subsubsection{Quadratic nonlinear dynamics}
Consider a nonlinear dynamical system
\begin{align}
\frac{d}{dt}=x^2.
\end{align}
There is no Taylor series that satisfies \eqref{eq:lie-eigenfunction}, except the trivial solution $\varphi=0$ for $\lambda=0$.
Instead, we assume a Laurent series:
\begin{align*}
\varphi(x) &= \cdots +c_{-3}x^{-3} + c_{-2}x^{-2} + c_{-1}x^{-1} +  c_0 +c_1x + c_2x^2 + c_3x^3 + \cdots.
\end{align*}
The gradient and directional derivatives are given by:
\begin{align*}
\nabla \varphi & =\cdots -3c_{-3}x^{-4} - 2c_{-2}x^{-3} -c_{-1}x^{-2} + c_1 + 2c_2x + 3 c_3x^2 + 4c_4x^3+\cdots\\
\nabla\varphi\cdot f & =\cdots -3c_{-3}x^{-2} - 2c_{-2}x^{-1} -c_{-1} + c_1x^2 + 2c_2x^3 + 3 c_3x^4 + 4c_4x^5+\cdots.
\end{align*}
Solving for the coefficients of the Laurent series that satisfy \eqref{eq:lie-eigenfunction}, we find that all coefficients with positive index are zero, i.e. $c_k=0$ for all $k\geq 1$.
However, the nonpositive index coefficients are given by
 the recursion $\lambda c_{k+1} = k c_{k}$, for negative $k\leq -1$.
Thus, the Laurent series is
\begin{align*}
\varphi(x) &= c_0\left(1 - \lambda x^{-1} + \frac{\lambda^2}{2}x^{-2} - \frac{\lambda^3}{3!}x^{-3}+\cdots\right)=c_0e^{-\lambda/x}.
\end{align*}
This holds for all values of $\lambda\in\mathbb{C}$.  There are also other Koopman eigenfunctions that can be identified from the Laurent series.

\subsubsection{Polynomial nonlinear dynamics}
For a more general nonlinear dynamical system
\begin{align}
\frac{d}{dt}=ax^n,
\end{align}
$\varphi(x)=e^{\frac{\lambda}{(1-n)a}x^{1-n}}$ is an eigenfunction for all $\lambda\in\mathbb{C}$.

As mentioned in \cref{sec:koopman-eigenfunctions}, it is also possible to generate new eigenfunctions by taking powers of these primitive eigenfunctions; the resulting eigenvalues generate a lattice in the complex plane.

\section{Dynamic mode decomposition}\label{Sec:DMD}

Dynamic mode decomposition, originally introduced by Schmid~\cite{schmid2008aps,schmid2010jfm} in the fluid dynamics community, has rapidly become the standard algorithm to approximate the Koopman operator from data~\cite{rowley2009jfm,tu2014jcd,kutz2016book}.
Rowley et al.~\cite{rowley2009jfm} established the first connection between DMD and the Koopman operator.
The DMD algorithm was originally developed to identify spatio-temporal coherent structures from high-dimensional time-series data, as are commonly found in fluid dynamics.
DMD is based on the computationally efficient singular value decomposition (SVD), also known as proper orthogonal decomposition (POD) in fluid dynamics, so that it provides scalable dimensionality reduction for high-dimensional data.
The SVD orders modes hierarchically based on how much of the variance of the original data is captured by each mode; these modes remain invariant even when the order of the data is shuffled in time.
In contrast, the DMD modes are linear combinations of the SVD modes that are chosen specifically to extract spatially correlated structures that have the same coherent linear behavior in time, given by oscillations at a fixed frequency with growth or decay.
Thus, DMD provides dimensionality reduction in terms of a low-dimensional set of spatial modes along with a linear model for how the amplitudes of these modes evolve in time.
In this way, DMD may be thought of as a combination of SVD/POD in space with the Fourier transform in time, combining the strengths of each approach~\cite{chen2012jns,kutz2016book}.

{Several leading DMD variants, especially DMD with control and delay DMD, are closely related to subspace system identification methods, many of which predate DMD by decades.
However, modern Koopman theory provides a new interpretation for these methods when applied to nonlinear systems.
For example, there are close connections between delay DMD and the eigensystem realization algorithm (ERA), although the classical theory of ERA is only valid for strictly linear systems, while delay DMD approaches may be applied more broadly to nonlinear, and even chaotic systems.
Similarly, DMD and Koopman based approaches tend to also have the perspective of applying to very high dimensional systems, including the numerical considerations that must be accounted for.}

There are a number of factors that have led to the widespread adoption of DMD as a workhorse algorithm for processing high-dimensional spatiotemporal data.
The DMD algorithm approximates the best-fit linear {matrix} operator that advances high-dimensional measurements of a system forward in time~\cite{tu2014jcd}.
Thus, DMD approximates the Koopman operator restricted to the measurement subspace given by direct measurements of the state of a system.
DMD is valid for both experimental and simulated data, as it is based entirely on measurement data and does not require knowledge of the governing equations.
In addition, DMD is highly extensible because of its simple formulation in terms of linear algebra, resulting in innovations related to control, compressed sensing, and multi-resolution, among others.
Because of these strengths, DMD has been applied to a wide range of diverse applications beyond fluid mechanics, including neuroscience, disease modeling, robotics, video processing, power grids, financial markets, and plasma physics.
Many of these extensions and applications will be discussed more here and in \cref{Sec:Observables}.

\subsection{The DMD algorithm}
The DMD algorithm seeks a best fit linear {matrix} operator $\mathbf{A}$ that approximately advances the state of a system, $\mathbf{x}\in\mathbb{R}^n$, forward in time according to the linear dynamical system
\begin{align}\label{Eq:DMD:Propagator}
    \mathbf{x}_{k+1} = \mathbf{A}\mathbf{x}_k,
\end{align}
where $\mathbf{x}_k=\mathbf{x}(k\Delta t)$, and $\Delta t$ denotes a fixed time step that is small enough to resolve the highest frequencies in the dynamics.  Thus, the {matrix} $\mathbf{A}$ is an approximation of the Koopman operator $\mathcal{K}$ restricted to a measurement subspace spanned by direct measurements of the state $\mathbf{x}$.

DMD is fundamentally a data-driven algorithm, and the {matrix} $\mathbf{A}$ is approximated from a collection of snapshot pairs of the system, $\{(\bx(t_k),\bx(t_k')\}_{k=1}^m$, where $t_k' = t_k + \Delta t$.
A snapshot is typically a measurement of the full state of the system, such as a fluid velocity field sampled at a large number of spatially discretized locations, that is reshaped into a column vector of high dimension.
The original formulation of Schmid~\cite{schmid2010jfm} required data from a single trajectory with uniform sampling in time, so that $t_k=k\Delta t$.
Here we present the \emph{exact} DMD algorithm of Tu et al.~\cite{tu2014jcd}, which works for irregularly spaced data and data concatenated from multiple different time series.
Thus, in exact DMD, the times $t_k$ need not be sequential or evenly spaced, but for each snapshot $\mathbf{x}(t_k)$ there is a corresponding snapshot $\mathbf{x}(t_k')$ one time step $\Delta t$ in the future.
These snapshots are arranged into two data matrices, $\bX$ and $\bX'$:
\begin{subequations}
\begin{align}
\bX &= \begin{bmatrix} \vline & \vline & & \vline \\
\bx(t_1) & \bx(t_2) & \cdots & \bx(t_m) \\
 \vline & \vline & & \vline
 \end{bmatrix}\\
 \bX' &= \begin{bmatrix} \vline & \vline & & \vline \\
\bx(t_1') & \bx(t_2') & \cdots & \bx(t_m') \\
 \vline & \vline & & \vline
 \end{bmatrix}.
\end{align}
\end{subequations}

Equation~\eqref{Eq:DMD:Propagator} may be written in terms of these data matrices as
\begin{align}
\bX' \approx \bA \bX. \label{eq:DMD-matrix}
\end{align}
The best fit {matrix} $\bA$ establishes a linear dynamical system that approximately advances snapshot measurements forward in time, which may be formulated as an optimization problem
\begin{align}
\bA = \argmin_{\bA} \|\bX' - \bA \bX\|_F = \bX'\bX^\dagger\label{Eq:DMD:Definition}
\end{align}
where $\|\cdot\|_F$ is the Frobenius norm and $^\dagger$ denotes the pseudo-inverse.  The pseudo-inverse may be computed using the SVD of $\mathbf{X}=\mathbf{U}\boldsymbol{\Sigma}\mathbf{V}^*$ as $\mathbf{X}^\dagger=\mathbf{V}\boldsymbol{\Sigma}^{-1}\mathbf{U}^*$.  The matrices $\mathbf{U}\in\mathbb{C}^{n\times n}$ and $\mathbf{V}^{m\times m}$ are unitary, so that $\mathbf{U}^*\mathbf{U}=\mathbf{I}$ and $\mathbf{V}^*\mathbf{V}=\mathbf{I}$, where $^*$ denotes complex-conjugate transpose.  The columns of $\mathbf{U}$ are known as POD modes~\cite{taira2017,taira2020}.

{Because $\bA$ is an approximate representation of the Koopman operator restricted to a finite-dimensional subspace of linear measurements, we are often interested in the eigenvectors $\bPhi$ and eigenvalues $\bLambda$ of $\bA$:
\begin{align}
    \bA \bPhi = \bPhi \bLambda.
\end{align}
However,} the matrix $\bA$ has $n^2$ elements, so for high-dimensional data it may be intractable to represent{, let alone compute its eigendecomposition}.
Instead, the DMD algorithm seeks the leading spectral decomposition (i.e., eigenvalues and eigenvectors) of $\bA$ without ever explicitly constructing it.
The data matrices $\mathbf{X}$ and $\mathbf{X}'$ typically have far more rows than columns, i.e. $m\ll n$, so that $\mathbf{A}$ will have at most $m$ nonzero eigenvalues and non-trivial eigenvectors.
In practice, the effective rank of the data matrices, and hence the {matrix} $\mathbf{A}$, may be even lower, given by $r<m$.
Instead of computing $\mathbf{A}$ in \eqref{Eq:DMD:Definition}, we may project $\mathbf{A}$ onto the first $r$ POD modes in $\mathbf{U}_r$ and approximate the pseudo-inverse using the rank-$r$ SVD approximation $\bX\approx \mathbf{U}_r\boldsymbol{\Sigma}_r\mathbf{V}_r^*$:
\begin{subequations}
\begin{align}
\bAtilde &= \mathbf{U}_r^*\mathbf{A}\mathbf{U}_r \\
&= \mathbf{U}_r^*\mathbf{X}'\mathbf{X}^\dagger\mathbf{U}_r\\
& = \mathbf{U}_r^*\mathbf{X}'\mathbf{V}_r\boldsymbol{\Sigma}_r^{-1}\mathbf{U}_r^*\mathbf{U}_r\\
& = \mathbf{U}_r^*\mathbf{X}'\mathbf{V}_r\boldsymbol{\Sigma}_r^{-1}.
\end{align}
\end{subequations}

The leading spectral decomposition of $\mathbf{A}$ may be approximated from the spectral decomposition of the much smaller $\bAtilde$:
\begin{align}
\bAtilde \bW = \bW \bLambda.
\end{align}
The diagonal matrix $\bLambda$ contains the \emph{DMD eigenvalues}, which correspond to eigenvalues of the high-dimensional matrix $\bA$.
The columns of $\bW$ are eigenvectors of $\bAtilde$, and provide a coordinate transformation that diagonalizes the matrix.  These columns may be thought of as linear combinations of POD mode amplitudes that behave linearly with a single temporal pattern given by the corresponding eigenvalue $\lambda$.

\begin{figure}
\begin{center}
\includegraphics[width=\textwidth]{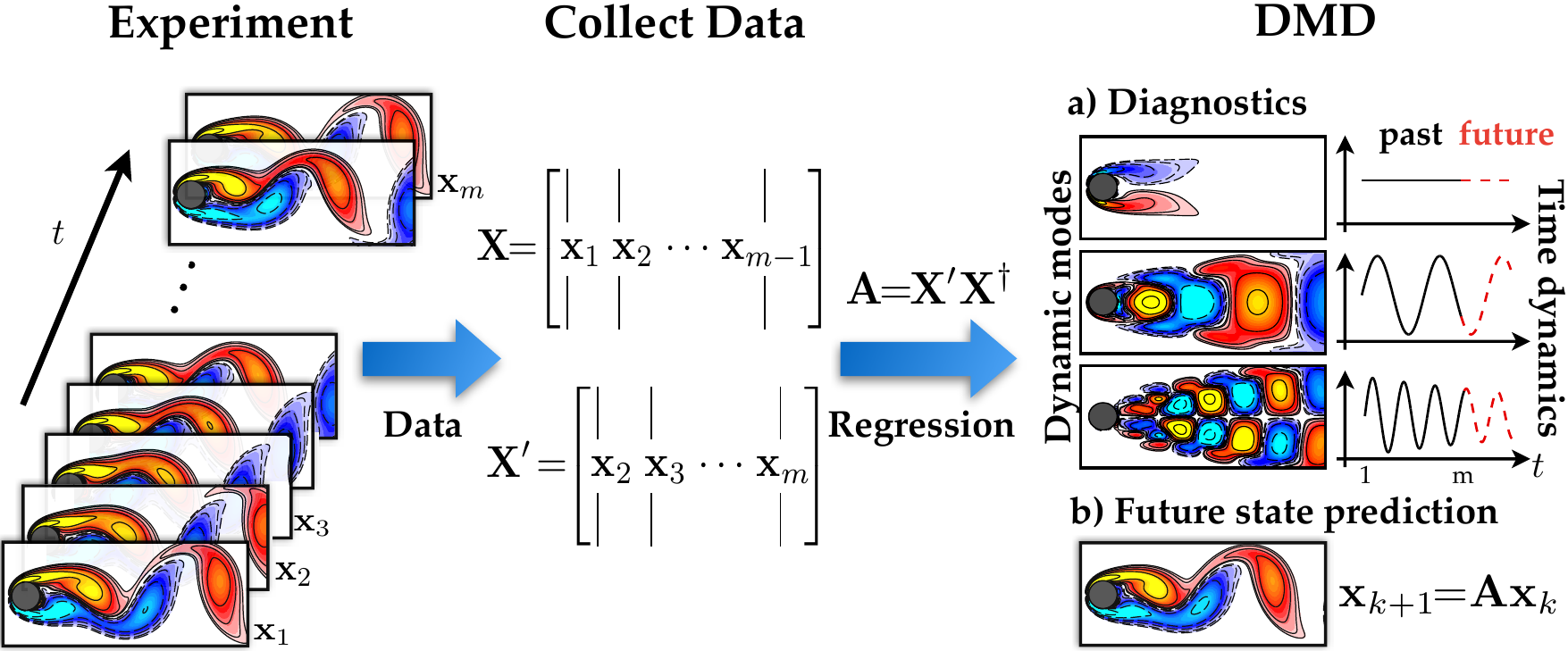}
\caption{Overview of DMD illustrated on the fluid flow past a circular cylinder at Reynolds number $100$.  \emph{Reproduced from Kutz et al.~\cite{kutz2016book}.}}\label{Fig:DMDOVERVIEW}
\end{center}
\end{figure}

The eigenvectors of $\bA$ are the \emph{DMD modes} $\bPhi$, and they are reconstructed using the eigenvectors $\bW$ of the reduced system and the time-shifted data matrix $\bX'$:
\begin{align}
\bPhi = \bX' \bVtilde \bSigmatilde^{-1}\bW.
\end{align}
Tu et al.~\cite{tu2014jcd} proved that these DMD modes are eigenvectors of the full $\bA$ matrix under certain conditions.
This approach is illustrated for a fluid flow in \cref{Fig:DMDOVERVIEW}.
There are also several open-source DMD implementations~\cite{kutz2016book,Demo18pydmd}.

\subsubsection{Spectral decomposition and the DMD expansion}
Once the DMD modes and eigenvalues are computed, it is possible to represent the system state in terms of the DMD expansion:
\begin{align}\label{Eq:DMD:SpectralExpansion}
\bx_k = \sum_{j=1}^r \bphi_j \lambda_j^{k-1} b_j = \bPhi \bLambda^{k-1} \bb,
\end{align}
where $\bphi_j$ are eigenvectors of $\bA$ (DMD modes), $\lambda_j$ are eigenvalues of $\bA$ (DMD eigenvalues), and $b_j$ are the mode amplitudes.
{The DMD expansion \eqref{Eq:DMD:SpectralExpansion} is directly analogous to the Koopman mode decomposition \eqref{Eq:KoopmanModeDecomposition}.
The DMD modes $\bphi_j$ approximate the Koopman modes $\bv_j$, the DMD eigenvalues $\lambda_j$ approximate the corresponding Koopman eigenvalues, and the mode amplitudes $b_j$ approximate the Koopman eigenfunctions evaluated at the initial condition $\varphi_j(\bx_0)$.
}

{To directly connect the DMD expansion to the Koopman mode decomposition \eqref{Eq:KoopmanModeDecomposition}, we may write the KMD explicitly in matrix form for the observable $\bg(\bx) = \bx$ as
\begin{subequations}
\begin{align}
    \bx_k \approx \sum_{j=1}^r\lambda_j^{k-1} \varphi_j(\bx_1) \bv_j & = \underbrace{
    \left[ \begin{array}{ccc} | & & | \\ \bv_1 & \cdots & \bv_r \\ | & & | \end{array}  \right]}_{\bPhi}
    \underbrace{\left[ \begin{array}{ccc} \lambda_1 &  & \\ & \ddots & \\ & & \lambda_r  \end{array} \right]}_{\bLambda}
    \underbrace{\left[ \begin{array}{c}\varphi_1(\bx_1) \\
    \vdots  \\ \varphi_r(\bx_1)  \end{array} \right]}_{\bb}.
\end{align}
\end{subequations}
Thus, comparing with the DMD expansion \eqref{Eq:DMD:SpectralExpansion}, the correspondence of terms is clear.
The KMD expansion may be written equivalently as
\begin{align}
    \bx_k \approx \sum_{j=1}^r\lambda_j^{k-1} \varphi_j(\bx_1) \bv_j =
    \left[ \begin{array}{ccc} | & & | \\ \bv_1 & \cdots & \bv_r \\ | & & | \end{array}  \right]
    \left[ \begin{array}{ccc} \varphi_1(\bx_1) &  & \\ & \ddots & \\ & & \varphi_r(\bx_1)  \end{array} \right]
    \left[ \begin{array}{c} \lambda_1  \\
    \vdots  \\ \lambda_r  \end{array} \right]
\end{align}
which makes it possible to express the data matrix $\bX$ as
\begin{align}
    \bX =
    \left[ \begin{array}{ccc} | & & | \\ \bphi_1 & \cdots & \bphi_r \\ | & & | \end{array}  \right]
    \left[ \begin{array}{ccc} b_1 &  & \\ & \ddots & \\ & & b_r  \end{array} \right]
    \left[ \begin{array}{ccc} \lambda_1 & \cdots & \lambda_1^{m-1} \\
    \vdots & \ddots & \vdots \\ \lambda_r & \cdots & \lambda_r^{m-1} \end{array} \right].
\end{align}
}

The amplitudes in $\bb$ are {often} given by
\begin{align}
\bb = \bPhi^{\dagger}\bx_1,\label{Eq:DMDModeAmplitude}
\end{align}
{using the first snapshot to determine the mixture of DMD mode amplitudes; note that this first snapshot $\bx_1$ from DMD is equivalent to the initial condition used to evaluate the Koopman eigenfunction in \eqref{Eq:KoopmanModeDecomposition}.}
Alternative approaches to compute $\bb$~\cite{chen2012jns,jovanovic2014pof,askham2017arxiv} will be discussed in \cref{Sec:DMD:AlternativeOpt}.

The spectral expansion in \eqref{Eq:DMD:SpectralExpansion} may be converted to continuous time by introducing the continuous eigenvalues $\omega = \log(\lambda)/\Delta t$:
\begin{align}
\bx(t) = \sum_{j=1}^r \bphi_j e^{\omega_j t} b_j = \bPhi \exp(\bOmega t)\bb,
\end{align}
where $\bOmega$ is a diagonal matrix containing the continuous-time eigenvalues $\omega_j$.
Thus, the data matrix $\mathbf{X}$ may be represented as
\begin{equation}
    \bX \approx \left[ \begin{array}{ccc} | & & | \\ \boldsymbol{\phi}_1 & \cdots & \boldsymbol{\phi}_r \\ | & & | \end{array}  \right]
    \left[ \begin{array}{ccc} b_1 &  & \\ & \ddots & \\ & & b_r  \end{array} \right]
    \left[ \begin{array}{ccc} e^{\omega_1 t_1} & \cdots & e^{\omega_1 t_m} \\
    \vdots & \ddots & \vdots \\ e^{\omega_r t_1} & \cdots & e^{\omega_r t_m} \end{array} \right]
    =  \bPhi \mbox{diag}(\bb) {\bf T}(\boldsymbol{\omega}).
\label{eq:dmd_opt}
\end{equation}

\subsubsection{Alternative optimizations for DMD}\label{Sec:DMD:AlternativeOpt}

The DMD algorithm is purely data-driven, and is thus equally applicable to experimental and numerical data.
When characterizing experimental data with DMD, the effects of sensor noise and stochastic disturbances must be accounted for.
Bagheri~\cite{bagheri2014} showed that DMD is particularly sensitive to the effects of noisy data, and it has been shown that significant and systematic biases are introduced to the eigenvalue distribution~\cite{duke2012error,bagheri2013jfm,dawson2016ef,hemati2017tcfd}.
Although increased sampling decreases the variance of the eigenvalue distribution, it does not remove the bias~\cite{hemati2017tcfd}.
This noise sensitivity has motivated several alternative optimization algorithms for DMD to improve the quality and performance of DMD over the standard optimization in \eqref{Eq:DMD:Definition}, which is a least-square fitting procedure involving the Frobenius norm.
These algorithms include the total least-squares DMD~\cite{hemati2017tcfd}, forward-backward DMD~\cite{dawson2016ef}, variable projection~\cite{askham2017arxiv}, and robust principal component analysis~\cite{Scherl2020prf}.

One of the simplest ways to remove the systematic bias of the DMD algorithm is by computing it both forward and backward in time and averaging the equivalent {matrices}, as proposed by Dawson et al.~\cite{dawson2016ef}.  Thus the two following approximations are considered
\begin{align}
\bX' \approx \bA_1 \bX  \hspace{0.5in} \mbox{and} \hspace{0.5in} \bX \approx \bA_2 \bX'
\end{align}
where $\bA_2^{-1}\approx\bA_1$ for noise-free data.  Thus the {matrix} $\bA_2$ is the inverse, or backward time-step, mapping the snapshots from $t_{k+1}$ to $t_{k}$.
The forward and backward time {matrices} are then averaged, removing the systematic bias from the measurement noise:
\begin{equation}
    \bA = \frac{1}{2} \left( \bA_1 + \bA_2^{-1} \right)
\end{equation}
where the optimization \eqref{Eq:DMD:Definition} can be used to compute both the forward and backward mapping $\bA_1$ and $\bA_2$.
This optimization can be formulated as
\begin{equation}
    \bA = \argmin_{\bA} \frac{1}{2} \left( \|\bX' - \bA \bX\|_F   +  \|\bX - \bA^{-1} \bX'\|_F \right),
    \label{eq:fb1}
\end{equation}
which is highly nonlinear and non-convex due to the inverse $\bA^{-1}$.
An improved optimization framework was developed by Azencot et al.~\cite{azencot2019consistent} which proposes
\begin{equation}
    \bA = \argmin_{\bA_1, \bA_2} \frac{1}{2} \left( \|\bX' - \bA_1 \bX\|_F   +  \|\bX - \bA_2 \bX'\|_F \right)
    \,\,\, \mbox{s.t.} \,\,\, \bA_1\bA_2={\bf I}, \,\, \bA_2\bA_1={\bf I},
    \label{eq:fb2}
\end{equation}
to circumvent some of the difficulties of the optimization in \eqref{eq:fb1}.

Hemati et al.~\cite{hemati2017tcfd} formulate another DMD algorithm, replacing the original least-squares regression with a total least-squares regression to account for the possibility of noisy measurements and disturbances to the state.
This work also provides an excellent discussion on the sources of noise and a comparison of various denoising algorithms.
The subspace DMD algorithm of Takeishi et al.~\cite{takeishi2017} compensates for measurement noise by computing an orthogonal projection of future snapshots onto the space of previous snapshots and then constructing a linear model.
Extensions that combine DMD with Bayesian approaches have also been developed~\cite{takeishi2017jcai}.

Good approximations for the  mode amplitudes ${\bf b}$ in \eqref{eq:dmd_opt} have also proven to be difficult to achieve, with and without noise.
Jovanovi\'c et al.~\cite{jovanovic2014pof} developed the first algorithm to improve the estimate of the modal amplitudes by promoting sparsity.
In this case, the underlying optimization algorithm is framed around improving the approximation \eqref{eq:dmd_opt} using the formulation
\begin{equation}
    \argmin_{{\bf b}} \left( \| \bX -   \bPhi \mbox{diag}(\bb) {\bf T}(\boldsymbol{\omega}) \|_F + \gamma \| {\bf b} \|_1 \right)
\end{equation}
where $\|\cdot\|_1$ denotes the $\ell_1$-norm penalization which promotes sparsity of the vector ${\bf b}$.
More recently, Askham and Kutz~\cite{askham2017arxiv} introduced the {\em optimized DMD} algorithm, which uses a variable projection method for nonlinear least squares to compute the DMD for unevenly timed samples, significantly mitigating the bias due to noise.
The optimized DMD algorithm solves the exponential fitting problem directly:
\begin{equation}
       \argmin_{ \boldsymbol{\omega}, \bPhi_{\bf b} } \| \bX -   \bPhi_{\bf b} {\bf T}(\boldsymbol{\omega}) \|_F.
\end{equation}
This has been shown to provide a superior decomposition due to its ability to optimally suppress bias and handle snapshots collected at arbitrary times.  The disadvantage of optimized DMD is that one must solve a nonlinear optimization problem.
{{However, by using statistical bagging methods, the optimized DMD algorithm can be stabilized and the {\em boosted optimized DMD} (BOP-DMD) method can not only improve performance of the decomposition, but also produce UQ metrics for the DMD eigenvalues and DMD eigenmodes~\cite{sashidhar2021bagging}. This provides a nearly optimal linear model for forecasting of dynamics.}}

{DMD is able to accurately identify an approximate linear model for dynamics that are linear, periodic, or quasi-periodic.
However, DMD is fundamentally unable to capture a linear dynamical system model with essential nonlinear features, such as multiple fixed points, unstable periodic orbits, or chaos~\cite{brunton2016plosone}.
As an example, DMD trained on data from the chaotic Lorenz system will fail to yield a reasonable linear model, and the resulting DMD matrix will also not capture important features of the linear portion of the Lorenz model.
The sparse identification of nonlinear dynamics (SINDy)~\cite{brunton2016pnas} is a related algorithm that identifies fully nonlinear dynamical systems models from data.
However, SINDy often faces scaling issues for high-dimensional systems that do not admit a low-dimensional subspace or sub-manifold.
In this case, the recent linear and nonlinear disambiguation optimization (LANDO) algorithm~\cite{baddoo2021kernel} leverages kernel methods to identify an implicit model for the full nonlinear dynamics, where it is then possible to extract a low-rank DMD approximation for the linear portion linearized about some specified operating condition.
In this way, the LANDO algorithm robustly extracts the linear DMD dynamics even from strongly nonlinear systems.
This work is part of a much larger effort to use kernels for learning dynamical systems and Koopman representations~\cite{williams2015jcd,fujii2019dynamic,das2020a,klus2020eigendecompositions,klus2020kernel,manohar2020kernel,baddoo2021kernel}.
}

\subsubsection{Krylov subspace perspective}\label{Sec:DMD:Krylov}
In the original formulation~\cite{schmid2010jfm,rowley2009jfm}, the matrices $\bX$ and $\bX'$ were formed from sequential snapshots, evenly spaced in time:
\begin{subequations}\label{Eq:DMD:SnapshotsSequential}
\begin{align}
\bX &= \begin{bmatrix} \vline & \vline & & \vline \\
\bx_1 & \bx_2 & \cdots & \bx_m \\
 \vline & \vline & & \vline
 \end{bmatrix}\\
 \bX' &= \begin{bmatrix} \vline & \vline & & \vline \\
\bx_2 & \bx_3 & \cdots & \bx_{m+1} \\
 \vline & \vline & & \vline
 \end{bmatrix}.
\end{align}
\end{subequations}
The columns of $\bX$ belong to a Krylov subspace generated by $\bA$ and $\bx_1$:
\begin{align}
\bX &\approx \begin{bmatrix} \vline & \vline & & \vline \\
\bx_1 & \bA\bx_1 & \cdots & \bA^{m-1}\bx_1 \\
 \vline & \vline & & \vline
 \end{bmatrix}.
\end{align}
Thus, DMD is related to Arnoldi iteration to find the dominant eigenvalues and eigenvectors of a matrix $\bA$.

The matrices $\bX$ and $\bX'$ are also related through the {shift matrix $\bS$}
\begin{align}
\bX' = \bX\bS,
\end{align}
{which is a finite-dimensional representation of the \emph{shift} operator.}
Thus, $\bS$ acts on columns of $\bX$, as opposed to $\bA$, which acts on rows of $\bX$.
The shift {matrix $\bS$ has the form of a companion matrix and} is given by
\begin{equation}
\bS =
\begin{bmatrix}
0 & 0 & 0 &\cdots & 0 & a_1 \\
1 & 0 & 0 &\cdots & 0 & a_2 \\
0 & 1 & 0 & \cdots & 0 & a_3\\
\vdots & \vdots & \vdots &\ddots & \vdots & \vdots \\
0 & 0 & 0 & \cdots & 1 & a_{m}
\end{bmatrix}.\label{eq:companion-matrix}
\end{equation}
In other words, the first $m-1$ columns of $\bX'$ are obtained by shifting the last $m-1$ columns of $\bX$, and the last column is obtained as a best-fit combination of the $m$ columns of $\bX$ that minimizes the residual.
The shift {matrix} may be viewed as a matrix representation of the Koopman operator, as it advances snapshots forward in time.
The $m\times m$ matrix $\bS$ has the same non-zero eigenvalues as $\bA$, so that  $\bS$ may be used to obtain dynamic modes and eigenvalues.
However, computations based on $\bS$ are not as numerically stable as the DMD algorithm presented above.

\subsection{Methodological extensions}
The DMD algorithm has been successful in large part because of its simple formulation in terms of linear regression and because it does not require knowledge of governing equations.
For these reasons, DMD has been extended to include several methodological innovations~\cite{kutz2016book} presented here.
Several of these extensions, including to nonlinear systems, delay measurements, and control, will be explored in much more detail in later sections.
Algorithms to handle noisy data were already discussed in \cref{Sec:DMD:AlternativeOpt}.

\subsubsection{Including inputs and control}
Often, the goal of obtaining reduced-order models is to eventually design more effective controllers to manipulate the behavior of the system.
Similarly, in many systems, such as the climate, there are external forcing variables that make it difficult to identify the underlying unforced dynamics.
Proctor et al.~\cite{proctor2016siads} introduced the DMD with control (DMDc) algorithm to disambiguate the natural unforced dynamics and the effect of forcing or actuation given by the variable $\bu$.
This algorithm is based on
\begin{align}
\bx_{k+1} \approx \bA \bx_k + \bB\bu_k,\label{Eq:DMDc:Propagator}
\end{align}
which results in another linear regression problem.
This algorithm was motivated by epidemiological systems (e.g., malaria or polio), where it is not possible to stop intervention efforts, such as vaccinations and bed nets, in order to characterize the unforced dynamics~\cite{proctor2014epj}.
DMDc will be explored extensively in \cref{Sec:Control:DMDc}.
{An input--output DMD~\cite{benner2018reduced} has also been formulated recently that fits into broader reduced-order modeling efforts~\cite{Benner2015siamreview,peherstorfer2016data,qian2020lift,benner2020operator}.}

\subsubsection{Compression and randomized linear algebra}
The DMD algorithm is fundamentally based on the assumption that there are dominant low-dimensional patterns even in high-dimensional data, such as fluid flow fields.
Randomized algorithms~\cite{halko2011siamreview} are designed to exploit these patterns to accelerate numerical linear algebra.
In the randomized DMD algorithm~\cite{bistrian2016ijnme,erichson2017arxiv} data is randomly project into a lower-dimensional subspace where computations may be performed more efficiently.
The existence of patterns also facilitate more efficient measurement strategies based on principles of \emph{sparsity} to reduce the number of measurements required in time~\cite{tu2014ef} and space~\cite{brunton2015jcd,gueniat2015pof,erichson2016jrtp}.
This has the broad potential to enable high-resolution characterization of systems from under-resolved measurements.
In 2014, Jovanovic et al.~\cite{jovanovic2014pof} used sparsity promoting optimization to identify the fewest DMD modes required to describe a data set; the alternative, testing and comparing all subsets of DMD modes, is computationally intractable.
Finally, libraries of DMD modes have also been used to identify dynamical regimes~\cite{kramer2015arxiv}, based on the sparse representation for classification~\cite{wright2009ieeetpami}, which was used earlier to identify dynamical regimes using libraries of POD modes~\cite{bright2013pof,brunton2014siads}.

\subsubsection{Nonlinear measurements and latent variables}
The connection between DMD and the Koopman operator~\cite{rowley2009jfm,tu2014jcd,kutz2016book} has motivated several extensions for strongly nonlinear systems.
The standard DMD algorithm is able to accurately characterize periodic and quasi-periodic behavior, even in nonlinear systems.
However, DMD models based on linear measurements of the system are generally not sufficient to characterize truly nonlinear phenomena, such as transients, intermittent phenomena, or broadband frequency cross-talk.
In Williams et al.~\cite{williams2015jnls,williams2015jcd}, DMD measurements were augmented to include nonlinear measurements of the system, enriching the basis used to represent the Koopman operator.
There are also important extensions of DMD to systems with latent variables.
Although DMD was developed for high-dimensional data, it is often desirable to characterize systems with incomplete measurements.
As an extreme example, consider a single measurement that oscillates as a sinusoid, $x(t) = \sin(\omega t)$.
Although this would appear to be a perfect candidate for DMD, the algorithm incorrectly identifies a real eigenvalue because the data does not have sufficient rank to extract a complex conjugate pair of eigenvalues $\pm i\omega$.
This paradox was first explored by Tu et al.~\cite{tu2014jcd}, where it was discovered that a solution is to stack delayed measurements into a larger matrix to augment the rank of the data and extract phase information.
Delay coordinates have also been used effectively to extract coherent patterns in neural recordings~\cite{brunton2016b}.
The connections between delay DMD and Koopman were subsequently investigated~\cite{brunton2017natcomm,arbabi2017,das2017arxiv,Kamb2020siads}.
Nonlinear measurements and latent variable will both be explored extensively in \cref{Sec:Observables}.

\subsubsection{Multiresolution}  DMD is often applied to complex, high-dimensional dynamical systems, such as fluid turbulence or epidemiological systems, that exhibit multiscale dynamics in both space and time.
Many multiscale systems exhibit transient or intermittent phenomena, such as the El Ni\~no observed in global climate data.
These transient dynamics are not captured accurately by DMD, which seeks spatio-temporal modes that are globally coherent across the entire time series of data.
To address this challenge, the multiresolution DMD (mrDMD) algorithm was introduced~\cite{kutz2016siads}, which effectively decomposes the dynamics into different timescales, isolating transient and intermittent patterns.
Multiresolution DMD modes were recently shown to be advantageous for sparse sensor placement by Manohar et al.~\cite{manohar2017arxiv}.

\subsubsection{Streaming and parallelized codes}  Because of the computational burden of computing the DMD on high-resolution data, several advances have been made to accelerate DMD in streaming applications and with parallelized algorithms.
DMD is often used in a streaming setting, where a moving window of snapshots are processed continuously, resulting in savings by eliminating redundant computations when new data becomes available.
Several algorithms exist for streaming DMD, based on the incremental SVD~\cite{hemati2014pof}, a streaming method of snapshots SVD~\cite{pendergrass2016arxiv}, and rank-one updates to the DMD matrix~\cite{zhang2019}.
The DMD algorithm is also readily parallelized, as it is based on the SVD.  Several parallelized codes are available, based on the QR~\cite{sayadi2016tcfd} and SVD~\cite{erichson2017arxiv,erichson2016arxiva,erichsonrandomized}.

{\subsubsection{Tensor formulations}  Most data used to compute DMD has additional spatial structure that is discarded when the data is reshaped into column vectors.
The tensor DMD extension of Klus et al.~\cite{klus2018tensor} performs DMD on a tensorial, rather than vectorized, representation of the data, retaining this additional structure.
In addition, this approach reduces the memory requirements and computational complexity for large-scale systems.
Extensions to this approach have been introduced based on reproducing kernel Hilbert spaces~\cite{fujii2019dynamic} and the extended DMD~\cite{nuske2019tensorbased}, and additional connections have recently been made between the Koopman mode decomposition and tensor factorizations~\cite{redman2021koopman}.
Tensor approaches to related methods, such as the sparse identification of nonlinear dynamics~\cite{brunton2016pnas}, have also been developed recently~\cite{Gelss2019mindy}.
}

\subsection{Domain applications}\label{Sec:Applications}
DMD has been widely applied to a diverse range of applications.
We will explore key applications in fluid dynamics, epidemiology, neuroscience, and video processing.
In addition, DMD has been used for robotics~\cite{berger2014ieee,abraham2019ieee,bruder2019proc,mamakoukas2019proc}, finance~\cite{mann2016qf}, power grids~\cite{sinha2020arxiv,susuki2011c,susuki2011jns,susuki2012nonlinear,korda2018arxiv}, and plasma physics~\cite{taylor2017arxiv,kaptanoglu2020pop}.

\begin{figure}
\vspace{-.1in}
\begin{center}
\begin{overpic}[width=.95\textwidth]{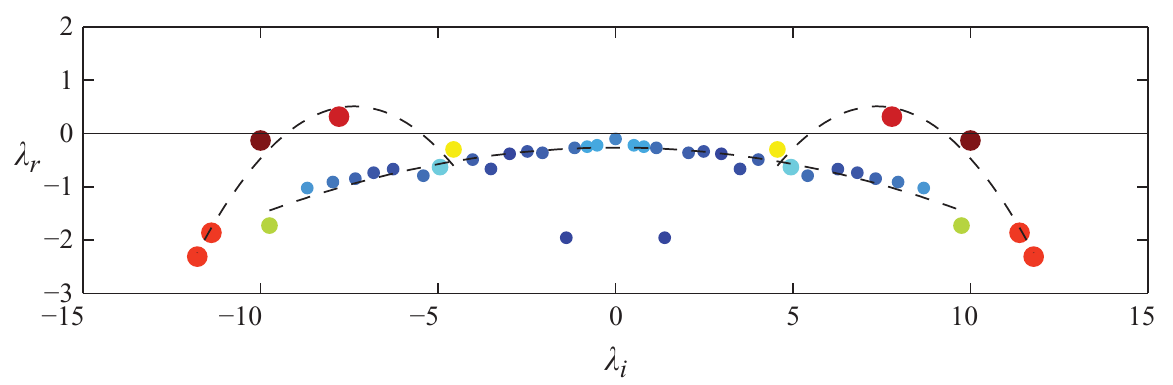}
\end{overpic}
\begin{overpic}[width=.85\textwidth]{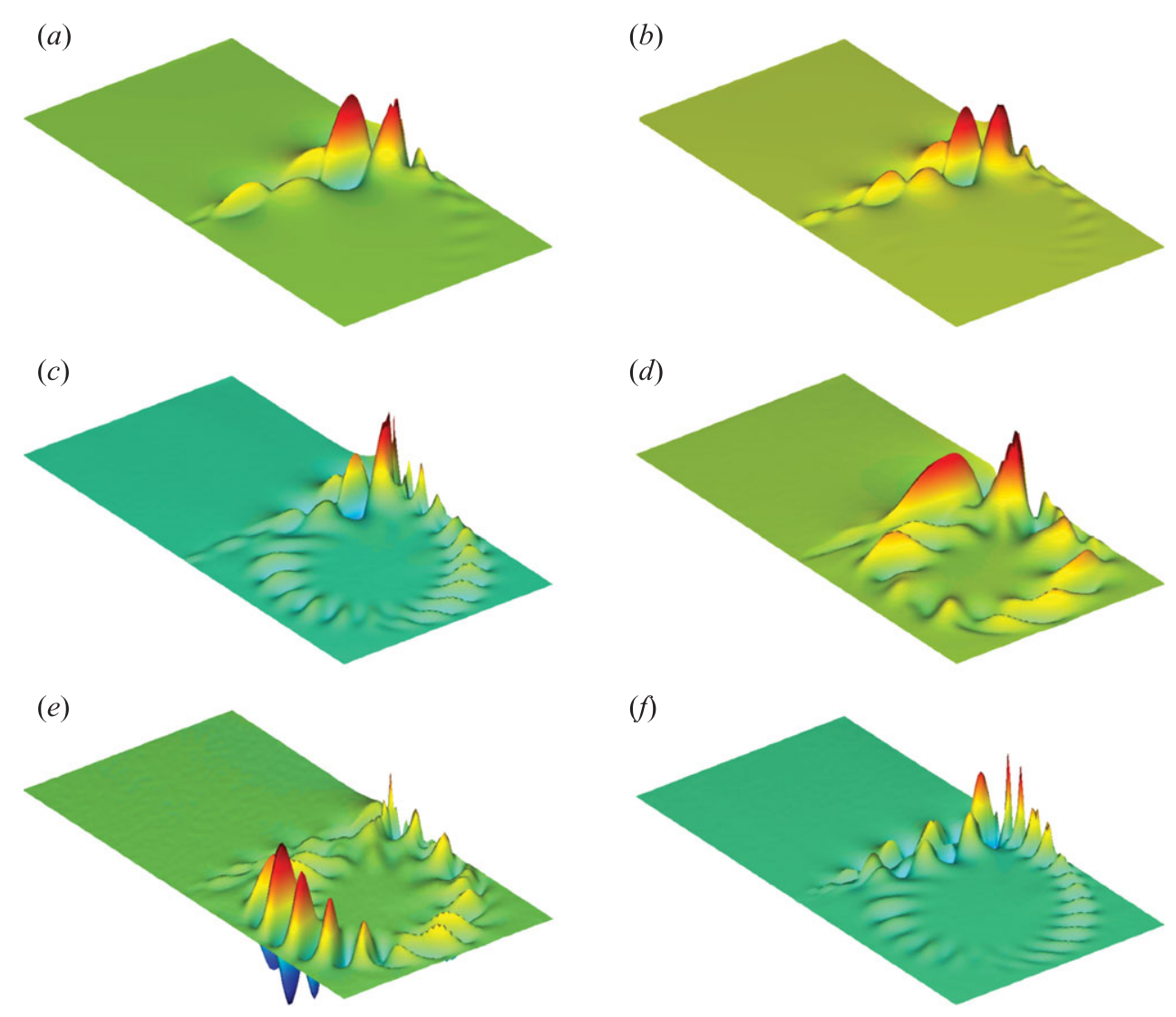}
\end{overpic}
\vspace{-.1in}
\caption{(top) DMD eigenvalue spectrum for the cavity flow at Reynolds number $Re=4500$.  (bottom) Corresponding DMD modes, visualized by the streamwise velocity.  (a) Most unstable dynamic mode, (b--d) DMD modes from the unstable branch, (e,f) DMD modes from the stable branch.  \emph{Modified with permission, from Schmid 2010 Journal of Fluid Mechanics}~\cite{schmid2010jfm}.}\label{Fig:Schmid1}
\end{center}
\vspace{-.1in}
\end{figure}

\subsubsection{Fluid dynamics}
\begin{figure}
\begin{center}
\begin{overpic}[width=\textwidth]{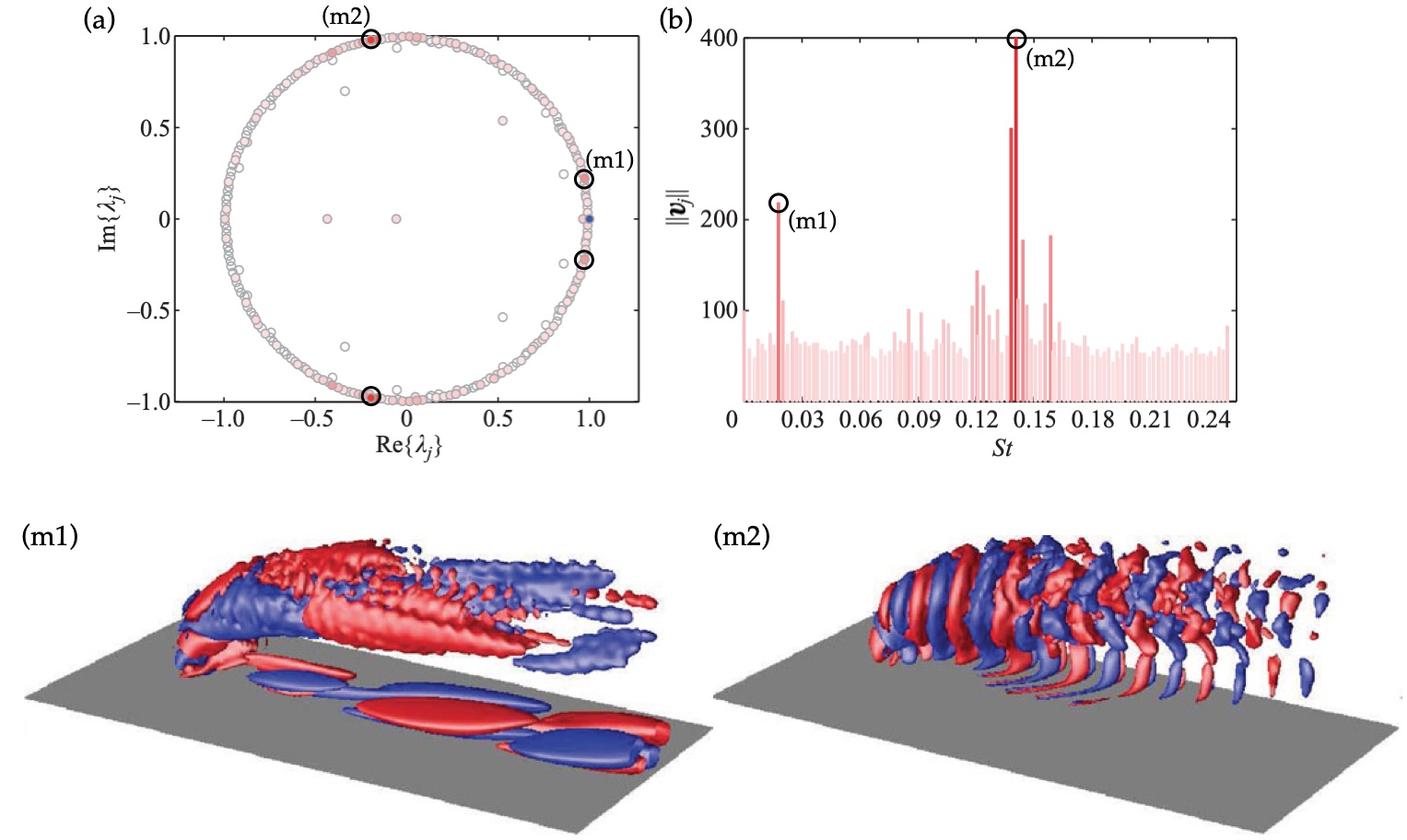}
\end{overpic}
\vspace{-.1in}
\caption{DMD eigenvalues for jet-in-crossflow plotted on the unit circle (a) and as a power spectrum versus Strouhal number (b).  Two DMD modes (m1) and (m2) are pictured below.  \emph{Modified with permission, from Rowley et al. 2009 Journal of Fluid Mechanics}~\cite{rowley2009jfm}.}\label{Fig:Rowley1}
\end{center}
\end{figure}

DMD originated in the fluid dynamics community~\cite{schmid2010jfm}, and has since been applied to a wide range of flow geometries (jets, cavity flow, wakes, channel flow, boundary layers, etc.), to study mixing, acoustics, and combustion, among other phenomena.
In Schmid~\cite{schmid2008aps,schmid2010jfm}, both a cavity flow and a jet were considered; the cavity flow example is shown in \cref{Fig:Schmid1}.
Rowley \emph{et al.}~\cite{rowley2009jfm} investigated a jet in cross-flow; modes are shown in \cref{Fig:Rowley1} with the corresponding eigenvalue spectrum.
Thus, it is no surprise that DMD has subsequently been used broadly in both cavity flows~\cite{schmid2010jfm,lusseyran2011flow,seena2011dynamic,basley2011experimental,basley2013space} and jets~\cite{bellani2011experimental,semeraro2012analysis,schmid:2012,schmid2011tcfd}.

DMD has also been applied to wake flows, including to investigate frequency lock-on~\cite{tu2011koopman}, the wake past a gurney flap~\cite{pan2011dynamical}, the cylinder wake~\cite{bagheri2013jfm}, and dynamic stall~\cite{dunne2015ef}.
Boundary layers have also been extensively studied with DMD~\cite{ostoich2013interaction,sayadi2014reduced,mizuno2011investigation}.
In acoustics, DMD has been used to capture the near-field and far-field acoustics that result from instabilities observed in shear flows~\cite{song2013global}.
In combustion, DMD has been used to understand the coherent heat release in turbulent swirl flames~\cite{moeck2013tomographic} and to analyze a rocket combustor~\cite{huang2013analysis}.
DMD has also been used to analyze non-normal growth mechanisms in thermoacoustic interactions in a Rijke tube.
DMD has been compared with POD for reacting flows~\cite{roy2015pre}.
DMD has also been used to analyze more exotic flows, including a simulated model of a high-speed train~\cite{muld2012flow}.
Shock turbulent boundary layer interaction (STBLI) has also been investigated with DMD to identify a pulsating separation bubble that is accompanied by shockwave motion~\cite{grilli:2012}.
DMD has also been used to study self-excited fluctuations in detonation waves~\cite{massa2012dynamic}.
Other problems include identifying hairpin vortices~\cite{tang2012dynamic}, decomposing the flow past a surface mounted cube~\cite{muld2012mode}, modeling shallow water equations~\cite{bistrian2015ijnmf}, studying nano fluids past a square cylinder~\cite{sarkar2013mixed}{, fluid-structure interaction~\cite{goza2018modal}} , and measuring the growth rate of instabilities in annular liquid sheets~\cite{duke2012experimental}.
{A modified \emph{recursive} DMD algorithm was also formulated by Noack et al.~\cite{Noack2016jfm} to provide an orthogonal basis for empirical Galerkin models in fluids.
The use of DMD in fluids fits into a broader effort to leverage machine learning for improved models and controllers~\cite{Brenner2019prf,bar2019learning,Brunton2020arfm, sanchez2020learning,li2020fourier,kochkov2021machine}, especially for turbulence closure modeling~\cite{Ling2016jfm,Kutz2017jfm,Duraisamy2019arfm,Maulik2019jfm,beetham2020formulating,beetham2021sparse}.}

\begin{figure}
\begin{center}
\begin{overpic}[width=.85\textwidth]{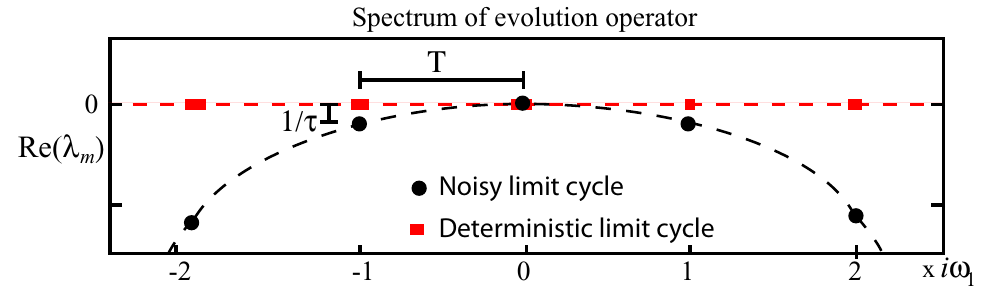}
\end{overpic}
\begin{overpic}[width=.9\textwidth]{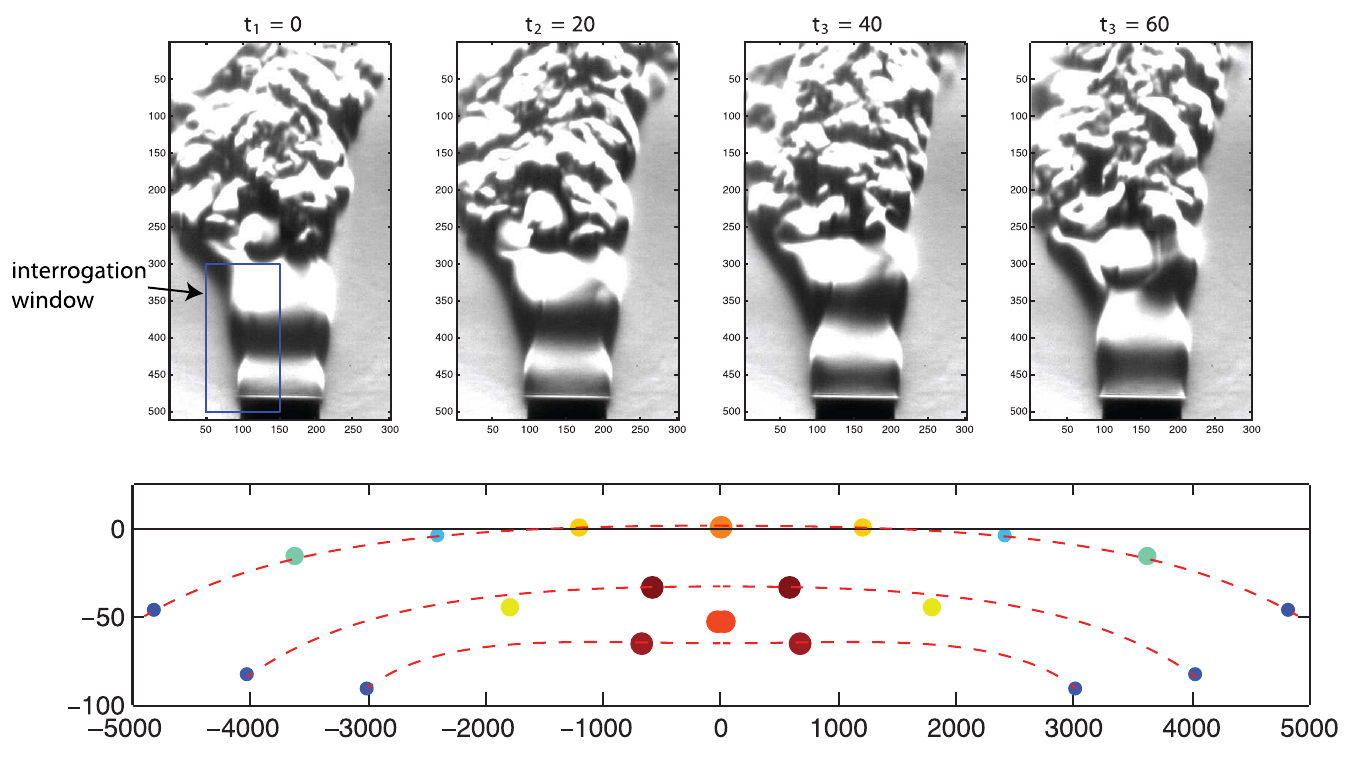}
\end{overpic}
\caption{(top) Koopman eigenvalue spectrum for oscillator with noise.    (middle, bottom) DMD spectrum for jet flow, exhibiting similar noise pattern.
\emph{Reproduced with permission, from Bagheri 2014 Physics of Fluids~\cite{bagheri2014}; (bottom, middle) panels originally from Schmid et al., 2011 Theoretical and Computational Fluid Dynamics~\cite{schmid2011tcfd}}.}
\label{Fig:Bagheri1}
\end{center}
\vspace{-.1in}
\end{figure}

It is interesting to note that the noise effects that were carefully analyzed by Bagheri~\cite{bagheri2014} explain the eigenvalue spectrum observed earlier by Schmid et al.~\cite{schmid2011tcfd} for a turbulent jet.
This comparison is shown in \cref{Fig:Bagheri1}.

\subsubsection{Epidemiology}
\begin{figure}
\begin{center}
\begin{overpic}[width=\textwidth]{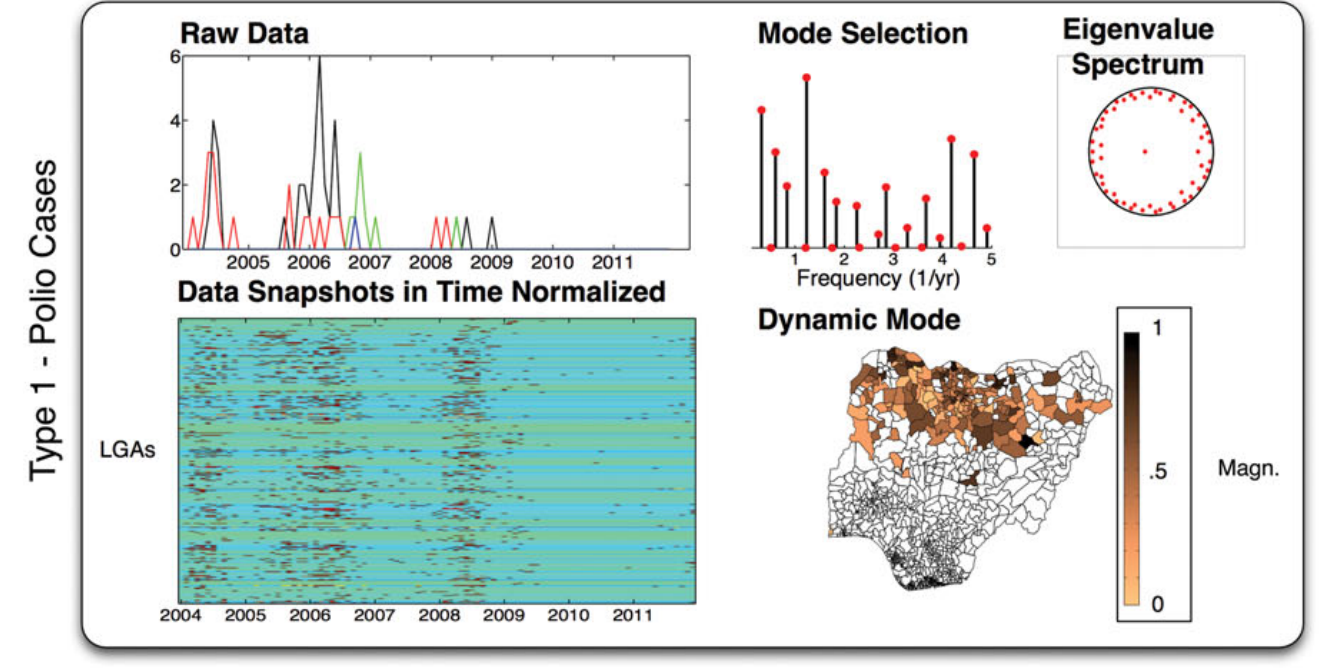}
\end{overpic}
\vspace{-.1in}
\caption{Results of DMD analysis for Polio cases in Nigeria. \emph{Reproduced with permission, from Proctor and Eckhoff 2014 International health~\cite{proctor2015ih}}.}\label{Fig:Proctor}
\end{center}
\end{figure}

Epidemiological data often consists of high-dimensional spatiotemporal time series measurements, such as the number of infections in a given neighborhood or city.
Thus, DMD provides particularly interpretable decompositions for these systems, as explored by  Proctor and Eckhoff~\cite{proctor2015ih} and illustrated in \cref{Fig:Proctor}.
Modal frequencies often correspond to yearly or seasonal fluctuations.
Moreover, the phase of DMD modes gives insight into how disease fronts propagate spatially, potentially informing future intervention efforts.
The application of DMD to disease systems also motivated the DMD with control algorithm~\cite{proctor2016siads}, since it is infeasible to stop vaccinations in order to identify the unforced dynamics.

\subsubsection{Neuroscience}
\begin{figure}
\begin{center}
\begin{overpic}[width=\textwidth]{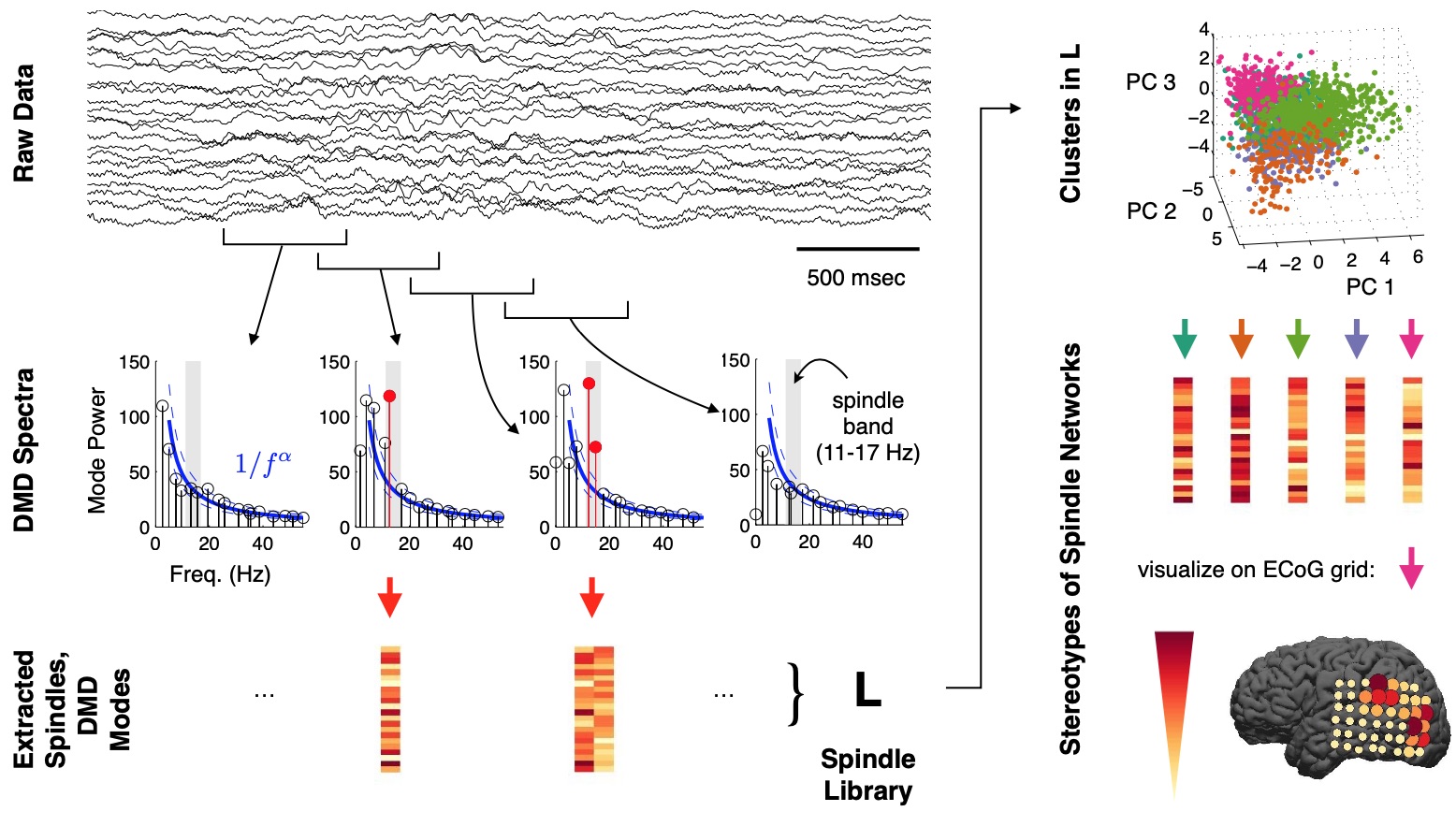}
\end{overpic}
\vspace{-.1in}
\caption{Illustration of DMD applied to human brain data from electroencephalogram (ECOG) measurements.  \emph{Reproduced with permission, from B. Brunton et al. 2016 Journal of Neuroscience Methods~\cite{brunton2016b}}.}\label{Fig:Brunton2016jns}
\end{center}
\end{figure}

 Complex signals from neural recordings are increasingly high-fidelity and high dimensional, with advances in hardware pushing the frontiers of data collection.
 DMD has the potential to transform the analysis of such neural recordings, as evidenced in a recent study by B. Brunton et al.~\cite{brunton2016b} that identified dynamically relevant features in ECOG data of sleeping patients, shown in \cref{Fig:Brunton2016jns}.
 Since then, several works have applied DMD to neural recordings or suggested possible implementation in hardware~\cite{agrawal2016book,broad2017proc,tirunagari2017mva,marrouch2020data}.

\subsubsection{Video processing}
\begin{figure}
\begin{center}
\begin{overpic}[width=\textwidth]{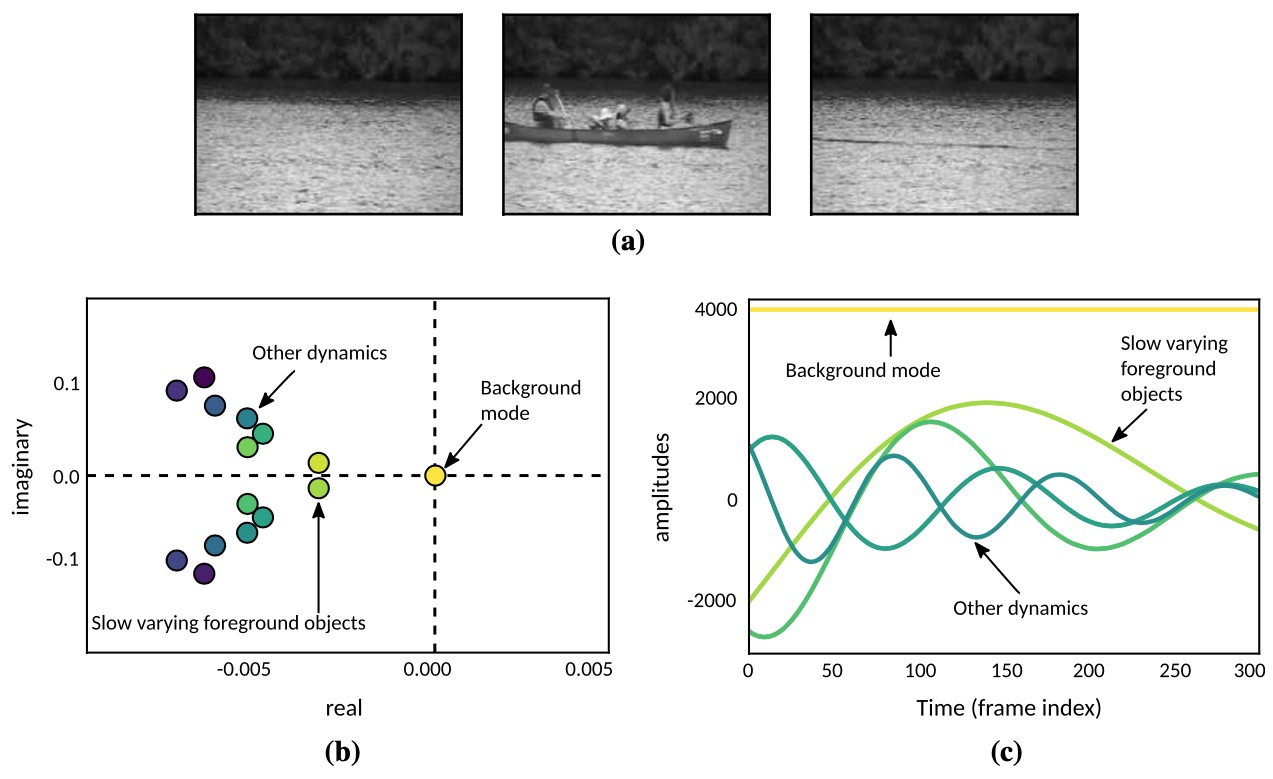}
\end{overpic}
\caption{DMD may be used to process videos, extracting dominant background modes with zero frequency. \emph{Reproduced with permission, from Erichson et al. 2016 Journal of Real Time Processing~\cite{erichson2016jrtp}}.}\label{Fig:ErichsonVideo}
\end{center}
\end{figure}

DMD has also been used for segmentation in video processing, as a video may be thought of as a time-series of high-dimensional snapshots (images) evolving in time.
Separating foreground and background objects in video is a common task in surveillance applications.
Real-time separation is a challenge that is only exacerbated by increasing video resolutions.
DMD provides a flexible approach for video separation, as the background may be approximated by a DMD mode with zero eigenvalue~\cite{grosek2014arxiv,erichson2016jrtp,pendergrass2016arxiv}, as illustrated in \cref{Fig:ErichsonVideo}.

\subsection{Connections to other methods}\label{Sec:DMD:Connections}
Linear regression models are common throughout the sciences, and so it is not surprising to find many methods that are closely related to DMD.
These related methods have been developed in a diverse set of scientific and mathematical disciplines, using a range of mathematical formulations.

\subsubsection{System identification}
DMD may be viewed as a form of system identification, which is a mature subfield of control theory.
Tu et al.~\cite{tu2014jcd} established a connection between DMD and the eigensystem realization algorithm (ERA)~\cite{juang1985jgcd}.
Subsequently, DMD with control~\cite{proctor2016siads} was also linked to ERA and the observer Kalman identification methods (OKID)~\cite{phan1992jas,phan1993jota}.
ERA and OKID construct input-output models using impulse response data and continuous disturbance data, respectively~\cite{juang1985jgcd,juang1991nasatm,phan1992jas,phan1993jota,juang1994book}.
Both of these methods were developed to construct input-output state-space models for aerospace applications where the systems tend to have higher rank than the number of sensor measurements~\cite{ho1965aac,juang1985jgcd}.
In contrast, DMD and DMDc were developed for systems with a large number of measurements and low-rank dynamics.
ERA and OKID have also been categorized as subspace identification methods~\cite{qin20061502}, which include the numerical algorithms for subspace system identification (N4SID), multivariable output error state space (MOESP), and canonical variate analysis (CVA)~\cite{vanoverschee199475,van.1996,katayama.2005,qin20061502}.
Algorithmically, these methods involve regression, model reduction, and parameter estimation steps,  similar to DMDc.
There are important contrasts regarding the projection scaling between all of these methods~\cite{proctor2016siads}, but the overall viewpoint is similar among these diverse methods.
{Recently, DMD has also been combined with an extended Kalman filter to simultaneously de-noise and estimate parameters~\cite{nonomura2018dynamic,nonomura2019extended}.
The problem of data assimilation, which arose in geosciences~\cite{law2015data,apte2008}, is closely related to the construction of the state observer in control theory.
  DMD- and POD-based data-driven model reduction was successfully combined with model-based Kalman filters to yield state and uncertainty estimates~\cite{iungo2015jp,mehta2018}.
  Furthermore, such models were shown to be highly effective even for model dimensions out of reach of classical particle filtering techniques~\cite{albarakati2021}.
}

{It is important to note that these subspace system identification approaches are designed to model linear systems, and it is unclear how to interpret the results when they are applied to nonlinear systems.
However, the strong connection between DMD and the Koopman operator provides a valuable new perspective for how to interpret these approaches, even when used to model data from nonlinear dynamical systems.
For example, we will see in \Cref{Sec:HAVOK} that ERA, which is valid for linear systems, may be applied to nonlinear systems with a modern Koopman interpretation.}

\subsubsection{ARIMA: Auto-regressive integrating moving average}
Autoregressive moving average (ARMA) models and the generalized autoregressive integrated moving average (ARIMA) models~\cite{arima} are commonly used in statistics and econometrics.
These models leverage time-series data to build models to forecast predictions into the future.
ARIMA models are often applied to data that are non-stationary.
Like DMD, ARMA and ARIMA models are characterized by a number of key parameters, one of them being the number of past time points used for forecasting a future point.
However, DMD correlates each time snapshot directly to the previous time snapshot.
Autoregressive models have a number of useful variants, including a generalization to a vector framework, i.e. the VARIMA  (vector ARIMA) model, and a generalization which includes seasonal effects, i.e. the SARIMA (seasonal ARIMA) model.
In the DMD architecture, seasonal variations are automatically included.
Moreover, if the mean of the data matrix $\bX$ is subtracted, then the companion matrix DMD formulation from \cref{Sec:DMD:Krylov} has been shown to be equivalent to a Fourier decomposition of the vector field in time~\cite{chen2012jns,hirsh2019}.
DMD can be thought of as taking advantage of both the vector nature of the data and any oscillatory (seasonal) variations.
Further, the real part of the DMD spectra allows one to automatically capture exponential trends in the data.

\subsubsection{LIM:  Linear inverse modeling}
Linear inverse models (LIMs) have been developed in the climate science community.
LIMs are essentially identical to the DMD architecture under certain assumptions~\cite{tu2014jcd}.
By construction, LIMS rely on low-rank modeling, like DMD, using low-rank truncation which are known as {\em empirical orthogonal functions} (EOFs)~\cite{eof1}.
Using EOFs, Penland~\cite{eof2} derived a method to compute a linear dynamical system that approximates data from a stochastic linear Markov system, which later came to be known as LIM~\cite{eof3}.
Under certain circumstances, DMD and LIM may be considered equivalent algorithms.
The equivalence of projected DMD and LIM gives us yet another way to interpret DMD analysis.
If the mean of the data is removed, then the low-rank map that generates the DMD eigenvalues and eigenvectors is simply the same map that yields the statistically most likely state in the future.
This is the case for both the exact and projected DMD algorithms, as both are built on the same low-order linear map.
LIM is typically performed for data where the mean has been subtracted, whereas DMD is valid with or without mean subtraction.
Regardless, there is a strong enough similarity between the two methods that the communities may find their particular innovations valuable to one another.

\subsubsection{PCR:  Principal component regression}
In the statistical sciences, {\em principal component regression} (PCR)~\cite{pcr} is a regression analysis technique that is once again based on the SVD (specifically PCA).
PCR regresses the outcome (also known as the response or, the dependent variable) on the principal components of a set of covariates (also known as predictors or, explanatory variables or, independent variables) based on a standard linear regression model.
Instead of regressing the dependent variable on the explanatory variables directly, PCR uses the principal components of the explanatory variables as regressors.
One typically uses only a subset of all the principal components for regression, thus making PCR a regularized procedure. Often the principal components with larger variances are selected as regressors; these principal components correspond to eigenvectors with larger eigenvalues of the sample variance-covariance matrix of the explanatory variables.
However, for the purpose of predicting the outcome, the principal components with low variances may also be important, in some cases even more important~\cite{pcr,witten}.
Unlike the DMD/LIM literature, which has been traditionally steeped in dynamics, the statistics literature is often concerned with regression of static data, mapping generic input data to target data.
Thus PCR is typically not specifically applied to time-series data, but is instead a general regression procedure that may or may not be applied to data from a dynamical system.
In some sense, the first two steps of the DMD algorithm may be viewed as performing PCR on snapshot data from a high-dimensional dynamical system.
However, PCR does not include the additional steps of eigendecomposition of the matrix $\tilde\bA$ or the reconstruction of high-dimensional coherent modes.
This last step is what relates DMD to the Koopman operator, connecting the data analysis to nonlinear dynamics.

\subsubsection{Resolvent analysis}
DMD and Koopman operator theory have also been connected to the resolvent analysis from fluid mechanics~\cite{sharma2016prf,Herrmann2020arxiv}.
Resolvent analysis seeks to find the most receptive states of a dynamical system that will be most amplified by forcing, along with the corresponding most responsive forcings~\cite{trefethen1993science,Jovanovic2005,McKeon2010b,Jovanovic2020}.
Sharma, Mezi\'c, and McKeon~\cite{sharma2016prf} established several important connections between DMD, Koopman theory, and the resolvent operator, including a generalization of DMD to enforce symmetries and traveling wave structures.
They also showed that the resolvent modes provide an optimal basis for the Koopman mode decomposition.
Typically, resolvent analysis is performed by linearizing the governing equations about a base state, often a turbulent mean flow.
However, this approach is \emph{invasive}, requiring a working Navier-Stokes solver.
Herrmann et al.~\cite{Herrmann2020arxiv} have recently developed a purely data-driven resolvent algorithm, based on DMD, that bypasses knowledge of the governing equations.
{DMD and resolvent analysis are also both closely related to the spectral POD~\cite{Towne2018jfm,schmidt2020guide,taira2017}, which is related to the classical POD of Lumley and provides time-harmonic modes at a set of discrete frequencies. }

\section{Koopman operator and modern nonlinear dynamics}\label{Sec:AdvancedKoopman}

\subsection{Eigenfunctions as nonlinear coordinate changes}
\label{sec:eigenfunctions-as-coordinates}

Koopman eigenfunctions provide an explicit coordinate transformation between the nonlinear dynamical system \(\dot \bx = \mathbf{f}(\bx)\) \eqref{Eq:ContinuousDynamicsGeneral} and a \emph{factor} of the Koopman (Lie) dynamical system \(\dot g  = \gen{} g\) \eqref{eq:koopman-dynamics-continuous}, where functions \(g\) are restricted to the span of a subset of Koopman eigenfunctions.
In this section, we provide an overview of results that exploit linearity of Koopman dynamics to arrive at conclusions about the original nonlinear dynamical system.

Consider a more general setting between two dynamical systems \(\dot \bx_{1} = \mathbf{f}_{1}(\bx_{1})\) for \(\bx_{1} \in \mathcal{X}_{1}\) and \(\dot \bx_{2} = \mathbf{f}_{2}(\bx_{2})\) for \(\bx_{2} \in \mathcal{X}_{2}\), respectively inducing flows \(\flow_{k} \colon \mathcal{X}_{k} \to \mathcal{X}_{k}\) given by \(\bx_{k}(t) = \flow_{k}^{t}(\bx_{k}(0)),\ k = 1,2\).
If there exists a homeomorphism (continuous function with a continuous inverse) \(h \colon \mathcal{X}_{1} \to \mathcal{X}_{2}\) such that its composition with the flows
\begin{equation}
  \label{eq:conjugacy}
  h \circ \flow_{1}^{t} \equiv \flow_{2}^{t} \circ h
\end{equation}
holds everywhere, the two dynamical systems are \emph{topologically conjugate}.\footnote{Restricting the domain to a strict subset of \(\mathcal{X}_{1,2}\) leads to \emph{local} conjugacy. Requiring higher regularity (differentiability, smoothness, analyticity) of \(h\) and \(h^{-1}\), or  extending \(h\) to also convert between the time domains of two systems leads to related concepts in the theory of differential equations and dynamical systems covered by standard textbooks, such as \cite{wiggins2003,perko2009,meiss2007}.}
The orbit structure of two topologically conjugate dynamical systems is qualitatively the same.
Therefore, if one identifies a pair of such systems where one of them is easier to analyze analytically, numerically, or is already familiar, it makes it possible to port results from the ``simpler'' to the more complicated system.
One seminal result of this sort is the Hartman--Grobman theorem~\cite[\S 19.12A]{wiggins2003}, which establishes that a nonlinear system \(\dot \bx = \mathbf{f}(\bx)\) in an open neighborhood of a hyperbolic fixed point \(\bx_{0}\) is locally  topologically conjugate to its linearization \(\dot{\bz} = \bA  \bz\), where \(\bA = D\mathbf{f}(\bx_{0})\),  formally justifying the use of the Jacobian in inferring stability of hyperbolic nonlinear equilibria.

Lan and Mezi\'{c}~\cite{lan2013physd} recognized that the conjugacy relationship resulting from the Hartman--Grobman theorem
\begin{equation}
  \label{eq:hartman-grobman-conjugacy}
  h( \mathbf{F}^{t}(\bx) ) = e^{\bA t} h(\bx),
\end{equation}
where \(e^{\bA t}\) is the flow of the linearization, implies the existence of a (vector-valued) Koopman eigenfunction \(\mathbf{h}\),  as it matches the definition \eqref{Eq:KoopmanEfun:Discrete}.
By diagonalization of \(\bA\), the components \(h_{k}\) of \(\mathbf{h}\) are precisely scalar-valued Koopman eigenfunctions.
Although H--G establishes that \eqref{eq:hartman-grobman-conjugacy} holds only in some neighborhood of the origin, \cite{lan2013physd} show that trajectories emanating from the boundary of that neighborhood (backward-in-time) can be used to extend the definition of eigenfunctions \(h_{k}\) up to the edge of the basin of attraction or up to the end of the interval of existence of trajectories, therefore proving an extension of the Hartman--Grobman theorem to the entire basin of attraction.
{Additional results by Eldering, Kvalheim, and Revzen further extend these results to normally-attracting invariant manifolds~\cite{eldering2018}.}
By analogous results for discrete-time dynamics, linearization of periodic and periodically-forced dynamical systems follows.

In \cite{bollt2018}, this idea is taken further, to construct the conjugacy \(h(\cdot)\) between two nonlinear systems that are known to be topologically conjugate. Explicit construction of conjugacies is generally a difficult task that is often eschewed when existence of conjugacy can be inferred from other properties of the system. It is demonstrated that while for one-dimensional systems the construction is straightforward, conjugacies on state spaces in \(\mathbb{R}^{n}\) require a more delicate approach, especially when dealing with eigenvalues with nontrivial geometric multiplicity.

In the setting of \emph{model/order reduction}, relations~\eqref{eq:conjugacy} and~\eqref{eq:hartman-grobman-conjugacy} are required to hold for non-invertible maps or on closed (or even singular) subsets of the domain. In this case the full orbit structure of systems is not typically qualitatively the same; however, much can be gained by studying a simpler, often significantly-lower dimensional, system, and transporting its properties onto the original dynamics.

\begin{figure}[t]
\vspace{-.1in}
  \centering
  \includegraphics[width=0.975\linewidth]{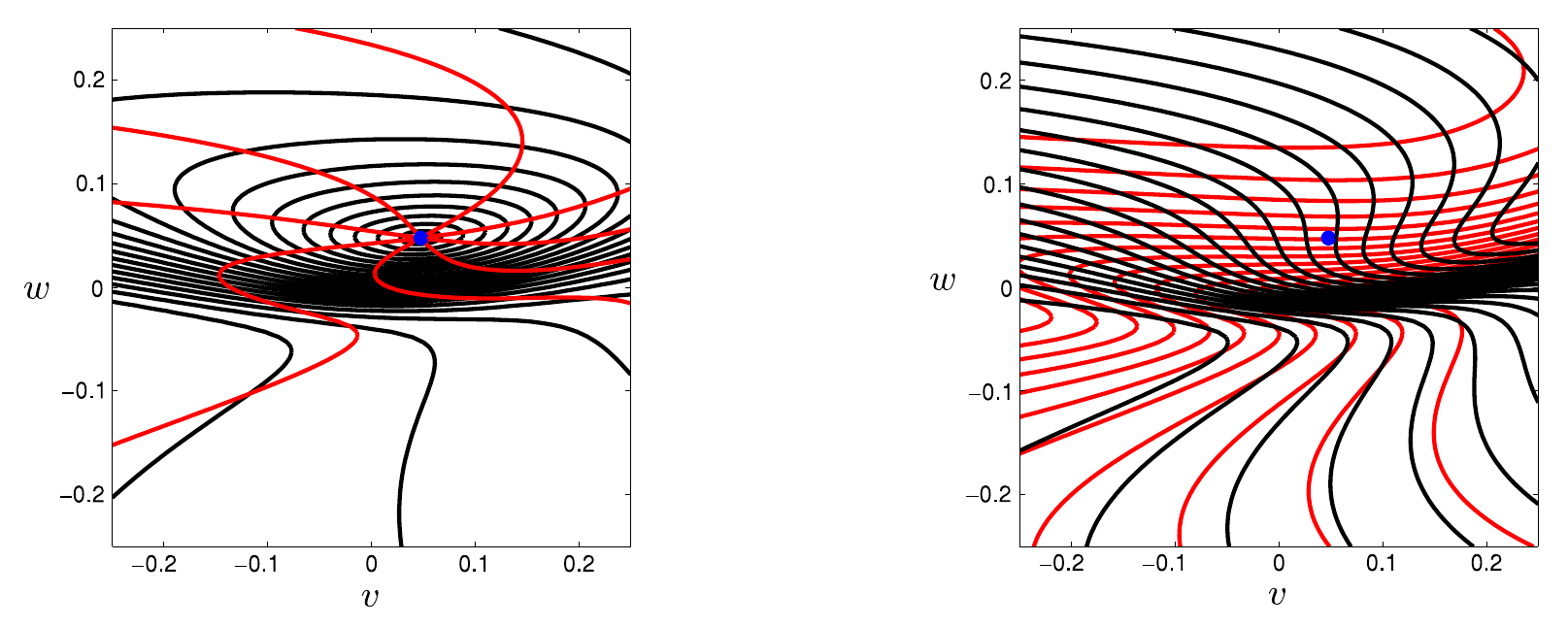}
  \vspace{-.1in}
  \caption[Isostables and isochrons]{Isostables and isochrons for the FitzHugh--Nagumo model acting as a deformed rectifiable coordinate system in vicinity of a focus (left) and node (right). \textit{Reproduced with permission, from Mauroy et al. 2013 Physica D~\cite{mauroy2013}.}}
  \label{fig:mauroy2013-isostables}
  \vspace{-.1in}
\end{figure}

In this context, Mauroy et al.~\cite{mauroy2012,mauroy2013} demonstrate that the coordinate transformations \(h_{k}\) can be numerically computed by forward integration of a trajectory and a Laplace average of an observable
\begin{equation}
  \label{eq:laplace-average}
  \tilde{g}_{\lambda}(\bx) \coloneqq  \lim_{T\to\infty} \frac{1}{T}\int_{0}^{T} g_{t}(\bx) e^{\bar{\lambda} t} dt,
\end{equation}
which is an extension of the harmonic Fourier average~\eqref{eq:harmonic-average}.
{ The existence and uniqueness of \(C^{k}\) Koopman eigenfunctions has been established around stable fixed points and periodic orbits by Kvalheim et al.~\cite{kvalheim2021,kvalheim2021a}, rigorously justifying their use as conjugacies.}

The isostables (level sets of the absolute value of eigenfunctions) and isochrons (level sets of the arguments of complex eigenfunctions) act either as rectifiable Cartesian coordinate systems for node equilibria, or as rectifiable polar (action-angle) coordinate systems for focus-type equilibria, as shown in \cref{fig:mauroy2013-isostables}, with clear generalizations to higher dimensional systems.

Further theoretical developments have led to extensions to nonlinear stability analysis and optimal control of dynamical systems~\cite{sootla2018,mauroy2017arxiv,mauroy2016ieeetac,sootla2016acc,mauroy2013cdc}. Papers \cite{wilson2019,wilson2016} apply this concept to synchronization of oscillators by extending the phase-response curves using isostables and isochrons computed as Koopman eigenfunctions. Notably, these developments demonstrate that Koopman eigenfunctions are a viable and practical path both in analytic~\cite{wilson2020} and in data-driven approaches~\cite{wilson2020a} to synchronization of oscillators.

\subsection{Phase portrait and symmetries}
\label{sec:symmetries-and-koopman}

Discussions of Hartman--Grobman theory and conjugacies are typically concerned with the behavior of two dynamical systems in the vicinity of an object of interest, such as a basin of attraction/repulsion of a fixed point or periodic orbit.
Here we describe how Koopman eigenfunctions can reflect the structure of the entire phase portrait of the dynamical system.

The phase portrait of a dynamical system on state space \(\mathcal{X}\) is a collection of all orbits \(\mathbf{o}(\bx) = \{ \bF^{t}(\bx) \}_{t \in \mathbb{R}}  \) emanating from each initial condition \(\bx\) by the flow map \(\bF^t\).
It is immediately clear that two initial conditions \(\mathbf{y},\mathbf{z}\) on the same orbit \(\mathbf{o}(\bx)\) generate the same orbit, \( \mathbf{o}(\bx) = \mathbf{o}(\by) = \mathbf{o}(\bz)\).
However, computationally pointwise equality between orbits may be challenging to verify when orbits are merely approximated by numerical integration, and when orbits contain thousands of points.
Instead, it may be practical to verify that \emph{orbital averages} of a subset of functions \(\mathcal{G}(\mathcal{X})\) are sufficiently close, resp.\ distant, to declare that two collections of points belong, resp.
do not belong, to the same orbit.

\begin{figure}[t]
\vspace{-.1in}
  \centering
  \begin{subfigure}[t]{0.45\linewidth}
\includegraphics[width=0.9\linewidth]{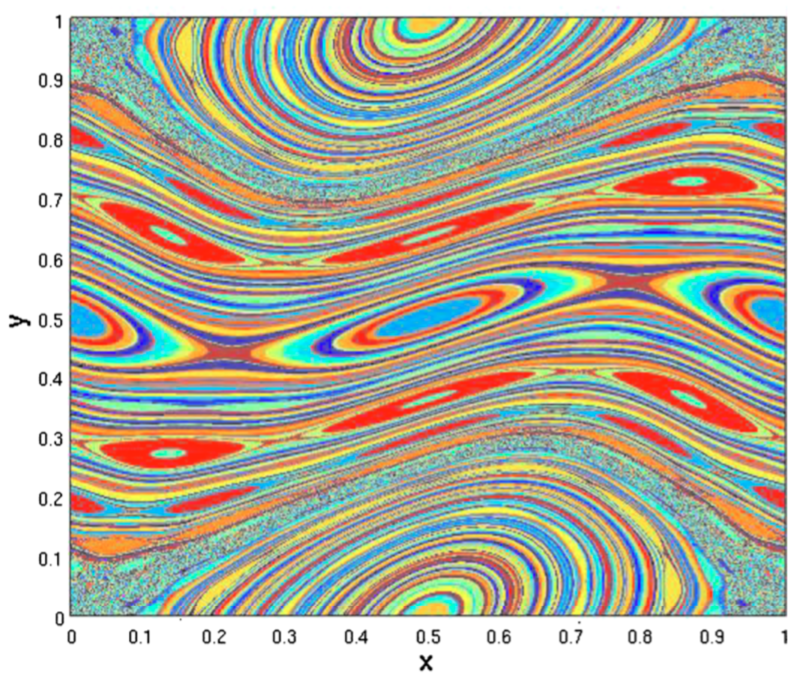}
  \caption{}\label{fig:ergodic-partition}
\end{subfigure}
  \begin{subfigure}[t]{0.45\linewidth}
\includegraphics[width=0.9\linewidth]{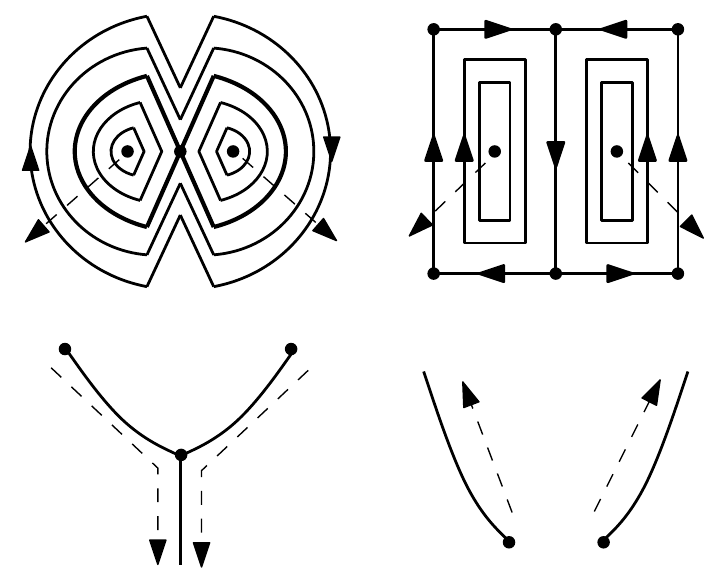}
  \caption{}\label{fig:quotient-coordinates}
\end{subfigure}
\vspace{-.2in}
\caption{(\subref{fig:ergodic-partition}) Approximate ergodic partition of the Chirikov Standard Map.  \textit{Reproduced with permission, from Levnaji\'c and Mezi\'c 2010 Chaos~\cite{levnajic2010}}. (\subref{fig:quotient-coordinates})  Sketch of ordering of orbits by ergodic quotient trajectories.  \textit{Reproduced with permission, from Budi{\v{s}}i{\'c} and I. Mezi\'c 2012 Physica D~\cite{budisic2012physd}}.}
\vspace{-.1in}
\end{figure}

As mentioned in \cref{sec:koopman-eigenfunctions}, the orbital average, or ergodic average, of any observable \(g \in \mathcal{G}(\mathcal{X})\) projects the observable on the invariant eigenspace of the Koopman operator.
Therefore the heuristic technique for comparing orbits just described can be interpreted by comparing whether points \(\by\) and \(\bz\) are mapped to the same point by a vector valued function \(\mathbf{P} \colon \mathcal{X} \to \mathbb{R}^{p}\), \(p \leq \infty\) whose components are all invariant Koopman eigenfunctions:
\begin{equation}
  \mathbf{Q}(\bx) =
  \begin{bmatrix}
    q_{1}(\bx) & q_{2}(\bx) & \dots & q_p(\bx)
  \end{bmatrix}, \text{ where } \koop q_{j} \equiv q_{j}.
\label{eq:erg-quot-embedding}
\end{equation}
Mezi\'{c}~\cite{mezic1999chaos} showed that when the number of averaged functions grows to infinity, this procedure manages to separate ergodic sets that are measure-theoretic counterparts to orbits.
This approach was applied to visualize the phase portrait~\cite{mezic1999chaos,levnajic2010,budisic2012physd}, where it is possible to assign pseudocolors to trajectories, generating visually appealing and detailed figures as in \cref{fig:ergodic-partition}.

Using the embedding function \eqref{eq:erg-quot-embedding} it is further possible to treat the ergodic quotient as a geometric object and compute its local coordinates.
Such coordinates parametrize the space of invariant (ergodic) sets; when continuous, they act as a way of ordering level sets of conserved quantities, even when no explicit (or global) formulas for conserved quantities is known.
The crucial step is to treat the elements \(q_{i}(\bx)\) in~\eqref{eq:erg-quot-embedding} as Fourier coefficients and define a metric using a Sobolev norm on the associated space.
The resulting geometric space can then be coordinatized using manifold learning algorithms.

The ergodic partition is a process of classifying global information, assuming access to many trajectories, into how they relate locally. The opposite direction, where local information is stitched into a global picture is equally as relevant.
Consider two disjoint \emph{invariant} sets \(A,B \subset \mathcal{X}\) of the dynamics \(\mathbf{F}^{t} \colon \mathcal{X} \to \mathcal{X}\).
A typical application of DMD starts from a single trajectory and constructs a finite-dimensional approximation of the Koopman operator.
{Depending on whether the initial condition is in \(A\) or in \(B\), any single-trajectory based computation can at best approximate the restricted Koopman operators \(\koop_{A}\) and \(\koop_{B}\) that do not have to have direct connections, except as restrictions of the Koopman operator associated with the dynamics on the entire space.}
Indeed, it was typical of early papers to either assume that the system has a quasiperiodic attractor, so that all trajectories quickly converge to this set, or to assume that the trajectory used for DMD is ergodic, so that it visits close to any other point in the state space, ruling out the existence of disjoint invariant sets \(A,B\) of positive measure.

If DMD is indeed computed twice, based on trajectories contained in disjoint ergodic sets, it is possible to ``stitch'' the two operators \(\koop_{A},\ \koop_{B}\) together by placing them into a block-diagonal stitched Koopman operator.
Assuming that the respective spaces of observables are \(L^{2}(A,\mu_{A})\), \(L^{2}(B,\mu_{B})\), the joint space over \(\mathcal{X} = A \cup B\) can be taken as
\begin{equation}
  \label{eq:joint-observable-space}
  \mathcal{G}(\mathcal{X}) = L^{2}(A \cup B,\mu_{A} + \mu_{B}) = L^{2}(A,\mu_{A}) \oplus L^{2}(B,\mu_{B}),
\end{equation}
since any function \(f \in \mathcal{G}(\mathcal{X}) \) can be decomposed into disjoint components owing to \(A,B\) being disjoint.
Reference~\cite{nandanoori2020} gives further theoretical backing to this process, as well an incremental and a non-incremental version of the data-driven DMD procedure built in this fashion.
For practical purposes, the assumption that \(A\) and \(B\) are truly dynamically separate is not needed; rather, it is simply sufficient to choose pairs of DMD trajectory samples in a trajectory from two disjoint sets in order to construct such restricted approximations that can be stitched~\cite{sinha2020}.

A particularly important application of stitching across invariant sets concerns phase spaces of systems with symmetries. When the orbit structure is symmetric with respect to a transformation, the analysis of the entire phase space can be dramatically simplified. Symmetry of the differential equations \eqref{Eq:ContinuousDynamicsGeneral} with respect to a symmetry group \(\Gamma\), or \emph{\(\Gamma\)-equivariance}, is defined as a conjugacy
\begin{equation}
  \label{eq:ode-equivariance}
  \mathbf{f}( \gamma \bx ) = \gamma \mathbf{f}(\bx), \forall \gamma \in \Gamma,
\end{equation}
where \(\gamma\) represent action by a group element on the state space. An analogous relationship holds for discrete dynamics~\eqref{Eq:DiscreteDynamics}.
In both cases, the implication is that given any orbit \(\{ \bx(t) \}_{t \in \mathbb{R}}\), there exist symmetry-related counterparts \(\{ \gamma \bx(t) \}_{t \in \mathbb{R}}\) generated by applying any \(\gamma \in \Gamma\) to the orbit.
Stability and asymptotic properties of the symmetry-related orbits are the same, which allows us to study just a portion of the state space in detail, and then export those results to other parts of the state space by symmetry.

Although symmetry-based arguments have long been used to simplify dynamical systems, especially in classification of bifurcations, explicit connections with the Koopman operator framework are fairly recent~\cite{kakubr2018arxiv,salova2019,mesbahi2019,sinha2020}.
In all cases, the central role is played by connections between the definition of the Koopman operator and the conjugacy in the definition of equivariance. The following two theorems appear as Thm.\ III.1 and its corrolary in~\cite{salova2019}, and Thm.\ 1 and Prop.\ 2 in \cite{sinha2020}.
\begin{theorem}
  For a \(\Gamma\)-equivariant dynamical system, \(\koop\) commutes with the action of all \(\gamma\in\Gamma\) for any observable \(g \in \mathcal{G}(\mathcal{X})\)
  \begin{equation}
    \label{eq:commutativity}
    [\gamma \circ (\koop g)](\bx) = [\koop (\gamma \circ g)](\bx).
  \end{equation}
\end{theorem}
\begin{proposition}
  Any eigenspace of \(\koop\) for a \(\Gamma\)-equivariant dynamical system is \(\Gamma\)-\emph{invariant}.
\end{proposition}
Proofs of both statements follow by manipulation of the definitions of an eigenfunction \eqref{Eq:KoopmanEfun:Discrete} and equivariance~\eqref{eq:ode-equivariance}, and we omit them here. For cyclic groups, i.e., finite groups \(\mathbb{Z}_{n} = \{1, \gamma, \gamma^{2},\cdots,\gamma^{n-1}\}\), a more detailed argument is given in \cite{mesbahi2019}, in particular, stating that Koopman modes associated with a particular eigenvalue are similarly symmetric.

The symmetry arguments fit well with the idea of ``stitching'' the global Koopman operator from ``local'' Koopman operators, because symmetry removes the need for simulating additional trajectories in order to explore various invariant sets.
{Consider two disjoint invariant sets \(A,B \subset \mathcal{X}\) that are related by a linear group action \(\tau_{\gamma}\), i.e., \(A = \tau_{\gamma} B\).
Choosing a common finite set (dictionary) of observables, Koopman operators may be given (approximate) representations on each of these sets by computing the approximations from trajectories in \(A\) and \(B\), yielding Koopman matrices \(\mathbf{K}_{A}, \mathbf{K}_{B}\).
These matrices are conjugate (similar)
\begin{equation}
  \label{eq:koopman-symmetry-conjugation}
  \mathbf{K}_{A} = \bT_{\gamma}^{-1} \mathbf{K}_{B} \bT_{\gamma}
\end{equation}
where the conjugating matrix \(\bT\) is the group action on the space of observables, as Sinha et al.~\cite{sinha2020} show.
As a result, assuming that  \(\mathbf{K}_{A}\) can be computed, the matrix \(\mathbf{K}_{B}\) can then be computed without additional simulations.}

\begin{figure}[t]
  \vspace{-.1in}
  \centering
  \includegraphics[width=0.7\linewidth]{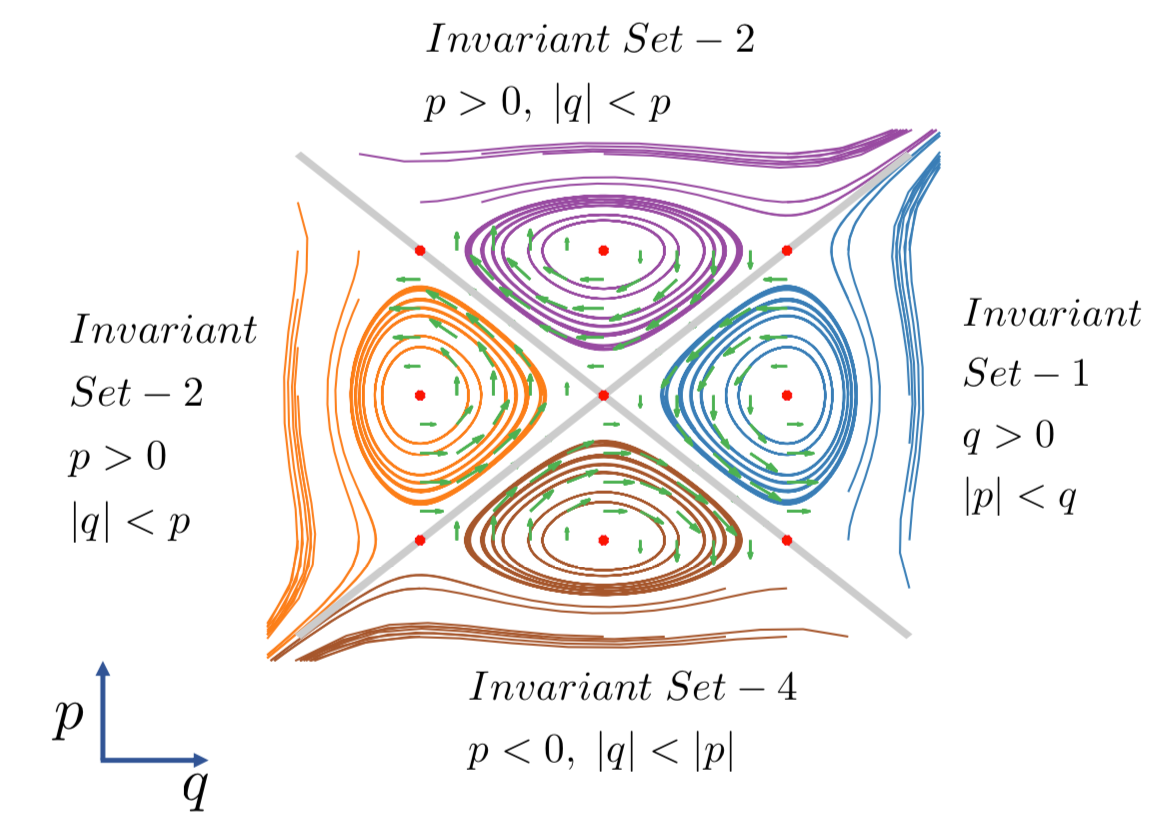}
  \vspace{-.1in}
  \caption{Invariant sets for a dynamical system with \(Z_{2} \times Z_{2}\) symmetry whose Hamiltonian is \(H(p,q) = p^{4}/4 - 9p^{2}/2 - q^{4}/4 + 9q^{2}/2\).  \textit{Reproduced with permission, from Sinha et al.~\cite{sinha2020}.}}
  \vspace{-.15in}
\end{figure}

Even in the case that the full Koopman operator approximation may be directly computed, converting it to a block-diagonal form simplifies the task of computing its spectral properties as eigendecomposition may be performed on individual blocks instead of the entire operator.
{Commutativity of two linear operators, the full Koopman operator and the multiplication by its approximate matrix representation, implies that they preserve each others eigenspaces}, further implying that even when the space of observables \(\mathcal{G}(\mathcal{X})\) is not \emph{constructed} as a direct product \eqref{eq:joint-observable-space}, it is possible to perform a change of basis based on the symmetry group such that \(\mathcal{G}(\mathcal{X})\)  decomposes into so-called \emph{isotypic} subspaces, invariant with respect to \(\koop\).
Consequently, this  block-diagonalizes a finite-dimensional representation of the Koopman operator.
In~\cite{salova2019}, this idea is taken as a starting point for a practical block-diagonalization of matrices appearing in DMD.

An example of isotypic decomposition of functions \(g \in \mathcal{G}(\mathcal{X})\) associated with the group \(\mathbb{Z}_{2}\), where \(\gamma(x) = -x\), is a decomposition into odd and even components
\begin{equation}
g(x) = \underbrace{\frac{g(x) + g(-x)}{2}}_{=: g_{e}(x)} + \underbrace{\frac{g(x) - g(-x)}{2}}_{=: g_{o}(x)}
\end{equation}

\begin{figure}[t]
  \centering
  \vspace{-.1in}
  \includegraphics[trim=0 0 0 64, clip,width=0.9\linewidth]{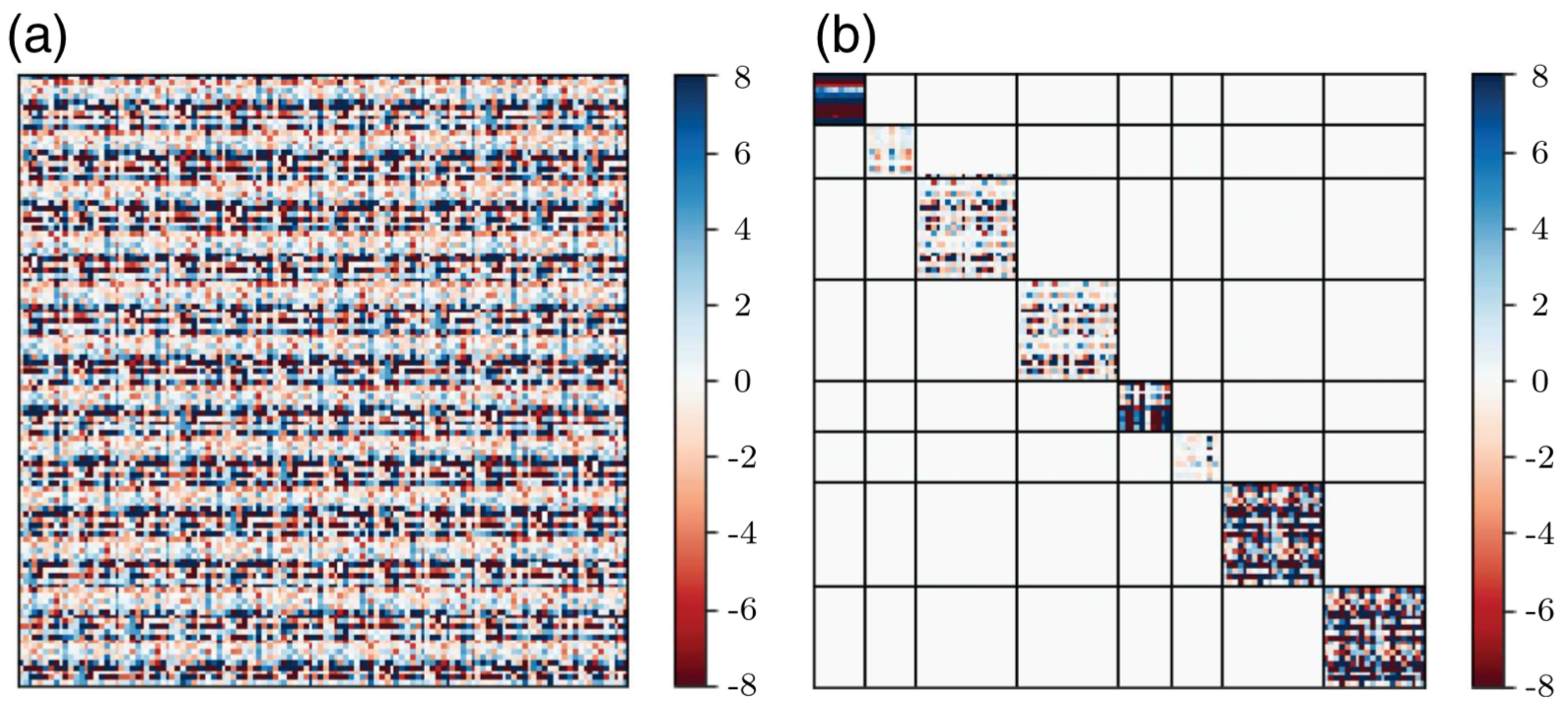}
  \vspace{-.1in}
  \caption{Sparsity structure of DMD matrices before (left) and after (right) block-diagonalization procedure applied to a coupled Duffing oscillator system with \(Z_{2} \times D_{3}\) symmetry.  \textit{Reproduced with permission, from Salova et al. 2019 Chaos~\cite{salova2019}}.}
  \vspace{-.2in}
\end{figure}

Generalization of this procedure to arbitrary finite groups is used as a preconditioning step in a modified eDMD (see \cref{Sec:EDMD}) procedure in \cite{salova2019}, resulting in a block-diagonalized form of the associated DMD matrix that approximates the Koopman operator.
Consequently, this reduces the computational effort needed to compute eigenvalues of the Koopman operator.

\subsection{Adjoint: The Perron--Frobenius operator}
\label{sec:perron-frobenius}

The Perron--Frobenius (PF) operator, also known as the Ruelle transfer operator, is a linear representation of nonlinear dynamics that traces its roots to the mathematics underpinning statistical and quantum physics, paralleling the development of the Koopman operator.
Instead of evolving measurement functions (observables) taking values from the domain of the dynamics, as the Koopman operator does, the PF operator evolves measures (distributions) supported on the domain of the dynamics.
As the PF and Koopman operators can be shown to be formally adjoint in appropriately defined function spaces, we summarize the basic concepts related to the PF operator here.
The fundamentals are well documented in textbooks and monographs~\cite{lasota1994,bollt2013,gaspard1998,chaosbook} and we point to them for a more precise and general introduction of these topics.

Let the domain of dynamics \(\mathcal{X}\) be given the structure of a Borel measurable space.
Given a probability measure \(\mu\), and any measurable set $A$, define the \emph{Perron--Frobenius operator} as
\begin{equation}
  \label{eq:perron-frobenius}
  \pf^{t} \mu (A) \coloneqq  \mu( \flow^{-t}(A) ),
\end{equation}
where \(\flow^{-t}(A) \coloneqq  \left\{ \bx \in \mathcal{X} \colon \flow^{t}(\bx) \in A \right\}\) is the pre-image of \(A\) through the dynamics.
Similar to the Koopman operator, the family \(\pf^{t}\) forms a monoid.
An alternative formulation by Lasota and Mackey~\cite{lasota1994} replaces the action on the space of probability measures with an action on a function space.
This assumes that the flow map \(\flow^{t}\) is nonsingular with respect to some ground measure \(m\), e.g., a Lebesgue measure, meaning
\begin{equation}
  m(A) \not = 0 \quad \Longrightarrow \quad m( \flow^{-t}(\mathcal{A}) ) \not = 0.
\end{equation}
Interpreting $g \in L^{1}(\mathcal{X},dm)$ as densities of measures, i.e., $d\mu = g dm$, it is possible to define the PF operator \(\pf^{t} \colon L^{1}(\mathcal{X},dm) \to L^{1}(\mathcal{X},dm)\)  as
\begin{equation}
  \label{eq:weak-perron-frobenius}
  \int_{A} \pf^{t} g(\bx) dm = \int_{\flow^{-t}(A)} g(\bx) dm,
\end{equation}
for any Borel set \(A\).
If the flow map is additionally smooth, this definition is equivalent to
\begin{equation}
  \label{eq:smooth-pf}
  \pf^{t} g(\bx) = \int_{\flow^{-t}(\bx)} \frac{g(\bs)}{\abs{\nabla \flow^{t}(\bs)} } dm(\bs),
\end{equation}
where \(\abs{\nabla \flow^{t}}\) indicates the determinant of the Jacobian matrix of derivatives of the flow map.

The two formulations \eqref{eq:perron-frobenius} and \eqref{eq:weak-perron-frobenius} are connected by interpreting \(g \in L^{1}\) as densities defining probability measures absolutely continuous with respect to \(dm\), i.e., \(d\mu =  g dm\).

Assuming that the dynamics can be restricted to a space of finite measure, it holds that \(L^{\infty}(\mathcal{X},dm) \subset L^{2}(\mathcal{X},dm) \subset L^{1}(\mathcal{X},dm)\). In this setting, the
Koopman and PF operators can be defined, respectively, on \(L^{\infty}(\mathcal{X},dm)\) and \(L^{1}(\mathcal{X},dm)\), or both on \(L^{2}(\mathcal{X},dm)\), and one can show that they are adjoint to each other:
\begin{equation}
  \label{eq:adjoint}
  \left\langle \pf^{t} f, g \right \rangle =   \left\langle f, \koop^{t} g \right \rangle, \text{ where } \left\langle f, g \right\rangle \coloneqq  \int_{\mathcal{X}} \bar f(\bx) g(\bx) dm.
\end{equation}
Since the proof  \cite[\S 3.3]{lasota1994} proceeds by the standard argument of approximation by simple functions, i.e., linear combinations of indicator functions, this relationship extends to a wider range of spaces.
Adjoint operators have the same spectrum, although their eigenfunctions do differ, as is the case in general for eigenvectors of matrices and their adjoints (transposes).
This connection is partially responsible for the parallel development of Koopman and PF techniques in various contexts of applied mathematics.
{For invertible non-singular dynamics that preserve the Lebesgue measure,\cite[Cor 3.2.1.]{lasota1994} shows that the Perron--Frobenius operator of time-forward map $\flow^{t}$ is precisely the Koopman operator of the time-backward map $\flow^{-t}$.
Several papers~\cite{koltai2019,klus2017data,klus2015numerical,bittracher2015,govindarajan2019} have treated both transfer operators, that is Koopman and Perron--Frobenius operators, using the same formalism, especially in the context of volume-preserving invertible dynamics and in the context of stochastic differential equations.
}

{ An early, but still used, technique to approximate the Perron--Frobenius operator by a stochastic matrix is termed \emph{Ulam's method} after a conjecture by S.~Ulam~\cite[\S IV.4]{ulam1960problems}, ultimately proved to be correct by Y.~Li~\cite{li1976}.}
First, assume that the flow \(\flow^{t}\) preserves a finite probability measure \(m\) on a bounded set \(\mathcal{X}\), partition \(\mathcal{X}\) into subsets \(\{S_{i}\}\), each with the associated characteristic function \(\chi_{i}\).
Choosing a fixed time \(T\), we then form the stochastic \emph{Ulam matrix}
\begin{equation}
  \mathbf{U} = (u_{ij}), \quad u_{ij} \coloneqq  m( \flow^{-T}(S_{i}) \cap S_{j}).
\label{eq:ulam-matrix}
\end{equation}
Entries \(u_{ij}\) can be approximately computed by seeding a random (large) collection of initial values in each set \(S_{j}\) and evolving by \(\flow^{T}\).
The proportion of resulting endpoints of trajectories that land in a set \(S_{i}\) is entered in \(u_{ij}\).

The described process amounts to a Monte Carlo integration of the integral \(m( \flow^{-T}(S_{i}) \cap S_{j}) = \int_{S_{j}} \chi_{(\flow^{-T}(S_{i}))}dm(\bx)\)~\cite{dellnitz1997a,dellnitz1997}, or more generally, a procedure to compute a Galerkin approximation of the PF operator using piecewise constant functions~\cite{dellnitz1999}.
A computationally efficient implementation as a code GAIO~\cite{dellnitz2001} was shown to be able to approximate the eigenvalue spectrum of the PF operator and its eigenfunctions~\cite{dellnitz1999} both in the \(L^{2}\)-space and in fractional Sobolev spaces~\cite{froyland2014}.
{These set-oriented methods have also been connected with more classical geometric structures, such as finite time Lyapunov exponents~\cite{tallapragada2013set}.}

{ An alternative to the Ulam approximation is the so-called periodic approximation of the transfer operators, where a permutation matrix is used instead of the more general stochastic matrix.
While this approach dates back to the classical theory of measure-preserving transformations in~\cite{halmos1944a,katok1967,lax1971} and \cite[P.IV]{cornfeld1982}, it has been computationally revisited recently by~\cite{govindarajan2019} to approximate the continuous spectrum of Koopman and Perron--Frobenious operators and compute their spectral projectors.}

Eigenfunctions of the PF operator correspond to complex-valued distributions of points in state space that evolve according to the associated eigenvalue
\begin{equation}
  \label{eq:evolution-pf}
  \pf^{t} \rho = \lambda^{t} \rho.
\end{equation}
When \(\lambda=1\), \(\rho\) are invariant densities, which can be used to estimate the sets containing dynamical attractors, Sinai--Ruelle--Bowen measures, and partition the dynamics into invariant sets.
Even away from the limit of the Galerkin approximation, the Ulam matrix \(\mathbf{U}\) can be interpreted as a Markov chain transition matrix on a directed graph, which allows for a definition of almost-invariant sets in the state space~\cite{froyland2003} and a variational formulation of the problem of locating invariant sets in the state space.
Eigenfunctions for eigenvalues with \(\abs{\lambda} \not = 1\) are associated with the mixing properties and escape rates  in the state space of dynamical systems~\cite{Froyland2010b}.

The infinitesimal generator for PF is the Liouville operator \(\mathcal{A}\), defined analogously to \eqref{Eq:Koopman:InfinitesimalGenerator}, which satisfies
\begin{equation}
  \pf^{t} = e^{t \mathcal{A}}.
\end{equation}
Similar to  the Lie generator of the Koopman operator~\eqref{eq:Lie-operator}, when dynamics are specified using a velocity field \(\dot \bx = \mathbf{f}(\bx)\) the Liouville operator can be shown to satisfy
\begin{equation}
  \label{eq:Liouville-operator}
  \mathcal{A} \rho = -\operatorname{div}( \rho \mathbf{f} ),
\end{equation}
further leading to a partial differential equations that eigenfunctions  of the Liouville operator and PF operator must satisfy
\begin{equation}
  \label{eq:Liouville-eigenfunction}
   \mu \rho + \operatorname{div}( \rho \mathbf{f} ) = 0,
 \end{equation}
 whenever \eqref{eq:evolution-pf} is satisfied with \(\lambda = e^{\mu} \in \mathbb{C}/\{0\}\).

 Approximating \(\mathcal{A}\) instead of \(\pf\) leads to so-called \emph{simulation-free} numerical  techniques that can be interpreted either as finite-volume methods for the advection of PDEs, or as spectral collocation methods~\cite{froyland2013,bittracher2015}.

 Invariant eigenfunctions of both the Koopman and Perron--Frobenius operators have been used to extract invariant sets in the state space of dynamical systems.
 In function spaces where these operators are dual, eigenfunctions of both operators theoretically contain the same information.
 However, in reality, the choice may be made based on practical constraints.
 For example, approximation of invariant sets via Koopman eigenfunctions in \cite{budisic2012physd,levnajic2010,levnajic2015} relies on long-duration trajectories, while Ulam's approximation of PF typically requires short bursts of trajectories but seeded densely in the domain.

\subsection{Spectrum beyond eigenvalues}
\label{sec:koopman-spectrum}

Spectral characterization of infinite dimensional operators requires a separate treatment of two concepts: of \emph{the spectrum}, which generalizes the set of eigenvalues of finite dimensional operators, and of integration against \emph{spectral measures}, which takes the role of the eigenvalue-weighted summation appearing in the spectral decomposition theorem for normal matrices.
Standard textbooks on functional analysis commonly provide an introductory treatment of these concepts; however, among them we highlight~\cite{reed1978,reed1980,lax2002}, which include examples relating to the Koopman operator.

{In linear algebra, that is for operators induced by matrix multiplication, the \emph{spectrum} \(\sigma(\bT)\) is synonymous with the set of eigenvalues \(\lambda\), that is scalars \(\lambda\) such that
\begin{equation}
  \label{eq:eigenvalue}
    \bT \bxi = \lambda \bxi,
    \text{ or }
    (\bT - \lambda \mathbf{I}) \bxi = 0,
  \end{equation}
  for some vector \(\bxi \in \mathcal{G}\), termed the eigenvector.
  To extend the concept of a spectrum to operators \(\mathcal{T} : \mathcal{G} \to \mathcal{G}\) on Banach spaces, we interpret~\eqref{eq:eigenvalue} as a statement that eigenvalues  \(\lambda\) are those scalars for which the operator \(\mathcal{T} - \lambda \mathcal{I}\) does not have a bounded inverse.
  The spectrum \(\sigma(\mathcal{T})\) can be further classified into subsets based on the reason for why \((\mathcal{T} - \lambda \mathcal{I})^{-1}\) (\emph{the resolvent}) fails to exist as a bounded operator.

  For \(\lambda \in \sigma_{p}(\mathcal{T})\) (\emph{the point spectrum}) \(\mathcal{T} - \lambda \mathcal{I}\) is non-injective; this coincides with eigenvalues and is equivalent to a finite-dimensional spectrum.

  For \(\lambda \in \sigma_{c}(\mathcal{T})\) (\emph{the continuous spectrum}) the range of \(\mathcal{T} - \lambda \mathcal{I}\) is not the whole codomain (non-surjective), however it is dense in the codomain.
  This amounts to showing that for \(\forall \varepsilon\), there exists an observable \(g_{\varepsilon}\) such that the analogue to the relationship~\eqref{eq:eigenvalue} holds approximately  \(\norm{\mathcal{T}g_{\varepsilon} - \lambda g_{\varepsilon}} < \varepsilon\).
  In the context of Koopman theory, this was studied classically by V.~Rokhlin~\cite{rokhlin1966,nadkarni1979}, with more examples in~\cite[\S 13]{cornfeld1982}

  For \(\lambda \in \sigma_{r}(\mathcal{T})\) (\emph{the residual spectrum}) the range of \(\mathcal{T} - \lambda \mathcal{I}\) is not the whole codomain (non-surjective), and it is not even dense in the codomain.
  There are standard examples of Koopman operators that have continuous and residual spectra, e.g., those equivalent to multiplication operators and shift operators.
}
  The spectral decomposition theorem for normal matrices, which have orthogonal eigenvectors, states that the action of the matrix can be represented by the decomposition
  \begin{equation}
    \label{eq:spectral-decomposition-matrices}
    \bT^{n} \mathbf{g} = \sum_{\lambda \in \sigma(\bT)} \lambda^{n} \bxi_{\lambda}\avg{\bxi_{\lambda},\mathbf{g}},
  \end{equation}
  where \(\bxi_{\lambda}\avg{\bxi_{\lambda},\cdot}\) are orthogonal spectral projections onto eigenvectors \(\bxi_{\lambda}\).
  Generalizing~\eqref{eq:spectral-decomposition-matrices} to (infinite dimensional) Koopman operators requires that the operators are normal, which holds for certain dynamical systems.
  For example, when the flow \(\flow\) is invertible and preserves a finite measure \(\mu\), e.g., when the associated velocity field is divergence-free on a bounded domain, working with the Hilbert space \(\mathcal{H}\) of square-integrable observables \(g \in  \mathcal{H} = L^{2}(\mathcal{X},\mu)\) results in a unitary, and therefore normal, Koopman operator~\cite{vonneumann1932a}.
Then, the classical spectral resolution theorem (due to Hilbert, Stone, and Hellinger) applies to \(\koop\)~\cite{mezic2004physicad,mezic2005nd,lax2002} as
\begin{equation}
  \label{eq:spectral-decomposition}
  \begin{aligned}[t]
    \koop^{n} g = \int_{-\pi}^{\pi} e^{in\omega} d[\mathbb{E}(\omega) g]
    = \underbrace{\sum_{k} e^{in\omega_{k}} \mathbb{P}_{k} g}_{\text{atomic}} + \underbrace{\int_{-\pi}^{\pi} e^{in\omega} d[\mathbb{E}_{c}(\omega) g]}_{\text{continuous}},
\end{aligned}
\end{equation}
where the operator-valued \emph{spectral measure} \(\mathbb{E}\) forms a partition of unity and can be separated into atomic projections \(\mathbb{P}_{k}\) and the continuous part \(\mathbb{E}_{c}\).
This setting covers a wide-range of steady-state dynamics~\cite{eckmann1985,lasota1994,cornfeld1982}.

The atomic part of the spectral measure \(\mathbb{E}(\omega)\) is supported on frequencies \(\omega_{k}\) that yield eigenvalues \(e^{i\omega_{k}}\) of \(\koop\) and are associated with regular dynamics.
As for matrices, when the eigenvalues are simple, the eigenspace projections \(\mathbb{P}_{k}\) can be written using eigenfunctions \(\varphi_{k}\) of \(\koop\) as
\begin{equation}
\mathbb{P}_{k} g = \avg{g, \varphi_{k}}\varphi_{k}.\label{eq:eigenprojection}
\end{equation}
The atomic part of the spectral decomposition therefore aligns with the framework described in \cref{Sec:Koopman}.

There is no counterpart to the continuous spectral measure \(\mathbb{E}_{c}(\omega)\) in finite-dimensional settings; therefore, to interpret \(\mathbb{E}_{c}(\omega)\) and its connections to the dynamics requires routes that do not involve eigenvalues.
In the remainder of this section we summarize
\begin{inparaenum}[(a)]
  \item how existence/non-existence of \(\mathbb{E}_{c}(\omega)\) is connected to asymptotic statistical properties of the dynamics,
  \item how local structure of \(\mathbb{E}_{c}(\omega)\) connects to evolution of Koopman on subspaces of observables, and
  \item how to modify the space of observables in order to convert the continuous spectrum into a continuum of eigenvalues.
  \end{inparaenum}

\subsubsection{Global properties}

An operator-valued measure \(\mathbb{E}_{c}(\omega)\) can be converted to a ``plain'' scalar measure by choosing a function (observable) \(g \in L^{2}(\mathcal{X},\mu)\) and studying its autocorrelation sequence \(\avg{\koop^{n} g, g}\). Applying \eqref{eq:spectral-decomposition} here yields
\begin{equation}
  \label{eq:fourier-measure}
  \avg{\koop^{n} g, g} = \int_{-\pi}^{\pi}e^{in\omega} d\underbrace{\avg{\mathbb{E}(\omega) g, g}}_{=: \sigma_{g}(\omega)},
\end{equation}
where \(\sigma_{g}(\omega)\) is the \emph{Fourier or spectral measure} for the evolution of \(g\).

If the space \(L^{2}(\mathcal{X},\mu)\) is defined with respect to an ergodic measure \(\mu\), the autocorrelation function of the observable can be computed along trajectory initialized at any \(\bx\)
\begin{equation}
  \label{eq:autocorrelation}
\corr(n) \coloneqq  \lim_{K \to \infty} \frac{1}{K} \sum_{k = 0}^{K-1} g(\bx_{k+n})g(\bx_{k}) = \lim_{K \to \infty}  \frac{1}{K} \sum_{k = 0}^{K-1} [\koop^{n} g](\bx_{k})g(\bx_{k})
\end{equation}
and is directly related to the Fourier coefficients of \(\sigma_{g}\)
\begin{align}
\int_{-\pi}^{\pi} e^{i n \omega} d\sigma_{g}(\omega) &=  \avg{\koop^n g, g} = \lim_{K \to \infty}  \frac{1}{K} \sum_{k = 0}^{K-1} [\koop^{n} g](\bx_{k})g(\bx_{k})  =  \corr(n).\label{eq:correlation-fourier-spectrum}
\end{align}
In other words, the Fourier measure \(\sigma_{g}(\omega)\) is the Fourier transform of the \emph{autocorrelation function} \(\corr(n)\), which allows for characterization of irregular dynamics in terms of the measure \(\sigma_{g}\)~\cite{koopman1932pnas}.

In general \(\mathbb{E}\), and therefore some \(\sigma_{g}\), contain all three components that are mutually singular (e.g., the Lebesgue decomposition of a measure~\cite{strichartz1990}):
\begin{itemize}
\item atomic component, supported on frequencies of eigenvalues;
\item absolutely continuous component, having a spectral density and corresponding to mixing (stochastic-like) behavior{~\cite{rokhlin1966,eisner2015}};
\item singularly continuous component, having a fractal structure, and corresponding to weakly-anomalous transport{~\cite{zaks2017,pikovsky1995}}.
\end{itemize}
If \(\sigma_{g}\) is absolutely continuous, it has a spectral density (the Radon--Nikodym derivative) and then \(\corr(n) \to 0\), so the time-separated samples \(g(\bx_{k})\) will be asymptotically equivalent to independent random variables. If the same holds for all non-constant observables \(g \in \mathbf{1}^{\perp}\), the dynamics is mixing.
If \(\sigma_{g}\) has neither atomic components nor a spectral density, it is \emph{singularly continuous}. In this case, samples in the time trace are correlated infinitely often no matter their separation in time, but the correlation occurs rarely enough that \emph{on average} they appear uncorrelated. In this case, the measure can be thought of as having a fractal structure as its support is neither a full interval, nor a collection of discrete points.
If this holds for all \(g \in \mathbf{1}^{\perp}\) then the dynamics is weakly mixing; this is the signature behavior of the so-called weakly-anomalous transport~\cite{zaslavsky2002,zaslavsky1991,knill1998,hof1998}.

\subsubsection{Local properties}
To interpret how the local variation of \(\sigma_{g}(\omega)\) is associated with the dynamics, we investigate integrals over spectral intervals \([a,b] \subset [-\pi,\pi]\) \(\int_{a}^{b}e^{in\omega} d\sigma_{g}(\omega)\). It can be shown that such restrictions on the spectral domain are equivalent to restrictions of \(\koop\)  to certain subspaces of observables~\cite[Prop. 2.7]{queffelec2010}.
More precisely, there exists an orthogonal projection \(P_{[a,b]}\) such that the following equality holds
  \begin{equation}
    \label{eq:localization}
    \int_{a}^{b}e^{in\omega} d\sigma_{g}(\omega) = \int_{-\pi}^{\pi} e^{in\omega }\mathbf{1}_{[a,b]} d\sigma_{g}(\omega) = \avg{ \koop^{n} P_{[a,b]} g, P_{[a,b]} g}.
  \end{equation}
  The range of the projection \(P_{[a,b]}\) is the invariant subspace
  \begin{equation}
    \mathcal{H}_{[a,b]} = \{ g \in L^{2}(\mathcal{X},d\mu) \colon \sigma_{g}(\mathbb{T} / [a,b]) = 0 \}. \label{eq:localized-subspace}
  \end{equation}
  In other words, localizing the spectral integral results in a compression of \(\koop^{t}\) to some dynamically-invariant subspace of the space of observables.
  This holds more generally not only for intervals, but for arbitrary measurable sets in the spectral domain.

  In applications, the selection of observables often comes before the analysis of dynamics.
  If the chosen observable happens to already be in some subspace \(\mathcal{H}_{[a,b]}\), the restricted integral would be equivalent to its full evolution \(\avg{\koop^{n} g,g} = \avg{\koop^{n} P_{[a,b]} g, P_{[a,b]} g} \).
  In other words, if one would try to infer the ``full'' evolution of the \(\koop\) from a single observable, an unintentional  choice of the observable \(g\) from an invariant subspace \(\mathcal{H}_{[a,b]}\) may result that instead of the full operator, \(\koop\), only its compression \(P_{[a,b]}^{\top}\koop P_{[a,b]}\)  may be reconstructed.
  On the other hand, there exists a subspace of observables for which associated Fourier measures  are ``maximal'', in the sense that any zero set for a maximal measure, i.e. spectral ``bandwidth'' exhibiting no dynamics, is a zero set for any other Fourier measure.
  This implies that a judicious choice of observable can make it possible to fully characterize statistical properties of the system from the autocorrelation timeseries of a single observable.

\begin{figure}[t]
  \centering
  \vspace{-.2in}
\includegraphics[width=0.45\linewidth]{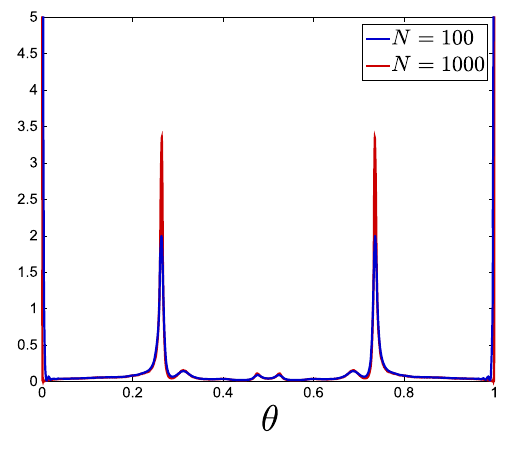}
  \includegraphics[width=0.45\linewidth]{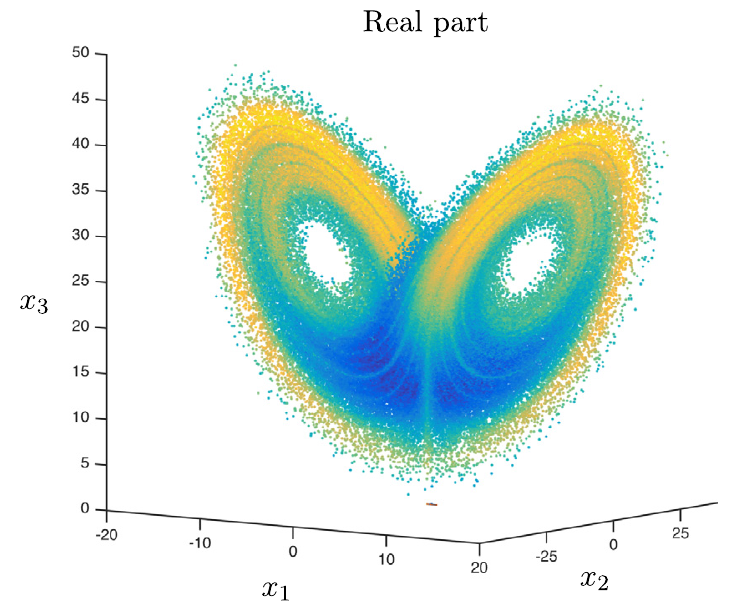}
  \vspace{-.15in}
  \caption{Approximation to the density of spectral measure for the Lorenz'63 model (left) and the real part of the projection of the evolution onto the cyclic vector of the subspace associated with spectral interval \([0.24,0.28]\) (right). \textit{Reproduced with permission, from Korda et al. 2020 Appl. Comput. Harmon. Anal.~\cite{korda2020}.}}
  \label{fig:koopman-continuous-spectral-measure}
  \vspace{-.25in}
\end{figure}

  If the content for many spectral intervals \([a,b] \subset \mathbb{T}\) is of interest, we may want to approximate the weak derivative \(d\sigma_{g}/d\theta\) and visualize it. In~\cite{korda2020} the Fourier coefficients (or trigonometric moments), computed by ergodic averages~\eqref{eq:correlation-fourier-spectrum}, are used to formulate the moment problem for the density~\(d\sigma_{g}/d\theta\), which is solved using a variational approach based on the Christoffel--Darboux kernel.
  Based on the approximated density, for any given interval \([a,b]\) in the spectral domain one can construct a type of an ergodic average that computes the associated projection \(P_{[a,b]}\) of the evolution, resulting in the analog of eigenvectors for the continuous spectral measure.
  Constructions involved relate to the HAVOK algorithm~\cite{brunton2017natcomm,arbabi2017} (see \cref{Sec:HAVOK}) due to the correlation function that connects time-delayed samples of the observable.
  \Cref{fig:koopman-continuous-spectral-measure} demonstrates this approach on an example of the Lorenz'63 dynamical system that is known to be mixing on the attractor~\cite{Luzzatto2005} by computing its spectral density and a spectral projection for a spectral interval containing significant dynamical content.

\subsubsection{Removing the continuous spectrum by ``rigging''}
To illustrate that existence of the continuous spectrum (as defined by non-surjectivity of \(\koop - \lambda \mathcal{I}\)) is not necessarily connected with irregular behavior, \cite{mezic2019} studies in detail the pendulum
\begin{equation}
  \label{eq:physical-pendulum}
  \dot \alpha = v, \quad  \dot v = -\sin\theta
\end{equation}
which can be converted to action-angle variables with \(I \geq 0 \) representing the conserved quantity, and a periodic variable  \(\theta \in \mathbb{S}^{1}\) the momentum,
\begin{equation}
  \label{eq:action-angle-oscillator}
  \dot I = 0,\quad \dot \theta = I.
\end{equation}
The exact solution is
\begin{equation}
  \label{eq:aa-solution}
  I(t) = I_{0},\quad \theta(t) = I_{0} t + \theta_{0} \pmod{2\pi}.
\end{equation}
Any observable \(g \colon \mathbb{R}_{+} \times \mathbb{S}^{1} \to \mathbb{R}\) without dependence on the angle coordinate \({g(I, \theta) = g(I)}\) is clearly an eigenfunction of \(\koop^{t}\) with \(\lambda = 1\), since
\begin{equation}
\koop^{t} g(I_{0}) = g( I(t) ) = g(I_{0}) \label{eq:eigenvalue-at-1},
\end{equation}
implying that the eigenvalue \(\lambda = 1\) has an infinite multiplicity. This rules out ergodicity, and therefore mixing, with respect to the \(L^{2}\) space of observables defined over any annulus \([I_{1},I_{2}] \times \mathbb{S}^{1}\) of positive area. In other words, the atomic spectrum detects regular behavior of dynamics.

At the same time, no function with a variation in the angular direction is an eigenfunction, as the only translation-invariant function on \(\mathbb{S}^{1}\) is a constant, despite all trajectories being clearly periodic. However, if observables are taken from a space that includes generalized functions (distributions),
then a Dirac-$\delta$ supported on a single level set of $I$
\begin{equation}
g_{c}(I,\theta) = \delta(I - c) e^{i\theta}\label{eq:distribution}
\end{equation}
would indeed be a (generalized) eigenfunction as
\begin{equation}
  \koop^{t}g_{c}(I,\theta) = \delta(I - c) e^{i I t + \theta} = e^{i c t} \delta(I - c)e^{i\theta} = e^{i c t} g_{c}(I,\theta).
\end{equation}
{ This example has also been revisited in~\cite{bollt2021}.}

The example above illustrates how a \emph{rigged}\footnote{This metaphor is intended to evoke a utilitarian \emph{rigging of a ship}, rather than a nefarious \emph{rigging of a match}.}, or \emph{equipped} Hilbert space can convert a continuous spectrum (containing no eigenvalues) to a continuum of eigenvalues~\cite{antoniou1997a,suchanecki1996,tasaki1993,madrid2002a}.
Instead of a single \(L^{2}\) space of observables, the approach employs a so-called Gelfand triple \(\Gamma \subset L^{2} \subset \Gamma^{\dagger}\), where ``the rigging'' consists of \(\Gamma\), a subset of judiciously-chosen test functions, and \(\Gamma^{\dagger}\) its dual. In the example above, \(\Gamma^{\dagger}\) contains generalized functions (distributions).
By enlarging the domain of the Koopman operator, surjectivity of \((\koop-\lambda \mathcal{I})\) can be mostly restored, resulting in shrinking of the continuous spectrum to a set of discrete values~\cite{slipantschuk2020}.

As a result, the continuous projection measure in the Hellinger--Stone spectral theorem~\eqref{eq:spectral-decomposition}
can be replaced by an integral against eigenvector-based projections, analogous to~\eqref{eq:spectral-decomposition-matrices}
  \begin{equation}
    \label{eq:rigged-spectral-decomposition}
    \avg{\rho,\koop g} = \sum_{\lambda} \lambda \avg{\rho, \psi_{\lambda}}\avg{\varphi_{\lambda},g},
  \end{equation}
  for an observable \(g \in \Gamma\) and a density \(\rho \in \Gamma^{\dagger}\), while \(\varphi_{\lambda}, \psi_{\lambda}\) are elements of a biorthonormal basis in \(L^{2}\).
  The extended spectrum \(\lambda\) now contains both the \(L^{2}\)-eigenvalues of \(\koop\) and Ruelle--Pollicott resonances \cite{Ruelle1986,Pollicott1985} associated with infinitesimal stochastic perturbations of the Koopman operator~\cite{Chekroun2014,Chekroun2020}.
  For example, while the map \(x \mapsto 2x\) on \(\mathbb{S}^{1}\) has only a constant eigenfunction and only the eigenvalue at 1, a representation in the rigged Hilbert space where \(\Gamma\) are analytic functions yields functions \(\varphi_{k}, \psi_{k}\) that are related to Bernoulli polynomials, and values \(\lambda_{k} = 2^{-k}\)~\cite{antoniou1997}.

  While the rigged Hilbert space framework has existed since the 1950s and Gelfand, the approach has been used to analyze the Koopman and Perron--Frobenius operators only since the 1990s, and for quantum theory a decade later \cite{madrid2005,madrid2002a}.
  Only in the past two years have the modern numerical approaches such as DMD, started to connect to this theory \cite{mezic2019,slipantschuk2020}, so we expect the growth of interest in this area in the coming years.

\subsection{Koopman operators for nonautonomous and stochastic dynamics}
\label{sec:koopman-stochastic}

The theory developed so far was based on the time-independent or autonomous dynamics~\eqref{Eq:ContinuousDynamicsGeneral} and~\eqref{Eq:DiscreteDynamics}.
This clearly does not cover most models used in practice.
Common sources of time variability include changes in system parameters, the presence of input forcing, stochastic noise, and control or actuation.
Time variability induced by feedback control is highly structured, and thus Koopman theory can be developed in more detail depending on the structure, as described in~\cref{Sec:Control}.
Even though the original justification for DMD-style algorithms was based on the autonomous Koopman framework, the algorithms were applied to data generated by nonautonomous dynamics, either by tacitly assuming that the time variation is negligible, or by employing various sliding window techniques.

\subsubsection{Sliding and multiresolution analysis} Consider the nonautonomous dynamics
\begin{equation}
  \label{eq:nonautonomous}
  \dot \bx = \bF(\bx, t),\quad \bx \in \mathcal{X},
\end{equation}
and assume that over a time window \(t \in [\tau, \tau+T]\) the function \(\bF\) remains approximately constant
\begin{equation}
  \label{eq:approx-constant}
  \bF(\cdot, t) \approx \bF(\cdot, \tau).
\end{equation}
Furthermore, assume that this holds over a continuous range of starting points \(\tau\), while maintaining a constant window size \(T\).

A \emph{sliding window} implies that the snapshots of observations generated by data collected over each time window \([\tau_{i}, \tau_{i}+T]\), for some \(i = 1,2,\dots\), are separately processed by a DMD algorithm to produce eigenvalues \(\lambda_{k}(\tau, T)\) and DMD modes \(\phi_{k}(\tau,T)\) that depend on the parameters of the time window.
This approach is neither new nor unique to Koopman analysis.
In signal processing, computing the discrete Fourier transform over a sliding window goes under several names, including the \emph{short time Fourier transform}, \emph{spectrogram}, \emph{sliding discrete Fourier transform}, or \emph{time-dependent Fourier transform}, and is a standard topic in many signal processing textbooks~\cite{Oppenheim2014, Stoica2005}.
If, in addition to the starting point \(\tau\), the length of the window \(T\) is systematically varied as well, this is known as \emph{multiresolution} or \emph{multiscale analysis}.

The first systematic treatment of multiresolution analysis in the context of DMD presented
the basic sliding strategy of computing DMD in several passes over the same set of snapshots~\cite{kutz2016book,kutz2016siads}.
Within each pass, the window size was kept constant, and the starting point of the window was moved in non-overlapping fashion.
Several additional strategies for sampling the data and moving the window were discussed, including connections to classical signal processing techniques, such as the Gabor transform~\cite{kutz2016siads,brunton2016b}.
The overlap between consecutive windows can be exploited to compute a more accurate global reconstruction~\cite{dylewsky2019}, and to reduce the computational effort required to compute the DMD~\cite{zhang2019}.
Sliding window strategies have been found useful even in autonomous systems~\cite{chen2012jns}; for example, analytic solutions of autonomous dynamics that follow heteroclinic connections between an unstable and a stable fixed point were analyzed, where a well-tuned sliding window is able to discern the difference in DMD spectra near each of the fixed points~\cite{page2019}.

\subsubsection{Process formulation for the nonautonomous Koopman operator}

The extension of definitions associated with the Koopman operator follow two paths for extending the general theory of dynamical systems: the so-called \emph{process} formulation, or the \emph{skew-product} formulation.

For an autonomous system, such as \(\dot \bx(t) = \bA \bx(t), \bx(t_{0}) = \bx_{0}\), the time dependence of the corresponding flow map \(\bx_{0} \mapsto \bx_{0}e^{\bA(t-t_{0})}\) can be written in terms of the duration \(t-t_{0}\) of the time interval over which the dynamics is evolving.
The semigroup property~\eqref{eq:semigroup-property} captures the consistency between advancing dynamics \(t \to t + \Delta t\) in one step, or in two steps \(t \to t+ \Delta t/2\) and then \(t + \Delta t/2 \to t + \Delta t\).

To illustrate the \emph{process} formulation of nonautonomous systems~\cite{kloeden_rasmussen_2011,caraballo2017,macesic2018siads}, consider the simple nonautonomous ODE
\begin{align}
  \label{eq:nonautonomous-example}
    \dot x(t) &= \cos(t) x, \ x(t_{0}) = x_{0} \\
    \shortintertext{solved by}
    x(t) &= \flow_{t_{0}}^{t}(x_{0}) = x_{0} e^{\sin(t_{0}) - \sin(t)}.
\end{align}
The time dependence of the flow map cannot be expressed simply in terms of the duration \(t-t_{0}\).
The consistency of advancing dynamics across adjoining time intervals is now captured by the \emph{cocycle} property
\begin{equation}
  \label{eq:cocycle}
  \flow^{t}_{t_{0}}(\bx) = \flow^{t}_{\tau}( \flow^{\tau}_{t_{0}}(\bx) ), \forall \bx,\ 0 \leq t_{0} \leq \tau \leq t.
\end{equation}

The Koopman operator can now naturally be formed as a composition operator with this two-parameter flow map~\cite{macesic2018siads}
\begin{equation}
  \label{eq:cocycle-koopman}
  \koop^{t}_{t_{0}} g = g \circ \flow^{t}_{\tau}.
\end{equation}
Its Lie generator, given by
\begin{equation}
\gen_{t_{0}} g \coloneqq  \lim_{t\to t_{0}} \frac{\koop^{t}_{t-t_{0}} g-g}{t},\label{eq:cocycle-koopman-generator}
\end{equation}
itself depends on time \(t_{0}\), in contrast to the autonomous case, and in turn eigenvalues and eigenvectors also depend on time.
Such a construction of the Koopman operator appears to match the sliding-window approaches to DMD discussed earlier.
However,~\cite{macesic2018siads} demonstrate that a sliding-window approach may make large mistakes, especially when the window overlaps the region of rapid change in system parameters; the same source describes an algorithm that is able to detect the local error of a DMD approximation and adjust accordingly, which is particularly effective for so-called hybrid systems, in which the time dependency of the equations is discontinuous.

A word of caution is needed here; in nonautonomous systems (finite- or infinite-dimensional) eigenvalues of the system matrix do not always correctly predict the stability of trajectories.
For example, it is possible to formulate a finite-dimensional dynamical system \(\dot \bx = \bA(t) \bx\) such that eigenvalues of \(\bA(t)\) (the generator of the flow) are time-independent and have negative real value, while the system admits a subspace of unstable solutions~\cite{markus1960,meiss2007}.
Furthermore, \cite{macesic2018siads} find that methods based on Krylov subspaces, e.g., snapshot-based DMD algorithms, result in substantial errors in the real parts of eigenvalues when the time-dependence of eigenvalues is pronounced, and suggest that the problem can be mitigated by a guided selection of observables.

\subsubsection{Skew-product formulation for the nonautonomous Koopman operator}

We turn now to the \emph{skew-product} formulation,  which can incorporate both deterministic and stochastic nonautonomous systems.
Consider again the ODE~\eqref{eq:nonautonomous-example} but now introduce an additional periodic variable \(y \in \mathbb{S}^{1}\)
\begin{equation}
  \label{eq:nonautonomous-example-skew}
  \begin{aligned}
  \dot x(t) &= \cos(y) x, \\
  \dot y(t) &= 1.
\end{aligned}
\end{equation}
The added variable plays the role of time, re-establishing the semigroup property of the associated flow map which now acts on the extended space \(\flow^{t} \colon \mathcal{X} \times \mathbb{S}^{1} \to  \mathcal{X} \times \mathbb{S}^{1}\).
The dynamics of \(y(t)\) is itself autonomous and is sometimes referred to as the \emph{driver} for the \emph{driven} dynamics of \(\bx(t)\). The skew-product construction appeared in classical ergodic theory, for example~\cite{cornfeld1982}.

General analysis of skew-products does not require that \(y(t)\) is as simple as in~\eqref{eq:nonautonomous-example-skew};
rather, an extension of~\eqref{Eq:DiscreteDynamics} to the more general skew form
\begin{equation}
  \label{eq:skew-dynamics}
  \begin{aligned}
  \bx_{k+1} &= \bF(\bx_{k},\by_{k}) \\
  \by_{k+1} &= \bG(\by_{k})
\end{aligned}
\end{equation}
can be treated in an analogous fashion.
Assuming that the driving dynamics evolves on a compact state space, or otherwise is measurable with a well-defined invariant measure, is sufficient to justify computation of eigenfunctions using ergodic averages~\eqref{eq:harmonic-average} for time-periodic systems and for systems driven by quasiperiodic dynamics~\cite{budisic2012physd,susuki2020}.

There are two possible formulations of the Koopman operator associated with \eqref{eq:skew-dynamics}.
The first formulation treats the skew-flow as an autonomous system on the extended state space and acts by composing an observable \(g \in \mathcal{G}(\mathcal{X} \times \mathcal{Y})\) with the flow
\begin{equation}
  \label{eq:skew-Koopman-auto}
  [\koop g](\bx,\by) \coloneqq g( \bF(\bx,\by), \bG(\by) ).
\end{equation}
Since any particular observable on the original state space \(h \colon \mathcal{X}\mapsto\mathbb{C}\) can be trivially extended to \(h(\bx,\by) = h(\bx)\), this formulation is sufficient for studying dynamics of a finite collection of observables.
However, since the space \(\mathcal{G}(\mathcal{X} \times \mathcal{Y})\) is larger than \(\mathcal{G}(\mathcal{X})\), representing \(\koop\) in a particular basis of observables requires working with a practically larger set, e.g., instead of monomials \(x^{k}, k=1,2,\dots\), one has to work with \(x^{k}y^{j}, k=1,2,\dots, j=1,2,\dots\). The problem is, of course, more acute the higher the dimension of \(\mathcal{Y}\).

\subsubsection{Stochastic Koopman operator}

An alternative formulation of the Koopman operator for \eqref{eq:skew-dynamics} acts on the observables in the original state variable only, \(\mathcal{G}(\mathcal{X})\), but retains a parametric dependence of the Koopman operator on the initial condition of the driving system.
The skew-flow map used in this case is the interpretation of the second argument in \(\bF(\bx, \by)\) as a parameter for the flow map \(\bF_{\by}(\bx) \coloneqq \bF(\bx, \by)\).
The replacement for the semigroup property of autonomous systems is then the skew-flow property
\begin{equation}
  \bF_{\by}^{t+s} = \bF_{\bG^{s}(\by)}^{t}\circ \bF_{\by}^{s}.
\end{equation}
The Koopman operator is then defined as the composition operator with respect to the flow
\begin{equation}
  \label{eq:skew-koopman}
  [\koop^{t}_{\by} g](\bx) \coloneqq g( \bF^{t}_{\by}(\bx) ).
\end{equation}
In contrast to the ``cocycle Koopman operator'' \eqref{eq:cocycle-koopman}, in which time-dependence is represented by an additional time-like parameter, in this ``skew Koopman operator'' the time-dependence is an added state-like parameter.

Depending on the properties of the driving system, this framework can encompass not only systems such as \eqref{eq:nonautonomous-example-skew}, but also so-called \emph{random dynamical systems} (RDS)~\cite{arnold1998,caraballo2017} for which the driving system on \(\mathcal{Y}\) is simply assumed to be a measure-preserving dynamical system~\cite{macesic2020, crnjaric-zic2020}.

While several special cases have been discussed in recent literature~\cite{crnjaric-zic2020}, we focus here on the case of Markovian RDS.
This is the case for dynamics generated by a nonlinear stochastic differential equation
\begin{equation}
  \label{eq:sde}
  d\bx(t) = \mathbf{f}(\bx) dt + \mathbf{\sigma}(\bx) d\mathbf{w}(t).
\end{equation}
It is then possible to define a \emph{stochastic} Koopman operator by computing the expectation of the skew-Koopman operator~\eqref{eq:skew-koopman} (with continuous time domain) over the invariant probability of the driving flow~\cite{arnold1998,mezic2004physicad}:
\begin{equation}
  \label{eq:stochastic-koopman}
  [\koop^{t}_{S} g](\bx) \coloneqq \mathbb{E}_{\by}\{ g( \bF^{t}_{\by}(\bx) ) \}.
\end{equation}
When the driving system is a stochastic system, this establishes the Koopman operator as the action of the dynamics, averaged over the distribution of the stochastic input.

Furthermore, assuming continuous and bounded functions on the right-hand side of the SDE~\eqref{eq:sde}, it can be shown that \eqref{eq:stochastic-koopman} is a strongly-continuous semigroup, which is a consequence of the Chapman--Kolmogorov equation, with a well-defined generator \(\gen_{S}\) acting on a space of twice-differentiable observables.
In this case, the analogue of the PDE formulation of the generator~\eqref{eq:Lie-operator} is given by
\begin{equation}
  \label{eq:stochastic-lie-pde}
  \gen_{S} g = \nabla g\ \mathbf{f} + \frac{1}{2} \operatorname{Tr}( \boldsymbol{\sigma} \nabla^{2} g \boldsymbol{\sigma}^{\top}).
\end{equation}
{The generator-based approach to the Koopman operator for SDEs has been pursued in~\cite{klus2020}, which includes several applications and discusses the model-predictive control (MPC) framework based on the Koopman generator.}

Several DMD algorithms have been adapted to the RDS setting for the Koopman operator~\cite{crnjaric-zic2020, takeishi2017}, with convergence assurances.
An explicit optimization-based approximation of the Koopman operator in the RDS context has been given in~\cite{sinha2020a}.

In summary, the theory behind DMD-style algorithms for nonautonomous and stochastic dynamical systems have seen a rapid development in recent years, bringing about both justification for applying such algorithms at first glance ``off-the-label'', and providing additional guidance for reduction of bias and errors in computation of eigenvalues and eigenmodes.

\subsection{Partial differential equations}
\label{sec:koopman-PDEs}
There already exists a clear connection of Koopman theory to PDEs~\cite{kutz2018}.
Just as with ordinary differential equations, the goal is to discover a linearizing transformation of the governing nonlinear PDE to a new PDE model which evolves linearly in the new coordinate system.
{{
Thus instead of a finite dimensional autonomous ordinary differential equation \eqref{Eq:ContinuousDynamicsGeneral} on a state space \(\mathcal{X}\subseteq \mathbb{R}^{n}\),
we instead generate a Koopman operator that maps functions to functions, i.e. infinite-dimensional spaces to infinite-dimensional spaces.  Thus the flow map operator, or time-\(t\) map,  \(\flow^{t} \colon \mathcal{X} \to \mathcal{X}\) advances initial conditions forward along the trajectory by a time $t$ where \(\mathcal{G}(\mathcal{X})\) is again a set of \emph{measurement functions}
\(g \colon \mathcal{X} \to \mathbb{C}\).  Such a transformation is no different than, for instance, a Fourier transform that maps a spatially dependent PDE to a wavenumber representation.  This is a common solution technique for linear PDEs since the mapping to the wavenumber representation also diagnolizes the evolution dynamics so that each Fourier mode satisfies an ODE.  Importantly, PDEs generically have both a discrete and continuous spectrum, unlike the discrete spectrum of finite-dimensional dynamical systems.  In practice, however, computing solutions of PDEs requires discretization of the solution, thus rendering the evolution dynamics finite-dimensional.

Linearizing transforms for nonlinear PDEs have historically been done}}
%
in a number of contexts, specifically the Cole--Hopf transformation for Burgers' equation with diffusive regularization and the inverse scattering transform (IST) for completely integrable PDEs.
{{For IST, for instance, a rigorous analytic formulation can be explicitly constructed that maps infinite-dimensional function spaces to new infinite-dimensional functions spaces whose evolution dynamics is given by a linear operator~\cite{ablowitz1974inverse}.}}
The connection of these two analytic transformations is considered below.
Such linearizing transformations have been difficult to achieve in practice.  However, data-driven methods have opened new pathways for constructing these transformations.
Neural networks, diffusion maps, and time-delay embeddings all allow for the data-driven construction of mappings capable of transforming nonlinear PDEs into linear PDEs whose Koopman operator can be constructed.
In this section, we consider the connection of some of the historically developed methods to Koopman theory.
In \cref{Sec:Observables}, we show how such linearizing embeddings are constructed with modern data-driven methods.

The dynamics of nonlinear PDEs evolve on manifolds which are often difficult to characterize and are rarely known analytically.
However, an appropriate choice of coordinates can in some cases, {\em linearize} the dynamics.  For instance, the nonlinear evolution governed by Burgers' PDE equation can be linearized by the Cole--Hopf transformation~\cite{hopf50,cole51}, thus providing a linear model that can trivially describe the evolution dynamics.
Such exact solutions to nonlinear PDEs are extremely rare and do not often exist in practice, with the inverse scattering transform for Korteweg--deVries, nonlinear Schr\"odinger, and other integrable PDEs being the notable exceptions~\cite{ist}.
As will be shown in \cref{sec:NNKoop}, neural networks provide a data-driven method to learn coordinate embeddings such as those analytically available from the Cole--Hopf transform and IST.

To demonstrate the construction of a specific and exact Koopman operator, we consider Burgers' equation with diffusive regularization and its associated Koopman embedding~\cite{kutz2018,page2019,page2018a,balabane2020koopman}.  The evolution is governed by diffusion with a nonlinear advection~\cite{burgers}:
\begin{equation}
  u_t + u u_x - \epsilon u_{xx} = 0  \,\,\,\,\,\,\,\, \epsilon>0, \,\,  x\in[-\infty,\infty].
  \label{burgers}
\end{equation}
When $\epsilon=0$, the evolution can lead to shock formation in finite time.  The presence of the diffusion term regularizes the PDE, ensuring continuous solutions for all time. In independent, seminal contributions,
Hopf~\cite{hopf50} and Cole~\cite{cole51} derived a transformation that linearizes the PDE.  The Cole--Hopf transformation
is defined as follows
\begin{equation}
  u=-2\epsilon v_x/v \, .
  \label{colehopf}
\end{equation}
The transformation to the new variable $v(x,t)$ replaces the nonlinear PDE (\ref{burgers}) with the linear, diffusion equation
\begin{equation}
  v_t = \epsilon v_{xx}
\end{equation}
where it is noted that $\epsilon>0$ in (\ref{burgers}) in order to produce a well-posed PDE.  Denoting $\hat{v}=\hat{v}(k,t)$ as the Fourier transform of $v(x,t)$ with wavenumber $k$ gives the analytic solution
\begin{equation}
  \hat{v}=\hat{v}_0 \exp( -\epsilon k^2 t )
  \label{bkoop}
\end{equation}
where $\hat{v}_0=\hat{v}(k,0)$ is the Fourier transform of the initial condition $v(x,0)$.  Thus to construct the Koopman operator, we can then combine the transform to the variable $v(x,t)$ from (\ref{colehopf})
with the Fourier transform to define the observable $ g(u) = \hat{v}$.
This gives the Koopman operator
\begin{equation}
   {\cal K} = \exp( -\epsilon k^2 t )  \, .
   \label{eq:burg_koo}
\end{equation}
This is one of the rare instances where an explicit expression for the Koopman operator and the
observables can be constructed analytically.  As such, Burgers' equation allows one to build explicit representations of Koopman operators that characterize its nonlinear evolution~\cite{kutz2018,page2018a}.

The inverse scattering transform~\cite{ist} for other canonical and integrable PDEs, such as the Korteweg--deVries and nonlinear Schr\"odinger equations,  also can lead to an explicit expression for the Koopman operator, but the scattering transform and its inversion are much more difficult to construct in practice.  Peter Lax developed a general mathematical framework that preceded IST theory and which provided a
general principle for associating nonlinear evolutions with linear operators so that the eigenvalues of the linear operator are integrals of the nonlinear equation~\cite{lax1968integrals}.   The scattering theory and its association with nonlinear evolution equations was then placed on more rigorous foundations by the seminal contribution of Ablowitz, Kaup, Newell and Segur known as the AKNS scheme~\cite{ablowitz1974inverse}.

{In brief, the method developed by Lax for constructing analytic solutions for nonlinear evolution equations involved constructing suitable linear operators whose compatibility condition was the evolution equation itself.  For a given nonlinear evolution equation
\begin{equation}
   u_t = N (u, u_x, u_{xx} , \cdots) ,
\end{equation}
the goal is to posit a {\em Lax pair} of operators
\begin{equation}\label{eq:laxpair}
\begin{aligned}
\mathcal{A}\phi &= \lambda \phi \\
\frac{d}{dt}\phi &= \mathcal{B} \phi
\end{aligned}
\end{equation}
where $\mathcal{A}$ and $\mathcal{B}$ are linear operators. Specifically, the operator $\mathcal{A}$ is a spatial operator that is self-adjoint and does not depend explicitly on $t$, while the operator $\mathcal{B}$ is a time-evolution linear operator.
In the context of Koopman theory, the time evolution \(\frac{d}{dt}\phi = \mathcal{B}\phi\) can be taken to be equivalent to
\( \frac{d}{dt}g = \gen{} g \)
in \eqref{eq:koopman-dynamics-continuous}.\footnote{It is typical to denote the \emph{spatial} Lax operators by \(\mathcal{L}\); in our case this may lead to the confusion with the Koopman generator.}
Importantly, $\mathcal{A}$, $\mathcal{B}$ and the evolution equation for $u(x,t)$ must be all self-consistent, or compatible, in order for the Lax theory to hold.  Self-consistency is achieved by taking the time-derivative of \eqref{eq:laxpair} with respect to time  and enforcing solutions that have an iso-spectral evolution with respect to these operators so that $\lambda_t=0$.  This then gives
\begin{equation}
   \frac{d}{dt}\mathcal{A}  + [\mathcal{A},\mathcal{B}] =0
\end{equation}
where $[\mathcal{A},\mathcal{B}]=\mathcal{A}\mathcal{B}-\mathcal{B}\mathcal{A}$ represents the commutator of the operators.  }Importantly, within this framework, the operators $\mathcal{A}$ and $\mathcal{B}$ are linear. Thus once found, the evolution dynamics in the transformed coordinate system is {\em linear}, much like what occurs in the Burgers' example.  Of course, such a general mathematical framework only holds for integrable PDEs~\cite{ablowitz1974inverse}.  However, it does show that the Koopman operator framework is directly associated with the Lax pairs, and in particular with the linear time evolution operator $\mathcal{B}=\gen{}$, connecting IST and Koopman theory explicitly.
Parker and Page~\cite{parker2019koopman} recently developed a detailed analysis of fronts and solitons in a variety of systems and explicitly connected their findings to the IST.  Moreover, they showed how {\em two} Koopman decompositions, upstream and downstream of the localized structure, can be used to derive a full Koopman decomposition that leverages the IST mathematical machinery.

{{In addition to the explicit Koopman representations provided by the Cole--Hopf transformation and IST for a class of PDEs, data-driven methods are providing a suite of techniques for constructing Koopman operators for PDEs. DMD provides an efficient algorithm to do so~\cite{kutz2018,page2019,page2018a,balabane2020koopman}.  And more recently, Nakao and Mezic~\cite{nakao2020}
introduce the concept of Koopman eigenfunctionals of PDEs and use the notion of conjugacy to develop spectral expansion of the Koopman operator. Their methodology is developed on a number of PDEs and demonstrates a promising direction of research for constructing Koopman PDE representations. Finally, as}}
will be shown in \cref{sec:NNKoop}, neural networks provide an ideal, data-driven mathematical construct to learn coordinate embeddings such as those analytically available from the Cole--Hopf transform and IST.
Indeed, this has been done even for the Kuramoto--Sivashinsky equation~\cite{gin2019deep}.
Additional connections between the Cole--Hopf transform and Koopman eigenfunctions in the context of ODEs are discussed in \cite{bollt2018}.

\section{Data-driven observable selection and embeddings}\label{Sec:Observables}
Linearizing a nonlinear dynamical system near fixed points or periodic orbits provides a \emph{locally} linear representation of the dynamics~\cite{guckenheimer1986}, enabling the limited use of linear techniques for prediction, estimation, and control~\cite{brunton2019data}.
The Koopman operator seeks \emph{globally} linear representations that are valid far from fixed points and periodic orbits.
In the data-driven era, this amounts to finding a coordinate system, or \emph{embedding}, defined by nonlinear observable functions that span a Koopman-invariant subspace.
Finding, approximating, and representing these observables and embeddings is still a central challenge in modern Koopman theory.
Dynamic mode decomposition~\cite{schmid2010jfm,rowley2009jfm,kutz2016book} from \cref{Sec:DMD} approximates the Koopman operator restricted to a space of linear measurements with a best-fit linear model advancing these measurements from one time to the next.
However, linear DMD alone is unable to capture many essential features of nonlinear systems, such as multiple fixed points and transients.
It is thus common to augment DMD with nonlinear functions of the measurements~\cite{williams2015jnls}, although there is no guarantee that these functions will form a closed subspace under the Koopman operator~\cite{brunton2016plosone}.
In this section, we cover several leading approaches to identify and represent Koopman embeddings from data, including the extended dynamic mode decomposition~\cite{williams2015jnls} and methods to directly identify eigenfunctions~\cite{kaiser2017arxiv}, the use of time-delay coordinates~\cite{brunton2017natcomm}, as well as machine learning approaches such as diffusion maps~\cite{giannakis2019} and deep neural networks~\cite{wehmeyer2017arxiv,mardt2017arxiv,Takeishi2017neurips,yeung2017arxiv,otto2017arxiv,lusch2017arxiv}.
{The convergence and robustness of these data-driven approximation algorithms is also discussed.}

\subsection{Extended DMD}\label{Sec:EDMD}
Although DMD~\cite{schmid2010jfm} has become a standard numerical approach to approximate the Koopman operator~\cite{rowley2009jfm,tu2014jcd,kutz2016book}, it is based on linear measurements of the system and is unable to identify nonlinear changes of coordinates necessary to approximate the Koopman operator for strongly nonlinear systems.
The extended dynamic mode decomposition (eDMD)~\cite{williams2015jnls,williams2015jcd,klus2015numerical} was introduced by Williams et al. to address this issue, so that the best-fit linear DMD regression is performed on an augmented vector containing nonlinear measurements of the state.
The eDMD approach was recently shown to be equivalent to the earlier variational approach of conformation dynamics (VAC)~\cite{noe2013variational,nuske2014jctc,nuske2016variational,wu2017variational,wu2020variational} from the molecular dynamics literature, as explored in the excellent review by Klus et al.~\cite{klus2017data}.

In eDMD, an augmented state ${\bz\in\mathbb{R}^p}$ is constructed from nonlinear measurements of the state $\bx$ given by the functions $g_k$:
\begin{equation}\label{Eq:eDMD:Basis}
\bz = \bg(\bx) = \begin{bmatrix}
g_1(\bx)\\
g_2(\bx)\\
\vdots\\
g_p(\bx)
\end{bmatrix}.
\end{equation}
The vector $\bz$ may contain the original state $\bx$ as well as nonlinear measurements, so often $p\gg n$.
Next, two data matrices are constructed, as in DMD:
\begin{subequations}
\begin{align}
\bZ &= \begin{bmatrix} \vline & \vline & & \vline \\
\bz_1 & \bz_2 & \cdots & \bz_m \\
 \vline & \vline & & \vline
 \end{bmatrix}, \quad\quad\quad\quad
 \bZ' = \begin{bmatrix} \vline & \vline & & \vline \\
\bz_2 & \bz_3 & \cdots & \bz_{m+1} \\
 \vline & \vline & & \vline
 \end{bmatrix}.
\end{align}
\end{subequations}
Here $\bz_k=\bg(\bx_k)=\bg(\bx(k\Delta t))$, where we assume data is sampled at regular intervals in time, for simplicity.
As in DMD, a best-fit linear {matrix} operator $\bA_{\bZ}$ is constructed that maps $\bZ$ into $\bZ'$:
\begin{align}\label{Eq:EDMD:Regression}
\bA_{\bZ} = \argmin_{\bA_{\bZ}} \| \bZ' - \bA_{\bZ}\bZ\|_F =  \bZ' \bZ^{\dagger}.
\end{align}
Because the augmented vector $\bz$ may be significantly larger than the state $\bx$, it is typically necessary to employ kernel methods to compute this regression~\cite{williams2015jcd}.
{The use of kernel methods to approximate the Koopman operator with DMD has become an important research topic in recent years~\cite{williams2015jcd,fujii2019dynamic,das2020a,klus2020eigendecompositions,klus2020kernel,manohar2020kernel,baddoo2021kernel}.}
In principle, the functions $\{g_k\}_{k=1}^p$ form an enriched basis in which to approximate the Koopman operator.
In the limit of infinite data, the extended DMD {matrix} converges to the Koopman operator projected onto the subspace spanned by these functions~\cite{korda2017arxiv}.
However, if these functions do not span a Koopman invariant subspace, then the projected operator will have spurious eigenvalues and eigenvectors that differ from the true Koopman operator~\cite{brunton2016plosone}.

For example, consider the trivial example of a diagonalized linear system with eigenvalues $\lambda \in \{1,2,5\}$, eigenvectors $\bxi_j$ in the coordinate directions $x_j$, and a naive measurement that mixes the first two eigenvectors:
\begin{align}
    \frac{d}{dt}\begin{bmatrix} x_1 \\ x_2 \\ x_3\end{bmatrix} = \begin{bmatrix} 1 & 0 & 0 \\ 0 & 2 & 0 \\ 0 & 0 & 5\end{bmatrix} \begin{bmatrix} x_1 \\ x_2 \\ x_3\end{bmatrix} \quad\text{with}\quad
    y = \begin{bmatrix} 1 & 1 & 0\end{bmatrix} \begin{bmatrix} x_1 \\ x_2 \\ x_3\end{bmatrix}.
\end{align}
In this case, DMD will predict a spurious eigenvalue of $3$, which is the sum of the first two eigenvalues $\lambda=1$ and $\lambda=2$, since the measurement is a sum of the first two eigenvectors.
Therefore, it is essential to use validation and cross-validation techniques to ensure that eDMD models are not overfit, as discussed below.

Eigenfunctions of the Koopman operator form a Koopman invariant subspace and provide an ideal basis, or coordinate system, in which to represent the dynamics.
However, Koopman eigenfunctions may not admit a finite representation in any standard basis.
Thus, these eigenfunctions may only be approximately represented in a given finite basis.
It is possible to approximate an eigenfunction $\varphi(\bx)$ as an expansion in terms of the set of candidate functions $\{g_k(\bx)\}_{k=1}^p$ from \eqref{Eq:eDMD:Basis} as:
\begin{equation}\label{Eq:eDMD:Eigenfunction}
\varphi(\bx) \approx \sum_{k=1}^p \xi_k g_k(\bx) = \bxi^T\bg(\bx).
\end{equation}
In discrete-time, a Koopman eigenfunction $\varphi(\bx)$ evaluated on a trajectory of snapshots $\{\bx_1,\cdots,\bx_m\}$ will satisfy:
\begin{align}\label{Eq:eDMDData}
\lambda\begin{bmatrix}
\varphi(\bx_1)&
\varphi(\bx_2)&
\cdots &
\varphi(\bx_{m})
\end{bmatrix} =
\begin{bmatrix}
\varphi(\bx_2)&
\varphi(\bx_3)&
\cdots &
\varphi(\bx_{m+1})
\end{bmatrix}.
\end{align}
Expanding the eigenfunction $\varphi$ using \eqref{Eq:eDMD:Eigenfunction} this equality becomes
\begin{align}
\lambda\begin{bmatrix}\bxi^T\bz_1 & \bxi^T\bz_2 & \cdots & \bxi^T\bz_m\end{bmatrix} = \begin{bmatrix}\bxi^T\bz_2 & \bxi^T\bz_3 & \cdots & \bxi^T\bz_{m+1}\end{bmatrix}
\end{align}
which is possible to write as a matrix system of equations in terms of the data matrices $\bZ$ and $\bZ'$:
\begin{equation}\label{Eq:DiscreteEDMD}
\lambda\bxi^T\bZ-\bxi^T\bZ' = \mathbf{0}.
\end{equation}
If we seek a \emph{least-squares} fit to~\eqref{Eq:DiscreteEDMD}, this reduces to extended DMD~\cite{williams2015jcd,williams2015jnls}:
\begin{equation}\label{Eqn:KoopmanEfunLSR}
\lambda\bxi^T =  \bxi^T\bZ'\bZ^{\dagger}.
\end{equation}
The eigenfunctions $\varphi(\bx)$ are formed from left eDMD eigenvectors $\bxi^T$ of $\bZ'\bZ^\dagger$ as ${\varphi\approx\bxi^T\bz}$ in the basis $\{g_k(\bx)\}_{k=1}^p$.
The right eigenvectors are the eDMD \emph{modes}{, similar to DMD modes}.

It is essential to confirm that predicted eigenfunctions actually behave linearly on trajectories, by comparing them with the predicted dynamics ${\varphi(\bx_{k+1}) = \lambda \varphi(\bx_k)}$, as the regression above will result in spurious eigenvalues and eigenvectors unless the basis elements $g_k$ span a Koopman invariant subspace~\cite{brunton2016plosone}.
It is common to include the original state $\bx$ in the augmented eDMD vector $\bz$.
However, it was shown that including the state $\bx$ in eDMD results in closure issues for systems with multiple fixed points, periodic orbits, or other attractors, because these systems cannot be topologically conjugate to a finite-dimensional linear eDMD system with a single fixed point~\cite{brunton2016plosone}.
For example, the Duffing oscillator in \cref{Fig:KoopmanOverview} has three fixed points, so no finite linear system can accurately evolve the state $\bx$ near all three of these fixed points.
However, eigenfunctions like the Hamiltonian, may be accurately expanded in a basis.
Thus, it is critical to sort out the accurate and spurious eigenfunctions in eDMD; often eigenfunctions corresponding to lightly damped eigenvalues can be better approximated, as they have a significant signature in the data.

One approach to prevent overfitting is to promote sparsity, as in the sparse identification of nonlinear dynamics (SINDy)~\cite{brunton2016pnas}.
This principle of parsimony may also be used to identify Koopman eigenfunctions by selecting only the few most important terms in the basis $\{g_k(\bx)\}_{k=1}^p$ needed to approximate $\varphi$~\cite{kaiser2017arxiv}.

As with standard DMD, the data in $\bZ$ does not need to be generated from a single trajectory, but can instead be sampled more efficiently, such as with latin hypercube sampling or sampling from a distribution over the phase space.
In this case, the data in $\bZ'$ must be obtained by advancing the data in $\bZ$ one time step forward.
Moreover, reproducing-kernel Hilbert spaces (RKHS) can be employed to { ensure that the regularity of $\varphi(\bx)$ allows one to extend the computations on samples \emph{locally}, to patches of state space}.

For continuous-time dynamics, the eigenfunction dynamics
\begin{align}
 \frac{d}{dt}\varphi(\bx) = \lambda \varphi(\bx)
\end{align}
may be written in terms of the approximation $\varphi(\bx)\approx\bxi^T\bg(\bx)$:
\begin{align}
    \frac{d}{dt}\bxi^T\bg(\bx) = \lambda\bxi^T\bg(\bx).
\end{align}
Applying the chain rule results in
\begin{align}
    \bxi^T \bGamma(\bx,\dot{\bx}) = \lambda \bxi^T \bg(\bx)
\end{align}
where $\bGamma$ is given by:
\begin{align}
    \bGamma(\bx,\dot{\bx}) = \begin{bmatrix}\nabla g_1(\bx)\cdot\dot{\bx} \\ \nabla g_2(\bx)\cdot\dot{\bx} \\ \vdots \\ \nabla g_p(\bx)\cdot\dot{\bx} \end{bmatrix}.
\end{align}
Each term is a directional derivative, representing the possible terms in $\nabla\varphi(\bx)\cdot\mathbf{f}(\bx)$ from \eqref{eq:lie-eigenfunction}.
 It is then possible to construct a data matrix $\bGamma$ evaluated on the trajectories from $\bX$ and $\dot{\bX}$:
\begin{align*}
\bGamma =
\begin{bmatrix}
\nabla g_1(\bx_1)\cdot\dot{\bx}_1 & \nabla g_1(\bx_2)\cdot\dot{\bx}_2 & \cdots & \nabla g_1(\bx_m)\cdot\dot{\bx}_m\\
\nabla g_2(\bx_1)\cdot\dot{\bx}_1 & \nabla g_2(\bx_2)\cdot\dot{\bx}_2 & \cdots & \nabla g_2(\bx_m)\cdot\dot{\bx}_m \\
\vdots & \vdots & \ddots & \vdots \\
\nabla g_p(\bx_1)\cdot\dot{\bx}_1 & \nabla g_p(\bx_2)\cdot\dot{\bx}_2 & \cdots & \nabla g_p(\bx_m)\cdot\dot{\bx}_m\end{bmatrix}.
\end{align*}
The Koopman eigenfunction equation then becomes:
\begin{equation}
\lambda\bxi^T\bZ - \bxi^T\bGamma = \mathbf{0}.\label{Eq:SparseKoopman}
\end{equation}
Note that here we use notation where $\bGamma$ is the transpose of the notation in Kaiser et al.~\cite{kaiser2017arxiv} to be consistent with the eDMD notation above.

\subsection{Time delay coordinates}\label{Sec:HAVOK}
The DMD and eDMD algorithms are based on the availability of full-state measurements, which are typically quite high-dimensional.
However, it is often the case that only partial observations of the system are available, so that there are hidden, or \emph{latent}, variables.
In this case, it is possible to use time-delayed measurements of the system to build an augmented state vector, resulting in an intrinsic coordinate system that forms a Koopman-invariant subspace~\cite{brunton2017natcomm}.
The use of time-delay coordinates as a Koopman coordinate system relies on the conditions of the Takens embedding theorem~\cite{takens1981lnm} being satisfied, so that the delay-embedded attractor is diffeomorphic to the attractor in the original full-state coordinates.

\begin{figure}
\begin{center}
\begin{overpic}[width=\textwidth]{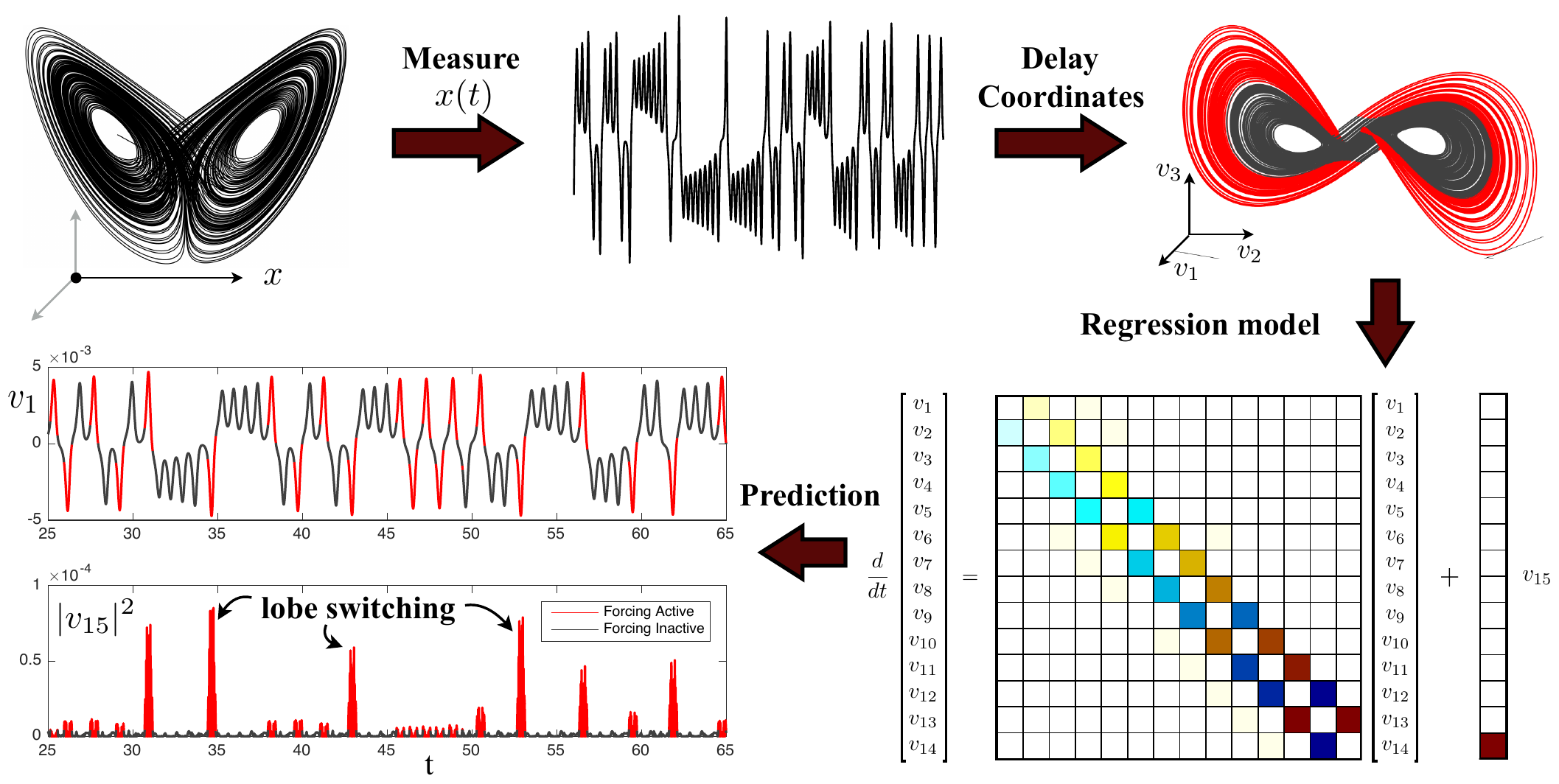}
\end{overpic}
\vspace{-.1in}
\caption{ Decomposition of chaos into a linear system with forcing.  A time series $x(t)$ is stacked into a Hankel matrix $\mathbf{H}$.  The SVD of $\mathbf{H}$ yields a hierarchy of  \emph{eigen} time series that produce a delay-embedded attractor.  A best-fit linear regression model is obtained on the delay coordinates $\mathbf{v}$; the linear fit for the first $r-1$ variables is excellent, but the last coordinate $v_r$ is not well-modeled as linear.  Instead, $v_r$ is an input that forces the first ${r-1}$ variables.  Rare forcing events correspond to lobe switching in the chaotic dynamics.  This architecture is called the Hankel alternative view of Koopman (HAVOK) analysis, from~\cite{brunton2017natcomm}.  \textit{Figure modified from Brunton et al.~\cite{brunton2017natcomm}.}}\label{Fig:HAVOK}
\end{center}
\end{figure}

The time-delay measurement scheme is illustrated schematically in \cref{Fig:HAVOK} on the Lorenz'63 system.
In this example, we have a single scalar measurement signal $x(t)$ from the original three-state Lorenz system.
It is possible to construct a Hankel matrix $\mathbf{H}$ from a time-series of this scalar measurement:
\begin{eqnarray}
\mathbf{H} = \begin{bmatrix}x(t_1) & x(t_2) &  \cdots & x(t_{p}) \\
x(t_2) & x(t_3) &\cdots & x(t_{p+1}) \\
\vdots & \vdots  &\ddots & \vdots \\
x(t_q) & x(t_{q+1})& \cdots & x(t_{m})
 \end{bmatrix}.\label{Eq:Hankel}
\end{eqnarray}
Each column of $\mathbf{H}$ may be obtained by advancing the previous column  forward in time by one time step.
Thus, we may re-write \eqref{Eq:Hankel} in terms of the Koopman operator $\mathcal{K}$:
\begin{eqnarray}
\mathbf{H} = \begin{bmatrix}x(t_1) & \mathcal{K} x(t_1) &  \cdots & \mathcal{K}^{p-1}x(t_{1}) \\
\mathcal{K} x(t_1) & \mathcal{K}^2 x(t_1)  & \cdots & \mathcal{K}^p x(t_{1}) \\
\vdots & \vdots  &\ddots & \vdots \\
\mathcal{K}^{q-1} x(t_1) & \mathcal{K}^q x(t_{1}) & \cdots & \mathcal{K}^{m-1} x(t_{1})
 \end{bmatrix}.\label{Eq:HAVOK}
\end{eqnarray}
For a sufficient volume of data, the system will converge to an attractor, so that the columns of $\mathbf{H}$ become approximately linearly dependent.
In this case, it is possible to obtain a Koopman-invariant subspace by computing the SVD of $\mathbf{H}$:
\begin{align}
    \mathbf{H} = \mathbf{U}\boldsymbol{\Sigma}\mathbf{V}^*.\label{Eq:Hankel:SVD}
\end{align}
The columns of $\mathbf{U}$ and $\mathbf{V}$ from the SVD are arranged hierarchically by their ability to model the columns and rows of $\mathbf{H}$, respectively.
The low-rank approximation in \eqref{Eq:Hankel:SVD} provides a \emph{data-driven} measurement system that is approximately invariant to the Koopman operator  for states on the attractor.
By definition, the dynamics map the attractor into itself, making it \emph{invariant} to the flow.
Often, $\mathbf{H}$ will admit a low-rank approximation by the first $r$ columns of $\mathbf{U}$ and $\mathbf{V}$, so that these columns approximate a Koopman-invariant subspace.
Thus, the columns of \eqref{Eq:Hankel} are well-approximated by the first $r$ columns of $\mathbf{U}$.
The first $r$ columns of $\mathbf{V}$ provide a time series of the magnitude of each of the columns of $\mathbf{U\Sigma}$ in the data.
By plotting the first three columns of $\mathbf{V}$, we obtain an embedded attractor for the Lorenz system, as in \cref{Fig:HAVOK}.

Because the columns of $\mathbf{H}$, and hence $\mathbf{U}$, form a Koopman-invariant subspace, it is possible to perform DMD on the two matrices formed from the first $p-1$ and last $p-1$ columns of $\mathbf{H}$.
In practice, we recommend performing DMD on a similar set of two matrices formed from the first $p-1$ and last $p-1$ columns of $\mathbf{V}$, since the columns of $\mathbf{V}$ are the coordinates of the system in the $\mathbf{U\Sigma}$ frame.
This results in a linear regression model on the variables in $\mathbf{V}$
\begin{eqnarray}
\frac{d}{dt}\mathbf{v}(t) = \mathbf{A}\mathbf{v}(t).\label{Eq:HAVOK1}
\end{eqnarray}
Champion et al.~\cite{champion2019discovery} showed that this linear system captures the dynamics of weakly nonlinear systems.
For chaotic systems, however, even with an approximately Koopman-invariant measurement system, there remain challenges to identifying a closed linear model.
A linear model, however detailed, cannot capture multiple fixed points or the unpredictable behavior characteristic of chaos with a positive Lyapunov exponent~\cite{brunton2016plosone}.
Instead of constructing a closed linear model for the first $r$ variables in $\mathbf{V}$, we build a linear model on the first $r-1$ variables and impose the last variable, $v_r$, as a forcing term~\cite{brunton2017natcomm}:
\begin{eqnarray}
\frac{d}{dt}\mathbf{v}(t) = \mathbf{A}\mathbf{v}(t) + \mathbf{B}{v}_r(t),\label{Eq:ChaosModel}
\end{eqnarray}
where $\mathbf{v}=\begin{bmatrix}v_1 & v_2 & \cdots & v_{r-1}\end{bmatrix}^T$ is a vector of the first ${r-1}$ eigen-time-delay coordinates.
In all of the chaotic examples explored~\cite{brunton2017natcomm}, the linear model on the first $r-1$ terms is accurate, while no linear model represents $v_r$.
Instead, $v_r$ is an input forcing to the linear dynamics in \eqref{Eq:ChaosModel}, which approximates the nonlinear dynamics.
The statistics of $v_r(t)$ are non-Gaussian, with long tails correspond to rare-event forcing that drives lobe switching in the Lorenz system; this is related to rare-event forcing distributions observed and modeled by others~\cite{majda2012nonlinearity,sapsis2013pnas,majda2014pnas}.

The Hankel matrix has been used for decades in system identification, for example in the eigensystem realization algorithm (ERA)~\cite{juang1985jgcd} and the singular spectrum analysis (SSA)~\cite{broomhead1989prsla}.
{These early algorithms were developed specifically for linear systems, and although they were often applied to weakly nonlinear systems, it was unclear how to interpret the resulting models and decompositions.
Modern Koopman operator theory has provided a valuable new perspective for how to interpret the results of these classical Hankel-based approaches when applied to nonlinear systems. }
Computing DMD on a Hankel matrix was first introduced by Tu et al.~\cite{tu2014jcd} and was used by B. Brunton et al.~\cite{brunton2016b} in the field of neuroscience.
The connection between the Hankel matrix and the Koopman operator, along with the linear regression models in \eqref{Eq:ChaosModel}, was established by Brunton et al.~\cite{brunton2017natcomm} in the Hankel alternative view of Koopman (HAVOK) framework.
Several subsequent works have provided additional theoretical foundations for this approach~\cite{arbabi2017,das2017arxiv,Kamb2020siads,champion2019discovery,hirsh2019}.
Hirsh et al.~\cite{hirsh2019} established connections between HAVOK and the Frenet-Serret frame from differential geometry, motivating a more accurate computational modeling approach.
The HAVOK approach is also often referred to as delay-DMD~\cite{tu2014jcd} or Hankel-DMD~\cite{arbabi2017}.
A connection between delay embeddings and the Koopman operator was established as early as 2004 by Mezi\'{c} and Banaszuk~\cite{mezic2004physicad}, where a stochastic Koopman operator is defined and a statistical Takens theorem is proven.  Other work has investigated the splitting of dynamics into deterministic linear, and chaotic stochastic dynamics~\cite{mezic2005nd}.
The use of delay coordinates may be especially important for systems with long term memory effects and where the Koopman approach has recently been shown to provide a successful analysis tool~\cite{stanton2016pre}.

\subsection{Diffusion maps for Koopman embeddings}\label{Sec:Diffusion}

Diffusion maps are a recently developed nonlinear dimensionality-reduction technique for embedding high-dimensional data on nonlinear manifolds~\cite{coifman2005pnas,coifman2006acha,coifman2008mmas,nadler2006acha}.
Diffusion maps leverage the underlying geometry, and in particular its local similarity structure, to create an organization of data that is heavily influenced by local structures.
The distance between data is measured by a kernel, for example the Gaussian kernel, which takes the form of the fundamental Green's function solution of the heat equation and is proportional to the connectivity between two data points.
The diffusion kernel is given by
\begin{equation}
    k(\bx_j,\bx_k) = \exp \left( - \frac{\|\bx_j-\bx_k\|}{\alpha} \right),
\end{equation}
where $\bx_j$ and $\bx_k$ are two data points and $\alpha$ determines the range of influence of the kernel.
Thus, points that are not sufficiently close have an approximately zero connectivity since the kernel decays as a Gaussian between data points.

The diffusion kernel plays the role of a normalized likelihood function.  It also has important properties when performing spectral analysis of the distance matrix constructed from the embedding.
These include symmetry, $k(\bx,\by)=k(\by,\bx)$, and positivity preserving, $k(\bx,\by)\geq 0$.
The basic diffusion mapping algorithm computes a kernel matrix ${\bf K}$ whose elements are given by $K_{j,k}=k(\bx_j,\bx_k)$.
After normalization of the rows of the kernel matrix, the eigenvalues and eigenvectors of ${\bf K}$ are computed and the data is projected onto the dominant $r$-modes.
This $r$-dimensional subspace is the low-dimensional embedding of the diffusion map.
Coifman and Lafon~\cite{coifman2006acha} demonstrated that this mapping gives a low dimensional parametrization of the geometry and density of the data. In the field of data analysis, this construction is known as the {\em normalized graph Laplacian}.

Diffusion maps thus provide a dimensionality reduction method that exploits the geometry and density of the data. The diffusion map can be directly used to construct a Koopman model by using a DMD regression on the time evolution in the diffusion coordinates.
The methodology can also be used for forecasting~\cite{manohar2020kernel}, for example leveraging time-delay embeddings to provide a nonparametric forecasting method for data generated by ergodic dynamical systems~\cite{giannakis2019}.
Such a representation is based upon the Koopman and Perron-Frobenius groups of unitary operators in a smooth orthonormal basis which is acquired from time-ordered data through the diffusion maps algorithm. Giannakis~\cite{giannakis2019} establishes in such a representation  a correspondence between Koopman operators and Laplace-Beltrami operators constructed from data in Takens delay-coordinate space, using this correspondence to provide an interpretation of diffusion-mapped delay coordinates for ergodic systems.

\subsection{Neural networks for Koopman embeddings}\label{sec:NNKoop}
Despite the promise of Koopman embeddings, obtaining tractable representations has remained a central challenge.
Even for relatively simple dynamical systems, the eigenfunctions of the Koopman operator may be arbitrarily complex and will only be approximately represented in a finite basis.
Deep learning is well-suited for representing such arbitrarily complex functions, and has recently shown tremendous promise for discovering and representing Koopman embeddings and Koopman forecasts~\cite{wehmeyer2017arxiv,mardt2017arxiv,Takeishi2017neurips,yeung2017arxiv,otto2017arxiv,li2017chaos,lusch2017arxiv,mardt2020deep,azencot2020forecasting,lange2020fourier,eivazi2021recurrent}.
{In addition to leveraging neural networks to learn Koopman embeddings, Koopman theory is also being applied to understand neural networks~\cite{manojlovic2020applications,dogra2020optimizing}, and algorithms more generally~\cite{dietrich2020koopman}.}

\begin{figure}[t]
\centering
\begin{overpic}[width=\linewidth]{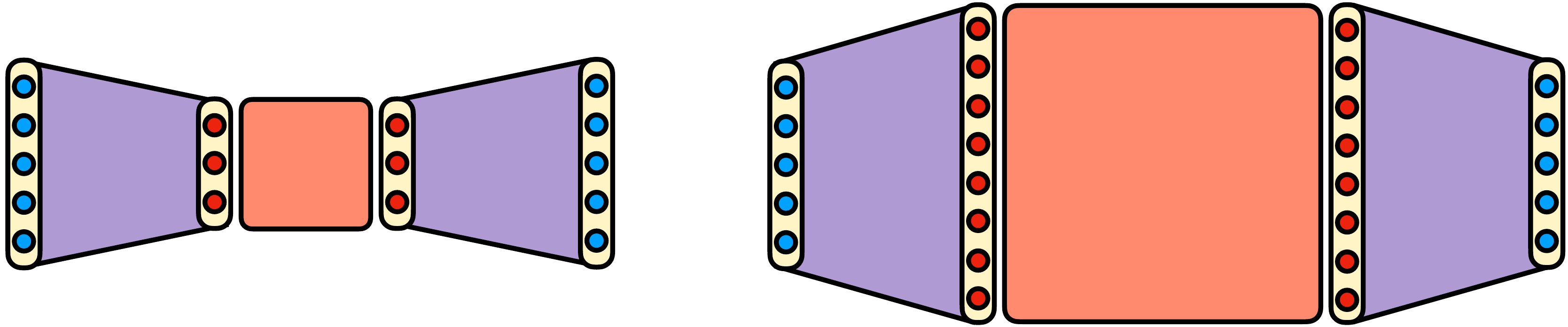}
\put(0,19){(a)}
\put(48.5,19){(b)}
\put(0.5,1.5){$\bx_k$}
\put(12.5,4){$\bz_k$}
\put(23.5,4){$\bz_{k+1}$}
\put(36.5,1.5){$\bx_{k+1}$}
\put(49,1.5){$\bx_k$}
\put(61.3,-2){$\bz_k$}
\put(84,-2){$\bz_{k+1}$}
\put(97.,1.5){$\bx_{k+1}$}
\put(6.5,9.5){$\boldsymbol{\varphi}$}
\put(18.2,9.5){$\mathbf{K}$}
\put(30.7,9.5){$\boldsymbol{\psi}$}
\put(55,9.5){$\boldsymbol{\varphi}$}
\put(73,9.5){$\mathbf{K}$}
\put(91.5,9.5){$\boldsymbol{\psi}$}
\end{overpic}
\vspace{-.1in}
\caption{Competing neural network architectures to approximate the Koopman operator.  (a) Key Koopman eigenfunctions are extracted with a deep auto-encoder network. (b) Alternatively, the system is lifted to a higher dimension where a linear model is identified.  In these architectures $\boldsymbol{\varphi}$ is the encoder and $\boldsymbol{\psi}$ is the decoder.
}
\label{fig:schema}
\vspace{-.15in}
\end{figure}
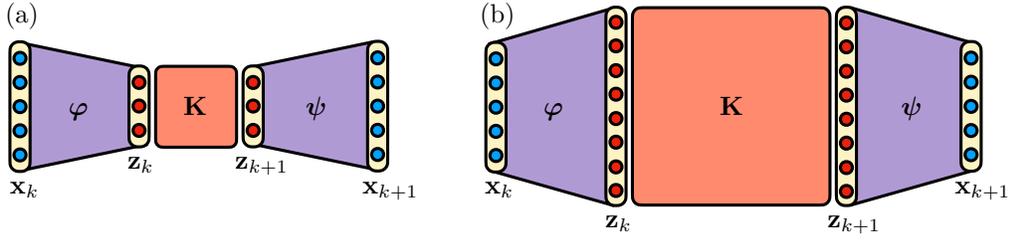

There are two leading deep neural network architectures that have been proposed for Koopman embeddings, shown in \cref{fig:schema}.
In the first architecture, the deep auto-encoder architecture extracts a few key latent variables $\bz=\boldsymbol{\varphi}(\bx)$ to parameterize the dynamics.
In the second architecture, the high-dimensional input data is lifted to an even higher dimension, where the evolution is approximately linear.
In either Koopman neural network, an additional constraint is enforced so that the dynamics must be linear on these latent variables, given by the {matrix} $\mathbf{K}$.
The constraint of linear dynamics is enforced by the loss function $\|\boldsymbol{\varphi}(\bx_{k+1})-\bK \boldsymbol{\varphi}(\bx_k)\|$, where $\bK$ is a matrix.
In general, linearity is enforced over multiple time steps, so that additional terms $\|\boldsymbol{\varphi}(\bx_{k+p})-\bK^p \boldsymbol{\varphi}(\bx_k)\|$ are added to the loss function.

Autoencoder networks have the advantage of a low-dimensional latent space, which may promote interpretable solutions.
Autoencoders are already widely used to model complex systems, for example in fluid mechanics~\cite{Brunton2020arfm}, and they may be viewed as nonlinear extensions of the singular value decomposition, which is central to the DMD algorithm.
In this way, a deep Koopman network based on an autoencoder may be viewed as a nonlinear generalization of DMD.
Similarly, if the matrix $\mathbf{K}$ is diagonalized, then the embedding functions $\boldsymbol{\varphi}$ correspond to Koopman eigenfunctions~\cite{lusch2017arxiv,gin2019deep}.
Variational autoencoders are also used for stochastic dynamical systems, such as molecular dynamics, where the map back to physical configuration space from the latent variables is probabilistic~\cite{wehmeyer2017arxiv,mardt2017arxiv}.
In contrast, the second paradigm, where measurements are lifted to a higher-dimensional space, is related to many results in machine learning, such as Cover's theorem~\cite{lange2020fourier}, where nonlinear problems tend to become more linear in higher-dimensional embeddings.

\begin{figure}[t]
\centering
\begin{overpic}[width=\linewidth]{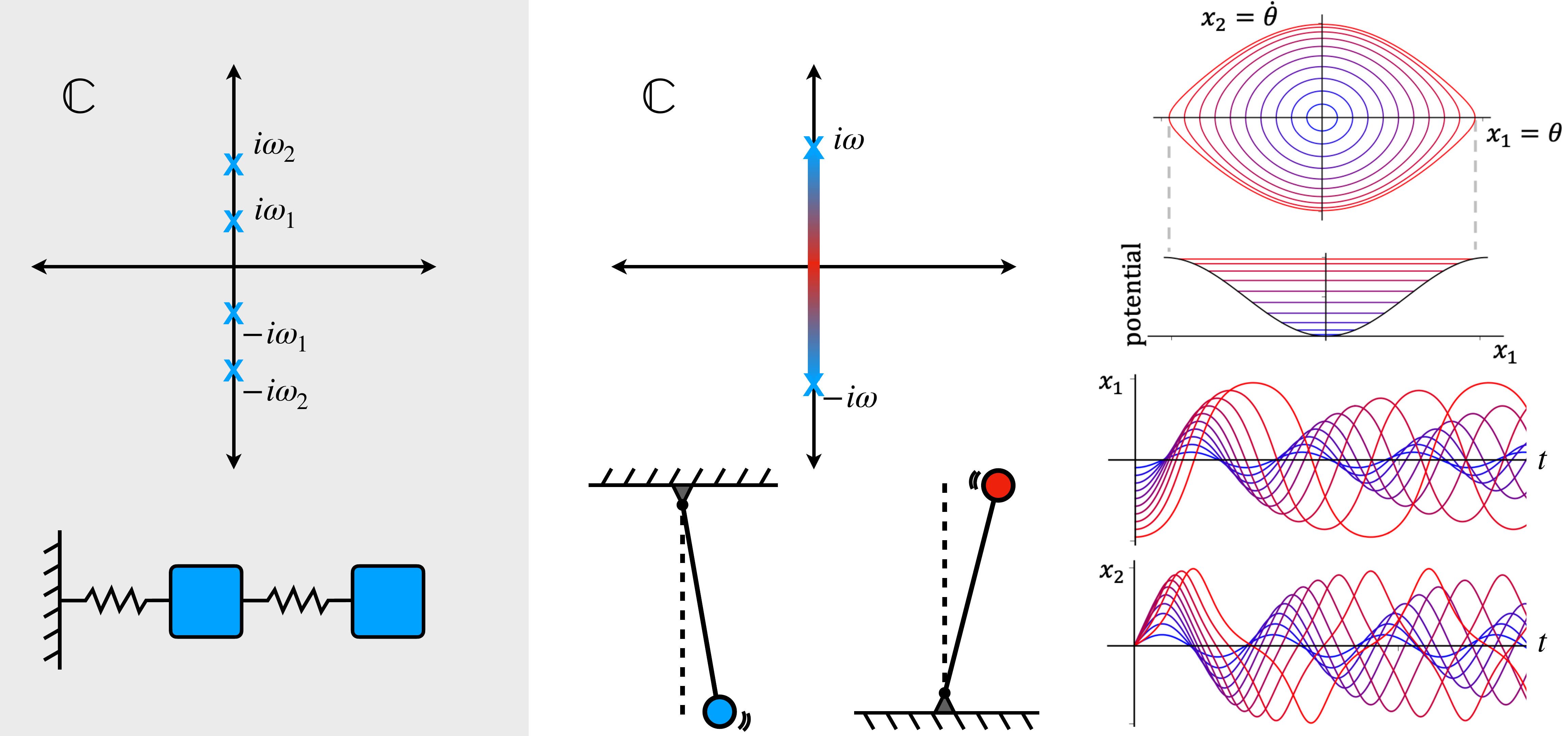}
\put(.5,44){(a) Discrete spectrum}
\put(37.25,44){(b) Continuous spectrum}
\end{overpic}
\vspace{-.2in}
\caption{Comparison of discrete vs. continuous spectrum dynamics. \textit{Right panel reproduced from Lusch et al.~\cite{lusch2017arxiv}.}}
\label{fig:ContinuousSpectrum}
\vspace{-.2in}
\end{figure}

For simple systems with a discrete eigenvalue spectrum, a compact representation may be obtained in terms of a few autoencoder variables.
However, dynamics with continuous eigenvalue spectra defy standard low-dimensional  Koopman representations, including the autoencoder network above.
Continuous spectrum dynamics are ubiquitous, ranging from the simple pendulum to nonlinear optics and broadband turbulence.
For example, the classical pendulum, given by
\begin{align}
\ddot{x} = -\sin(\omega x)
\end{align}
exhibits a continuous range of frequencies, from $\omega$ to $0$, as the amplitude of the pendulum oscillation is increased, as illustrated in \cref{fig:ContinuousSpectrum}.
Thus, the continuous spectrum confounds a simple description in terms of a few Koopman eigenfunctions.
Indeed, away from the linear regime, an infinite Fourier sum is required to approximate the continuous shift in frequency, which may explain why the high-dimensional lifting approach has been widely used in Koopman neural networks.

In a recent work by Lusch et al.~\cite{lusch2017arxiv}, an auxiliary  network is used to parameterize the continuously varying eigenvalue, enabling a network structure that is both parsimonious and interpretable.
In contrast to other network structures, which require a large autoencoder layer to encode the continuous frequency shift with an asymptotic expansion in terms of harmonics of the natural frequency, the parameterized network is able to identify a single complex conjugate pair of eigenfunctions with a varying imaginary eigenvalue pair.
If this explicit frequency dependence is unaccounted for, then a high-dimensional network is necessary to account for the shifting frequency and eigenvalues.
Recently, this framework has been generalized to identify linearizing coordinate transformations for PDE systems~\cite{gin2019deep}, such as the Cole--Hopf transform of the nonlinear Burgers' equation into the linear heat equation.
Related work has been developed to identify the analogues of Green's functions for nonlinear systems~\cite{gin2020deepgreen}.

{
\subsection{Convergence of data-driven Koopman approximations}
\label{sec:approximations}

The dynamic mode decomposition has been seen, almost from its beginning, as an algorithm for computational Koopman analysis, giving it strong connections to the theory of dynamical systems~\cite{rowley2009jfm,mezic2005nd}.
The computationally effective DMD explicitly avoids approximating the Koopman operator, but instead enables the approximation of its eigenvalues and eigenmodes.
In this section, we summarize research efforts to establish the convergence and effectiveness of DMD and other techniques in approximating the Koopman operator and its spectral decomposition~\eqref{Eq:KoopmanModeDecomposition}.
At the time of writing, the study of the quality of Koopman approximations remains a vigorous area of research, with many contributions available only in preprint format.

The original snapshot DMD was interpreted as an Arnoldi iteration of the Koopman operator on the set of observables (see ~\Cref{Sec:DMD:Krylov}).
In this formulation, the problem of accurate computation of eigenvalues of the Koopman operator is replaced by the problem of computing eigenvalues of the companion matrix~\eqref{eq:companion-matrix}.
Computations involving the companion matrix and its diagonalizing transform, the Vandermonde matrix, are ill-conditioned; an excellent summary of the delicateness of working with companion matrices from the perspective of numerical linear algebra is given in a sequence of papers by Drma\v{c} et al.~\cite{drmac2018,drmac2019,drmac2020}.
Despite numerical fragility, the Vandermonde matrix interpretation leads to valuable insights into connections between DMD, the discrete Fourier transform~\cite{chen2012jns}, and Prony analysis~\cite{susuki2015,zhen2021,arbabi2017}.
The optimized DMD~\cite{chen2012jns}, and more recently the stabilization of the Vandermonde interpretation of DMD~\cite{drmac2020,drmac2019}, have made the Arnoldi interpretation a practically viable approach to computing DMD.

The convergence and robustness of DMD-based algorithms has benefited from the interpretation of  exact DMD~\cite{schmid2010jfm,tu2014jcd} as computing spectral data of a low-rank approximation of the Koopman operator via regression on a subspace of observables.
The extended DMD (eDMD)~\cite{williams2015jnls} (see ~\cref{Sec:EDMD}) enables this argument by augmenting the state vector, amounting to the vector-valued observable \(\bg(\bx) = \bx\), with nonlinear transformations of the input data, typically chosen from a basis of functions spanning the space of observables.
This in turn allows one to interpret the DMD regression~\eqref{Eq:DMD:Definition} as a regression on a subspace whose dimension can now be tuned to include more, or less, of the full function space.
Initial convergence arguments by Williams et al.~\cite{williams2015jnls} and Klus et al.~\cite{klus2015numerical} showed that eDMD is equivalent to a Galerkin projection of the Koopman operator onto a subspace of the \(L^{2}\) space with respect to an ergodic measure.
This was further refined and improved upon by Korda and Mezi\'{c}~\cite{korda2017arxiv}, showing that the eDMD approximation converges in operator norm to the compression of the Koopman operator to a chosen subspace of observables for ergodic dynamics in the limit of long trajectories.
When the sequence of observable subspaces limits to the entire \(L^{2}\) space, the projections converge in the strong operator topology (SOT), i.e., pointwise with respect to any observable, to the Koopman operator.
While SOT is not enough to guarantee convergence of the spectrum, the eigenvalues of the Koopman operator are among accumulation points of eDMD eigenvalues~\cite{korda2017arxiv}.

Instead of targeting the Koopman operator directly, several approaches have sought to approximate its infinitesimal generator~\eqref{Eq:Koopman:InfinitesimalGenerator}, which satisfies the evolution equation~\eqref{eq:Lie-derivative}: \(\dot g_{t} = \gen{} g_{t}\).
Since many models of physics are stated as differential equations, the generator may have a direct link to the physics or may be more easily exploited for control~\cite{klus2020,sechi2021,mauroy2020tac}.
Additionally, in certain contexts it is possible to establish stronger convergence results when working with the generator instead of the Koopman operator itself~\cite{rosenfeld2021}.
Similar arguments apply to the approximation of the Perron--Frobenius operator and its generator~\cite{froyland2016,froyland2013}, discussed in \cref{sec:perron-frobenius}.
The Koopman generator can be approximated by exploiting any of its connections to the Koopman operator.
One approach is to first approximate the Koopman operator family, and then use a finite-difference approximation to compute the Lie derivative of the Koopman operator~\cite{giannakis2020a,sechi2021}.
Alternatively, the relationship between the Koopman operator and its generator
\begin{equation}
  \label{eq:lie-exponential}
  \koop{}^{t} = \exp[t \gen{}]
\end{equation}
can be inverted, and the generator computed by approximating the matrix logarithm of the Koopman operator~\cite{mauroy2020tac}.
Another approach is to apply the eDMD regression to~\eqref{eq:Lie-derivative} by computing time derivatives of the basis of observables~\cite{klus2020,klus2020kernel,rosenfeld2021,rosenfeld2021a}.
Finally, Giannakis and Das~\cite{das2017arxiv,giannakis2019,giannakis2020} approach the problem of approximation of the Koopman and its generator as a manifold-learning problem on a space-time manifold.
This gives rise to an approximation problem for the generators of the evolution on the manifold, which is successfully resolved for ergodic dynamical systems, such as those evolving on a chaotic attractor.

In certain contexts it is possible to study the quality of an approximation to the Koopman spectrum and eigenfunctions directly, without appeal to the convergence of the operator approximation.
Early approximations of spectral properties in this manner date back to Neumann, Wiener, and Wintner~\cite{wiener1941,neumann1932pnas} and appear computationally in the modern context in~\cite{mezic1999chaos,mezic2004physicad}.
For the Koopman operator defined on an \(L^{p}\) space of functions with respect to an ergodic measure, harmonic/Fourier averages~
\eqref{eq:harmonic-average} converge to a projection of an observable onto an eigenspace associated with a chosen eigenvalue \(\omega\).
Spatially, the convergence of such averages is assured with respect to the \(L^{p}\) norm, and the \(p=2\) case corresponds to the orthogonal projection onto the eigenspace associated with \(e^{i\omega}\); if \(e^{i\omega}\) is not an eigenvalue, the averages converge to the value zero almost everywhere.
Temporally, the rate of convergence can be arbitrarily slow; however, in case of regular discrete-time dynamics it scales as \(\mathcal{O}(N^{-1})\) with the number of iterates \(N\), and for mixing dynamics \(\mathcal{O}(N^{-\sfrac{1}{2}})\), with more detailed results available~\cite{assani2004,kachurovskii1996,mezic2020}.
While it is possible to numerically investigate such convergence starting with an arbitrary continuous observable, as was done in \cite{levnajic2010,levnajic2015,budisic2012chaos}, pointwise evaluations of functions in \(L^{p}\) spaces are ill-defined, and therefore insufficient for theoretical guarantees.
In the reproducing-kernel Hilbert space formulation of the Koopman operator~\cite{das2020a}, harmonic averaging can be used to compute the RKHS norm of the projection that can be used to detect the eigenvalue spectrum.
The harmonic averaging can, in principle, be established for eigenvalues off the unit circle, although such an algorithm is numerically unfavorable~\cite{mohr2014arxiv,mohr2016a,kvalheim2021,kvalheim2021a}.
Furthermore, successive projections using ergodic averages, termed \emph{generalized Laplace analysis} (GLA), can be used to incrementally account for all (quasi-)regular components of the evolution of observables~\cite{mohr2014arxiv}.
While the GLA process is ill-conditioned if the infinite ergodic averages are truncated, it is possible to provide an analogous process for the case of finite-length data~\cite{drmac2019}.

  As discussed in the \cref{sec:koopman-spectrum}, the spectral measure of the Koopman operator~\eqref{eq:spectral-decomposition} can have a continuous component, in addition to the atomic (eigenvalue) spectrum.
  Approximation of the non-eigenvalue spectrum was studied in the context of the DMD on time delayed observables~\cite{brunton2017natcomm,arbabi2017}.
  Whereas eDMD builds up the space of observables using an explicit choice of basis functions, the delay DMD (a.k.a. HAVOK or Hankel DMD from \cref{Sec:HAVOK}) instead uses delayed copies of state observations, which connects the DMD formalism with auto- and cross-correlations of observables~\eqref{eq:autocorrelation}.
   The distribution of eigenvalues of delay DMD~\cite{arbabi2017} cannot be used to approximate the continuous spectral density, as eigenvalues will form accumulations only around Koopman eigenvalues, but in the limit will distribute uniformly~\cite{korda2020}.

   Instead, the approximation of the continuous spectral measure can be posed as a truncated moment problem~\cite{korda2020}.
   The moments of the spectral measure for ergodic systems, computed using autocorrelations of trajectories, are used as inputs into the algorithm with the reconstruction of the atomic and absolutely-continuous components achieved using Christoffel--Darboux kernels.
   Additional results concerning convergence of the delay DMD can be found in~\cite{zhen2021}.
   Koopman modes were initially conceived as projections of a chosen set of observables onto individual eigenfunctions.

   While the non-atomic spectrum does not have eigenfunctions associated with it, it is possible to compute the mode-analogues for the continuous spectrum via spectral projections~\cite{korda2020}, ultimately leading to approximations of the operator itself.
   For measure-preserving dynamics, spectral projectors can be computed using the periodic approximation~\cite{govindarajan2019}, related to the Ulam approximation already discussed in \Cref{sec:perron-frobenius}.
   An alternative interpretation of modes for the non-atomic spectrum is offered by Giannakis~\cite{giannakis2020a}, which searches for approximately-coherent observables over a chosen window of observation, demonstrating their construction even for systems where no true coherent observables, that is eigenfunctions, exist.

Koopman operator theory was initially developed for observables chosen from a Lebesgue space \(L^{p}\).
Unfortunately, elements of \(L^{p}\) spaces are equivalence classes that are specified only up to a measure-zero set; consequently, computations based on finite pointwise samples of a function cannot be naively linked to strong convergence results based on \(L^{p}\) arguments.
Inferring the behavior of an observable from a finite pointwise sample requires working with a more restrictive set of observables, typically continuous functions with additional regularity properties. These properties are formally incorporated through the language of reproducing kernel Hilbert spaces (RKHS), which afford both a sufficient regularity and the inner-product structure.
An added benefit of working with RKHS spaces is that the observables do not need to take values from a linear space; for example, angle-valued observables taking values from \(\mathbb{S}^{1}\) or direction-valued observables taking values from projective spaces can be handled in this setting~\cite{klus2020eigendecompositions,das2020a}.
Kernel-based DMD, introduced by Williams et al. in~\cite{williams2015jcd} and, in parallel, by Kawahara~\cite{kawahara2016neurips} were the first computational efforts in this direction.
More recent results~\cite{klus2020eigendecompositions} demonstrate decompositions that apply both to Koopman and Perron--Frobenius operators.
Establishing the existence of the spectrum for Koopman operators on Hardy-type spaces~\cite{mohr2014arxiv} was the first theoretical step in pursuit of RKHS-based Koopman theory, with a more broad and rigorous theoretical basis presented by Mezi\'{c} in~\cite{mezic2019}.
As mentioned, computational approximation of harmonic averages is set on a firmer ground in RKHS~\cite{das2020a}.
More recent results pursue the convergence of approximations to Koopman operators and its generators in this setting~\cite{das2021,rosenfeld2021}.
}

\section{Koopman theory for control}\label{Sec:Control}
The Koopman operator framework is especially relevant for engineering applications in control~\cite{mauroy2020book,otto2021koopman}, for which it offers new opportunities for the control of nonlinear systems by circumventing theoretical and computational limitations due to nonlinearity.
Nonlinear control methods, such as feedback linearization and sliding mode control, overcome some of these limitations.
However, these approaches often do not generalize beyond a narrow class of systems, and deriving stability and robustness conditions, for instance, can become a demanding exercise.
Koopman-based methods provide a linear framework that can exploit mature theoretical and computational methods, with successes already demonstrated in a wide range of challenging applications, including fluid dynamics~\cite{arbabi2018cdc,peitz2018feedback}, robotics~\cite{abraham2017conf,abraham2019ieee,bruder2019proc,mamakoukas2019proc}, power grid~\cite{korda2018arxiv,netto2018tps}, traffic~\cite{ling2018ieee},
biology~\cite{hasnain2019arxivb}, logistics~\cite{hogg2019plosone}, and chemical processes~\cite{narasingam2019}.
Koopman analysis achieves this by representing the nonlinear dynamics in a globally linear framework, without linearization.
Thus Koopman analysis is able to generalize the Hartman-Grobman theorem to the entire basin of attraction of a stable or unstable equilibrium or periodic point~\cite{lan2013physd}.
Further, as the Koopman operator acts on observables, it is amenable to data-driven (model free) approaches which have been extensively developed in recent years~\cite{proctor2016siads,williams2016ifac,korda2016arxiv,proctor2017siads,surana2016cdc,kaiser2017arxiv,kakubr2018arxiv,peitz2017arxiv,abraham2019ieee}.
The resulting models have been shown to reveal insights into global stability properties~\cite{sootla2016acc,mauroy2016ieeetac}, observability/controllability~\cite{vaidya2007cdc,goswami2017cdc,yeung2018acc}, and sensor/actuator placement ~\cite{sinha2016jmaa,sharma2018arxiv_a} for the underlying nonlinear system.

Koopman theory is closely related to Carleman linearization~\cite{carleman1932am}, which also embeds finite-dimensional dynamics into infinite-dimensional linear systems. Carleman linearization has been used for decades to obtain truncated linear (and bilinear) state estimators~\cite{krener1974,brockett1976automatica,amini2019ifac} and to examine stability, observability, and controllability of the underlying nonlinear system~\cite{loparo1978tac,banks1992infinite,mozyrska2006,mozyrska2008}.
However, the applicability is restricted to polynomial (or analytical) systems.
In contrast, the Koopman operator framework does not rely on the analyticity of the vector field,
but applies to general nonlinear systems, including systems with discontinuities.
Extending Koopman operator theory for actuated systems was first noted in~\cite{mezic2004physicad}, which interpreted the stochastic forcing in random dynamical systems as actuation.
The first Koopman-based control schemes were published more than a decade later, powered by the algorithmic development of DMD~\cite{proctor2016siads}.
More recently, Koopman models have been increasingly used in combination with LQR~\cite{brunton2016plosone,mamakoukas2019proc,mamakoukas2020arxiv}, state-dependent LQR~\cite{kaiser2017arxiv}, and model predictive control (MPC)~\cite{korda2016arxiv,kaiser2017arxivb}. Other directions include optimal control for switching control problems~\cite{peitz2017arxiv,peitz2018feedback}, Lyapunov-based stabilization~\cite{huang2018cdc,huang2019arxiv}, eigenstructure assignment~\cite{hemati2017aiaa}, and, more recently, active learning~\cite{abraham2019ieee}.
MPC~\cite{garriga2010model,lee2011springer,mayne1997proc,morari1999model,qin1997proc,garcia1989model,camacho2013model,allgower2004nonlinear,eren2017jgcd} stands out as a main facilitator for the success of Koopman-based control, with applications including power grids~\cite{korda2018arxiv}, high-dimensional fluid flows~\cite{arbabi2018cdc,morton2019neurips}, and electrical drives~\cite{hanke2018arxiv}.

In this section, we review the mathematical formulation for a Koopman operator control framework, beginning with model-based control and moving to data-driven methods.
We will describe several approaches for identifying control-oriented models including dynamic mode decomposition with control (DMDc), extended DMD with control (eDMDc), extensions based on SINDy, and the use of delay coordinates.
Further, we compare these approaches on numerical examples and
discuss the use of the Koopman operator for analyzing important properties, such as stability, observability, and controllability, of the underlying system.

\subsection{Model-Based Control}\label{Sec:KoopmanControl:Math}
Beginning with model-based control theory, consider the non-affine nonlinear control system
\begin{subequations}\label{Eq:NonlinearControlSystem}
\begin{align}
\dot{\bx} &= {\bf f}(\bx,\bu),\quad \bx(0)=\bx_0\\
\by &= {\bf c}(\bx,\bu),
\end{align}
\end{subequations}
where $\bx\in \mathcal{X}\subseteq\mathbb{R}^n$ is the state vector, $\bu\in \mathcal{U}\subseteq\mathbb{R}^q$ is the vector of control inputs, and $\by\in \mathcal{Y}\subseteq\mathbb{R}^p$ is the output.
Unless noted otherwise, we assume full-state measurements $\by = \bx$.
Equation~\eqref{Eq:NonlinearControlSystem} represents a non-autonomous dynamical system, where the input $\bu$ may be interpreted as a perturbed parameter or control actuation.
In the context of control, we seek typically to determine a feedback control law ${\bf d}:\mathcal{Y}\rightarrow \mathcal{U}$,
\begin{equation}\label{Eq:FeedbackControlLaw}
	\bu = {\bf d}(\by),
\end{equation}
that maps measurements $\by$ to control inputs $\bu$ to modify the behavior of the closed-loop system.

\subsubsection{Model predictive control}\label{Sec:MPC}
MPC is one of the most successful model-based control schemes~\cite{garriga2010model,lee2011springer,mayne1997proc,morari1999model,qin1997proc,garcia1989model,camacho2013model,allgower2004nonlinear,eren2017jgcd}.
Over the last two decades, MPC has gained increasing popularity due to its success in a wide range of applications, its ability to incorporate customized cost functions and constraints, and extensions to nonlinear systems.
In particular, it has become the {\em de-facto} standard advanced control method in process industries~\cite{qin1997proc} and has gained considerable traction in the aerospace industry due to its versatility~\cite{eren2017jgcd}.
Many of the successes of Koopman-based controls leverage the MPC framework,
with its adaptive nature allowing one to compensate for modeling discrepancies and to account for disturbances.

\begin{figure}[tb]
    \centering
    \includegraphics[width=0.8\textwidth]{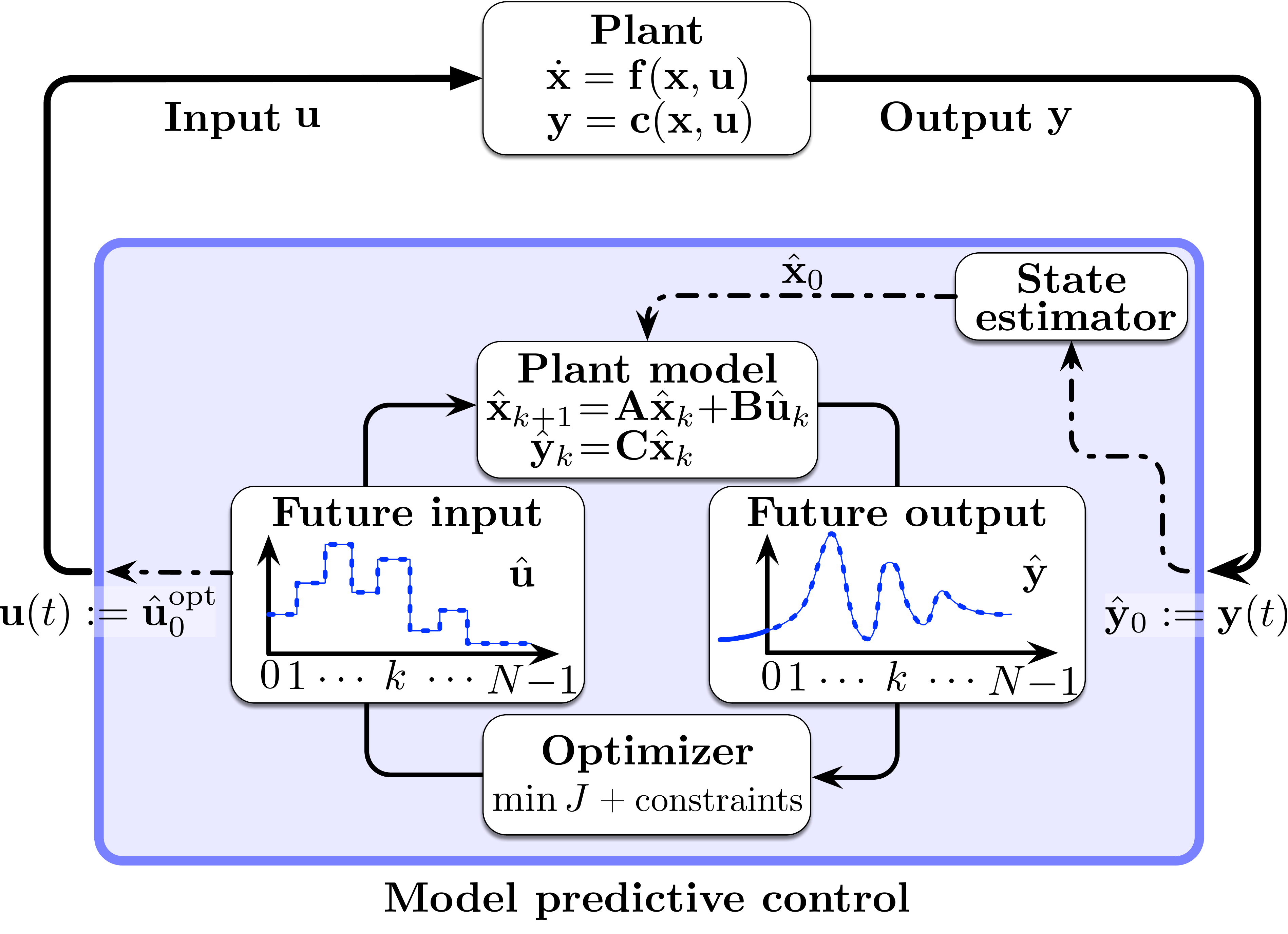}
    \caption{Schematic of the receding horizon model predictive control framework.}
    \label{fig:MPC}
\end{figure}

MPC is an adaptive control procedure solving an open-loop optimization problem over a receding horizon (see schematic in \cref{fig:MPC}).
The optimization problem aims to solve for a sequence of control inputs $\{\hat{\bu}_0,\hat{\bu}_2,\ldots,\hat{\bu}_{N-1}\}$ over the time horizon $T=N\Delta t$ that minimizes a pre-defined objective function $J$. Typically, only the first control input $\hat{\bu}_0^{\mathrm{opt}}$ is applied and then a new measurement is collected.
At each time instant a new measurement is collected, the optimization problem is re-initialized
and, thus, adaptively determines optimal control actions adjusting to model inaccuracies and changing conditions in the environment.
The most critical part of MPC is the identification of a dynamical
model that accurately and efficiently represents the system behavior in the presence of actuation.
If the model is linear, minimization of a quadratic cost functional subject to linear constraints results in a tractable convex problem.
There exist several variants of MPC, including nonlinear MPC, robust MPC, and explicit MPC, which are, however, more computationally expensive and thus limit their real-time applicability.
Combining Koopman-based models with linear MPC has the potential to significantly extend the reach of linear MPC for nonlinear systems.

The receding-horizon optimization problem can be stated as follows.
Linear MPC aims to minimize the following quadratic objective function
\begin{align}\label{Eq:CostFunction}
    \operatornamewithlimits{\min}\limits_{\hat\bu(\cdot\vert \by)\in\mathbb{U}}
  J
  =
   \operatornamewithlimits{\min}\limits_{\hat\bu(\cdot\vert \by)\in\mathbb{U}}
  \operatornamewithlimits{\sum}\limits_{k=0}^{N-1} \vert\vert\hat{\by}_{k}-{\bf r}_{k}\vert\vert_{\bQ}^2
  + \vert\vert\hat{\bu}_{k}\vert\vert_{\bR}^2
  + \vert\vert\Delta\hat{\bu}_{k}\vert\vert_{\bR_{\Delta}}^2
\end{align}
subject to discrete-time, linear system dynamics
\begin{subequations}\label{Eq:DiscreteTimeLinearControlSystem}
\begin{align}
\hat\bx_{k+1} &= \bA \hat\bx_{k} + \bB\hat\bu_{k},\\
\hat\by_{k} &= \bC \hat\bx_{k},
\end{align}
\end{subequations}
and state and input constraints
\begin{subequations}
\begin{align}
\by_{min} &\leq \hat{\by}_k \leq\by_{max},\\
\bu_{min} &\leq \hat{\bu}_k \leq \bu_{max},
\end{align}
\end{subequations}
where $\Delta\hat{\bu}_k \coloneqq  \hat{\bu}_k-\hat{\bu}_{k-1}$ is the control input rate.
Each term in the cost function~\eqref{Eq:CostFunction} is computed as the weighted norm of a vector, i.e. $\vert\vert\by\vert\vert_{\bQ}^2 \coloneqq  \by^T\bQ\by$.
In the model~\eqref{Eq:DiscreteTimeLinearControlSystem},
$\bA:\mathcal{X}\rightarrow \mathcal{X}$ is the state transition matrix,
$\bB:\mathcal{U}\rightarrow \mathcal{X}$ is the control matrix, and
$\bC:\mathcal{X}\rightarrow \mathcal{Y}$ the measurement matrix.
The weight matrices
$\bR\in\mathbb{R}^{q\times q}$,
$\bR_{\Delta}\in\mathbb{R}^{q\times q}$, and
$\bQ\in\mathbb{R}^{n\times n}$
are positive semi-definite and penalize the inputs, input rates, and deviations of the predicted output $\hat{\by}$ along a trajectory ${\bf r}$, respectively, and set their relative importance.
We define the control sequence to be solved over the receding horizon as $\hat\bu(0,\ldots,N-1\vert \by)\coloneqq  \{\hat{\bu}_0,\hat{\bu}_2,\ldots,\hat{\bu}_{N-1}\}$ given the measurement $\by$. The measurement $\by$ is the current output of the plant, whose dynamics are governed by a generally nonlinear system~\eqref{Eq:NonlinearControlSystem}, and is used to estimate the initial condition $\hat{\bx}_0$ for the optimization problem.
The general feedback control law~\eqref{Eq:FeedbackControlLaw} is then:
\begin{equation}
    {\bf d}(\by) = \hat\bu^{\mathrm{opt}}(0\vert \by)  = \hat{\bu}_0^{\mathrm{opt}}
\end{equation}
given a specific $\by$ and selecting the first entry of the optimized control sequence.

\begin{figure}[tb]
	\centering
	\includegraphics[width=0.8\textwidth]{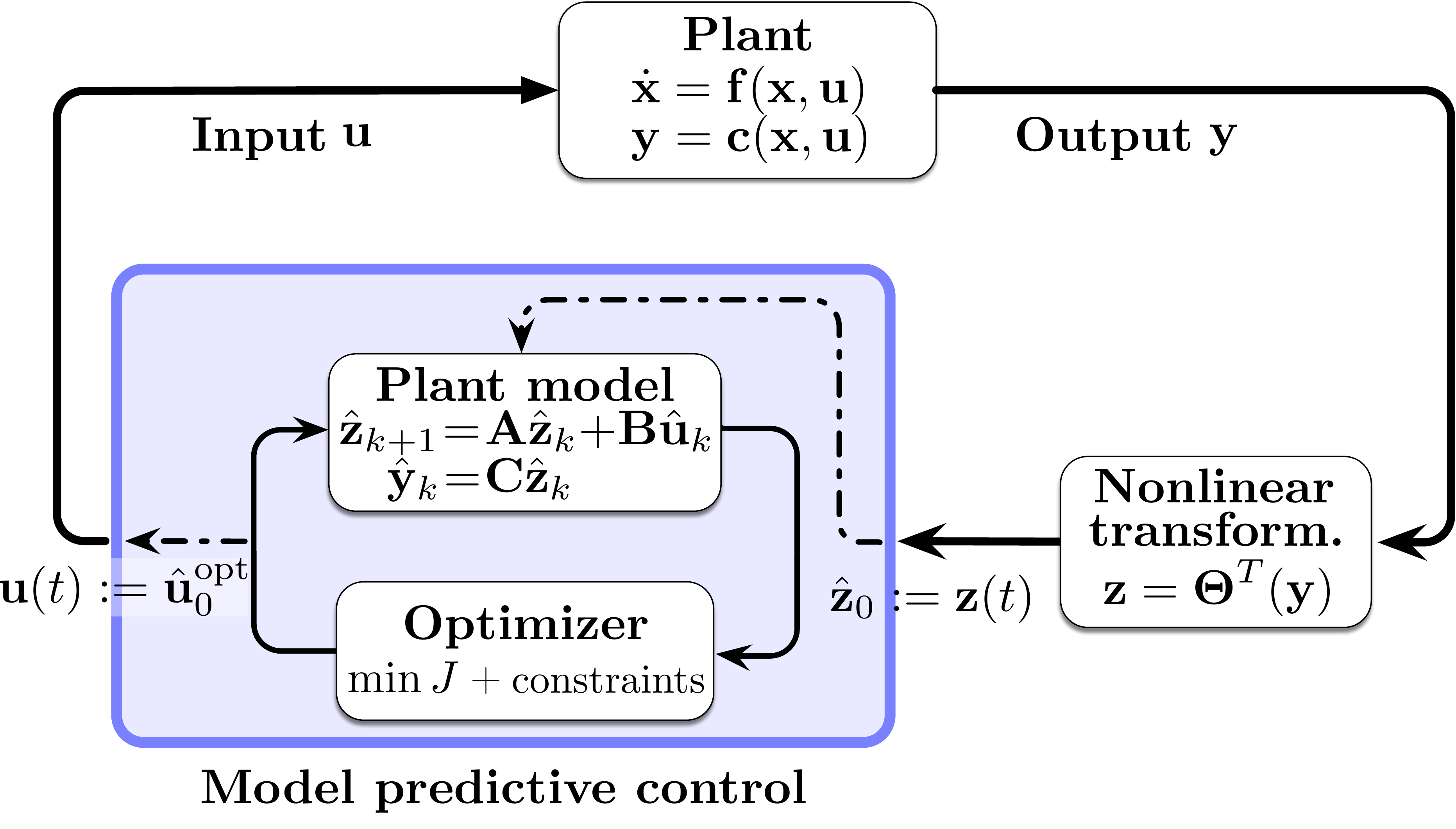}
	\caption{Schematic of the model predictive control framework incorporating a model based on the Koopman operator.}
	\label{fig:MPC-KO}
\end{figure}

Two primary research thrusts, which can be roughly categorized into approaches for either discrete or continuous inputs, have integrated Koopman theory and MPC. For the latter, the Koopman-MPC framework is schematically depicted in \cref{fig:MPC-KO}.
Output measurements are lifted into a higher-dimensional space using a nonlinear transformation. Dynamics are modeled in the lifted space, typically by solving a linear least-squares regression problem, and the resulting model is employed in the MPC optimization procedure. Besides the goal of achieving increased predictive power via a Koopman-based model, this approach further provides the possibility to readily incorporate nonlinear cost functions and constraints in a linear fashion by incorporating these directly in the set of observables~\cite{williams2016ifac}.

\subsubsection{Koopman operator theory for control systems}
\label{Sec:KoopmanForControl}
Koopman theory for control requires disambiguating the unforced dynamics from the effect of actuation.
The first Koopman-based approaches were developed for discrete-time systems, which  are more general and form a superset containing those induced by continuous-time dynamics.  Such discrete-time dynamics are often more consistent with experimental measurements and actuation, and may be preferred for numerical analysis.  However, there has been increasing effort in developing formulations building on the Lie operator, i.e.\ the infinitesimal generator of the semigroup of Koopman operators, for system identification and control.

Consider the non-affine, continuous-time control system given by (\ref{Eq:NonlinearControlSystem}a)
that is fully observable, i.e.\ $\by = \bx$ in (\ref{Eq:NonlinearControlSystem}b).
For every initial condition $\bx\in \mathcal{X}$ and control function $\bu\in \mathcal{U}$, there exists a unique solution $\bF^t(\bx,\bu)$ at time $t$ with initial condition $\bF^0(\bx,\bu) = \bx$.
System~\eqref{Eq:NonlinearControlSystem} represents a family of differential equations parameterized by the control functions $\bu$.
In order to analyze this from a Koopman perspective, it is convenient to introduce the \emph{control flow}, which defines~\eqref{Eq:NonlinearControlSystem} as a single dynamical system (for more details about the dynamical systems perspective of control theory we refer to Colonius et al.~\cite{colonius2000book}).
The control flow $\tilde{\bF}^t(\bx,\bu):\mathbb{R}\times \mathcal{X}\times \mathcal{U}\rightarrow \mathcal{X}\times \mathcal{U}$ associated with~\eqref{Eq:NonlinearControlSystem} is given by the map:
\begin{equation}
\tilde{\bF}^t(\bx, \bu) = (\bF^t(\bx,\bu), \boldsymbol\Theta^t(\bu)),
\end{equation}
where $\boldsymbol\Theta^t(\bu)$ is the shift on $\mathcal{U}$, so that $\boldsymbol\Theta^t(\bu)(s)=\bu(s+t)$, $s\in\mathbb{R}$.
Skew-product flows, such as $\tilde{\bF}^t(\bx,\bu)$, arise in topological dynamics to study non-autonomous systems, e.g. with explicit time dependency or parameter dependency.
Actuation renders the dynamical system and its associated Koopman family (or its generator) non-autonomous.
By defining the Koopman operator on the extended state $\tilde{\bx}\coloneqq [\bx^T,\bu^T]^T$, the Koopman operator becomes autonomous and is equivalent to the Koopman operator associated with the unforced dynamics.
Under further special considerations, standard numerical schemes for autonomous systems become readily applicable for system identification.
Beyond considering the extended state $\tilde{\bx}$, we also make the dependency on $\bx$ and $\bu$ explicit to disambiguate the state and inputs.

Let $g(\bx,\bu):\mathcal{X}\times \mathcal{U}\rightarrow \mathbb{C}$ be a scalar observable function of the extended state space. Each observable is an element of an infinite-dimensional Hilbert space and the semigroup of Koopman operators $\mathcal{K}^t:\mathcal{G}(\mathcal{X},\mathcal{U}) \rightarrow \mathcal{G}(\mathcal{X},\mathcal{U})$ acts on these observables according to:
\begin{equation}\label{Eq:NonlinearContinuuousTimeControlSystem:KoopmanOperator}
g(\bx(t),\bu(t)) = \mathcal{K}^tg(\bx_0,\bu_0)
= g(\tilde{\bF}^t(\bx_0,\bu_0)).\\
\end{equation}
Here, it is assumed that the Koopman operator acts on the extended state space in the same manner as the Koopman operator associated with the unforced, autonomous dynamical system.
A Koopman eigenfunction $\varphi(\bx,\bu)$ corresponding to eigenvalue $\lambda$ then satisfies
\begin{equation}\label{Eq:NonlinearContinuuousTimeControlSystem:KoopmanEigenfunction}
\varphi(\bx(t),\bu(t)) = \mathcal{K}^t \varphi(\bx_0,\bu_0) = \lambda^t \varphi(\bx_0,\bu_0).
\end{equation}
Further, a vector-valued observable
\begin{equation}\label{Eqn:VectorValuedObservableWithControl}
{\bf g}(\bx,\bu)
\coloneqq  \begin{bmatrix}
g_1(\bx,\bu)\\
\vdots\\
g_{p}(\bx,\bu)
\end{bmatrix}
\end{equation}
can be written in terms of the infinite  Koopman expansion as
\begin{equation}
\bg(\bx(t),\bu(t))
= \mathcal{K}^t\bg(\bx(0),\bu(0))
= \sum\limits_{j=1}^{\infty} \lambda_j^t\varphi_j(\bx_0,\bu_0)\bv_j,
\end{equation}
where $\bv_j = [\langle \varphi_j,g_1\rangle,\ldots,\langle \varphi_j,g_p\rangle]$.
This representation encompasses dynamics on $\bu$ itself, which may appear due to external perturbations when $\bu$ is interpreted as a perturbed parameter to the system.
While the actuation dynamics are typically known or set for both open-loop and closed-loop control,
it provides a convenient starting point for system identification.
Indeed, it is a useful representation for data-driven approaches that identify the underlying system dynamics and control simultaneously.  Depending on the choice of observable functions, further simplifications are possible to identify a model for the state dynamics by incorporating the effect of control, which is discussed in \cref{Sec:KoopmanControl:SystemID}.

We consider special cases in the discrete-time and continuous-time settings which are modeled within the {\em Koopman with inputs and control} (KIC) framework~\cite{proctor2017siads}.
The time-varying actuation input may evolve dynamically according to $\dot{\bu} = {\bf h}(\bu)$ or it may be governed by a state feedback control law~\eqref{Eq:FeedbackControlLaw}, $\bu = \bd (\bx)$, as in closed-loop control applications.

\paragraph{Discrete-time formulation}
For a discrete-time control system
\begin{subequations}\label{Eq:NonlinearDiscrete-TimeControlSystem}
	\begin{align}
	\bx_{k+1} &= {\bf F}(\bx_k,\bu_k),
	\end{align}
\end{subequations}
with initial condition $\bx_0$, the Koopman operator advances measurement functions according to
\begin{equation}\label{Eq:NonlinearDiscreteTimeControlSystem:KoopmanOperator}
\mathcal{K}g(\bx_{k},\bu_{k})
= g(\bF(\bx_{k},\bu_k),\bu_{k+1})
= g(\bx_{k+1},\bu_{k+1}).\\
\end{equation}
Koopman eigenpairs $(\varphi,\lambda)$ associated with~\eqref{Eq:NonlinearDiscreteTimeControlSystem:KoopmanOperator} satisfy:
\begin{equation}
\mathcal{K}\varphi(\bx_k,\bu_k)
= \varphi(\bF(\bx_{k},\bu_k),\bu_{k+1})
= \varphi(\bx_{k+1},\bu_{k+1})
= \lambda\varphi(\bx_k,\bu_k).
\end{equation}
By defining the control flow $\tilde{\bF}$ on the extended state space $\tilde{\bx} \coloneqq  [\bx^T,\bu^T]^T$, the dynamics become autonomous and can be written as
\begin{equation}\label{Eq:DiscreteTimeKoopmanControl:ExtendedState}
\mathcal{K}g(\tilde\bx_k)
= g(\tilde{\bF}(\tilde{\bx}_k))
= g(\tilde{\bx}_{k+1}).
\end{equation}
If $\bu_{k+1} = \bH(\bu_k)$, then $g(\bF(\bx_k), \bH(\bu_k)) = g(\bx_{k+1}, \bu_{k+1})$, thus allowing for a suitable choice of observable functions that can simultaneously model and identify dynamics for $\bx$ and $\bu$.
In many situations, we are not interested in the dynamics of $\bu$ itself but only of $\bx$, e.g. in control where $\bu$ is usually a design variable.
For instance, if the dynamics of $\bu$ is prescribed by a specific state-feedback law ${\bu_k} = \bD(\bx_k)$, then
$\mathcal{K}g(\bx_k,\bu_k)
= g(\bF(\bx_k,\bu_k),\bD(\bx_k))
= g(\bx_{k+1},\bD(\bx_{k+1}))$.
By defining $\bF^{\bD}(\bx)\coloneqq \bF(\bx,\bD(\bx))$ and restricting the observable to be a function solely of the state, the Koopman operator is associated with the autonomous dynamics $\bF^{\bD}$ for a given control law $\bD$:
$\mathcal{K}g(\bx_k)
= g(\bF^{\bD}(\bx_k))
= g(\bx_{k+1})$.
If instead, we consider a constant exogenous forcing or discrete control action $\bar{\bu}\in\mathcal{U}$, where $\mathcal{U}$ is the set of discrete inputs.  The Koopman operator may be defined for each discrete input separately:
$\mathcal{K}^{\bar{\bu}}g(\bx_k,\bar{\bu})
= g(\bF({\bx}_{k},\bar{\bu}),\bar{\bu})
= g(\bx_{k+1},\bar{\bu})$.
In general, the Koopman operator and its associated eigenfunctions are parameterized by the discrete control input $\bar{\bu}$ and the dynamics are autonomous for each $\bar{\bu}$.
The flow map $\bF^{\bar{\bu}}(\bx)\coloneqq \bF(\bx,\bar{\bu})$ is then defined for each $\bar{\bu}$.  Considering only the reduced dynamics on observable functions of the state, we then obtain
\begin{equation}\label{Eqn:KoopmanControl:DiscreteTimeLPV}
\mathcal{K}^{\bar\bu}g(\bx_k)
= g(\bF^{\bar{\bu}}(\bx_k))
= g(\bx_{k+1}).
\end{equation}
By switching from continuous to discrete inputs, a single model is replaced with a family of models, so that the specific control dependency of the state does not have to be captured. Instead of optimizing the input itself, one may then optimize the switching times between inputs, as in Peitz et al.~\cite{peitz2017arxiv}.

\paragraph{Lie operator formulation}
There has been increasing interest in the control formulation for the infinitesimal generator of the Koopman operator family.
It can be shown~\cite{lasota1994} that if observables $g$ are continuously differentiable with compact support, then they satisfy the first-order partial differential equation~\eqref{eq:lie-eigenfunction}.
The Lie operator $\gen{}g = \nabla_{\tilde{\bx}} g\cdot\tilde{\bf f}$ associated with the dynamics of the control system~\eqref{Eq:NonlinearControlSystem} induces the dynamics
\begin{equation}\label{Eq:GeneratorControl:ObservableDynamics}
\frac{d}{dt}g(\bx,\bu) = \gen{}g(\bx,\bu).
\end{equation}
Similarly, smooth eigenfunctions corresponding to the eigenvalue $\mu$ satisfy
\begin{equation}
\frac{d}{dt}\varphi(\bx,\bu) = \gen{}\varphi(\bx,\bu)  = \mu \varphi(\bx,\bu).
\end{equation}
Note that the smooth eigenfunction of the generator is also a Koopman eigenfunction~\eqref{Eq:NonlinearContinuuousTimeControlSystem:KoopmanEigenfunction} and their eigenvalues are connected via $\mu = \log(\lambda)$ (see also \cref{Sec:Koopman}).

We may rewrite the Lie operator in~\eqref{Eq:GeneratorControl:ObservableDynamics} explicitly using the chain rule
\begin{equation}\label{Eq:GeneratorControl:ChainRule}
\gen{}g(\bx,\bu)
= \nabla_{\bx}g(\bx,\bu)\cdot\dot{\bx} + \nabla_{\bu}g(\bx,\bu)\cdot\dot{\bu}.
\end{equation}
The derivatives $\dot{\bx} = {\bf f}$ and $\dot{\bu} = {\bf h}$ are both velocity vectors (state or action change per time unit) and generally depend on $\bx$ and $\bu$.
These describe local changes in the observable function due to local changes in $\bx$ and external forcing via $\bu$.
For a dynamically evolving $\bu$ it may be possible to approximate the Koopman operator on the extended state using the system identification framework based on the generator formulation~\cite{mauroy2020tac}.
It is also possible to consider $\dot{\bu}$ as the new input to the system~\cite{moore1967}, while applying more traditional methods.
Equation~\eqref{Eq:GeneratorControl:ChainRule} represents the adjoint equation of the controlled Liouville equation (see e.g.~\cite{brockett2007ams}), which describes how a density function evolves in state space under the effect of external forcing. Using the analogy to the scattering of particles, these transport operators have been used to model and derive optimal control for agent-based systems~\cite{kwee2001}.
For the state-feedback law $\bu = \bd(\bx)$, then the right-hand side of~\eqref{Eq:GeneratorControl:ObservableDynamics} becomes
$\gen{}g(\bx,\bu)
= \left[ \nabla_{\bx}g(\bx,{\bf d}(\bx)) + \nabla_{\bf d}g(\bx,{\bf d}(\bx))\cdot\nabla_{\bx}{\bf d}(\bx) \right] \cdot{\bf f}(\bx,{\bf d}(\bx))$.
Since this is solely a function of the state, the (autonomous) Koopman operator associated with a particular state-feedback law $\bd(\bx)$ can be defined on the reduced state $\bx$:
$\gen{}g(\bx)
= \nabla_{\bx}g(\bx)\cdot{\bf f}(\bx,\bd(\bx))$.

Constant exogenous forcing or discrete control actions $\bar{\bu}\in\mathcal{U}$ render the system autonomous for each $\bar{\bu}$ so that
$\gen{}^{\bar\bu}g(\bx,\bar{\bu})
= \nabla_{\bx}g(\bx,\bar{\bu})\cdot{\bf f}(\bx,\bar\bu)$
with $\dot{{\bu}}={\bf 0}$.
In order to apply most control methods a model for the reduced state $\bx$ is required. Thus, by restricting the space of observables defined on the reduced state $\bx$, we obtain
\begin{equation}
\gen{}g(\bx)
= \nabla_{\bx}g(\bx)\cdot{\bf f}(\bx,\bar\bu).
\end{equation}
This representation is the starting point for data-driven formulations, e.g., to estimate the associated Koopman operator by assuming the zero-order hold for the input across consecutive snapshots, or Koopman operators parameterized by the discrete control input so that a gain-scheduled or interpolated controller may be enacted.

\paragraph{Bilinearization} Many applications are modeled by control-affine systems
\begin{equation}\label{Eqn:ControlAffineSystem}
\dot{\bx} = {\bf f}_0(\bx) + \sum_{j=1}^q {\bf f}_j(\bx) u_j.
\end{equation}
Using Carleman linearization, it is possible to transform this system into a (generally infinite-dimensional) bilinear model using multivariable monomials of the state as observables, which is then truncated at a suitable order.
However, this may nevertheless lead to systems with an undesirably high dimensionality for the required accuracy.
A general nonlinear system $\dot{\bx}=\mathbf{f}(\bx,\bu)$ with analytic $\mathbf{f}$ can also be transformed into an infinite-dimensional bilinear system of the form $\dot{\bz} = \bA\bz + \dot{\bu}\bB\bz $, where components of $\bz$ consist of multivariable monomials $x_i^k u_j^l$.
Note the similarity to the Koopman generator PDE~\eqref{Eq:GeneratorControl:ChainRule}: considering a smooth vector-valued observable $\bg(\bx,\bu)$, interpreting $\dot{\bu}$ as new input, and assuming that $\nabla_{\bu}\bg(\bx,\bu)$ lies in the span of $\bg(\bx,\bu)$, this too leads to a bilinear equation and is equivalent to Carleman linearization.

In general, any vector-valued smooth observable with compact support that is solely a function of the state $\bx$, which evolves according to~\eqref{Eqn:ControlAffineSystem}, satisfies
\begin{equation}
\frac{d}{dt}\bg({\bx}) = \mathcal{L}^{\bu}\bg(\bx) = \nabla_{\bx}\bg(\bx)\cdot {\bf f}_0(\bx) +  \nabla_{\bx}\bg(\bx)\cdot\left(\sum_{j=1}^q {\bf f}_j(\bx) u_j\right),
\end{equation}
where $\mathcal{L}^{\bu}$ denotes the non-autonomous Lie operator due to the external control input $\bu$.
If $\nabla_{\bx}\bg(\bx)\cdot {\bf f}_j(\bx)$ lies in the span of ${\bf g}(\bx)$ for all $j=1,\ldots,q$, then there exists a set of constant matrices $\bB_j$ so that
\begin{equation}\label{Eqn:BilinearEquationForObservables}
\frac{d}{dt}\bg({\bx}) = \bA\bg(\bx) +  \sum_{j=1}^q u_j\bB_j\bg(\bx).
\end{equation}
leading to a bilinear equation for the observable $\bg$.
Here, $\bA$ and $\bB_j$ decompose the {matrix representation of the} Lie operator into an unforced term $\bA$ and forcing terms $\bB_j$ (compare also \cref{Fig:KoopmanControl:MethodsOverview}).
The following theorem appears as Thm. 2 in~\cite{goswami2017cdc} and states the bilinearizability condition for Eqns.~\eqref{Eqn:ControlAffineSystem} and~\eqref{Eqn:BilinearEquationForObservables}:
\begin{theorem}
  Suppose there is a finite set of Koopman eigenfunctions associated with the unforced vector field ${\bf f}_0$, $\varphi_j(\bx)$ for $j=1,\ldots,n$, and these form an invariant subspace of the Lie derivatives $\mathcal{L}_{{\bf f}_j}$, $j=1,\ldots,q$. Then there exists a finite-dimensional bilinear representation for~\eqref{Eqn:ControlAffineSystem}, and in turn also for~\eqref{Eqn:BilinearEquationForObservables}.
\end{theorem}
A special class of observables are Koopman eigenfunctions (or eigenfunctions of the Lie operator), which behave linearly in time by definition.
Specifically, smooth eigenfunctions of the Koopman operator associated with the uncontrolled system satisfy
\begin{equation}\label{Eqn:BilinearEquationForEigenfunctions}
\frac{d}{dt}\varphi({\bx}) = \mu \varphi(\bx) +  \nabla_{\bx}\varphi(\bx)\cdot\left(\sum_{j=1}^q {\bf f}_j(\bx) u_j\right),
\end{equation}
with unforced dynamics determined by its associated eigenvalue $\mu$ and
where the second terms reflects explicitly how the control input affects these eigenfunctions.
Koopman eigenfunctions are associated with the global behavior of the system and represent a Koopman-invariant subspace.
Equation~\eqref{Eqn:BilinearEquationForEigenfunctions} may also represent a bilinear equation (for further details see \cref{Sec:Control:EigenfunctionSubspace}).
Importantly, for systems with low-rank structure it is possible to describe the dynamics in terms of a small number of eigenfunctions which sufficiently well describe the global behavior and avoid the explosion of observables of the Carleman linearization.
Even though the system may not be completely bilinearizable, approximate bilinearizability may be sufficient for accurate prediction~\cite{goswami2017cdc}.
Furthermore, there exists a vast literature of diagnostic and control techniques tailored to bilinear systems that can be utilized.

\subsection{Data-Driven Control}
\label{Sec:KoopmanControl:SystemID}
\begin{figure}[tb]
	\centering
	\includegraphics[width=0.85\textwidth, trim = 0 0 0 0cm, clip=true]{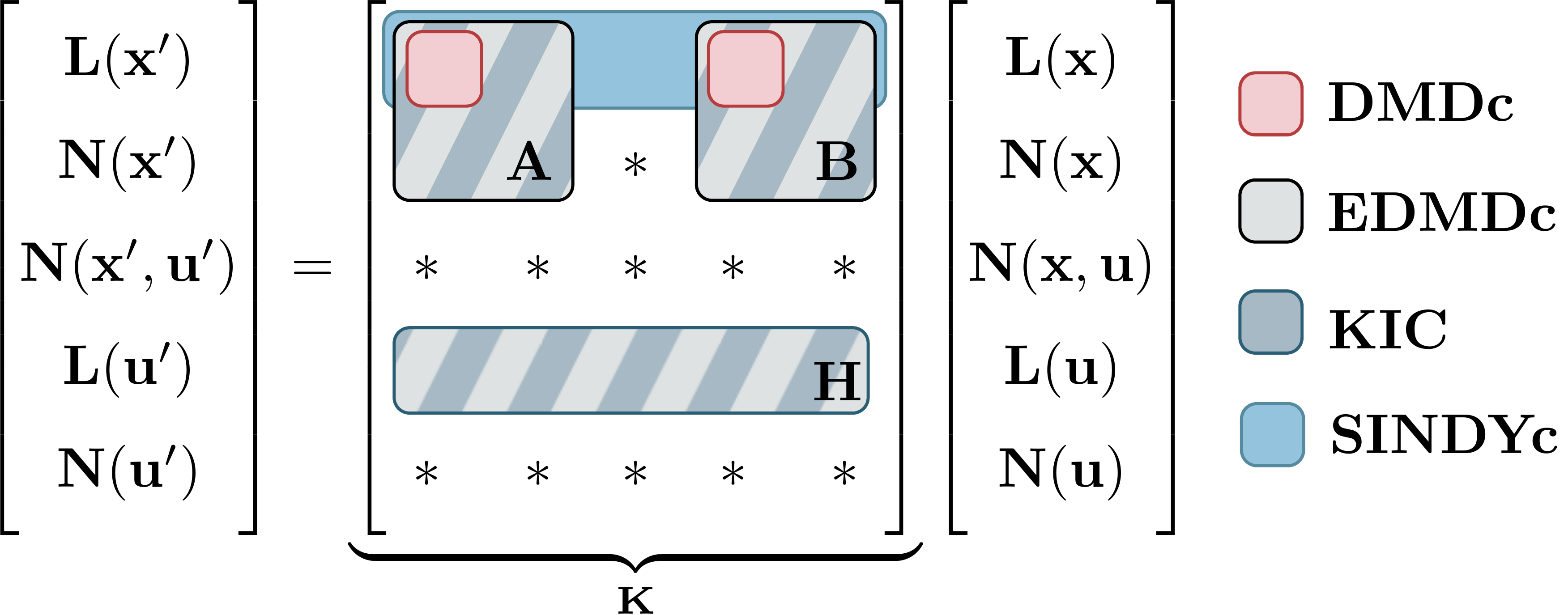}
	\caption{Data-driven finite-dimensional approximation of the Koopman operator defined on the extended state space $(\bx,\bu)$ in practice. As with eDMDc it is also possible to identify evolution dynamics $\bu'=\bH(\bx,\bu)$ for the control input within the KIC framework. A state-feedback law $\bu = \bD(\bx)$ may be hidden within $\bA$ and can be identified with knowledge of $\bu$.  Vector-valued observables are splitted into a linear ${\bf L}$ and nonlinear ${\bf N}$ part, e.g. ${\bg}(\bx)\coloneqq  [{\bf L}^T(\bx)\; {\bf N}^T(\bx)]^T$.}
	\label{Fig:KoopmanControl:MethodsOverview}
\end{figure}

\phantom{}
In contrast to model-based control, there are many emerging
advances in equation-free, data-driven system identification that leverage Koopman theory.
The relationships between some of the more prominent methods and their connection to the Koopman operator are shown in \cref{Fig:KoopmanControl:MethodsOverview}.
%

\subsubsection{Dynamic mode decomposition with control}
\label{Sec:Control:DMDc}
A major strength of DMD is the ability to describe complex and high-dimensional dynamical systems in terms of a small number of dominant modes, which represent spatio-temporal coherent structures.
Reducing the dimensionality of the system from $n$ (often millions or billions) to $r$ (tens or hundreds) enables faster and lower-latency prediction and estimation, which
generally translates directly into controllers with higher performance and robustness.
Thus, compact and efficient representations of complex systems such as fluid flows have been long-sought, resulting in the field of reduced-order modeling.
However, the original DMD algorithm was designed to characterize naturally evolving systems, without accounting for the effect of actuation and control.

The {\em dynamic mode decomposition with control} (DMDc) algorithm by
Proctor et al.~\cite{proctor2016siads} extends DMD to disambiguate between the natural unforced dynamics and the effect of actuation.
This generalizes the DMD regression to the control form
\begin{align}
\bx_{k+1} \approx \bA \bx_k + \bB\bu_k . \label{Eq:DMDc:Propagator-repeat}
\end{align}
For DMDc, snapshot pairs $\{\bx_j,\bx_{j+1},\bu_j \}_{j=1}^m$ of the state variable $\bx$ and actuation command $\bu$ are collected and organized into the following data matrices
\begin{equation}\label{Eq:DMDc:Data}
\bX  = \begin{bmatrix}\vert &       & \vert\\
\bx_1  &\cdots &\bx_{m-1}\\
\vert  &       & \vert\end{bmatrix},\;
\bX' = \begin{bmatrix}\vert &       & \vert\\
\bx_2  &\cdots &\bx_{m}\\
\vert  &       & \vert\end{bmatrix},\;
\boldsymbol{\Upsilon}  = \begin{bmatrix}\vert &       & \vert\\
\bu_1  &\cdots &\bu_{m-1}\\
\vert  &       & \vert\end{bmatrix},\;
\end{equation}
where $\bX'$ is a time-shifted version of $\bX$.
Then the system matrices $\bA$ and $\bB$ are determined jointly by solving the following least-squares optimization problem
\begin{align}
\operatornamewithlimits{\min}\limits_{\begin{bmatrix}
	\bA & \bB
	\end{bmatrix}}
\norm[\Big]{\bX'-\begin{bmatrix}
	\bA & \bB
	\end{bmatrix}
	\begin{bmatrix}
	\bX\\
	\boldsymbol{\Upsilon}
	\end{bmatrix}}^2_2,
\end{align}
where the solution is given as $[\bA \; \bB] = \bX'[\begin{smallmatrix}\bX\\ \boldsymbol{\Upsilon} \end{smallmatrix}]^{\dagger}$.

Originally motivated by intervention efforts in epidemiology~\cite{proctor2014epj}, DMDc based on linear and nonlinear measurements of the system has since been used with MPC for enhanced control of nonlinear systems by Korda and Mezi\'c~\cite{korda2016arxiv} and by Kaiser et al.~\cite{kaiser2017arxivb}, with the DMDc method performing surprisingly well, even for strongly nonlinear systems.  DMDc allows for a flexible regression framework, which can accommodate undersampled measurements of an actuated system and identify an accurate and efficient low-order model using compressive system identification~\cite{bai2017aiaa}, which is closely related to the eigensystem realization algorithm (ERA)~\cite{juang1985jgcd}.

\subsubsection{Extended dynamic mode decomposition with control}
\label{Sec:Control:EDMDc}
The DMDc framework provides a simple and efficient numerical approach for system identification,
in which the Koopman operator is approximated with a best-fit linear model advancing linear observables and linear actuation variables.
As discussed in \cref{Sec:EDMD}, eDMD is an equivalent approach to DMD that builds on nonlinear observables. An extension of eDMD for controlled systems was first introduced in Williams et al.~\cite{williams2016ifac} to approximate the Koopman operator associated with the unforced system and to correct for inputs affecting the system dynamics and data.
Inputs are handled as time-varying system parameters and the Koopman operator is modeled as a parameter-varying operator drawing inspiration from {\em linear parameter varying} (LPV) models.
This approach has been generalized in Korda et al.~\cite{korda2016arxiv} to identify the system matrices $\bA$ and $\bB$ in the higher-dimensional observable space which disambiguates the unforced dynamics and control on observables.
In particular, the Koopman operator is defined as an autonomous operator on the extended state $\tilde{\bx}\coloneqq [\bx^T,\bu^T]^T$ as in~\eqref{Eq:DiscreteTimeKoopmanControl:ExtendedState} and observables can be nonlinear functions of the state and input, i.e.\ $g(\bx,\bu)$.
However, in practice, simplifications are employed to allow for convex formulations of the control problem.

Assuming that the observables are nonlinear functions of the state and linear functions of the control input, i.e.\
${\bg(\bx, \bu) \coloneqq  [\theta_1(\bx),\ldots, \theta_p(\bx), u_1,\ldots,u_q]^T\in\mathbb{R}^{p+q}}$,
and by restricting the dynamics of interest to the state observables $\theta(\bx)$ themselves, the linear evolution equation to be determined is
\begin{align}
\bz_{k+1} \approx \bA \bz_k + \bB\bu_k,\label{Eq:EDMDc:Propagator}
\end{align}
where $\bz\in\mathbb{R}^p$ is the vector-valued observable as in \cref{Sec:EDMD} defined as
\begin{align}
\bz \coloneqq  \bTheta^T(\bx) = \begin{bmatrix}
\theta_1(\bx)\\
\theta_2(\bx)\\
\vdots\\
\theta_p(\bx)\\
\end{bmatrix}.
\end{align}
Analogous to DMDc~\eqref{Eq:DMDc:Data}, the (time-shifted) data matrices in the lifted space, $\bZ = \bTheta^T(\bX)$ and $\bZ' = \bTheta^T(\bX')$, are evaluated given data $\bX',\bX, \boldsymbol{\Upsilon}$.
The system matrices $\bA,\bB$ are determined from the least-squares regression problem
\begin{align}
\operatornamewithlimits{\min}\limits_{\begin{bmatrix}
	\bA & \bB
	\end{bmatrix}}
\norm[\Big]{
	\bZ'-
	\begin{bmatrix}
	\bA & \bB
	\end{bmatrix}
	\begin{bmatrix}
	\bZ\\
	\boldsymbol{\Upsilon}
	\end{bmatrix}}^2_2
\end{align}
where the solution is given as $[\bA\;\bB] = \bZ'[\begin{smallmatrix}\bZ\\ \boldsymbol{\Upsilon} \end{smallmatrix}]^{\dagger}
= \bTheta^T(\bX')[\begin{smallmatrix}\bTheta^T(\bX)\\ \boldsymbol{\Upsilon} \end{smallmatrix}]^{\dagger}$.
The state $\bx$ is often included in the basis or space of observables, e.g.\
$\bz = [\bx^T,\; \bTheta(\bx)]^T$,
and can then be estimated by selecting the appropriate elements of the observable vector $\bz$, so that $\bx = \bC \bz$ with the measurement matrix $\bC = [\bI_{n\times n}\;{\bf 0}]$.
If the state vector is not included as observable, one may approximate the measurement matrix by solving an equivalent least-squares problem:
\begin{align}\label{Eqn:KoopmanControl:EDMDc:MeasurementMatrix}
\operatornamewithlimits{\min}\limits_{\bC}
\norm[\Big]{
	\bX-
	\bC\bZ}^2_2
\end{align}
Alternatively, it has been shown that the state can be estimated using, e.g., multi-dimensional scaling~\cite{kawahara2016neurips}.
If full-state information is not available, but only input-output data is available, the observable vector $\bz$ often needs to be augmented with past input and output values, as is typically done in system identification~\cite{ljung1999book}, following attractor reconstruction methods, such as the Takens embedding theorem~\cite{takens1981lnm}.

As eDMD is known to suffer from overfitting, it is important to regularize the problem based on, for instance, the $L_{1,2}$-norm for group sparsity~\cite{williams2016ifac} or the $L_1$-norm~\cite{kaiser2017arxiv}.
eDMDc yields a linear model for the controlled dynamics in the space of observables through a nonlinear transformation of the state. As with DMDc, it can be combined with any model-based control approach.
In  particular, when combined with MPC for a quadratic cost function, the resulting MPC optimization problem can be shown~\cite{korda2016arxiv} {to lead to a convex quadratic programming problem, whose computational cost remains comparable to that of linear MPC and is independent of the number of observables}, and to outperform MPC with a model based on a local linearization or Carleman linearization.

{The advantage of a linear representation can also represent a restriction for some systems in the sense that only linear control methods can be applied. Bruder et al.~\cite{bruder2020arxiva} have recently demonstrated that the accuracy of approximate bilinear approximations improves with increasing number of basis functions, while that of linear approximations do not necessarily limiting their control performance. For instance, all control-affine systems can be reformulated as an infinite-dimensional bilinear system (see also Sec.~\ref{Sec:KoopmanForControl}), but not necessarily into an equivalent linear one.
As a consequence, bilinear Koopman control can represent a balancing middle ground outperforming linear Koopman-based MPC performance-wise and nonlinear MPC computationally-wise.}

\subsubsection{Generalizations of DMDc and eDMDc}
\label{Sec:Control:GeneralizationsDMDc}

The simple formulation of DMDc and eDMDc as a linear regression problem allows for a number of generalizations that can be easily exploited for integration into existing linear model-based control methods.  A number of these generalizations are highlighted here.

\paragraph{Mixed observables}
eDMDc allows for generalization by using a suitable choice of observable functions.  It has also been proposed~\cite{korda2016arxiv,proctor2017siads} to include nonlinear observables of the input, e.g $g(u)=u^2$, or mixed observables of the state and input, e.g $g(x,u)=xu$.
The advantage of incorporating these measurements has yet to be demonstrated and a subsequent control optimization is not straightforward.
However, for systems forced by external inputs or parameters
this may provide new opportunities to identify their underlying dynamics, e.g. $\dot{\bu}={\bf h}(\bx,\bu)$ or a state-feedback law $\bu = \bd(\bx)$ (see also \cref{Fig:KoopmanControl:MethodsOverview}).

\paragraph{Input and output spaces}
Within the KIC framework~\cite{proctor2017siads}, different domain and output spaces of the Koopman operator approximation have been examined.
This develops more formally what has been implicitly carried out in eDMDc-like methods to obtain models for only the state itself.
If only observables that are functions of the state and not of the input are considered,
the output space of the Koopman operator can be restricted to a subspace of the Hilbert space spanned by observables which are solely functions of the state and not of the extended state space:
\begin{equation}\label{KoopmanControl:DifferentDomainOutputSpaces}
\bz_{\bx,k+1} =
\bK_{\bx}\bz_k
= \bK_{\bx}
\begin{bmatrix}
\bg(\bx_k)\\ \bg(\bx_k,\bu_k)\\ \bg(\bu_k)
\end{bmatrix}
\end{equation}
with vector-valued observables
$\bz_{\bx}\coloneqq  \bg(\bx)$ and
$\bz \coloneqq
[ \bg^T(\bx),\; \bg^T(\bx,\bu),\; \bg^T(\bu)]^T$.
This is reasonable as the prediction of the future actuation input $\bu$ is not of interest for the purpose of control.
Analogously, Koopman eigenfunction expansions of observables must be distinguished based on the domain and output spaces.

\paragraph{Time delay coordinates}
Time delay coordinates provide an important and universal class of measurement functions for systems which display long-term memory effects.
They are especially important in real-world applications where limited access to full state measurements are available (see~\cref{Sec:HAVOK}), and they have demonstrated superior performance for control, compared with models using monomials as observables~\cite{kaiser2017arxivb}.
For simplicity, we will consider a sequence of scalar input-output measurements $u(t)$ and $x(t)$.
From these, we may construct a delay vector of inputs $\bu_k \coloneqq  [u_k,\; u_{k-1},\; u_{k-2},\;\ldots,u_{k-m_u}]^T$ and outputs $\bz_k \coloneqq  \bg(\bx_k) = [x_k,\; x_{k-1},\; x_{k-2},\;\ldots,x_{k-m_x}]^T$, respectively.
Here, it is assumed $m=m_x=m_u$ for simplicity. The dynamics may then be represented as
\begin{subequations}
	\begin{align}
	\bz_{k+1} &= \bA \bz_k + \bB \bu_k,  \label{Eqn:KoopmanControl:TimeDelayCoordinates} \\
	y_k &= \begin{bmatrix} 1 & 0 & \ldots & 0\end{bmatrix}\bz_k = x_k, \label{Eqn:KoopmanControl:TimeDelayCoordinatesMeasurement}
	\end{align}
\end{subequations}
where the current state $x$ is recovered from the first component of $\bz_k$.
Both the system matrix $\bA$ and the control matrix $\bB$ must satisfy an upper triangular structure to not violate causality; otherwise, current states will depend on future states and inputs.
Here, $\bu_k$ is constructed from past inputs, while effectively the system has only a single input.
Thus, it is recommended to augment $\bz$ with past inputs, i.e.\ $\bz_k\coloneqq \bg(\bx_k,\bu_{k-1})=[x_k, (x,u)_{k-1}, (x,u)_{k-2}, \ldots, (x,u)_{k-m}]$, so that the current actuation value $u_k$ appears as a single input to the system~\cite{korda2016arxiv,Kaiser2020chapter}:
\begin{subequations}
	\begin{align}
	\bz_{k+1} &= \hat{\bA} \bz_k + \hat{\bb}\, u_k,  \label{Eqn:KoopmanControl:TimeDelayCoordinates:AugmentedState} \\
	y_k &= \begin{bmatrix} 1 & 0 & \ldots & 0\end{bmatrix}\bz_k = x_k. \label{Eqn:KoopmanControl:TimeDelayCoordinatesMeasurement:AugmentedState}
	\end{align}
\end{subequations}
This delay observable $\bz_k=\bg(\bx_k,\bu_{k-1})$ has also been used in eDMDc if only input-output data is available~\cite{korda2016arxiv}. Without the delay information, eDMDc may fail despite the lifting into a higher-dimensional observable space (as for the examples in \cref{Sec:KoopmanControlStrategies}).
By combining~\eqref{Eqn:KoopmanControl:TimeDelayCoordinates} and ~\eqref{Eqn:KoopmanControl:TimeDelayCoordinatesMeasurement} into a single equation, we obtain an equivalent description to~\eqref{KoopmanControl:DifferentDomainOutputSpaces} with different domain and output spaces:
\begin{equation}
z_{1,k+1} = y_{k+1} =
\begin{bmatrix}
\bC\bA & \bC \bB
\end{bmatrix}
\begin{bmatrix}
\bz_k \\ \bu_k
\end{bmatrix}.
\end{equation}
This formulation is analogous to autoregressive-exogenous (ARX) models in linear system identification~\cite{ljung1999book} , where the current output is represented as a linear superposition of past measurements and inputs.

\paragraph{Parametrized models}
The control input can also be restricted to a finite set of discrete values.
The dynamics are then modeled as a linear parameter-varying model, so that each discrete input is associated with an autonomous Koopman operator~\eqref{Eqn:KoopmanControl:DiscreteTimeLPV}, which is approximated using either DMD or eDMD~\cite{peitz2017arxiv}. The control input is thus not modeled explicitly, but implicitly, by switching the model whenever the actuation command switches. The advantage of this formulation is that each model associated with a discrete control input remains autonomous.
Recently, this approach has been extended to continuous control input by interpolating between the finite set of control actions, so that the system dynamics are predicted using interpolated models based on the infinitesimal generator of the Koopman operator, leading to bilinear models~\cite{peitz2020arxiv} similar to~\eqref{Eqn:BilinearEquationForObservables}.
This also alleviates performance issues associated with the time discretization of the models when optimizing their switching times~\cite{peitz2020arxiv,klus2020}.

\paragraph{Deep Koopman models for control}
Finding a suitable basis or feature space that facilitates a Koopman-invariant subspace remains challenging.
Due to their flexibility and rich expressivity~\cite{raghu2017pmlr}, there is increased interest in using neural networks to learn complex Koopman representations (see also \cref{sec:NNKoop}).
Extending these architectures to incorporate the effect of control requires careful handling of the control variable to ensure tractability of the subsequent control optimization, which a nonlinear transformation of the control variable generally impedes.
Building on their work on deep DMD~\cite{yeung2017arxiv} and Koopman Gramians~\cite{yeung2018acc}, Liu et al.~\cite{liu2017arxiv}
propose to decompose the space of observables as three independent neural networks, associated with parameters $\boldsymbol{\theta}_x$, $\boldsymbol{\theta}_{xu}$, and $\boldsymbol{\theta}_u$, respectively, to which the input is either the state variable alone, the control variable alone, or both:
\begin{equation}\label{Eqn:KoopmanControl:DeepDMDc}
\bz_{\bx,k+1} 
=
\bK\bz_k
=
\begin{bmatrix}
\bK_x & \bK_{xu} & \bK_u
\end{bmatrix}
\begin{bmatrix}
\bz_{\bx,k}\\
\bz_{\bx\bu,k}\\
\bz_{\bu, k}
\end{bmatrix}
\quad\text{with}\quad
\begin{matrix}
\bz_{\bx} \coloneqq  \bg(\bx,\boldsymbol{\theta}_x),\\
\hspace{0.45cm}\bz_{\bx\bu}\coloneqq  \bg(\bx,\bu,\boldsymbol{\theta}_{xu}),\\
\bz_{\bu} \coloneqq  \bg(\bu,\boldsymbol{\theta}_u).
\end{matrix}
\end{equation}
The unknown Koopman matrices $\bK_{\bullet}$ and neural network parameters $\boldsymbol{\theta}_{\bullet}$, where $\bullet = x$, $xu$, or $u$, can then be optimized jointly~\cite{liu2017arxiv} or separately by switching alternately between model update and parameter update, as in~\cite{han2020arxiv,li2019arxiv}.
While the model is linear in the latent space, it is not tractable for prediction and control of the original system.
In order to map back to the original state for predictions and to solve the control optimization problem, the states and control inputs are also directly included as observables, so that
$\bz_{\bx} \coloneqq  [\bx,\bg(\bx,\boldsymbol{\theta}_x)]^T$ and
$\bz_{\bu} \coloneqq  [\bu,\bg(\bu,\boldsymbol{\theta}_u)]^T$.
While $\bz_{\bx\bu} \coloneqq  \bg^T(\boldsymbol{\theta}_{xu})$ is kept the same, it is also possible to neglect $\bK_{xu}$ under certain conditions. In subsequent work~\cite{you2018ifac,hasnain2019arxivb}, slightly modified approaches are studied while maintaining the tractability of the control synthesis problem.

Analogous to eDMDc, a multi-layer perceptron neural network can be used solely to lift the state space, while the effect of the control input in the lifted space remains linear~\cite{han2020arxiv}, i.e.
$\bz_{\bx,k+1} = \bA \bz_{\bx,k} + \bB \bu$,
where $\bK_x = \bA$, $\bK_u = \bB$, $\bK_{xu}=0$, $\bz_{\bx}\coloneqq \bg(\bx,\boldsymbol{\theta}_x)$, and
$\bz_{\bu}\coloneqq \bu$. The transformation from the latent space to the original state is determined via a least-squares optimization problem analogous to~\eqref{Eqn:KoopmanControl:EDMDc:MeasurementMatrix} where $\bZ$ then represents the output of the neural network $\bz_{\bx}$.
Autoencoders provide a natural framework for learning Koopman representations~\cite{lusch2017arxiv} (see also \cref{sec:NNKoop}) combining the state lifting operation and inverse transformation. Recently, this approach has been extended for control by incorporating the control term~\cite{han2020arxiv}, so that only the state is encoded into $\bz_{\bx,k}$ and the shifted feature state $\bz_{\bx,k+1}$ is decoded.
These formulations facilitate immediate application of standard linear control methods.
Graph neural networks (GNN)~\cite{li2019arxiv} have also been exploited to encode and decode the compositional structure of the system, modeling similar physical interactions with the same parameters, resulting in a block-structured, scalable Koopman matrix.
The GNNs yield a nonlinear transformation of the state, while the effect of the control input is assumed to be linear in the transformed latent space as above, resulting in a quadratic program for control synthesis.

%
\subsubsection{Control in intrinsic Koopman eigenfunction coordinates}
\label{Sec:Control:EigenfunctionSubspace}

Eigenfunctions of the Koopman operator represent a natural set of observables,
as their temporal behavior is linear by design.
Furthermore, these \emph{intrinsic} observables are associated with global properties of the underlying system.
Koopman eigenfunctions have been used for observer design within the Koopman canonical transform (KCT)~\cite{surana2016cdc,surana2016linear} and within the Koopman reduced-order nonlinear identification and control (KRONIC) framework~\cite{kaiser2017arxiv}, which both typically yield a global bilinear representation of the underlying system.
The Koopman eigenfunctions that are used to construct a reduced description are those associated with the point spectrum of the Koopman operator and are restricted to the unactuated, autonomous dynamics.
In particular, we consider  control-affine systems of the form
\begin{align}
\dot{\bx} = {\bf f}(\bx)+\sum_{j=1}^l {\bf h}_j(\bx)u_j,\label{Eq:KRONIC:ControlAffineSystem}
\end{align}
where ${\bf h}_j(\bx)$ are the control vector fields.

The observable vector is defined as
\begin{align}
\bz \coloneqq  \mathcal{T}(\bx) = \begin{bmatrix}
\varphi_1(\bx)\\
\varphi_2(\bx)\\
\vdots\\
\varphi_p(\bx)\\
\end{bmatrix},
\end{align}
where
$\mathcal{T}:\mathbb{R}^n\rightarrow \mathbb{C}$ represents a nonlinear transformation of the state $\bx$ into eigenfunction coordinates $\varphi(\bx)$, which are associated with the unforced dynamics $\dot{\bx} = {\bf f}(\bx)$.
If the $\varphi$'s are differentiable at $\bx$, their evolution equation can be written in terms of the Lie operator~\eqref{eq:lie-to-koopman-eigenfunctions}, i.e.\ $\dot{\varphi}(\bx)=\mathcal{L}_{{\bf f}}\varphi(\bx)={\bf f}(\bx)\cdot\nabla\varphi(\bx)=\mu\varphi(\bx)$.
Hence, we obtain
\begin{equation}
\mathcal{L}_{{\bf f}}\mathcal{T}(\bx) = \boldsymbol{\Lambda}\mathcal{T}(\bx),
\end{equation}
where $\boldsymbol{\Lambda}=\mathrm{diag}(\lambda_1,\ldots,\lambda_p)$  is a matrix, with diagonal elements consisting of the eigenvalues $\lambda_i$ associated with the eigenfunctions $\varphi_j$.
Then, the dynamics for observables $\bz$ for the controlled system satisfy (see also~\eqref{Eqn:BilinearEquationForEigenfunctions})
\begin{subequations}\label{Eqn:KoopmanControl:BilinearSystem}
	\begin{align}
	\dot{\bz}
	&=  \left({\bf f}(\bx)+\sum_{j=1}^l {\bf h}_j(\bx)u_j\right) \cdot \nabla\mathcal{T} (\bx)
	=  \mathcal{L}_{{\bf f}}\mathcal{T}(\bx)  + \sum_{j=1}^l\mathcal{L}_{{\bf h}_j}\mathcal{T}(\bx) u_j \\
	&= \bLambda \bz + \sum_{j=1}^l \bb_j (\bz)u_j
	= \bLambda \bz + \bB(\bz)\bu ,\label{Eq:KRONIC:Propagator}
	\end{align}
\end{subequations}
where $\bb_j(\bz)\coloneqq   \mathcal{L}_{{\bf h}_j}\mathcal{T}(\bx)\vert_{\bx=\bC\bz}$ and $\bx=\bC\bz$ estimates the state $\bx$ from observables $\bz$.
The system~\eqref{Eqn:KoopmanControl:BilinearSystem} is bilinearizable under certain conditions, e.g. if $\mathcal{L}_{{\bf h}_j}\mathcal{T}(\bx)$ lies in the span of the set of eigenfunctions for all $j$~\cite{huang2018cdc}.  This admits a representation of the Lie operators of the control vector fields in terms of the Koopman eigenfunctions~\cite{surana2016cdc}.

For $\bu = {\bf 0}$ the dynamics~\eqref{Eq:KRONIC:Propagator} are represented in a Koopman-invariant subspace associated with the unforced system. For $\bu {\not =} {\bf 0}$ the control term $\bB(\bz)$ describes how these eigenfunctions are affected by $\bu$. The advantage of this formulation is that the dimension of the system scales with the number of eigenfunctions, and often a few dominant eigenfunctions may be sufficient to capture the principal behavior.
In general, the eigenfunctions can be identified using eDMD,
kernel-DMD or other variants, and the model~\eqref{Eq:KRONIC:Propagator} may be well represented as long as their span contains ${\bf f}(\bx)$, ${\bf h}_j(\bx)$, and $\bx$~\cite{surana2017cdc}.
Alternatively, the KRONIC framework seeks to learn eigenfunctions directly, which are sparsely  represented  in a dictionary of basis functions.
The control term $\bB(\bz)$ can either be determined from the gradient of the identified eigenfunctions for known ${\bf h}_j(\bx)$ or identified directly~\cite{kakubr2018arxiv} for unknown
${\bf h}_j(\bx)$ by similarly expanding it in terms of a dictionary.
Multiple works have since been published focusing on the direct identification of Koopman eigenfunctions~\cite{korda2020ieee,haseli2020b,Pan2020arxiv}.
In any case, validating the discovered eigenfunctions, i.e.\ ensuring that these behavior linearly as predicted by their associated eigenvalue, is critical for prediction tasks.

\subsection{Koopman-based control strategies}
\label{Sec:KoopmanControlStrategies}
In contrast to linear systems,
the control of nonlinear systems remains an engineering grand challenge.
Many real-world systems represent a significant hurdle to nonlinear control design, leading to linear approximations
which produce suboptimal solutions.
Linear Koopman-based models aim to capture the dynamics of a nonlinear system, and thus leverage linear control methods.
\begin{figure}[tb]
    \centering
    \includegraphics[width=1\textwidth]{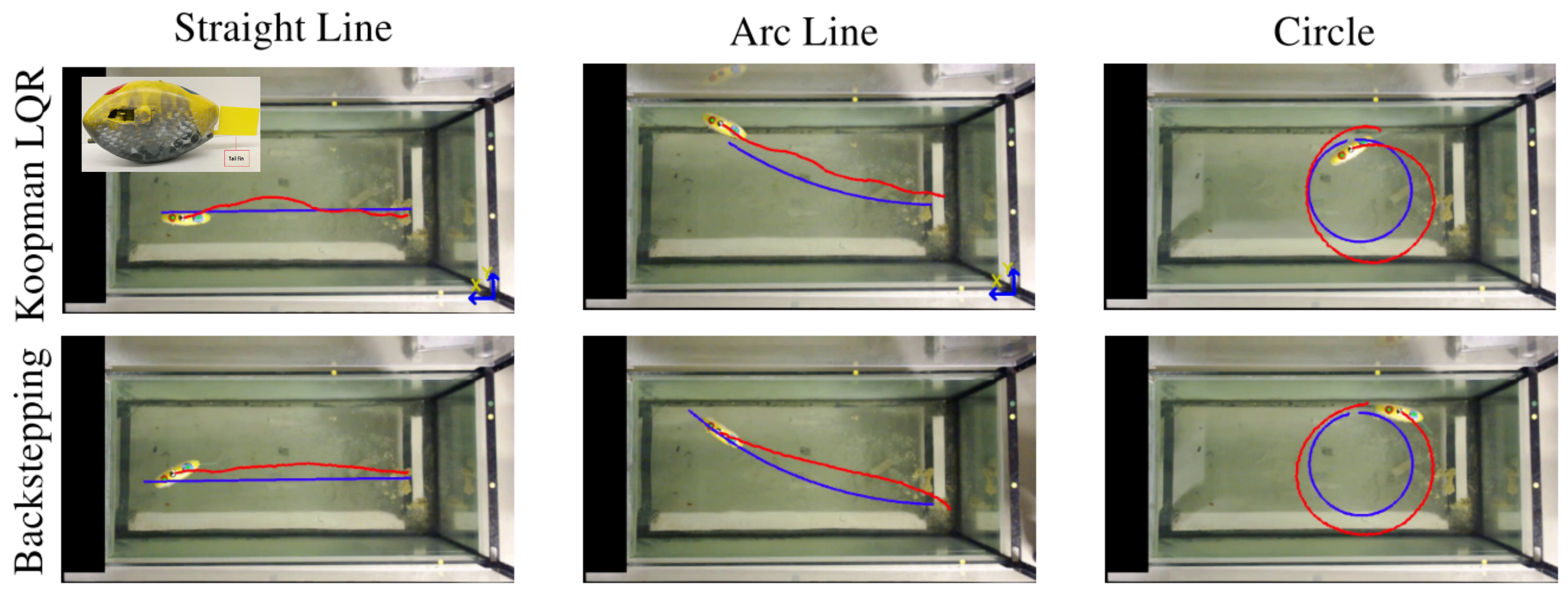}
    \caption{Tail-actuated robotic fish (top left corner) developed by the Smart Microsystems Lab at Michigan State University for trajectory tracking using Koopman LQR (top row) and backstepping (bottom row). Reference trajectory is displayed in blue, controlled trajectory in red.  \textit{Reproduced with permission, from Mamakoukas et al. 2019 Robotics: science and systems~\cite{mamakoukas2019proc}}.}
    \label{Fig:KoopmanControl:Results:FishRobot}
\end{figure}
The predictive power of linear models is increased within the Koopman operator framework
by replacing (or augmenting) the state space with nonlinear observables.
Nonlinear observables are not completely new in the context of control theory,
but appear also as the (control) Lyapunov function for stability, the value function in optimal control,
or as measurement functions in input-output linearization.
Normally the state itself is included in the set of observables.
While this is problematic from the perspective of finding a linear, closed form approximation to the Koopman operator associated with the underlying nonlinear system~\cite{brunton2016plosone}, it allows one to easily measure the state given the model output and the cost function associated with the original control problem.
If the state is not included, a measurement matrix or function mapping from the nonlinear observables to the state vector needs to estimated along with the system matrices, typically introducing additional errors.

{Koopman-based frameworks are amenable to the application of standard linear estimation and control theory and have been increasingly used in combination with optimal control~\cite{brunton2016plosone,mamakoukas2019proc,abraham2019ieee,kaiser2017arxiv,korda2016arxiv,kakubr2018arxiv,arbabi2018cdc,peitz2017arxiv,abraham2017conf}.}
In the simplest case, a linear quadratic regulator may be used, which has been successfully demonstrated in experimental robotics~\cite{mamakoukas2019proc,mamakoukas2020arxiv}, as in~\cref{Fig:KoopmanControl:Results:FishRobot}.
In addition to modifying the eigenvalues of the closed-loop system as in pole placement, the shape of Koopman eigenfunctions may be modified directly using eigenstructure assignment~\cite{hemati2017aiaa}.
Several other approaches build on a global bilinearization~\cite{surana2016cdc,kaiser2017arxiv,huang2018cdc} of the underlying system in terms of Koopman eigenfunctions (see \cref{Sec:Control:EigenfunctionSubspace}).
Under certain conditions, this may allow one to stabilize the underlying dynamical system with feedback linearization~\cite{goswami2017cdc}.
Feedback stabilization for nonlinear systems can also be achieved via a control Lyapunov function~\cite{surana2016cdc} expressed in Koopman eigenfunctions, for which the search can be formulated as a convex optimization~\cite{huang2018cdc,huang2019arxiv}.
Stabilizing solutions have also been determined by
constructing a Lyapunov function from stable Koopman models obtained by imposing stability constraints on the Koopman matrix~\cite{mamakoukas2020arxivb}.
Sequential action control~\cite{Ansari2016ieee} is an alternative to MPC in experimental robotics~\cite{abraham2017conf,abraham2019ieee}, which optimizes both the control input and the control duration. This is a convex approach applicable to control-affine systems, and as such, may also be combined with Koopman bilinearizations.

\begin{figure}[tb]
    \centering
    \includegraphics[width=1\textwidth]{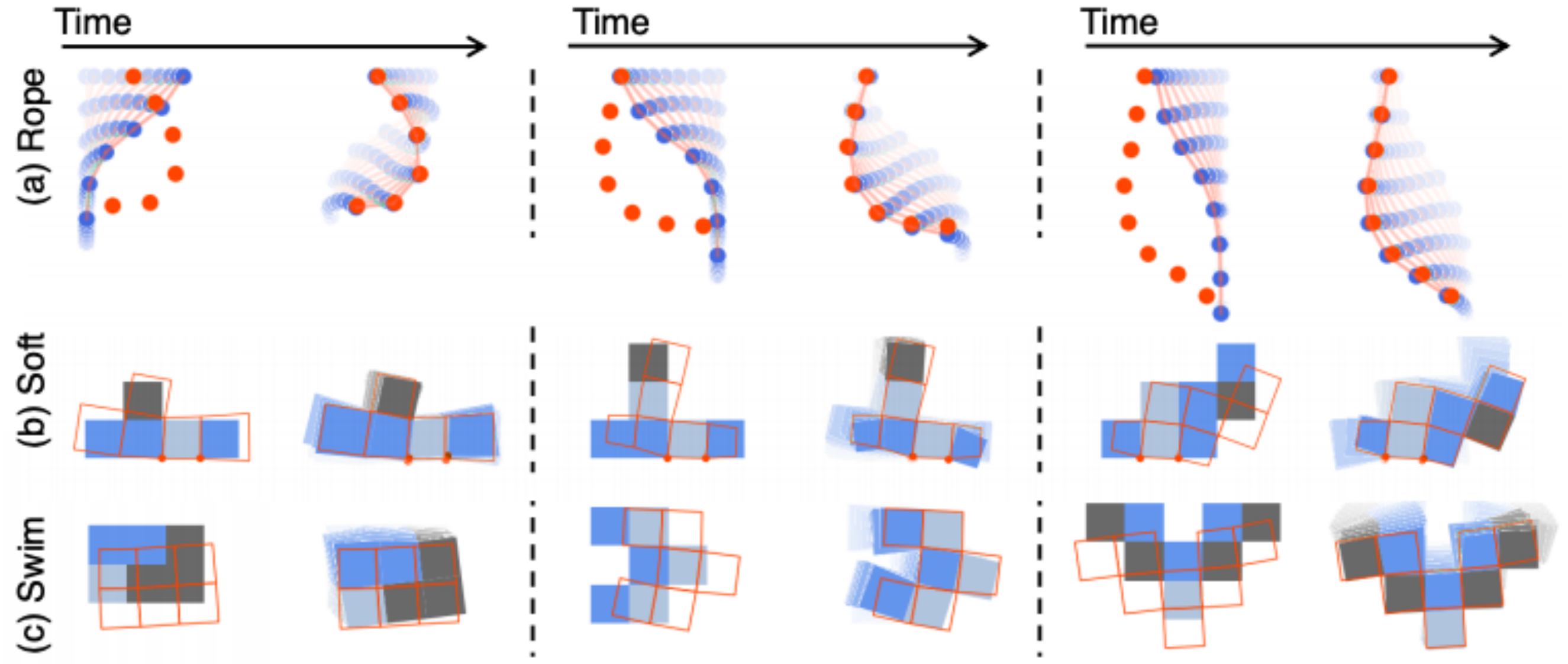}
    \caption{Model predictive control of compositional systems by combining Koopman operators with graph neural networks for an object-centric embedding. \textit{Reproduced with permission from Li et al. 2019 ICLR~\cite{li2019arxiv}}.}
    \label{Fig:KoopmanControl:Results:GNN-KO}
\end{figure}
Although the predictive power of Koopman models can be sensitive to the particular choice of observables and the training data, in particular receding-horizon methods such as MPC (see \cref{Sec:MPC}) provide a robust control framework that systematically compensate for model uncertainties by continually updating the optimization problem and taking into account new measurements.
{Incorporating Koopman-based models into MPC was first introduced in Korda \& Mezi\'c~\cite{korda2016arxiv} using eDMDc, and has subsequently been used for control in Koopman eigenfunction coordinates~\cite{kaiser2017arxiv,kakubr2018arxiv} and interpolated Koopman models~\cite{peitz2018arxiv,peitz2020arxiv} extending to continuous control inputs.}
The latter approach combines the data efficiency of estimating autonomous Koopman models associated with discrete control inputs with continuous control, which is achieved by relaxation of a switching control approach.
Convergence can be proven, albeit assuming an infinite-data limit and an infinite number of basis functions.
In general, guarantees on optimality, stability, and robustness of the controlled dynamical system remain limited. Of note are Koopman Lyapunov-based MPC that guarantees closed-loop stability and controller feasability~\cite{narasingam2019,son_narasingam2020}.
More recently, deep learning architectures have been increasingly used to represent the nonlinear mapping into the observable space.
These have also been combined with MPC with promising results~\cite{li2019arxiv} (see also \cref{Fig:KoopmanControl:Results:GNN-KO}).

The success of Koopman-based MPC belies the practical difficulties of approximating the Koopman operator.
On the one hand, there exist only a few systems with a known Koopman-invariant subspace and eigenfunctions, on which models could be analyzed and evaluated.
On the other hand, even the linearity property of the Koopman eigenfunctions associated with the model is rarely validated and MPC exhibits an impressive robustness when handling models with very limited predictive power.
Thus, there have been impressive practical successes, but there are open gaps in understanding how well the Koopman operator is actually approximated and what possible error bounds can be given.
In addition to theoretical advances, Koopman-based control has been increasingly applied in real-world problems such as power grids~\cite{korda2018arxiv,netto2018tps}, 
high-dimensional fluid flows~\cite{ma2011tcfd,arbabi2018cdc,peitz2018feedback,peitz2020arxiv},
biological systems~\cite{hasnain2019arxivb}, chemical processes~\cite{narasingam2019,son_narasingam2020}, human-machine systems~\cite{broad2017proc}, and experimental robotic systems~\cite{abraham2019ieee,bruder2019proc,mamakoukas2019proc,bruder2020arxiv,folkestad2020arxiv}.

\begin{figure}[tb]
    \centering
    \begin{overpic}[width=1\textwidth]{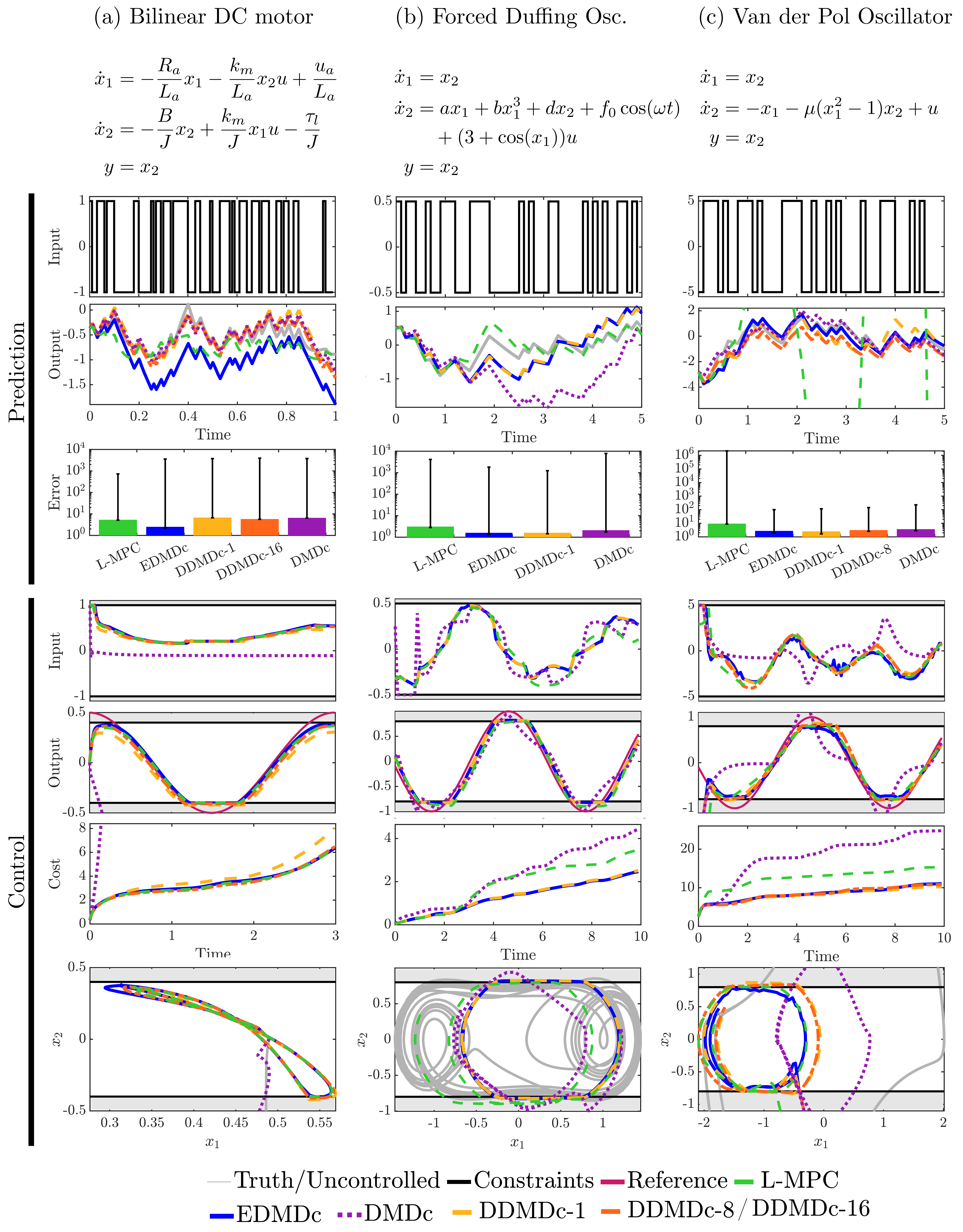}
        \put(3.5,31){{\rotatebox{90}{$J$}}}
    \end{overpic}
    \caption{Reference tracking using MPC combined with different Koopman operator models illustrated for three dynamical systems.}
    \label{Fig:KoopmanControl:Results:NumericalExamples}
    \vspace{-.2in}
\end{figure}

\begin{figure}[tb]
    \centering
    \includegraphics[width=1\textwidth]{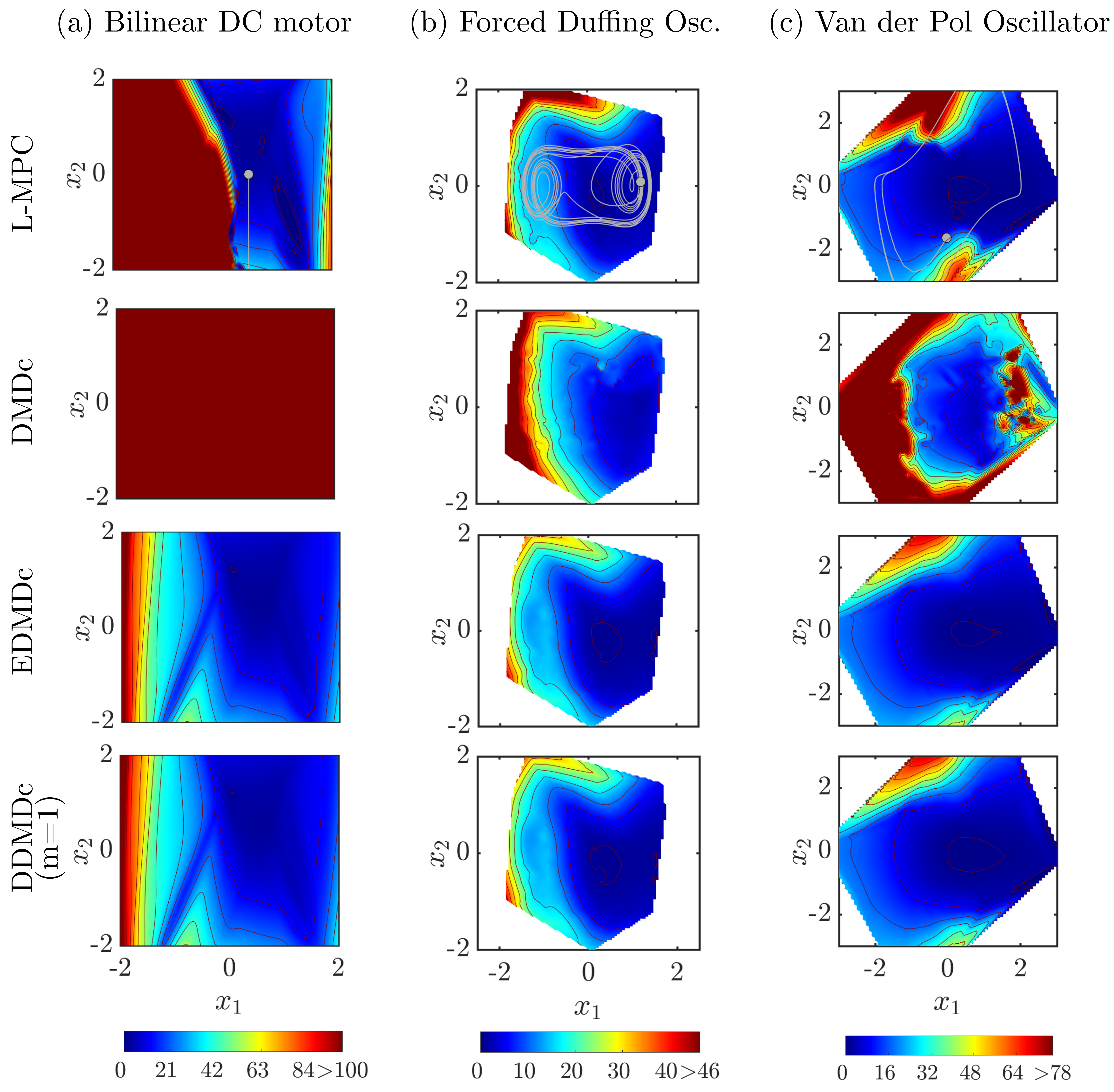}
    \caption{N-step controllable set for the reference tracking problems in \cref{Fig:KoopmanControl:Results:NumericalExamples}. The trajectory of the unforced system is displayed in gray in the top row. The same initial condition (gray circle) is used in \cref{Fig:KoopmanControl:Results:NumericalExamples}.}
    \label{Fig:KoopmanControl:Results:ControlDomain}
    \vspace{-.2in}
\end{figure}

Multiple variants of Koopman-based control have been demonstrated in numerical and experimental settings.
Here, we illustrate model predictive control (MPC, see \cref{Sec:MPC} for details) combined with different Koopman-based models for several nonlinear systems including a bilinear DC motor (see also~\cite{korda2016arxiv}), the forced Duffing oscillator, and the van der Pol oscillator.
An overview of the parameters, models and controls are presented in Tab.~\ref{Tab:KoopmanControl:Results:NumericalExamples} with results summarized in \cref{Fig:KoopmanControl:Results:NumericalExamples} and~\cref{Fig:KoopmanControl:Results:ControlDomain}.
\begin{table}[htb]
\centering
\begin{tabular}{ lrL{2.9cm}L{2.7cm}L{2.7cm}  }
\toprule
& & \multicolumn{3}{c}{Examples} \\
\cmidrule(r){3-5}
 & & \textbf{Bilinear} & \textbf{Forced Duffing} & \textbf{Van der Pol}\\
 & & \textbf{DC motor} & \textbf{oscillator} & \textbf{oscillator}\\
 & & $L_a=0.314$, & $a=1$, $b=-1$, & $\mu=2$\\
 & & $R_a=12.345$, & $d=-0.3$, $f_0 = 0.5$, & \\
 & & $k_m=0.253$, & $\omega = 1.2$ & \\
 & & $J=0.00441$, $B=0.00732$, $\tau_l=1.47$, $u_a=60$ & & \\
\midrule
\multirow{3}{*}{\rotatebox{90}{Training}} &  Timestep & $\Delta t = 0.01$ & $\Delta t = 0.1$ & $\Delta t = 0.1$\\
& Domain & $\bx\in[-1,1]^2$ & $\bx\in[-1.5,-1]^2$ & $\bx\in[-2,2]^2$ \\
& Input & $u\in[-1,1]$ & $u\in[-0.5,0.5]$ & $u\in[-5,5]$ \\
\midrule
\multirow{8}{*}{\rotatebox{90}{Model}} & DMDc & $\bz=\bx$ & $\bz=\bx$ & $\bz=\bx$\\
& & $p=3$ & $p=3$ & $p=3$\\
& DDMDc & $\bz = \bg^1(y,u)$ & $\bz = \bg^1(y,u)$ & $\bz = \bg^1(y,u)$\\
& & $p=3$ & $p=3$ & $p=3$\\
& & $\bz = \bg^{16}(y,u)$ & --- & $\bz = \bg^{8}(y,u)$\\
& & $p=33$ & & $p=17$\\
& eDMDc & $\bz = [\bg^1,\boldsymbol{\xi}]$ & $\bz = [\bg^1,\boldsymbol{\xi}]$ & $\bz = [\bg^1,\boldsymbol{\xi}]$\\
&       & $p=103$ & $p=103$ & $p=103$\\ 
 \hline
\multirow{5}{*}{\rotatebox{90}{Control}} & Reference & $r(t)=\frac{1}{2}\cos(\frac{2}{3}\pi t)$ & $r(t)=-\sin(t)$ & $r(t)=-\sin(t)$\\
& Horizon & $T=0.1$ & $T=1$ & $T=3$\\
& Weights & $Q=1$, $R=0.1$ & $Q=1$, $R=0.1$ & $Q=1$, $R=0.01$\\
& Constraints & $-0.4\leq x_2 \leq 0.4$ & $-0.8\leq x_2\leq 0.8$ & $-0.8\leq x_2\leq 0.8$\\
&  & $-1\leq u\leq 1$ & $-0.5\leq u\leq 0.5$ & $-5\leq u\leq 5$\\
\bottomrule
\end{tabular}
\caption{Parameters for numerical control examples.}
\label{Tab:KoopmanControl:Results:NumericalExamples}
\end{table}
Consider the single-input single-output (SISO) system
\begin{subequations}
	\begin{align}
    \dot{\bx} &= {\bf f}(\bx,u,t)\\
    y &= \bC \bx
    \end{align}
\end{subequations}
to be controlled using MPC to track a time-varying reference signal $r(t)$.
The cost function is defined as
\begin{equation}
    J = ||y_N-r_N||^2_{\bQ} +
    \sum\limits_{k=0}^{N-1} ||y_k-r_k||^2_{\bQ} + R||u_k||^2,
\end{equation}\label{Eqn:Control:Numerical:CostFun}
where the first term represents the terminal cost and the subscript $k$ corresponds to the $k$th timestep $t_k\coloneqq k\Delta t$.
The control sequence $\{u_0,u_1,\ldots,u_{N-1} \}$ over the $N-$step horizon is determined by minimizing the cost function subject to the model dynamics and state
and control constraints, $y_{min}\leq y\leq y_{max}$ and $u_{min}\leq u\leq u_{max}$, respectively.
These results are compared with a local linearization of the system dynamics (referred to as LMPC in the following).

While DMDc relies on full-state measurements,
only the output $y=\bC\bx$ is evaluated in the cost function~\eqref{Eqn:Control:Numerical:CostFun}.
For the remaining models, time delay coordinates are employed and we define the time-delayed observable vector as
\begin{equation}
    \bg^m(y,u)=[y_k,(y,u)_{k-1},\ldots,(y,u)_{k-m}]^T
\end{equation}
with $(y,u)_{k}\coloneqq  [y_k,u_k]$ and where $m$ denotes the number of delays.
The DDMDc model, which is short for DMDc applied to delay coordinates, is based on time delay coordinates discussed in \cref{Sec:Control:GeneralizationsDMDc}. However, past control values are also included in the observable vector $\bg^m(y,u)$.
The observable vector associated with the eDMDc models is defined as
$\bz = [\bg^m(y,u),\boldsymbol{\xi}(y,u)]^T$ (analogous to~\cite{korda2016arxiv}), so that the output $y$ itself is a coordinate of the observable.  The variable $\boldsymbol{\xi}(y,u)$ represents a particular choice of lifted nonlinear observables. Here, we employ thin plate spline radial basis functions, $\xi_j = r^2\ln ||\bg^m(y,u)-\bc_j||$ with $r=1$ and $100$ randomly sampled centers within the considered domain.
If the output is not included as an observable, an estimate of the mapping from $\boldsymbol{\xi}(y,u)$ to the output $y$ is required in order to evaluate the cost function. By incorporating the output in the model this can be avoided; however, this also restricts the model~\cite{brunton2016plosone}.

Training data for all examples is collected from $10^2$ initial conditions uniformly distributed within the examined domain stated in Tab.~\ref{Tab:KoopmanControl:Results:NumericalExamples}.
A different random, binary actuation signal within the range of the input constraints is applied to sample each trajectory in the training set. Each trajectory consists of $10^3$ samples, so that in total the training set consists of $10^5$ samples.
Data efficiency in system identification is an important topic, as only limited measurements may be available to train models and each model may have different requirements. For instance, as little as $10$ samples may be sufficient to train a DMDc model that successfully stabilizes an unstable fixed point using MPC (see, e.g., Fig.7 in~\cite{kaiser2017arxivb}).

The prediction and control performance is summarized in \cref{Fig:KoopmanControl:Results:NumericalExamples}.
The first part illustrates the prediction accuracy for a single trajectory and shows the mean relative prediction error $e = \frac{1}{N} \sum_{k=1}^N\,(y_k-\hat{y}_k)/y_k$ with $N=T/\delta t$,
where $y$ and $\hat{y}$ are the true and predicted values, respectively.
The bar displays the median value of $e$ for each model evaluated over $300$ random initial conditions, and the black lines indicate extremal values.
The second part shows the control performance for the same initial condition. The cost is evaluated by~\eqref{Eqn:Control:Numerical:CostFun} {and shown in the sixth row of Fig.~\ref{Fig:KoopmanControl:Results:NumericalExamples}}.
We note that the model based on a local linearization (LMPC) can be used to successfully control the bilinear DC motor without constraint violations or infeasibility issues, in contrast to~\cite{korda2016arxiv}, by using a smaller prediction horizon for MPC. A significantly larger horizon would exceed the predictability of the model and result in infeasibility and performance issues.
While the prediction accuracy may be comparable with other models, DMDc-MPC fails for the DC motor and is inferior in the other examples for tracking a time-varying reference signal. However, it may be possible to stabilize an unstable fixed point using DMDc in certain cases, despite having limited predictive power~\cite{kaiser2017arxivb}.
DDMDc$-m$ models, where $m$ denotes the number of delays, achieves comparable predictive and control capabilities as eDMDc; intriguingly, even with only one delay $m=1$ for the Duffing and van der Pol systems.
MPC is surprisingly robust with respect to model errors. Even in cases with limited predictive power (compare, e.g., the early divergence of LMPC for the van der Pol oscillator), the system is successfully controlled with MPC.
The robustness of MPC w.r.t. modeling inaccuracies may hide the fact that the model may actually be a poor representation of the underlying dynamics or its Koopman operator.

In \cref{Fig:KoopmanControl:Results:ControlDomain}, the control performance is assessed over the domain of the training data, as the prediction performance can vary significantly depending on the initial condition.
Each initial condition is color-coded with the cost  $J$~\eqref{Eqn:Control:Numerical:CostFun}, which is evaluated
for a duration of $1$, $4$, and $3$ time units for the three examples (left to right), respectively.
DDMDc$-1$ and eDMDc, both using $m=1$ delays, are nearly indistinguishable, while the model dimension is $p=3$ for the first and $p=103$ for the latter.
This raises the important, and so far relatively unanswered question, of which set of basis functions should be employed to train a Koopman model. Time delay coordinates perform extremely well for prediction and control in a variety of systems and have the benefit of well-developed theoretical foundations. More sophisticated basis functions may perform similarly well, but it is not clear that they warrant the increased design effort given the comparable performance.

\subsection{Stability}
Stability is central in the analysis and control of real-world systems.
The Lyapunov function plays an essential role in global stability analysis and control design of nonlinear systems.
Finding or constructing a Lyapunov function poses a significant challenge for general nonlinear systems.
Importantly, spectral properties of the Koopman operator have been associated  with geometric properties such as Lyapunov functions and measures, the contracting metric, and isostables which are used to study global stability of the underlying nonlinear system.
For instance, building on~\cite{vaidya2006ccc}, it has been shown~\cite{mauroy2013,mauroy2013cdc} that Lyapunov's second method implicitly uses an operator-theoretic approach and Lyapunov functions represent special observables advanced through the Koopman operator.

Consider a dynamical system $\dot{\bx} = {\bf f}(\bx)$ with fixed point $\bx^{*}$.  A Lyapunov function ${V:\mathbb{R}^n\rightarrow \mathbb{R}}$ must satisfy
\begin{equation}\label{Eq:LyapunovFunction}
\dot{V}(\bx) = \nabla V\cdot {\bf f}(\bx)\leq 0\;\text{for all}\;\bx{\not =}\bx^{*}.
\end{equation}
If $\dot{V}(\bx) \leq 0$ the system is asymptotically stable in the sense of Lyapunov.
There is also an equivalent definition of the Lyapunov function for discrete-time systems.
The dynamics of the Lyapunov equation~\eqref{Eq:LyapunovFunction}
can be formulated directly in terms of the generator of the Koopman operator family acting on a nonnegative observable~\eqref{Eq:GeneratorControl:ChainRule} (here without input $\bu={\bf 0}$).
Thus, the Lyapunov function decays everywhere under action of the Lie
operator and is related to its spectral properties and those of the Koopman operator semigroup.
The global stability of a fixed point can be established through the existence of a set of $C^1$ eigenfunctions of the Koopman operator associated with the eigenvalues of the Jacobian of the vector field with the Koopman eigenfunctions being used to define a Lyapunov function and contracting metric~\cite{mauroy2013cdc}.
A corresponding numerical scheme based on a Taylor expansion may then be used to compute stability properties, including the domain of attraction. These methods have been extended to characterize the global stability properties of hyperbolic fixed points, limit cycles, and non-analytical eigenfunctions~\cite{mauroy2016ieeetac}.

\subsection{Observability and controllability}
Observability and controllability play a crucial role in the design of sensor-based estimation and control.
However, they have limited validity when applied to linearized nonlinear systems.
There exist analogous (local) controllability and observability criteria for nonlinear systems using Lie derivatives~\cite{hermann1977ieee,khalil1996book}.
However, these criteria are typically restricted to a specific class of systems and their evaluation remains challenging for even low-dimensional systems, becoming computationally intractable for high-dimensional systems.
Operator-theoretic approaches provide new opportunities to assess observability and controllability of high-dimensional nonlinear systems using linear techniques and obtaining corresponding estimates from given data.
Generally, two states are indistinguishable if their future output sequences are identical. Thus, a nonlinear system is considered nonlinearly observable if any pair of states is distinguishable.
Common criteria may be divided into two classes: rank conditions and Gramians~\cite{brunton2019data}.
Rank conditions yield a binary decision, while Gramians are used to examine the degree of observability and controllability. The latter is also used for balanced truncation to obtain reduced-order models that balance a joint observability and controllability measure.
In contrast, corresponding metrics for nonlinear systems rely on the derivation of its linear equivalent using Lie derivatives, which can be challenging.
Reformulating a nonlinear system into a linear system allows one to apply linear criteria in a straightforward manner.  Thus, Koopman operator-theoretic techniques have been increasingly studied in the context of observability and controllability of nonlinear systems.

An immediate extension of linear rank conditions to nonlinear systems is achieved by applying the linear criteria to the representation in the lifted state space.
In~\cite{surana2016linear}, nonlinear observability is evaluated in a Koopman-based observer, which provides a linear representation of the underlying system in Koopman eigenfunction coordinates (see also \cref{Sec:Control:EigenfunctionSubspace}).
The underlying nonlinear system is then considered nonlinearly observable if the pair $(\bA,\bC_H)$ is observable, which can be determined via the rank condition of the corresponding observability matrix. These ideas have been applied to study pedestrian crowd flow~\cite{benosman2017ifac},
extended further to input-output nonlinear systems
resulting in bilinear or Lipschitz formulations~\cite{surana2016cdc}, and used
to compute controllability and reachability~\cite{goswami2017cdc}.
The observability and controllability Gramians can also be computed in the lifted observable space~\cite{surana2016linear}.
In this case, the observability/controllability of the underlying system is related to the observability/controllability of the observables. The underlying assumption is that the state, input, and output are representable in terms of a few Koopman eigenfunctions.
More recently, a deep learning DMD framework~\cite{yeung2017arxiv} has been extended to incorporate Koopman Gramians to maximize internal subsystem observability and disturbance rejection from unwanted noise from other subsystems~\cite{liu2018acc}.

As shown in the following example, if a system can be represented in a Koopman-invariant subspace, linear criteria applied to the Koopman representation can be equivalent to the corresponding nonlinear criteria applied to analyze the underlying the nonlinear system.
However, a Koopman-invariant subspace is rarely obtained and open questions remain, e.g., whether an approximate Koopman-invariant subspace may be sufficient and how observability and controllability estimates depend on the specific choice of observable functions.

\paragraph{Example: Nonlinear system with single fixed point and a slow manifold}
We examine the relationship between controllability of the original nonlinear system and the Koopman system
using the corresponding controllability criteria.
We consider the system examined in \cref{Sec:SlowManifold}, with the addition of control:
\begin{equation}\label{Eqn:SlowManifold:Control}
\dot{\bx}
= {\bf f}_0(\bx) + {\bf f}_1(\bx) u
= \begin{bmatrix}
 \mu x_1\\
 \lambda(x_2-x_1^2)
\end{bmatrix}
+ \begin{bmatrix}
 0\\ 1
\end{bmatrix}u.
\end{equation}
The additional control input $u$ actuates only the second state $x_2$.
The system becomes unstable for either $\lambda>0$ or $\mu>0$.
Choosing observables ${\bf y} = (x_1,x_2,x_1^2)^T$, the system admits a fully linear, closed description in a three-dimensional Koopman invariant subspace:
\begin{equation}\label{Eqn:SlowManifold:Control:Koopman}
	\frac{d}{dt}\begin{bmatrix} y_1\\ y_2\\ y_3\end{bmatrix} =
	{\bf A}\begin{bmatrix} y_1\\ y_2\\ y_3\end{bmatrix} + {\bf B} u
	=
	\begin{bmatrix} \mu & 0 & 0 \\ 0& \lambda& - \lambda\\ 0 & 0 & 2\mu\end{bmatrix}
	\begin{bmatrix}y_1\\ y_2\\ y_3\end{bmatrix}+
	\begin{bmatrix}0\\1\\0\end{bmatrix}u
\end{equation}
with constant control vector ${\bf B}$.

The controllability for a control-affine system $\dot{\bx}=\mathbf{f}_0(\bx) + \sum_{j=1}^q \mathbf{f}_j(\bx)u_j$ can be assessed through the following theorem.
We define the Lie derivative of $\mathbf{g}$ with respect to $\mathbf{f}$ as $\bL_{\mathbf{f}} \mathbf{g}\coloneqq \nabla \mathbf{g}\cdot \mathbf{f}$, the Lie bracket as $[\mathbf{f},\mathbf{g}]\coloneqq \nabla \mathbf{g}\cdot \mathbf{f}-\nabla \mathbf{f}\cdot \mathbf{g}$, and the recursive Lie bracket as $\operatorname{ad}_{\mathbf{f}} \mathbf{g} = [\mathbf{f},\mathbf{g}] = \nabla \mathbf{g}\cdot \mathbf{f} - \nabla \mathbf{f}\cdot \mathbf{g}$.
Hunt's theorem (1982)~\cite{hunt1982ieee} states that a nonlinear system is (locally) controllable if there exists an index $k$ such that
${\mathcal{C} = [\mathbf{f}_0,\mathbf{f}_1,\ldots,\mathbf{f}_q, \operatorname{ad}_{\mathbf{f}_0}\mathbf{f}_1,\ldots, \operatorname{ad}_{\mathbf{f}_0}\mathbf{f}_q, \ldots, \operatorname{ad}_{\mathbf{f}_0^k}\mathbf{f}_1,\ldots, \operatorname{ad}_{\mathbf{f}_0^k}\mathbf{f}_q]}$
has $n$ linearly independent columns.
The controllability matrix for the nonlinear system~\eqref{Eqn:SlowManifold:Control} is given by
\begin{equation}
\mathcal{C} =
\begin{bmatrix}
{\bf f}_0& \operatorname{ad}_{\mathbf{f}_0}\mathbf{f}_1 & \ldots & \operatorname{ad}_{\mathbf{f}_0^k}\mathbf{f}_1
\end{bmatrix} =
\begin{bmatrix}
\begin{bmatrix} 0\\1 \end{bmatrix} &
\begin{bmatrix} 0\\-\lambda \end{bmatrix} &
\begin{bmatrix} 0\\\lambda^2 \end{bmatrix} &
\ldots &
\begin{bmatrix} 0\\(-1)^k \lambda^k \end{bmatrix}
\end{bmatrix},
\end{equation}
which has rank 1 for any $k$. The system is uncontrollable in $x_1$ direction.

Analogously, we can construct the linear controllability matrix for the Koopman system~\eqref{Eqn:SlowManifold:Control:Koopman}:
\begin{equation}
\mathcal{C} = \begin{bmatrix}
\bB & \bA \bB & \bA^2\bB
\end{bmatrix}
=
\begin{bmatrix}
0 & 0 & 0\\
1 & \lambda & \lambda^2\\
0 & 0 & 0
\end{bmatrix},
\end{equation}
which is also of rank 1.
Another very useful test is the Popov-Belevitch-Hautus (PBH) test connecting controllability to a relationship between the columns of $\bB$ and the eigenspace of $\bA$. The PBH test states that the pair $(\bA,\bB)$ is controllable iff the column rank of $[(\bA - \alpha\bI)\;\bB]$ is equal to $n$ for all eigenvalues $\alpha(\bA)\in\mathbb{C}$.
This test confirms that the rank is $n=3$ for the eigenvalue $\lambda$ and 1 for both eigenvalues $\mu$ and $2\mu$ associated with~\eqref{Eqn:SlowManifold:Control:Koopman}.
Note that~\eqref{Eqn:SlowManifold:Control}  can still be stabilized as long as $\mu$ is stable.

\subsection{Sensor and actuator placement}
Optimizing sensor and actuator locations for data collection and decision-making is a crucial and challenging task in any real-world application.
Optimal placement of sensors and actuators amounts to an intractable brute force search, as these represent combinatorial problems suffering from the curse of dimensionality.
Rigorous optimization remains an open challenge for even linear problems; thus, approaches generally rely on heuristics.
In this realm, sparse sampling and greedy selection techniques, such as gappy POD~\cite{everson1995karhunen} and DEIM~\cite{drmac2016siam}, and more recently
sparsity-promoting algorithms, such as compressed sensing~\cite{candes2006bieeetit,donoho2006ieeetit},  have played an increasingly
important role in the context of sensor/actuator selection~\cite{manohar2017csm}.
These methods rely on exploiting the ubiquitous low-rank structure of data typically found in high-dimensional systems.
Operator-theoretic methods fit well into this perspective, as they are able to provide a tailored feature basis capturing global, dynamically persistent structures from data.  This can be combined with existing heuristic selection/sampling methods which have been demonstrated on a range of practical applications, such as environmental controls in buildings~\cite{vaidya2012jmaa,fontanini2016be,sharma2017ibpsa}, fluid dynamics~\cite{kaiser2018jcp,manohar2017arxiv}, and biology~\cite{hasnain2019arxiv}.

Importantly, operator-theoretic methods generalize to nonlinear systems, e.g. for estimating nonlinear observability and controllability,
providing practical means to systematically exploit nonlinear system properties for sensor/actuator placement within a linear framework.
As already noted, controllability and observability Gramians can be generalized for nonlinear systems based on the
Koopman and Perron-Frobenius operators (or Liouville and adjoint Liouville operators as their continuous-time counterpart) and subsequently used for sensor and actuator placement by maximizing the support (or the L2 norm) of the finite-time Gramians.
In~\cite{vaidya2012jmaa,sinha2016jmaa}, set-oriented methods have been utilized to approximate the (adjoint) Lie operators, i.e.\ the domain is discretized into small cells and few cells are selected as sensors/actuators, and the location can be optimized by solving a convex optimization problem. A greedy heuristic approach based on these ideas is proposed in~\cite{fontanini2016be}, which further investigates different criteria such as maximizing the sensing volume (sensor coverage), response time and accuracy (relative measure transported to the sensor in finite time) and incorporating spatial constraints.
The framework has been further extended to incorporate uncertainty~\cite{sharma2017ibpsa}.
More recently, observability Gramians based on the Koopman operator have been used to inform optimal sensor placement in the presence of noisy and temporally sparse data~\cite{hasnain2019arxivb}.
Sensor placement is here facilitated by maximizing the finite-time output energy in the lifted space.
More broadly, balanced truncation and model reduction has recently been used for efficient placement of sensors and actuators to simultaneously maximize the controllability and observability Gramians~\cite{manohar2018arxivb},
which may be promising direction to extend to nonlinear systems by using Gramians in the lifted observable space.

\section{Discussion and outlook}\label{Sec:Discussion}

In this review, we have explored the use of Koopman operator theory to characterize and approximate nonlinear dynamical systems in a linear framework.
Finding linear representations of nonlinear systems has a broad appeal, because of the potential to enable advanced prediction, estimation, and control using powerful and standardized linear techniques.
However, there appears to be a \emph{conservation of difficulty} with Koopman theory, where nonlinear dynamics are traded for linear but infinite dimensional dynamics, with their own associated challenges.
Thus, a central focus of modern Koopman analysis is to find a finite set of nonlinear measurement functions, or coordinate transformations, in which the dynamics appear linear, and the span of which may be used to approximate other quantities of interest.
In this way, Koopman theory follows in the tradition of centuries of work in mathematical physics to find effective coordinate systems to simplify the dynamics.
Now, these approaches are augmented with \emph{data-driven} modeling, leveraging the wealth of measurement data available for modern systems of interest, along with high-performance computing and machine learning architectures to process this data.

This review has also described several leading data-driven algorithms to approximate the Koopman operator, based on the dynamic mode decomposition (DMD) and variants.
DMD has several advantages that have made it widely adopted in a variety of disciplines.
First, DMD is purely data-driven, as it does not require governing equations, making it equally applicable to fluid dynamics and neuroscience.
In addition, the DMD algorithm is formulated in terms of simple linear algebra, so that it is highly extensible, for example to include control.
For these reasons, DMD has been applied broadly in several fields.
Although it is possible to apply DMD to most time-series data, researchers are still developing diagnostics to assess when the DMD approximation is valid and to what extent it is useful for prediction and control.

Despite the widespread use of DMD, there are still considerable challenges associated with applying DMD to strongly nonlinear systems, as linear measurements are often insufficient to span a Koopman-invariant subspace.
 Although significant progress has been made connecting DMD to nonlinear systems~\cite{williams2015jcd}, choosing nonlinear measurements to augment the DMD regression is still not an exact science.
Identifying measurement subspaces that remain closed under the Koopman operator is an ongoing challenge~\cite{brunton2016plosone}.
In the past decade, several approaches have been proposed to extend DMD, including with nonlinear measurements, time-delayed observations, and deep neural networks to learn nonlinear coordinate transformations.
These approaches have had varying degrees of success for systems characterized by transients and intermittent phenomena, as well as systems with broadband frequency content associated with a continuous eigenvalue spectrum.
It is expected that neural network representations of dynamical systems, and Koopman embeddings in particular, will remain a growing area of interest in data-driven dynamics.
Combining the representational power of deep learning with the elegance and simplicity of Koopman embeddings has the potential to transform the analysis and control of complex systems.

One of the most exciting areas of development in modern Koopman theory is around its use for the control of nonlinear systems~\cite{otto2021koopman}.
Koopman-based control is an area of intense focus, as even a small amount of nonlinearity often makes control quite challenging, and an alternative linear representation may enable dramatic improvements to robust control performance with relatively standard linear control techniques.
Model predictive control using DMD and Koopman approximations is particularly interesting, and has been applied to several challenging systems in recent years.
However, there are still open questions about how much of this impressive performance is due to the incredible robustness of MPC as opposed to the improved predictive capabilities of approximate Koopman models.
The goal of more effective nonlinear control will likely continue to drive applied Koopman research.

Despite the incredible promise of Koopman operator theory, there are still several significant challenges that face the field.
In a sense, these challenges guarantee that there will be interesting and important work in this field for decades to come.
There are still open questions about how the choice of observables impacts what can be observed in the spectrum.
Similarly, there is no \emph{return path} from a Koopman representation back to the governing nonlinear equations.
More generally, there is still little known about how properties about the nonlinear system $\dot{\bx}=\mathbf{f}(\bx)$ carry over to the Koopman operator, and vice versa.
For example, relationships between symmetries in the nonlinear dynamics and how they manifest in the Koopman spectrum is at the cusp of current knowledge.
A better understanding about how to factor out symmetries and generalize the notion of conserved quantities and symmetries to \emph{near-invariances} may provide insights into considerably more complex systems.
The community is still just now beginning to rectify the local and global perspectives on Koopman.
In addition, most theory has been formulated for ODEs, and connections to spatiotemporal systems are still being developed.
Finally, there are open questions around whether or not there is a Heisenberg uncertainty principle for the Koopman operator.

In addition to these theoretical open questions, there are several areas of applied research that are under active development.
Applied Koopman theory, driven by the dynamic mode decomposition, has largely been applied within the fluid dynamics community.
Koopman's original theory was a critical piece in resolving Boltzmann's ergodic hypothesis for gas dynamics.
It is likely that Koopman theory will provide similar enabling technology for the characterization of fully turbulent fluid dynamics, which have defied macroscopic coarse-graining for over a century.
Unifying algorithmic innovations is also an ongoing challenge, as it is not always obvious which algorithm should be used for a particular problem.
Several open source software libraries are being developed to ease this burden, including
\vspace{.1in}
\begin{itemize}
    \item PyDMD (\url{https://github.com/mathLab/PyDMD});
    \item PyKoopman (\url{https://github.com/dynamicslab/pykoopman}); {
    \item Data-driven dynamical systems toolbox (\url{https://github.com/sklus/d3s});
    \item deeptime (\url{http://https://github.com/deeptime-ml/deeptime}).}
\end{itemize}
\vspace{.1in}
As a parting thought, it is important to note that nonlinearity is one of the most fascinating features of dynamical systems, providing a wealth of rich phenomena.
In a sense, nonlinearity provides an amazing way to parameterize dynamics in an extremely compact and efficient manner.
For example, with a small amount of cubic nonlinearity in the Duffing equation, it is possible to parameterize continuous frequency shifts and harmonics, whereas this may require a comparatively complex Koopman parameterization.
Realistically, the future of dynamical systems will utilize the strengths of both traditional nonlinear, geometric representations along with emerging linear, operator theoretic representations.

\section*{Acknowledgements}
We would like to thank Igor Mezi\'{c} for early discussions and encouragement about this review.
We also thank Shervin Bagheri, Bing Brunton, Bethany Lusch, { Ryan Mohr}, Frank Noe, Josh Proctor, Clancy Rowley, and Peter Schmid for many fruitful discussions on DMD, Koopman theory, and control.

\newcommand{\noopsort}[1]{}


\begin{thebibliography}{100}

\bibitem{ablowitz1974inverse}
{\sc M.~J. Ablowitz, D.~J. Kaup, A.~C. Newell, and H.~Segur}, {\em The inverse
  scattering transform-{F}ourier analysis for nonlinear problems}, Studies in
  Applied Mathematics, 53 (1974), pp.~249--315.

\bibitem{ist}
{\sc M.~J. Ablowitz and H.~Segur}, {\em Solitons and the inverse scattering
  transform}, vol.~4, Siam, 1981.

\bibitem{abraham2017conf}
{\sc I.~Abraham, G.~de~la Torre, and T.~Murphey}, {\em Model-based control
  using {K}oopman operators}, in Proceedings of Robotics: Science and Systems,
  Cambridge, Massachusetts, July 2017,
  \url{https://doi.org/10.15607/RSS.2017.XIII.052}.

\bibitem{abraham2019ieee}
{\sc I.~Abraham and T.~D. Murphey}, {\em Active learning of dynamics for
  data-driven control using {K}oopman operators}, IEEE Transactions on
  Robotics, 35 (2019), pp.~1071--1083.

\bibitem{marsdenmtaa}
{\sc R.~Abraham, J.~E. Marsden, and T.~Ratiu}, {\em Manifolds, Tensor Analysis,
  and Applications}, vol.~75 of Applied Mathematical Sciences, Springer-Verlag,
  1988.

\bibitem{agrawal2016book}
{\sc M.~Agrawal, S.~Vidyashankar, and K.~Huang}, {\em On-chip implementation of
  {ECoG} signal data decoding in brain-computer interface}, in Mixed-Signal
  Testing Workshop (IMSTW), 2016 IEEE 21st International, IEEE, 2016, pp.~1--6.

\bibitem{albarakati2021}
{\sc A.~Albarakati, M.~Budi{\v s}i{\'c}, R.~Crocker, J.~{Glass-Klaiber},
  S.~Iams, J.~Maclean, N.~Marshall, C.~Roberts, and E.~S. Van~Vleck}, {\em
  Model and data reduction for data assimilation: {{Particle}} filters
  employing projected forecasts and data with application to a shallow water
  model}, Computers \& Mathematics with Applications,  (2021),
  \url{https://doi.org/10.1016/j.camwa.2021.05.026}.

\bibitem{alfatlawi2019arxiv}
{\sc M.~Alfatlawi and V.~Srivastava}, {\em An incremental approach to online
  dynamic mode decomposition for time-varying systems with applications to
  {EEG} data modeling}, Journal of Computational Dynamics, 7 (2020),
  pp.~209--241.

\bibitem{allgower2004nonlinear}
{\sc F.~Allg{\"o}wer, R.~Findeisen, and Z.~K. Nagy}, {\em Nonlinear model
  predictive control: From theory to application}, J. Chin. Inst. Chem. Engrs,
  35 (2004), pp.~299--315.

\bibitem{amini2019ifac}
{\sc A.~Amini, Q.~Sun, and N.~Motee}, {\em Carleman state feedback control
  design of a class of nonlinear control systems}, IFAC PapersOnLine, 52-20
  (2019), pp.~229--234.

\bibitem{Ansari2016ieee}
{\sc A.~R. Ansari and T.~D. Murphey}, {\em Sequential action control: Closed
  form optimal control for nonlinear and nonsmooth systems}, IEEE Transactions
  on Robotics, 32 (2016), pp.~1196--1214.

\bibitem{antoniou1997}
{\sc I.~Antoniou, L.~Dmitrieva, Y.~Kuperin, and Y.~Melnikov}, {\em Resonances
  and the extension of dynamics to rigged {{Hilbert}} space}, Computers \&
  Mathematics with Applications, 34 (1997), pp.~399--425,
  \url{https://doi.org/10.1016/S0898-1221(97)00148-X}.

\bibitem{antoniou1997a}
{\sc I.~Antoniou, B.~Qiao, and Z.~Suchanecki}, {\em Generalized spectral
  decomposition and intrinsic irreversibility of the {{Arnold Cat Map}}},
  Chaos, Solitons \& Fractals, 8 (1997), pp.~77--90,
  \url{https://doi.org/10.1016/S0960-0779(96)00056-2}.

\bibitem{apte2008}
{\sc A.~Apte, C.~K. R.~T. Jones, A.~M. Stuart, and J.~Voss}, {\em Data
  assimilation: {{Mathematical}} and statistical perspectives}, International
  Journal for Numerical Methods in Fluids, 56 (2008), pp.~1033--1046,
  \url{https://doi.org/10.1002/fld.1698}.

\bibitem{arbabi2018cdc}
{\sc H.~Arbabi, M.~Korda, and I.~Mezi\'c}, {\em A data-driven {K}oopman model
  predictive control framework for nonlinear partial differential equations},
  in Proceedings of the 2018 IEEE Conference on Decision and Control (CDC),
  IEEE, 2018.

\bibitem{arbabi2017}
{\sc H.~Arbabi and I.~Mezi{\'c}}, {\em Ergodic {{Theory}}, {{Dynamic Mode
  Decomposition}}, and {{Computation}} of {{Spectral Properties}} of the
  {{Koopman Operator}}}, SIAM Journal on Applied Dynamical Systems, 16 (2017),
  pp.~2096--2126, \url{https://doi.org/10.1137/17M1125236}.

\bibitem{arnold1998}
{\sc L.~Arnold}, {\em Random {{Dynamical Systems}}}, Springer {{Monographs}} in
  {{Mathematics}}, {Springer-Verlag}, {Berlin Heidelberg}, 1998,
  \url{https://doi.org/10.1007/978-3-662-12878-7}.

\bibitem{askham2017arxiv}
{\sc T.~Askham and J.~N. Kutz}, {\em Variable projection methods for an
  optimized dynamic mode decomposition}, SIAM Journal on Applied Dynamical
  Systems, 17 (2018), pp.~380--416.

\bibitem{assani2004}
{\sc I.~Assani}, {\em Spectral characterization of
  {{Wiener}}\textendash{{Wintner}} dynamical systems}, Ergodic Theory and
  Dynamical Systems, 24 (2004), pp.~347--365,
  \url{https://doi.org/10.1017/S0143385703000324}.

\bibitem{azencot2020forecasting}
{\sc O.~Azencot, N.~B. Erichson, V.~Lin, and M.~Mahoney}, {\em Forecasting
  sequential data using consistent koopman autoencoders}, in International
  Conference on Machine Learning, PMLR, 2020, pp.~475--485.

\bibitem{azencot2019consistent}
{\sc O.~Azencot, W.~Yin, and A.~Bertozzi}, {\em Consistent dynamic mode
  decomposition}, SIAM Journal on Applied Dynamical Systems, 18 (2019),
  pp.~1565--1585.

\bibitem{baddoo2021kernel}
{\sc P.~J. Baddoo, B.~Herrmann, B.~J. McKeon, and S.~L. Brunton}, {\em Kernel
  learning for robust dynamic mode decomposition: Linear and nonlinear
  disambiguation optimization ({LANDO})}, arXiv preprint arXiv:2106.01510,
  (2021).

\bibitem{bagheri2013jfm}
{\sc S.~Bagheri}, {\em Koopman-mode decomposition of the cylinder wake},
  Journal of Fluid Mechanics, 726 (2013), pp.~596--623.

\bibitem{bagheri2014}
{\sc S.~Bagheri}, {\em Effects of weak noise on oscillating flows: {{Linking}}
  quality factor, {{Floquet}} modes, and {{Koopman}} spectrum}, Physics of
  Fluids, 26 (2014), p.~094104, \url{https://doi.org/10.1063/1.4895898}.

\bibitem{bai2017aiaa}
{\sc Z.~Bai, E.~Kaiser, J.~L. Proctor, J.~N. Kutz, and S.~L. Brunton}, {\em
  Dynamic mode decomposition for compressive system identification}, AIAA
  Journal, 58 (2020).

\bibitem{balabane2020koopman}
{\sc M.~Balabane, M.~A. Mendez, and S.~Najem}, {\em Koopman operator for
  {B}urgers's equation}, Physical Review Fluids, 6 (2021), p.~064401.

\bibitem{banks1992infinite}
{\sc S.~Banks}, {\em Infinite-dimensional {C}arleman linearization, the {L}ie
  series and optimal control of non-linear partial differential equations},
  International journal of systems science, 23 (1992), pp.~663--675.

\bibitem{bar2019learning}
{\sc Y.~Bar-Sinai, S.~Hoyer, J.~Hickey, and M.~P. Brenner}, {\em Learning
  data-driven discretizations for partial differential equations}, Proceedings
  of the National Academy of Sciences, 116 (2019), pp.~15344--15349.

\bibitem{basley2013space}
{\sc J.~Basley, L.~R. Pastur, N.~Delprat, and F.~Lusseyran}, {\em Space-time
  aspects of a three-dimensional multi-modulated open cavity flow}, Physics of
  Fluids (1994-present), 25 (2013), p.~064105.

\bibitem{basley2011experimental}
{\sc J.~Basley, L.~R. Pastur, F.~Lusseyran, T.~M. Faure, and N.~Delprat}, {\em
  Experimental investigation of global structures in an incompressible cavity
  flow using time-resolved {PIV}}, Experiments in Fluids, 50 (2011),
  pp.~905--918.

\bibitem{battaglia2018relational}
{\sc P.~W. Battaglia, J.~B. Hamrick, V.~Bapst, A.~Sanchez-Gonzalez,
  V.~Zambaldi, M.~Malinowski, A.~Tacchetti, D.~Raposo, A.~Santoro, R.~Faulkner,
  et~al.}, {\em Relational inductive biases, deep learning, and graph
  networks}, arXiv preprint arXiv:1806.01261,  (2018).

\bibitem{beetham2020formulating}
{\sc S.~Beetham and J.~Capecelatro}, {\em Formulating turbulence closures using
  sparse regression with embedded form invariance}, Physical Review Fluids, 5
  (2020), p.~084611.

\bibitem{beetham2021sparse}
{\sc S.~Beetham, R.~O. Fox, and J.~Capecelatro}, {\em Sparse identification of
  multiphase turbulence closures for coupled fluid--particle flows}, Journal of
  Fluid Mechanics, 914 (2021).

\bibitem{bellani2011experimental}
{\sc G.~Bellani}, {\em Experimental studies of complex flows through
  image-based techniques},  (2011).

\bibitem{benner2020operator}
{\sc P.~Benner, P.~Goyal, B.~Kramer, B.~Peherstorfer, and K.~Willcox}, {\em
  Operator inference for non-intrusive model reduction of systems with
  non-polynomial nonlinear terms}, Computer Methods in Applied Mechanics and
  Engineering, 372 (2020), p.~113433.

\bibitem{Benner2015siamreview}
{\sc P.~Benner, S.~Gugercin, and K.~Willcox}, {\em A survey of projection-based
  model reduction methods for parametric dynamical systems}, SIAM review, 57
  (2015), pp.~483--531.

\bibitem{benner2018reduced}
{\sc P.~Benner, C.~Himpe, and T.~Mitchell}, {\em On reduced input-output
  dynamic mode decomposition}, Advances in Computational Mathematics, 44
  (2018), pp.~1751--1768.

\bibitem{benosman2017ifac}
{\sc M.~Benosman, H.~Mansour, and V.~Huroyan}, {\em Koopman-operator
  observer-based estimation of pedestrian crowd flows}, IFAC-PapersOnLine, 50
  (2017), pp.~14028---14033.

\bibitem{berger2014ieee}
{\sc E.~Berger, M.~Sastuba, D.~Vogt, B.~Jung, and H.~B. Amor}, {\em Estimation
  of perturbations in robotic behavior using dynamic mode decomposition},
  Journal of Advanced Robotics, 29 (2015), pp.~331--343.

\bibitem{birkhoff1931pnas}
{\sc G.~D. Birkhoff}, {\em Proof of the ergodic theorem}, Proceedings of the
  National Academy of Sciences, 17 (1931), pp.~656--660.

\bibitem{birkhoff1932pnas}
{\sc G.~D. Birkhoff and B.~O. Koopman}, {\em Recent contributions to the
  ergodic theory}, Proceedings of the National Academy of Sciences, 18 (1932),
  pp.~279--282.

\bibitem{bistrian2016ijnme}
{\sc D.~Bistrian and I.~Navon}, {\em Randomized dynamic mode decomposition for
  non-intrusive reduced order modelling}, International Journal for Numerical
  Methods in Engineering,  (2016).

\bibitem{bistrian2015ijnmf}
{\sc D.~A. Bistrian and I.~M. Navon}, {\em An improved algorithm for the
  shallow water equations model reduction: Dynamic mode decomposition vs
  {POD}}, International Journal for Numerical Methods in Fluids,  (2015).

\bibitem{bittracher2015}
{\sc A.~Bittracher, P.~Koltai, and O.~Junge}, {\em Pseudogenerators of
  {{Spatial Transfer Operators}}}, SIAM Journal on Applied Dynamical Systems,
  14 (2015), pp.~1478--1517, \url{https://doi.org/10.1137/14099872X}.

\bibitem{bollt2021}
{\sc E.~Bollt}, {\em Geometric {{Considerations}} of a {{Good Dictionary}} for
  {{Koopman Analysis}} of {{Dynamical Systems}}: {{Cardinality}}, '{{Primary
  Eigenfunction}},' and {{Efficient Representation}}}, arXiv:1912.09570 [cs,
  math],  (2021), \url{https://arxiv.org/abs/1912.09570}.

\bibitem{bollt2018}
{\sc E.~M. Bollt, Q.~Li, F.~Dietrich, and I.~Kevrekidis}, {\em On {{Matching}},
  and {{Even Rectifying}}, {{Dynamical Systems}} through {{Koopman Operator
  Eigenfunctions}}}, SIAM Journal on Applied Dynamical Systems, 17 (2018),
  pp.~1925--1960, \url{https://doi.org/10.1137/17M116207X}.

\bibitem{bollt2013}
{\sc E.~M. Bollt and N.~Santitissadeekorn}, {\em Applied and Computational
  Measurable Dynamics}, vol.~18 of Mathematical {{Modeling}} and
  {{Computation}}, {Society for Industrial and Applied Mathematics (SIAM),
  Philadelphia, PA}, 2013.

\bibitem{arima}
{\sc G.~E.~P. Box, G.~M. Jenkins, and G.~C. Reinsel}, {\em \textsl{Time series
  analysis: Forecasting and control, 3rd Ed.}}, Prentice Hall, Englewood
  Cliffs, N.J., 1994.

\bibitem{boyce2017elementary}
{\sc W.~E. Boyce, R.~C. DiPrima, and D.~B. Meade}, {\em Elementary differential
  equations}, John Wiley \& Sons, 2017.

\bibitem{breiman2001statistical}
{\sc L.~Breiman et~al.}, {\em Statistical modeling: The two cultures (with
  comments and a rejoinder by the author)}, Statistical science, 16 (2001),
  pp.~199--231.

\bibitem{Brenner2019prf}
{\sc M.~Brenner, J.~Eldredge, and J.~Freund}, {\em Perspective on machine
  learning for advancing fluid mechanics}, Physical Review Fluids, 4 (2019),
  p.~100501.

\bibitem{bright2013pof}
{\sc I.~Bright, G.~Lin, and J.~N. Kutz}, {\em Compressive sensing and machine
  learning strategies for characterizing the flow around a cylinder with
  limited pressure measurements}, Physics of Fluids, 25 (2013), pp.~1--15.

\bibitem{broad2017proc}
{\sc A.~Broad, T.~Murphey, and B.~Argall}, {\em Learning models for shared
  control of human-machine systems with unknown dynamics}, Robotics: Science
  and Systems Proceedings,  (2017).

\bibitem{brockett1976automatica}
{\sc R.~W. Brockett}, {\em Volterra series and geometric control theory},
  Automatica, 12 (1976), pp.~167--176.

\bibitem{brockett2007ams}
{\sc R.~W. Brockett}, {\em Optimal control of the {L}iouville equation}, AMS IP
  Studies in Advanced Mathematics, 39 (2007), p.~23.

\bibitem{broomhead1989prsla}
{\sc D.~Broomhead and R.~Jones}, {\em Time-series analysis}, in Proceedings of
  the Royal Society of London A: Mathematical, Physical and Engineering
  Sciences, vol.~423, The Royal Society, 1989, pp.~103--121.

\bibitem{bruder2020arxiv}
{\sc D.~Bruder, X.~Fu, R.~B. Gillespie, C.~D. Remy, and R.~Vasudevan}, {\em
  {Koopman-based Control of a Soft Continuum Manipulator Under Variable Loading
  Conditions}}, IEEE Robotics and Automation Letters, 6 (2021), pp.~6852--6859.

\bibitem{bruder2020arxiva}
{\sc D.~Bruder, X.~Fu, and R.~Vasudevan}, {\em Advantages of bilinear {K}oopman
  realizations for the modeling and control of systems with unknown dynamics},
  arXiv preprint arXiv:arXiv:2010.09961v3,  (2020),
  \url{https://arxiv.org/abs/https://arxiv.org/abs/2010.09961}.

\bibitem{bruder2019proc}
{\sc D.~Bruder, B.~Gillespie, C.~David~Remy, and R.~Vasudevan}, {\em Modeling
  and control of soft robots using the {K}oopman operator and model predictive
  control}, in Robotics: Science and Systems, 2019.

\bibitem{brunton2016b}
{\sc B.~W. Brunton, L.~A. Johnson, J.~G. Ojemann, and J.~N. Kutz}, {\em
  Extracting spatial\textendash temporal coherent patterns in large-scale
  neural recordings using dynamic mode decomposition}, Journal of Neuroscience
  Methods, 258 (2016), pp.~1--15,
  \url{https://doi.org/10.1016/j.jneumeth.2015.10.010}.

\bibitem{brunton2017natcomm}
{\sc S.~L. Brunton, B.~W. Brunton, J.~L. Proctor, E.~Kaiser, and J.~N. Kutz},
  {\em Chaos as an intermittently forced linear system}, Nature Communications,
  8 (2017), pp.~1--9.

\bibitem{brunton2016plosone}
{\sc S.~L. Brunton, B.~W. Brunton, J.~L. Proctor, and J.~N. Kutz}, {\em Koopman
  invariant subspaces and finite linear representations of nonlinear dynamical
  systems for control}, PLoS ONE, 11 (2016), p.~e0150171.

\bibitem{brunton2019data}
{\sc S.~L. Brunton and J.~N. Kutz}, {\em Data-driven science and engineering:
  Machine learning, dynamical systems, and control}, Cambridge University
  Press, 2019.

\bibitem{Brunton2020arfm}
{\sc S.~L. Brunton, B.~R. Noack, and P.~Koumoutsakos}, {\em Machine learning
  for fluid mechanics}, Annual Review of Fluid Mechanics, 52 (2020),
  pp.~477--508.

\bibitem{brunton2016pnas}
{\sc S.~L. Brunton, J.~L. Proctor, and J.~N. Kutz}, {\em Discovering governing
  equations from data by sparse identification of nonlinear dynamical systems},
  Proceedings of the National Academy of Sciences, 113 (2016), pp.~3932--3937.

\bibitem{brunton2015jcd}
{\sc S.~L. Brunton, J.~L. Proctor, J.~H. Tu, and J.~N. Kutz}, {\em Compressed
  sensing and dynamic mode decomposition}, Journal of Computational Dynamics, 2
  (2015), pp.~165--191.

\bibitem{brunton2014siads}
{\sc S.~L. Brunton, J.~H. Tu, I.~Bright, and J.~N. Kutz}, {\em Compressive
  sensing and low-rank libraries for classification of bifurcation regimes in
  nonlinear dynamical systems}, SIAM Journal on Applied Dynamical Systems, 13
  (2014), pp.~1716--1732.

\bibitem{budisic2009cdc}
{\sc M.~Budi{\v{s}}i{\'c} and I.~Mezi{\'c}}, {\em An approximate
  parametrization of the ergodic partition using time averaged observables}, in
  Decision and Control, 2009 held jointly with the 2009 28th Chinese Control
  Conference. CDC/CCC 2009. Proceedings of the 48th IEEE Conference on, IEEE,
  2009, pp.~3162--3168.

\bibitem{budisic2012physd}
{\sc M.~Budi{\v s}i{\'c} and I.~Mezi{\'c}}, {\em Geometry of the ergodic
  quotient reveals coherent structures in flows}, Physica D: Nonlinear
  Phenomena, 241 (2012), pp.~1255--1269,
  \url{https://doi.org/10.1016/j.physd.2012.04.006}.

\bibitem{budisic2012chaos}
{\sc M.~Budi{\v{s}}i{\'c}, R.~Mohr, and I.~Mezi{\'c}}, {\em Applied
  {K}oopmanism}, Chaos: An Interdisciplinary Journal of Nonlinear Science, 22
  (2012), p.~047510.

\bibitem{burgers}
{\sc J.~M. Burgers}, {\em A mathematical model illustrating the theory of
  turbulence}, Advances in applied mechanics, 1 (1948), pp.~171--199.

\bibitem{manohar2020kernel}
{\sc D.~Burov, D.~Giannakis, K.~Manohar, and A.~Stuart}, {\em Kernel analog
  forecasting: Multiscale test problems}, Multiscale Modeling \& Simulation, 19
  (2021), pp.~1011--1040.

\bibitem{camacho2013model}
{\sc E.~F. Camacho and C.~B. Alba}, {\em Model predictive control}, Springer
  Science \& Business Media, 2013.

\bibitem{candes2006bieeetit}
{\sc E.~J. Cand\`es, J.~Romberg, and T.~Tao}, {\em Robust uncertainty
  principles: exact signal reconstruction from highly incomplete frequency
  information}, IEEE Transactions on Information Theory, 52 (2006),
  pp.~489--509.

\bibitem{caraballo2017}
{\sc T.~Caraballo and X.~Han}, {\em Applied {{Nonautonomous}} and {{Random
  Dynamical Systems}}: {{Applied Dynamical Systems}}}, {Springer}, Jan. 2017.

\bibitem{carleman1932am}
{\sc T.~Carleman}, {\em Application de la th{\'e}orie des {\'e}quations
  int{\'e}grales lin{\'e}aires aux syst{\`e}mes d'{\'e}quations
  diff{\'e}rentielles non lin{\'e}aires}, Acta Mathematica, 59 (1932),
  pp.~63--87.

\bibitem{carleman1933theorie}
{\sc T.~Carleman}, {\em Sur la th{\'e}orie de l'{\'e}quation
  int{\'e}grodiff{\'e}rentielle de boltzmann}, Acta Mathematica, 60 (1933),
  pp.~91--146.

\bibitem{carleman1933systemes}
{\sc T.~Carleman}, {\em Sur les systemes lineaires aux d{\'e}riv{\'e}es
  partielles du premier ordrea deux variables}, CR Acad. Sci. Paris, 197
  (1933), pp.~471--474.

\bibitem{champion2019discovery}
{\sc K.~P. Champion, S.~L. Brunton, and J.~N. Kutz}, {\em Discovery of
  nonlinear multiscale systems: Sampling strategies and embeddings}, SIAM
  Journal on Applied Dynamical Systems, 18 (2019), pp.~312--333.

\bibitem{Chekroun2014}
{\sc M.~D. Chekroun, J.~D. Neelin, D.~Kondrashov, J.~C. McWilliams, and
  M.~Ghil}, {\em Rough parameter dependence in climate models and the role of
  {{Ruelle}}-{{Pollicott}} resonances}, Proceedings of the National Academy of
  Sciences, 111 (2014), pp.~1684--1690,
  \url{https://doi.org/10.1073/pnas.1321816111}.

\bibitem{Chekroun2020}
{\sc M.~D. Chekroun, A.~Tantet, H.~A. Dijkstra, and J.~D. Neelin}, {\em
  Ruelle\textendash{{Pollicott Resonances}} of {{Stochastic Systems}} in
  {{Reduced State Space}}. {{Part I}}: {{Theory}}}, Journal of Statistical
  Physics, 179 (2020), pp.~1366--1402,
  \url{https://doi.org/10.1007/s10955-020-02535-x}.

\bibitem{chen2012jns}
{\sc K.~K. Chen, J.~H. Tu, and C.~W. Rowley}, {\em Variants of dynamic mode
  decomposition: Boundary condition, {K}oopman, and {F}ourier analyses},
  Journal of Nonlinear Science, 22 (2012), pp.~887--915.

\bibitem{chicone1999evolution}
{\sc C.~Chicone and Y.~Latushkin}, {\em Evolution semigroups in dynamical
  systems and differential equations}, no.~70, American Mathematical Soc.,
  1999.

\bibitem{coifman2008mmas}
{\sc R.~R. Coifman, I.~G. Kevrekidis, S.~Lafon, M.~Maggioni, and B.~Nadler},
  {\em Diffusion maps, reduction coordinates, and low dimensional
  representation of stochastic systems}, Multiscale Modeling \& Simulation, 7
  (2008), pp.~842--864.

\bibitem{coifman2006acha}
{\sc R.~R. Coifman and S.~Lafon}, {\em Diffusion maps}, Applied and
  computational harmonic analysis, 21 (2006), pp.~5--30.

\bibitem{coifman2005pnas}
{\sc R.~R. Coifman, S.~Lafon, A.~B. Lee, M.~Maggioni, B.~Nadler, F.~Warner, and
  S.~W. Zucker}, {\em Geometric diffusions as a tool for harmonic analysis and
  structure definition of data: Diffusion maps}, Proceedings of the National
  Academy of Sciences of the United States of America, 102 (2005),
  pp.~7426--7431.

\bibitem{cole51}
{\sc J.~D. Cole}, {\em On a quasi-linear parabolic equation occurring in
  aerodynamics}, Quart. Appl. Math., 9 (1951), pp.~225--236.

\bibitem{colonius2000book}
{\sc F.~Colonius and W.~Kliemann}, {\em The Dynamics of Control}, Birkh\"auser
  Boston, 2000.

\bibitem{cornfeld1982}
{\sc I.~P. Cornfeld, S.~V. Fomin, and Y.~G. Sinai}, {\em Ergodic {{Theory}}},
  vol.~245 of Grundlehren Der Mathematischen {{Wissenschaften}}, {Springer New
  York}, {New York, NY}, 1982.

\bibitem{cranmer2020lagrangian}
{\sc M.~Cranmer, S.~Greydanus, S.~Hoyer, P.~Battaglia, D.~Spergel, and S.~Ho},
  {\em Lagrangian neural networks}, arXiv preprint arXiv:2003.04630,  (2020).

\bibitem{cranmer2019learning}
{\sc M.~D. Cranmer, R.~Xu, P.~Battaglia, and S.~Ho}, {\em Learning symbolic
  physics with graph networks}, arXiv preprint arXiv:1909.05862,  (2019).

\bibitem{crnjaric-zic2020}
{\sc N.~{\v C}rnjari{\'c}-{\v Z}ic, S.~Ma{\'c}e{\v s}i{\'c}, and I.~Mezi{\'c}},
  {\em Koopman {{Operator Spectrum}} for {{Random Dynamical Systems}}}, Journal
  of Nonlinear Science, 30 (2020), pp.~2007--2056,
  \url{https://doi.org/10.1007/s00332-019-09582-z}.

\bibitem{chaosbook}
{\sc P.~Cvitanovi{\'c}, R.~Artuso, R.~Mainieri, G.~Tanner, and G.~Vattay}, {\em
  Chaos: {{Classical}} and Quantum}, {Niels Bohr Inst.}, {Copenhagen}, 2016.

\bibitem{mozyrska2006}
{\sc M.~D and B.~Z}, {\em Dualities for linear control differential systems
  with infinite matrices}, Control Cybern, 35 (2006), pp.~887--904.

\bibitem{mozyrska2008}
{\sc M.~D and B.~Z}, {\em Carleman linearization of linearly observable
  polynomial systems}, in Sarychev A, Shiryaev A, Guerra M, Grossinho MdR (eds)
  Mathematical control theory and finance, Springer, Berlin, 2008,
  pp.~311--323.

\bibitem{das2017arxiv}
{\sc S.~Das and D.~Giannakis}, {\em Delay-coordinate maps and the spectra of
  {K}oopman operators}, Journal of Statistical Physics, 175 (2019),
  pp.~1107--1145.

\bibitem{das2020a}
{\sc S.~Das and D.~Giannakis}, {\em Koopman spectra in reproducing kernel
  {{Hilbert}} spaces}, Applied and Computational Harmonic Analysis, 49 (2020),
  pp.~573--607, \url{https://doi.org/10.1016/j.acha.2020.05.008}.

\bibitem{das2021}
{\sc S.~Das, D.~Giannakis, and J.~Slawinska}, {\em Reproducing kernel
  {{Hilbert}} space compactification of unitary evolution groups}, Applied and
  Computational Harmonic Analysis, 54 (2021), pp.~75--136,
  \url{https://doi.org/10.1016/j.acha.2021.02.004}.

\bibitem{dawson2016ef}
{\sc S.~T. Dawson, M.~S. Hemati, M.~O. Williams, and C.~W. Rowley}, {\em
  Characterizing and correcting for the effect of sensor noise in the dynamic
  mode decomposition}, Experiments in Fluids, 57 (2016), pp.~1--19.

\bibitem{dellnitz2001}
{\sc M.~Dellnitz, G.~Froyland, and O.~Junge}, {\em The algorithms behind
  {{GAIO}}-set oriented numerical methods for dynamical systems}, in Ergodic
  Theory, Analysis, and Efficient Simulation of Dynamical Systems, {Springer,
  Berlin}, 2001, pp.~145--174, 805--807.

\bibitem{dellnitz1997}
{\sc M.~Dellnitz and A.~Hohmann}, {\em A subdivision algorithm for the
  computation of unstable manifolds and global attractors}, Numerische
  Mathematik, 75 (1997), pp.~293--317,
  \url{https://doi.org/10.1007/s002110050240}.

\bibitem{dellnitz1997a}
{\sc M.~Dellnitz, A.~Hohmann, O.~Junge, and M.~Rumpf}, {\em Exploring invariant
  sets and invariant measures}, Chaos: An Interdisciplinary Journal of
  Nonlinear Science, 7 (1997), pp.~221--228,
  \url{https://doi.org/10.1063/1.166223}.

\bibitem{dellnitz1999}
{\sc M.~Dellnitz and O.~Junge}, {\em On the approximation of complicated
  dynamical behavior}, SIAM Journal on Numerical Analysis, 36 (1999),
  pp.~491--515, \url{https://doi.org/10.1137/S0036142996313002}.

\bibitem{Demo18pydmd}
{\sc N.~Demo, M.~Tezzele, and G.~Rozza}, {\em {PyDMD}: {P}ython dynamic mode
  decomposition}, The Journal of Open Source Software, 3 (2018), p.~530.

\bibitem{dietrich2020koopman}
{\sc F.~Dietrich, T.~N. Thiem, and I.~G. Kevrekidis}, {\em On the koopman
  operator of algorithms}, SIAM Journal on Applied Dynamical Systems, 19
  (2020), pp.~860--885.

\bibitem{dogra2020optimizing}
{\sc A.~S. Dogra and W.~T. Redman}, {\em Optimizing neural networks via
  {K}oopman operator theory}, arXiv preprint arXiv:2006.02361,  (2020).

\bibitem{donoho201750}
{\sc D.~Donoho}, {\em 50 years of data science}, Journal of Computational and
  Graphical Statistics, 26 (2017), pp.~745--766.

\bibitem{donoho2006ieeetit}
{\sc D.~L. Donoho}, {\em Compressed sensing}, IEEE Transactions on Information
  Theory, 52 (2006), pp.~1289--1306.

\bibitem{drmac2016siam}
{\sc Z.~Drma{\v{c}} and S.~Gugercin}, {\em A new selection operator for the
  discrete empirical interpolation method---improved a priori error bound and
  extensions}, SIAM Journal on Scientific Computing, 38 (2016), pp.~A631--A648.

\bibitem{drmac2018}
{\sc Z.~Drma{\v{c}}, I.~Mezi{\'{c}}, and R.~Mohr}, {\em {Data driven modal
  decompositions: Analysis and enhancements}}, SIAM J. Sci. Comput., 40 (2018),
  pp.~A2253--A2285, \url{https://doi.org/10.1137/17M1144155}.

\bibitem{drmac2019}
{\sc Z.~Drma{\v c}, I.~Mezi{\'c}, and R.~Mohr}, {\em Data {{Driven Koopman
  Spectral Analysis}} in {{Vandermonde}}--{{Cauchy Form}} via the {{DFT}}:
  {{Numerical Method}} and {{Theoretical Insights}}}, SIAM Journal on
  Scientific Computing, 41 (2019), pp.~A3118--A3151,
  \url{https://doi.org/10.1137/18M1227688}.

\bibitem{drmac2020}
{\sc Z.~Drma{\v c}, I.~Mezi{\'c}, and R.~Mohr}, {\em On {{Least Squares
  Problems}} with {{Certain Vandermonde}}--{{Khatri}}--{{Rao Structure}} with
  {{Applications}} to {{DMD}}}, SIAM Journal on Scientific Computing, 42
  (2020), pp.~A3250--A3284, \url{https://doi.org/10.1137/19M1288474}.

\bibitem{duke2012experimental}
{\sc D.~Duke, D.~Honnery, and J.~Soria}, {\em Experimental investigation of
  nonlinear instabilities in annular liquid sheets}, Journal of Fluid
  Mechanics, 691 (2012), pp.~594--604.

\bibitem{duke2012error}
{\sc D.~Duke, J.~Soria, and D.~Honnery}, {\em An error analysis of the dynamic
  mode decomposition}, Experiments in fluids, 52 (2012), pp.~529--542.

\bibitem{dunne2015ef}
{\sc R.~Dunne and B.~J. McKeon}, {\em Dynamic stall on a pitching and surging
  airfoil}, Experiments in Fluids, 56 (2015), pp.~1--15.

\bibitem{Duraisamy2019arfm}
{\sc K.~Duraisamy, G.~Iaccarino, and H.~Xiao}, {\em Turbulence modeling in the
  age of data}, Annual Reviews of Fluid Mechanics, 51 (2019), pp.~357--377.

\bibitem{dylewsky2019}
{\sc D.~Dylewsky, M.~Tao, and J.~N. Kutz}, {\em Dynamic mode decomposition for
  multiscale nonlinear physics}, Physical Review E, 99 (2019), p.~063311,
  \url{https://doi.org/10.1103/PhysRevE.99.063311}.

\bibitem{eckmann1985}
{\sc J.-P. Eckmann and D.~P. Ruelle}, {\em Ergodic theory of chaos and strange
  attractors}, Reviews of modern physics, 57 (1985), pp.~617--656,
  \url{https://doi.org/10.1103/RevModPhys.57.617}.

\bibitem{eisner2015}
{\sc T.~Eisner, B.~Farkas, M.~Haase, and R.~Nagel}, {\em Operator Theoretic
  Aspects of Ergodic Theory}, vol.~272, {Springer}, 2015.

\bibitem{eivazi2021recurrent}
{\sc H.~Eivazi, L.~Guastoni, P.~Schlatter, H.~Azizpour, and R.~Vinuesa}, {\em
  Recurrent neural networks and koopman-based frameworks for temporal
  predictions in a low-order model of turbulence}, International Journal of
  Heat and Fluid Flow, 90 (2021), p.~108816.

\bibitem{eldering2018}
{\sc J.~Eldering, M.~Kvalheim, and S.~Revzen}, {\em Global linearization and
  fiber bundle structure of invariant manifolds}, Nonlinearity, 31 (2018),
  pp.~4202--4245, \url{https://doi.org/10.1088/1361-6544/aaca8d}.

\bibitem{eren2017jgcd}
{\sc U.~Eren, A.~Prach, B.~B. Ko{\c{c}}er, S.~V. Rakovi{\'c}, E.~Kayacan, and
  B.~A{\c{c}}{\i}kme{\c{s}}e}, {\em Model predictive control in aerospace
  systems: Current state and opportunities}, Journal of Guidance, Control, and
  Dynamics,  (2017).

\bibitem{erichson2016jrtp}
{\sc N.~B. Erichson, S.~L. Brunton, and J.~N. Kutz}, {\em Compressed dynamic
  mode decomposition for real-time object detection}, Journal of Real-Time
  Image Processing,  (2016).

\bibitem{erichsonrandomized}
{\sc N.~B. Erichson, K.~Manohar, S.~L. Brunton, and J.~N. Kutz}, {\em
  Randomized cp tensor decomposition}, Machine Learning: Science and
  Technology, 1 (2020), p.~025012.

\bibitem{erichson2017arxiv}
{\sc N.~B. Erichson, L.~Mathelin, J.~N. Kutz, and S.~L. Brunton}, {\em
  Randomized dynamic mode decomposition}, SIAM Journal on Applied Dynamical
  Systems, 18 (2019), pp.~1867--1891.

\bibitem{erichson2016arxiva}
{\sc N.~B. Erichson, S.~Voronin, S.~L. Brunton, and J.~N. Kutz}, {\em
  Randomized matrix decompositions using {R}}, Journal of Statistical Softwar,
  89 (2019), pp.~1--48.

\bibitem{evensen2009data}
{\sc G.~Evensen}, {\em Data assimilation: the ensemble Kalman filter}, Springer
  Science \& Business Media, 2009.

\bibitem{everson1995karhunen}
{\sc R.~Everson and L.~Sirovich}, {\em Karhunen--{L}oeve procedure for gappy
  data}, JOSA A, 12 (1995), pp.~1657--1664.

\bibitem{Farazmand2012chaos}
{\sc M.~Farazmand and G.~Haller}, {\em Computing {L}agrangian coherent
  structures from their variational theory}, Chaos, 22 (2012),
  pp.~013128--1--013128--12.

\bibitem{folkestad2020arxiv}
{\sc C.~Folkestad, D.~Pastor, and J.~W. Burdick}, {\em Episodic koopman
  learning of nonlinear robot dynamics with application to fast multirotor
  landing},  (2020), pp.~9216--9222.

\bibitem{fontanini2016be}
{\sc A.~D. Fontanini, U.~Vaidya, and B.~Ganapathysubramanian}, {\em A
  methodology for optimal placement of sensors in enclosed environments: A
  dynamical systems approach}, Building and Environment, 100 (2016),
  pp.~145--161.

\bibitem{froyland2003}
{\sc G.~Froyland and M.~Dellnitz}, {\em Detecting and locating near-optimal
  almost-invariant sets and cycles}, SIAM Journal on Scientific Computing, 24
  (2003), pp.~1839--1863 (electronic),
  \url{https://doi.org/10.1137/S106482750238911X}.

\bibitem{froyland2014}
{\sc G.~Froyland, C.~{Gonz{\'a}lez-Tokman}, and A.~Quas}, {\em Detecting
  isolated spectrum of transfer and {{Koopman}} operators with {{Fourier}}
  analytic tools}, Journal of Computational Dynamics, 1 (2014), pp.~249--278,
  \url{https://doi.org/10.3934/jcd.2014.1.249}.

\bibitem{froyland2016}
{\sc G.~Froyland, G.~A. Gottwald, and A.~Hammerlindl}, {\em A trajectory-free
  framework for analysing multiscale systems}, Physica D: Nonlinear Phenomena,
  328-329 (2016), pp.~34--43,
  \url{https://doi.org/10.1016/j.physd.2016.04.010}.

\bibitem{froyland2013}
{\sc G.~Froyland, O.~Junge, and P.~Koltai}, {\em Estimating {{Long}}-{{Term
  Behavior}} of {{Flows}} without {{Trajectory Integration}}: {{The
  Infinitesimal Generator Approach}}}, SIAM Journal on Numerical Analysis, 51
  (2013), pp.~223--247, \url{https://doi.org/10.1137/110819986}.

\bibitem{Froyland2010b}
{\sc G.~Froyland and O.~Stancevic}, {\em Escape rates and
  {{Perron}}-{{Frobenius}} operators: {{Open}} and closed dynamical systems},
  Discrete \& Continuous Dynamical Systems - B, 14 (2010), p.~457,
  \url{https://doi.org/10.3934/dcdsb.2010.14.457}.

\bibitem{fujii2019dynamic}
{\sc K.~Fujii and Y.~Kawahara}, {\em Dynamic mode decomposition in
  vector-valued reproducing kernel hilbert spaces for extracting dynamical
  structure among observables}, Neural Networks, 117 (2019), pp.~94--103.

\bibitem{mamakoukas2020arxivb}
{\sc I.~A. G~Mamakoukas and T.~Murphey}, {\em Learning data-driven stable
  {K}oopman operators}, arXiv preprint arXiv:2005.04291,  (2020).

\bibitem{mamakoukas2019proc}
{\sc X.~T. G~Mamakoukas, M~Castano and T.~Murphey}, {\em Local {K}oopman
  operators for data-driven control of robotic systems}, in Proceedings of
  ``Robotics: Science and Systems 2019'', Freiburg im Breisgau, June 22-26,
  2019, IEEE, 2019.

\bibitem{garcia1989model}
{\sc C.~E. Garcia, D.~M. Prett, and M.~Morari}, {\em Model predictive control:
  theory and practice---a survey}, Automatica, 25 (1989), pp.~335--348.

\bibitem{garriga2010model}
{\sc J.~L. Garriga and M.~Soroush}, {\em Model predictive control tuning
  methods: A review}, Industrial \& Engineering Chemistry Research, 49 (2010),
  pp.~3505--3515.

\bibitem{gaspard1998}
{\sc P.~Gaspard}, {\em Chaos, Scattering and Statistical Mechanics}, vol.~9 of
  Cambridge {{Nonlinear Science Series}}, {Cambridge University Press,
  Cambridge}, 1998.

\bibitem{gaspard1995}
{\sc P.~Gaspard, G.~Nicolis, A.~Provata, and S.~Tasaki}, {\em Spectral
  signature of the pitchfork bifurcation: {{Liouville}} equation approach},
  Physical Review E, 51 (1995), pp.~74--94,
  \url{https://doi.org/10.1103/PhysRevE.51.74}.

\bibitem{Gelss2019mindy}
{\sc P.~Gel{\ss}, S.~Klus, J.~Eisert, and C.~Sch{\"u}tte}, {\em
  Multidimensional approximation of nonlinear dynamical systems}, Journal of
  Computational and Nonlinear Dynamics, 14 (2019).

\bibitem{giannakis2019}
{\sc D.~Giannakis}, {\em Data-driven spectral decomposition and forecasting of
  ergodic dynamical systems}, Applied and Computational Harmonic Analysis, 47
  (2019), pp.~338--396, \url{https://doi.org/10.1016/j.acha.2017.09.001}.

\bibitem{giannakis2020a}
{\sc D.~Giannakis}, {\em Delay-coordinate maps, coherence, and approximate
  spectra of evolution operators}, arXiv:2007.02195 [nlin, physics:physics],
  (2020), \url{https://arxiv.org/abs/2007.02195}.

\bibitem{giannakis2020}
{\sc D.~Giannakis and S.~Das}, {\em Extraction and prediction of coherent
  patterns in incompressible flows through space\textendash{}time {{Koopman}}
  analysis}, Physica D: Nonlinear Phenomena, 402 (2020), p.~132211,
  \url{https://doi.org/10.1016/j.physd.2019.132211}.

\bibitem{gin2019deep}
{\sc C.~Gin, B.~Lusch, S.~L. Brunton, and J.~N. Kutz}, {\em Deep learning
  models for global coordinate transformations that linearise {PDE}s}, European
  Journal of Applied Mathematics,  (2020), pp.~1--25.

\bibitem{gin2020deepgreen}
{\sc C.~R. Gin, D.~E. Shea, S.~L. Brunton, and J.~N. Kutz}, {\em {DeepGreen}:
  Deep learning of {G}reen's functions for nonlinear boundary value problems},
  arXiv preprint arXiv:2101.07206,  (2020).

\bibitem{goswami2017cdc}
{\sc D.~Goswami and D.~A. Paley}, {\em Global bilinearization and
  controllability of control-affine nonlinear systems: A {K}oopman spectral
  approach}, in Decision and Control (CDC), 2017 IEEE 56th Annual Conference
  on, IEEE, 2017, pp.~6107--6112.

\bibitem{govindarajan2019}
{\sc N.~Govindarajan, R.~Mohr, S.~Chandrasekaran, and I.~Mezic}, {\em On the
  {{Approximation}} of {{Koopman Spectra}} for {{Measure Preserving
  Transformations}}}, SIAM Journal on Applied Dynamical Systems, 18 (2019),
  pp.~1454--1497, \url{https://doi.org/10.1137/18M1175094}.

\bibitem{goza2018modal}
{\sc A.~Goza and T.~Colonius}, {\em Modal decomposition of fluid--structure
  interaction with application to flag flapping}, Journal of Fluids and
  Structures, 81 (2018), pp.~728--737.

\bibitem{grilli:2012}
{\sc M.~Grilli, P.~J. Schmid, S.~Hickel, and N.~A. Adams}, {\em Analysis of
  unsteady behaviour in shockwave turbulent boundary layer interaction},
  Journal of Fluid Mechanics, 700 (2012), pp.~16--28.

\bibitem{grosek2014arxiv}
{\sc J.~Grosek and J.~N. Kutz}, {\em Dynamic mode decomposition for real-time
  background/foreground separation in video}, arXiv preprint arXiv:1404.7592,
  (2014).

\bibitem{guckenheimer1986}
{\sc J.~Guckenheimer and P.~Holmes}, {\em Nonlinear Oscillations, Dynamical
  Systems, and Bifurcations of Vector Fields}, {New york,}, 1986, 1983.

\bibitem{gueniat2015pof}
{\sc F.~Gueniat, L.~Mathelin, and L.~Pastur}, {\em A dynamic mode decomposition
  approach for large and arbitrarily sampled systems}, Physics of Fluids, 27
  (2015), p.~025113.

\bibitem{halko2011siamreview}
{\sc N.~Halko, P.-G. Martinsson, and J.~A. Tropp}, {\em Finding structure with
  randomness: Probabilistic algorithms for constructing approximate matrix
  decompositions}, SIAM review, 53 (2011), pp.~217--288.

\bibitem{Haller2002pof}
{\sc G.~Haller}, {\em {Lagrangian} coherent structures from approximate
  velocity data}, Physics of Fluids, 14 (2002), pp.~1851--1861.

\bibitem{Haller2015arfm}
{\sc G.~Haller}, {\em Lagrangian coherent structures}, Annual Review of Fluid
  Mechanics, 47 (2015), pp.~137--162.

\bibitem{halmos1944a}
{\sc P.~R. Halmos}, {\em Approximation {{Theories}} for {{Measure Preserving
  Transformations}}}, Transactions of the American Mathematical Society, 55
  (1944), pp.~1--18, \url{https://doi.org/10.2307/1990137}.

\bibitem{han2020arxiv}
{\sc Y.~Han, W.~Hao, and U.~Vaidya}, {\em Deep learning of {K}oopman
  representation for control},  (2020), pp.~1890--1895.

\bibitem{hanke2018arxiv}
{\sc S.~Hanke, S.~Peitz, O.~Wallscheid, S.~Klus, J.~B\"ocker, and M.~Dellnitz},
  {\em Koopman operator based finite-set model predictive control for
  electrical drives}, arXiv:1804.00854,  (2018).

\bibitem{haseli2020b}
{\sc M.~Haseli and J.~Cortes}, {\em {Efficient Identification of Linear
  Evolutions in Nonlinear Vector Fields: {K}oopman Invariant Subspaces}},
  (2020), pp.~1746--1751, \url{https://doi.org/10.1109/cdc40024.2019.9029955}.

\bibitem{hasnain2019arxivb}
{\sc A.~Hasnain, N.~Boddupalli, S.~Balakrishnan, and E.~Yeung}, {\em Steady
  state programming of controlled nonlinear systems via deep dynamic mode
  decomposition},  (2020), pp.~4245--4251.

\bibitem{hasnain2019arxiv}
{\sc A.~Hasnain, N.~Boddupalli, and E.~Yeung}, {\em Optimal reporter placement
  in sparsely measured genetic networks using the {K}oopman operator},  (2019),
  pp.~19--24.

\bibitem{hemati2017aiaa}
{\sc M.~Hemati and H.~Yao}, {\em Dynamic mode shaping for fluid flow control:
  New strategies for transient growth suppression}, in 8th AIAA Theoretical
  Fluid Mechanics Conference, 2017, p.~3160.

\bibitem{hemati2017tcfd}
{\sc M.~S. Hemati, C.~W. Rowley, E.~A. Deem, and L.~N. Cattafesta}, {\em
  De-biasing the dynamic mode decomposition for applied {K}oopman spectral
  analysis}, Theoretical and Computational Fluid Dynamics, 31 (2017),
  pp.~349--368.

\bibitem{hemati2014pof}
{\sc M.~S. Hemati, M.~O. Williams, and C.~W. Rowley}, {\em Dynamic mode
  decomposition for large and streaming datasets}, Physics of Fluids
  (1994-present), 26 (2014), p.~111701.

\bibitem{hermann1977ieee}
{\sc R.~Hermann and A.~J. Krener}, {\em Nonlinear controllability and
  observability}, IEEE Transactions on Automatic Control, 22 (1977),
  pp.~728--740.

\bibitem{Herrmann2020arxiv}
{\sc B.~Herrmann, P.~J. Baddoo, R.~Semaan, S.~L. Brunton, and B.~J. McKeon},
  {\em Data-driven resolvent analysis}, Journal of Fluid Mechanics, 918 (2021).

\bibitem{hey2009fourth}
{\sc T.~Hey, S.~Tansley, K.~M. Tolle, et~al.}, {\em The fourth paradigm:
  data-intensive scientific discovery}, vol.~1, Microsoft research Redmond, WA,
  2009.

\bibitem{hirsh2019}
{\sc S.~M. Hirsh, K.~D. Harris, J.~N. Kutz, and B.~W. Brunton}, {\em Centering
  data improves the dynamic mode decomposition}, SIAM Journal on Applied
  Dynamical Systems, 19 (2020), pp.~1920--1955.

\bibitem{ho1965aac}
{\sc B.~L. Ho and R.~E. Kalman}, {\em Effective construction of linear
  state-variable models from input/output data}, in Proceedings of the 3rd
  Annual Allerton Conference on Circuit and System Theory, 1965, pp.~449--459.

\bibitem{hof1998}
{\sc A.~Hof and O.~Knill}, {\em Zero-dimensional singular continuous spectrum
  for smooth differential equations on the torus}, Ergodic Theory and Dynamical
  Systems, 18 (1998), pp.~879--888, \url{https://doi.org/null}.

\bibitem{hogg2019plosone}
{\sc J.~Hogg, M.~Fonoberova, I.~Mezi\'c, and R.~Mohr}, {\em Koopman mode
  analysis of agent-based models of logistics processes}, PloS one, 14 (2019),
  p.~e0222023.

\bibitem{HLBR_turb}
{\sc P.~J. Holmes, J.~L. Lumley, G.~Berkooz, and C.~W. Rowley}, {\em
  Turbulence, coherent structures, dynamical systems and symmetry}, Cambridge
  Monographs in Mechanics, Cambridge University Press, Cambridge, England,
  2nd~ed., 2012.

\bibitem{hopf50}
{\sc E.~Hopf}, {\em The partial differential equation $u_t + uu_x = \mu
  u_{xx}$}, Comm. Pure App. Math., 3 (1950), pp.~201--230.

\bibitem{huang2018cdc}
{\sc B.~Huang, X.~Ma, and U.~Vaidya}, {\em Feedback stabilization using
  {K}oopman operator}, in Decision and Control (CDC), 2018 IEEE 58th Annual
  Conference on, IEEE, 2018.

\bibitem{huang2019arxiv}
{\sc B.~Huang, X.~Ma, and U.~Vaidya}, {\em Data-driven nonlinear stabilization
  using {K}oopman operator}, The Koopman Operator in Systems and Control,
  (2020), pp.~313--334.

\bibitem{huang2013analysis}
{\sc C.~Huang, W.~E. Anderson, M.~E. Harvazinski, and V.~Sankaran}, {\em
  Analysis of self-excited combustion instabilities using decomposition
  techniques}, in 51st AIAA Aerospace Sciences Meeting, 2013, pp.~1--18.

\bibitem{hunt1982ieee}
{\sc L.~Hunt}, {\em Sufficient conditions for controllability}, IEEE
  Transactions on Circuits and Systems, 29 (1982), pp.~285--288.

\bibitem{iungo2015jp}
{\sc G.~V. Iungo, C.~Santoni-Ortiz, M.~Abkar, F.~Port{\'e}-Agel, M.~A. Rotea,
  and S.~Leonardi}, {\em Data-driven reduced order model for prediction of wind
  turbine wakes}, in Journal of Physics: Conference Series, vol.~625, IOP
  Publishing, 2015, p.~012009.

\bibitem{morton2019neurips}
{\sc M.~K. J~Morton, A~Jameson and F.~Witherden}, {\em Deep dynamical modeling
  and control of unsteady fluid flows}, in Advances in Neural Information
  Processing Systems 31 (NeurIPS 2018), 2019.

\bibitem{witten}
{\sc G.~James, D.~Witten, T.~Hastie, and R.~Tibshirani}, {\em \textsl{An
  Introduction to Statistical Learning}}, Springer, 2013.

\bibitem{pcr}
{\sc I.~T. Jolliffe}, {\em \textsl{A note on the Use of Principal Components in
  Regression}}, J. Roy. Stat. Soc. C, 31 (1982), pp.~300--303.

\bibitem{Jovanovic2020}
{\sc M.~R. Jovanovi{\'{c}}}, {\em {From bypass transition to flow control and
  data-driven turbulence modeling: An input-output viewpoint}}, Annu. Rev.
  Fluid. Mech., 53 (2021).

\bibitem{Jovanovic2005}
{\sc M.~R. Jovanovi{\'{c}} and B.~Bamieh}, {\em {Componentwise energy
  amplification in channel flows}}, J. Fluid Mech., 534 (2005), pp.~145--183.

\bibitem{jovanovic2014pof}
{\sc M.~R. Jovanovi{\'c}, P.~J. Schmid, and J.~W. Nichols}, {\em
  Sparsity-promoting dynamic mode decomposition}, Physics of Fluids, 26 (2014),
  p.~024103.

\bibitem{juang1994book}
{\sc J.~N. Juang}, {\em Applied System Identification}, Prentice Hall PTR,
  Upper Saddle River, New Jersey, 1994.

\bibitem{juang1985jgcd}
{\sc J.~N. Juang and R.~S. Pappa}, {\em An eigensystem realization algorithm
  for modal parameter identification and model reduction}, Journal of Guidance,
  Control, and Dynamics, 8 (1985), pp.~620--627.

\bibitem{juang1991nasatm}
{\sc J.~N. Juang, M.~Phan, L.~G. Horta, and R.~W. Longman}, {\em Identification
  of observer/{Kalman} filter {Markov} parameters: Theory and experiments},
  Technical Memorandum 104069, 1991.

\bibitem{kachurovskii1996}
{\sc A.~G. Kachurovski{\u \i}}, {\em Rates of convergence in ergodic theorems},
  Rossi\textbackslash u\i{} skaya Akademiya Nauk. Moskovskoe Matematicheskoe
  Obshchestvo. Uspekhi Matematicheskikh Nauk, 51 (1996), pp.~73--124,
  \url{https://doi.org/10.1070/RM1996v051n04ABEH002964}.

\bibitem{kakubr2018arxiv}
{\sc E.~Kaiser, J.~N. Kutz, and S.~L. Brunton}, {\em Discovering conservation
  laws from data for control}, in 2018 IEEE Conference on Decision and Control
  (CDC), IEEE, 2018, pp.~6415--6421.

\bibitem{kaiser2017arxivb}
{\sc E.~Kaiser, J.~N. Kutz, and S.~L. Brunton}, {\em Sparse identification of
  nonlinear dynamics for model predictive control in the low-data limit},
  Proceedings of the Royal Society of London A, 474 (2018).

\bibitem{Kaiser2020chapter}
{\sc E.~Kaiser, J.~N. Kutz, and S.~L. Brunton}, {\em Data-driven approximations
  of dynamical systems operators for control}, in The {{Koopman Operator}} in
  {{Systems}} and {{Control}}: {{Concepts}}, {{Methodologies}}, and
  {{Applications}}, A.~Mauroy, I.~Mezi{\'c}, and Y.~Susuki, eds., Lecture
  {{Notes}} in {{Control}} and {{Information Sciences}}, {Springer
  International Publishing}, {Cham}, 2020, pp.~197--234,
  \url{https://doi.org/10.1007/978-3-030-35713-9_8}.

\bibitem{kaiser2017arxiv}
{\sc E.~Kaiser, J.~N. Kutz, and S.~L. Brunton}, {\em Data-driven discovery of
  {K}oopman eigenfunctions for control}, Machine Learning: Science and
  Technology, 2 (2021), p.~035023.

\bibitem{kaiser2018jcp}
{\sc E.~Kaiser, M.~Morzy{\'n}ski, G.~Daviller, J.~N. Kutz, B.~W. Brunton, and
  S.~L. Brunton}, {\em Sparsity enabled cluster reduced-order models for
  control}, Journal of Computational Physics, 352 (2018), pp.~388--409.

\bibitem{kalman1960general}
{\sc R.~E. Kalman}, {\em On the general theory of control systems}, in
  Proceedings First International Conference on Automatic Control, Moscow,
  USSR, 1960.

\bibitem{kalman1963mathematical}
{\sc R.~E. Kalman}, {\em Mathematical description of linear dynamical systems},
  Journal of the Society for Industrial and Applied Mathematics, Series A:
  Control, 1 (1963), pp.~152--192.

\bibitem{Kamb2020siads}
{\sc M.~Kamb, E.~Kaiser, S.~L. Brunton, and J.~N. Kutz}, {\em Time-delay
  observables for {K}oopman: Theory and applications}, SIAM J. Appl. Dyn.
  Syst., 19 (2020), pp.~886--917.

\bibitem{kaptanoglu2020pop}
{\sc A.~A. Kaptanoglu, K.~D. Morgan, C.~J. Hansen, and S.~L. Brunton}, {\em
  Characterizing magnetized plasmas with dynamic mode decomposition}, Physics
  of Plasmas, 27 (2020), p.~032108.

\bibitem{katayama.2005}
{\sc T.~Katayama}, {\em Subspace Methods for System Identification},
  Springer-Verlag London, 2005.

\bibitem{katok1967}
{\sc A.~B. Katok and A.~M. Stepin}, {\em Approximations {{In Ergodic Theory}}},
  Russian Mathematical Surveys, 22 (1967), p.~77,
  \url{https://doi.org/10.1070/RM1967v022n05ABEH001227}.

\bibitem{kawahara2016neurips}
{\sc Y.~Kawahara}, {\em Dynamic mode decomposition with reproducing kernels for
  {K}oopman spectral analysis}, in Advances in Neural Information Processing
  Systems, 2016, pp.~911--919.

\bibitem{khalil1996book}
{\sc H.~K. Khalil}, {\em Noninear Systems}, Prentice-Hall, New Jersey, 1996.

\bibitem{kloeden_rasmussen_2011}
{\sc P.~E. Kloeden and M.~Rasmussen}, {\em Nonautonomous Dynamical Systems},
  Mathematical Surveys and Monographs, {American Mathematical Society},
  {Providence, R.I}, 2011.

\bibitem{klus2018tensor}
{\sc S.~Klus, P.~Gel{\ss}, S.~Peitz, and C.~Sch{\"u}tte}, {\em Tensor-based
  dynamic mode decomposition}, Nonlinearity, 31 (2018), p.~3359.

\bibitem{klus2015numerical}
{\sc S.~Klus, P.~Koltai, and C.~Sch{\"u}tte}, {\em On the numerical
  approximation of the {P}erron-{F}robenius and {K}oopman operator}, Journal of
  Computational Dynamics, 3 (2016), pp.~51--79.

\bibitem{klus2020kernel}
{\sc S.~Klus, F.~N{\"u}ske, and B.~Hamzi}, {\em Kernel-based approximation of
  the koopman generator and schr{\"o}dinger operator}, Entropy, 22 (2020),
  p.~722.

\bibitem{klus2017data}
{\sc S.~Klus, F.~N{\"u}ske, P.~Koltai, H.~Wu, I.~Kevrekidis, C.~Sch{\"u}tte,
  and F.~No{\'e}}, {\em Data-driven model reduction and transfer operator
  approximation}, Journal of Nonlinear Science, 28 (2018), pp.~985--1010.

\bibitem{klus2020}
{\sc S.~Klus, F.~N{\"{u}}ske, S.~Peitz, J.~H. Niemann, C.~Clementi, and
  C.~Sch{\"{u}}tte}, {\em {Data-driven approximation of the {K}oopman
  generator: Model reduction, system identification, and control}}, Phys. D
  Nonlinear Phenom., 406 (2020), pp.~1--32,
  \url{https://doi.org/10.1016/j.physd.2020.132416}.

\bibitem{klus2020eigendecompositions}
{\sc S.~Klus, I.~Schuster, and K.~Muandet}, {\em Eigendecompositions of
  transfer operators in reproducing kernel hilbert spaces}, Journal of
  Nonlinear Science, 30 (2020), pp.~283--315.

\bibitem{knill1998}
{\sc O.~Knill}, {\em Singular continuous spectrum and quantitative rates of
  weak mixing}, Discrete and Continuous Dynamical Systems, 4 (1998),
  pp.~33--42, \url{https://doi.org/10.3934/dcds.1998.4.33}.

\bibitem{kochkov2021machine}
{\sc D.~Kochkov, J.~A. Smith, A.~Alieva, Q.~Wang, M.~P. Brenner, and S.~Hoyer},
  {\em Machine learning accelerated computational fluid dynamics}, arXiv
  preprint arXiv:2102.01010,  (2021).

\bibitem{koltai2019}
{\sc P.~Koltai, H.~C. Lie, and M.~Plonka}, {\em Fr{\'e}chet differentiable
  drift dependence of {{Perron}}\textendash{{Frobenius}} and {{Koopman}}
  operators for non-deterministic dynamics}, Nonlinearity, 32 (2019),
  pp.~4232--4257, \url{https://doi.org/10.1088/1361-6544/ab1f2a}.

\bibitem{koopman1931pnas}
{\sc B.~O. Koopman}, {\em Hamiltonian systems and transformation in {H}ilbert
  space}, Proceedings of the National Academy of Sciences, 17 (1931),
  pp.~315--318.

\bibitem{koopman1932pnas}
{\sc B.~O. Koopman and J.-v. Neumann}, {\em Dynamical systems of continuous
  spectra}, Proceedings of the National Academy of Sciences, 18 (1932), p.~255.

\bibitem{korda2016arxiv}
{\sc M.~Korda and I.~Mezi{\'c}}, {\em Linear predictors for nonlinear dynamical
  systems: {K}oopman operator meets model predictive control}, Automatica, 93
  (2018), pp.~149--160.

\bibitem{korda2017arxiv}
{\sc M.~Korda and I.~Mezi{\'c}}, {\em On convergence of extended dynamic mode
  decomposition to the {K}oopman operator}, Journal of Nonlinear Science, 28
  (2018), pp.~687--710.

\bibitem{korda2020ieee}
{\sc M.~Korda and I.~Mezi\'c}, {\em {Optimal construction of {K}oopman
  eigenfunctions for prediction and control}}, EEE Transactions on Automatic
  Control,  (2020), p.~1.

\bibitem{korda2020}
{\sc M.~Korda, M.~Putinar, and I.~Mezi{\'c}}, {\em Data-driven spectral
  analysis of the {{Koopman}} operator}, Applied and Computational Harmonic
  Analysis, 48 (2020), pp.~599--629,
  \url{https://doi.org/10.1016/j.acha.2018.08.002}.

\bibitem{korda2018arxiv}
{\sc M.~Korda, Y.~Susuki, and I.~Mezi{\'c}}, {\em Power grid transient
  stabilization using {K}oopman model predictive control}, IFAC-PapersOnLine,
  51 (2018), pp.~297--302.

\bibitem{kowalski1991nonlinear}
{\sc K.~Kowalski, W.-H. Steeb, and K.~Kowalksi}, {\em Nonlinear dynamical
  systems and {C}arleman linearization}, World Scientific, 1991.

\bibitem{kramer2015arxiv}
{\sc B.~Kramer, P.~Grover, P.~Boufounos, M.~Benosman, and S.~Nabi}, {\em Sparse
  sensing and {DMD} based identification of flow regimes and bifurcations in
  complex flows}, SIAM Journal on Applied Dynamical Systems, 16 (2017),
  pp.~1164--1196.

\bibitem{krener1974}
{\sc A.~Krener}, {\em Linearization and bilinearization of control systems}, in
  Proc. 1974 Allerton Conf. on Circuit and System Theory, Monticello, vol.~834,
  Springer, Berlin, 1974.

\bibitem{Kutz2017jfm}
{\sc J.~N. Kutz}, {\em Deep learning in fluid dynamics}, Journal of Fluid
  Mechanics, 814 (2017), pp.~1--4.

\bibitem{kutz2016book}
{\sc J.~N. Kutz, S.~L. Brunton, B.~W. Brunton, and J.~L. Proctor}, {\em Dynamic
  Mode Decomposition: Data-Driven Modeling of Complex Systems}, SIAM, 2016.

\bibitem{kutz2016siads}
{\sc J.~N. Kutz, X.~Fu, and S.~L. Brunton}, {\em Multi-resolution dynamic mode
  decomposition}, SIAM Journal on Applied Dynamical Systems, 15 (2016),
  pp.~713--735.

\bibitem{kutz2018}
{\sc J.~N. Kutz, J.~L. Proctor, and S.~L. Brunton}, {\em Applied {{Koopman
  Theory}} for {{Partial Differential Equations}} and {{Data}}-{{Driven
  Modeling}} of {{Spatio}}-{{Temporal Systems}}}, Complexity,  (2018),
  p.~e6010634, \url{https://doi.org/10.1155/2018/6010634}.

\bibitem{kvalheim2021a}
{\sc M.~D. Kvalheim, D.~Hong, and S.~Revzen}, {\em Generic {{Properties}} of
  {{Koopman Eigenfunctions}} for {{Stable Fixed Points}} and {{Periodic
  Orbits}}}, IFAC-PapersOnLine, 54 (2021), pp.~267--272,
  \url{https://doi.org/10.1016/j.ifacol.2021.06.150}.

\bibitem{kvalheim2021}
{\sc M.~D. Kvalheim and S.~Revzen}, {\em Existence and uniqueness of global
  {{Koopman}} eigenfunctions for stable fixed points and periodic orbits},
  Physica D: Nonlinear Phenomena, 425 (2021), p.~132959,
  \url{https://doi.org/10.1016/j.physd.2021.132959}.

\bibitem{kwee2001}
{\sc I.~Kwee and J.~Schmidhuber}, {\em Optimal control using the transport
  equation: the {L}iouville machine}, Adaptive Behavior, 9 (2001).

\bibitem{lan2013physd}
{\sc Y.~Lan and I.~Mezi{\'c}}, {\em Linearization in the large of nonlinear
  systems and {K}oopman operator spectrum}, Physica D: Nonlinear Phenomena, 242
  (2013), pp.~42--53.

\bibitem{lange2020fourier}
{\sc H.~Lange, S.~L. Brunton, and J.~N. Kutz}, {\em From {F}ourier to
  {K}oopman: Spectral methods for long-term time series prediction.}, J. Mach.
  Learn. Res., 22 (2021), pp.~1--38.

\bibitem{lasota1994}
{\sc A.~Lasota and M.~C. Mackey}, {\em Chaos, {{Fractals}}, and {{Noise}}:
  {{Stochastic Aspects}} of {{Dynamics}}}, vol.~97 of Applied {{Mathematical
  Sciences}}, {Springer-Verlag}, {New York}, second~ed., 1994.

\bibitem{law2015data}
{\sc K.~Law, A.~Stuart, and K.~Zygalakis}, {\em Data assimilation}, Cham,
  Switzerland: Springer,  (2015).

\bibitem{lax2002}
{\sc P.~Lax}, {\em Functional {{Analysis}}}, {Wiley and Sons}, 2002.

\bibitem{lax1968integrals}
{\sc P.~D. Lax}, {\em Integrals of nonlinear equations of evolution and
  solitary waves}, Communications on pure and applied mathematics, 21 (1968),
  pp.~467--490.

\bibitem{lax1971}
{\sc P.~D. Lax}, {\em Approximation of measure preserving transformations},
  Communications on Pure and Applied Mathematics, 24 (1971), pp.~133--135,
  \url{https://doi.org/10.1002/cpa.3160240204}.

\bibitem{lee2011springer}
{\sc J.~H. Lee}, {\em Model predictive control: Review of the three decades of
  development}, International Journal of Control, Automation and Systems, 9
  (2011), pp.~415--424.

\bibitem{lee2020model}
{\sc K.~Lee and K.~T. Carlberg}, {\em Model reduction of dynamical systems on
  nonlinear manifolds using deep convolutional autoencoders}, Journal of
  Computational Physics, 404 (2020), p.~108973.

\bibitem{levnajic2010}
{\sc Z.~Levnaji{\'c} and I.~Mezi{\'c}}, {\em Ergodic theory and visualization.
  {{I}}. {{Mesochronic}} plots for visualization of ergodic partition and
  invariant sets}, Chaos: An Interdisciplinary Journal of Nonlinear Science, 20
  (2010), pp.~--, \url{https://doi.org/10.1063/1.3458896}.

\bibitem{levnajic2015}
{\sc Z.~Levnaji{\'c} and I.~Mezi{\'c}}, {\em Ergodic theory and visualization.
  {{II}}. {{Fourier}} mesochronic plots visualize (quasi)periodic sets}, Chaos:
  An Interdisciplinary Journal of Nonlinear Science, 25 (2015), p.~053105,
  \url{https://doi.org/10.1063/1.4919767}.

\bibitem{li2017chaos}
{\sc Q.~Li, F.~Dietrich, E.~M. Bollt, and I.~G. Kevrekidis}, {\em Extended
  dynamic mode decomposition with dictionary learning: A data-driven adaptive
  spectral decomposition of the {K}oopman operator}, Chaos: An
  Interdisciplinary Journal of Nonlinear Science, 27 (2017), p.~103111.

\bibitem{li1976}
{\sc T.~Y. Li}, {\em Finite approximation for the {{Frobenius}}-{{Perron}}
  operator. {{A}} solution to {{Ulam}}'s conjecture}, Journal of Approximation
  Theory, 17 (1976), pp.~177--186.

\bibitem{li2019arxiv}
{\sc Y.~Li, H.~He, J.~Wu, D.~Katabi, and A.~Torralba}, {\em Learning
  compositional {K}oopman operators for model-based control}, ICLR Eighth
  International Conference on Learning Representations,  (2019).

\bibitem{li2020fourier}
{\sc Z.~Li, N.~Kovachki, K.~Azizzadenesheli, B.~Liu, K.~Bhattacharya,
  A.~Stuart, and A.~Anandkumar}, {\em Fourier neural operator for parametric
  partial differential equations}, arXiv preprint arXiv:2010.08895,  (2020).

\bibitem{li2020multipole}
{\sc Z.~Li, N.~Kovachki, K.~Azizzadenesheli, B.~Liu, K.~Bhattacharya,
  A.~Stuart, and A.~Anandkumar}, {\em Multipole graph neural operator for
  parametric partial differential equations}, arXiv preprint arXiv:2006.09535,
  (2020).

\bibitem{li2020neural}
{\sc Z.~Li, N.~Kovachki, K.~Azizzadenesheli, B.~Liu, K.~Bhattacharya,
  A.~Stuart, and A.~Anandkumar}, {\em Neural operator: Graph kernel network for
  partial differential equations}, arXiv preprint arXiv:2003.03485,  (2020).

\bibitem{ling2018ieee}
{\sc E.~Ling, L.~Ratliff, and S.~Coogan}, {\em Koopman operator approach for
  instability detection and mitigation in signalized traffic}, in 21st
  International Conference on Intelligent Transportation Systems (ITSC), IEEE,
  2018, pp.~1297--1302.

\bibitem{Ling2016jfm}
{\sc J.~Ling, A.~Kurzawski, and J.~Templeton}, {\em Reynolds averaged
  turbulence modelling using deep neural networks with embedded invariance},
  Journal of Fluid Mechanics, 807 (2016), pp.~155--166.

\bibitem{liu2017arxiv}
{\sc Z.~Liu, S.~Kundu, L.~Chen, and E.~Yeung}, {\em Decomposition of nonlinear
  dynamical systems using {K}oopman gramians},  (2018), pp.~4811--4818.

\bibitem{ljung1999book}
{\sc L.~Ljung}, {\em System Identification: Theory for the User}, Pearson,
  2nd~ed., 1999.

\bibitem{loparo1978tac}
{\sc K.~Loparo and G.~Blankenship}, {\em Estimating the domain of attraction of
  nonlinear feedback systems}, IEEE Transactions on Automatic Control, 23
  (1978), pp.~602--608.

\bibitem{eof1}
{\sc E.~N. Lorenz}, {\em \textsl{Empirical orthogonal functions and statistical
  weather prediction}}, Technical report, Massachusetts Institute of
  Technology, Dec. (1956).

\bibitem{lu2021learning}
{\sc L.~Lu, P.~Jin, G.~Pang, Z.~Zhang, and G.~E. Karniadakis}, {\em Learning
  nonlinear operators via deeponet based on the universal approximation theorem
  of operators}, Nature Machine Intelligence, 3 (2021), pp.~218--229.

\bibitem{lusch2017arxiv}
{\sc B.~Lusch, J.~N. Kutz, and S.~L. Brunton}, {\em Deep learning for universal
  linear embeddings of nonlinear dynamics}, Nature Communications, 9 (2018),
  p.~4950.

\bibitem{lusseyran2011flow}
{\sc F.~Lusseyran, F.~Gueniat, J.~Basley, C.~L. Douay, L.~R. Pastur, T.~M.
  Faure, and P.~J. Schmid}, {\em Flow coherent structures and frequency
  signature: application of the dynamic modes decomposition to open cavity
  flow}, in Journal of Physics: Conference Series, vol.~318, IOP Publishing,
  2011, p.~042036.

\bibitem{Luzzatto2005}
{\sc S.~Luzzatto, I.~Melbourne, and F.~Paccaut}, {\em The {{Lorenz Attractor}}
  is {{Mixing}}}, Communications in Mathematical Physics, 260 (2005),
  pp.~393--401, \url{https://doi.org/10.1007/s00220-005-1411-9}.

\bibitem{ma2011tcfd}
{\sc Z.~Ma, S.~Ahuja, and C.~W. Rowley}, {\em Reduced order models for control
  of fluids using the eigensystem realization algorithm}, Theor. Comput. Fluid
  Dyn., 25 (2011), pp.~233--247.

\bibitem{macesic2020}
{\sc S.~Ma{\'c}e{\v s}i{\'c} and N.~{\v C}rnjari{\'c}-{\v Z}ic}, {\em Koopman
  {{Operator Theory}} for {{Nonautonomous}} and {{Stochastic Systems}}}, in The
  {{Koopman Operator}} in {{Systems}} and {{Control}}: {{Concepts}},
  {{Methodologies}}, and {{Applications}}, A.~Mauroy, I.~Mezi{\'c}, and
  Y.~Susuki, eds., Lecture {{Notes}} in {{Control}} and {{Information
  Sciences}}, {Springer International Publishing}, {Cham}, 2020, pp.~131--160,
  \url{https://doi.org/10.1007/978-3-030-35713-9_6}.

\bibitem{macesic2018siads}
{\sc S.~Ma{\'c}e{\v{s}}i{\'c}, N.~{\v{C}}rnjari\'c-{\v{Z}}ic, and
  I.~Mezi{\'c}}, {\em Koopman operator family spectrum for nonautonomous
  systems}, SIAM Journal on Applied Dynamical Systems, 17 (2018),
  pp.~2478--2515.

\bibitem{madrid2005}
{\sc R.~{\noopsort{madrid}}de~la Madrid}, {\em The role of the rigged
  {{Hilbert}} space in quantum mechanics}, European Journal of Physics, 26
  (2005), p.~287, \url{https://doi.org/10.1088/0143-0807/26/2/008}.

\bibitem{madrid2002a}
{\sc R.~{\noopsort{madrid}}de~la Madrid, A.~Bohm, and M.~Gadella}, {\em {Rigged
  Hilbert Space Treatment of Continuous Spectrum}}, Fortschritte der Physik, 50
  (2002), pp.~185--216,
  \url{https://doi.org/10.1002/1521-3978(200203)50:2<185::AID-PROP185>3.0.CO;2-S}.

\bibitem{majda2012nonlinearity}
{\sc A.~J. Majda and J.~Harlim}, {\em Physics constrained nonlinear regression
  models for time series}, Nonlinearity, 26 (2012), p.~201.

\bibitem{majda2014pnas}
{\sc A.~J. Majda and Y.~Lee}, {\em Conceptual dynamical models for turbulence},
  Proceedings of the National Academy of Sciences, 111 (2014), pp.~6548--6553.

\bibitem{mamakoukas2020arxiv}
{\sc G.~Mamakoukas, M.~L. Castano, X.~Tan, and T.~D. Murphey}, {\em
  Derivative-based koopman operators for real-time control of robotic systems},
  IEEE Transactions on Robotics,  (2021).

\bibitem{mann2016qf}
{\sc J.~Mann and J.~N. Kutz}, {\em Dynamic mode decomposition for financial
  trading strategies}, Quantitative Finance,  (2016), pp.~1--13.

\bibitem{manohar2017csm}
{\sc K.~Manohar, B.~W. Brunton, J.~N. Kutz, and S.~L. Brunton}, {\em
  Data-driven sparse sensor placement}, IEEE Control Systems Magazine, 38
  (2018), pp.~63--86.

\bibitem{manohar2017arxiv}
{\sc K.~Manohar, E.~Kaiser, S.~L. Brunton, and J.~N. Kutz}, {\em Optimized
  sampling for multiscale dynamics}, SIAM Multiscale modeling and simulation,
  17 (2019), pp.~117--136.

\bibitem{manohar2018arxivb}
{\sc K.~Manohar, J.~N. Kutz, and S.~L. Brunton}, {\em Optimal sensor and
  actuator placement using balanced model reduction}, arXiv preprint arXiv:
  1812.01574,  (2018).

\bibitem{manojlovic2020applications}
{\sc I.~Manojlovi{\'c}, M.~Fonoberova, R.~Mohr, A.~Andrej{\v{c}}uk,
  Z.~Drma{\v{c}}, Y.~Kevrekidis, and I.~Mezi{\'c}}, {\em Applications of
  koopman mode analysis to neural networks}, arXiv preprint arXiv:2006.11765,
  (2020).

\bibitem{mardt2020deep}
{\sc A.~Mardt, L.~Pasquali, F.~No{\'e}, and H.~Wu}, {\em Deep learning markov
  and koopman models with physical constraints}, in Mathematical and Scientific
  Machine Learning, PMLR, 2020, pp.~451--475.

\bibitem{mardt2017arxiv}
{\sc A.~Mardt, L.~Pasquali, H.~Wu, and F.~No{\'e}}, {\em {VAMP}nets: Deep
  learning of molecular kinetics}, Nature Communications, 9 (2018).

\bibitem{markus1960}
{\sc L.~Markus and H.~Yamabe}, {\em Global stability criteria for differential
  systems}, Osaka Mathematical Journal, 12 (1960), pp.~305--317.

\bibitem{marrouch2020data}
{\sc N.~Marrouch, J.~Slawinska, D.~Giannakis, and H.~L. Read}, {\em Data-driven
  koopman operator approach for computational neuroscience}, Annals of
  Mathematics and Artificial Intelligence, 88 (2020), pp.~1155--1173.

\bibitem{massa2012dynamic}
{\sc L.~Massa, R.~Kumar, and P.~Ravindran}, {\em Dynamic mode decomposition
  analysis of detonation waves}, Physics of Fluids (1994-present), 24 (2012),
  p.~066101.

\bibitem{Maulik2019jfm}
{\sc R.~Maulik, O.~San, A.~Rasheed, and P.~Vedula}, {\em Subgrid modelling for
  two-dimensional turbulence using neural networks}, Journal of Fluid
  Mechanics, 858 (2019), pp.~122--144.

\bibitem{mauroy2020tac}
{\sc A.~Mauroy}, {\em Koopman-based lifting techniques for nonlinear systems
  identification}, IEEE Transactions on Automatic Control, 65 (2020).

\bibitem{mauroy2017arxiv}
{\sc A.~Mauroy and J.~Goncalves}, {\em Koopman-based lifting techniques for
  nonlinear systems identification}, IEEE Transactions on Automatic Control, 65
  (2019), pp.~2550--2565.

\bibitem{mauroy2012}
{\sc A.~Mauroy and I.~Mezi{\'c}}, {\em On the use of {{Fourier}} averages to
  compute the global isochrons of (quasi)periodic dynamics}, Chaos: An
  Interdisciplinary Journal of Nonlinear Science, 22 (2012), p.~033112,
  \url{https://doi.org/10.1063/1.4736859}.

\bibitem{mauroy2013cdc}
{\sc A.~Mauroy and I.~Mezic}, {\em A spectral operator-theoretic framework for
  global stability}, in Decision and Control (CDC), 2013 IEEE 52nd Annual
  Conference on, IEEE, 2013, pp.~5234--5239.

\bibitem{mauroy2016ieeetac}
{\sc A.~Mauroy and I.~Mezi{\'c}}, {\em Global stability analysis using the
  eigenfunctions of the {K}oopman operator}, IEEE Transactions on Automatic
  Control, 61 (2016), pp.~3356--3369.

\bibitem{mauroy2013}
{\sc A.~Mauroy, I.~Mezic, and J.~Moehlis}, {\em Isostables, isochrons, and
  {{Koopman}} spectrum for the action-angle representation of stable fixed
  point dynamics}, Physica D: Nonlinear Phenomena, 261 (2013), pp.~19--30,
  \url{https://doi.org/10.1016/j.physd.2013.06.004}.

\bibitem{mauroy2020book}
{\sc A.~Mauroy, I.~Mezi\'c, and Y.~Susuki}, eds., {\em The {K}oopman Operator
  in Systems and Control: Concepts, Methodologies, and Applications}, Springer,
  2020.

\bibitem{mayne1997proc}
{\sc D.~Q. Mayne}, {\em Nonlinear model predictive control: An assessment}, in
  Proc. Chemical Process Control, V.~J.~C. Kantor, C.~E. Garcia, and
  B.~Carnahan, eds., AIchE, 1997, pp.~217--231.

\bibitem{McKeon2010b}
{\sc B.~J. McKeon and A.~S. Sharma}, {\em {A critical-layer framework for
  turbulent pipe flow}}, J. Fluid Mech., 658 (2010), pp.~336--382.

\bibitem{mehta2018}
{\sc P.~M. Mehta and R.~Linares}, {\em A {{New Transformative Framework}} for
  {{Data Assimilation}} and {{Calibration}} of {{Physical
  Ionosphere}}-{{Thermosphere Models}}}, Space Weather, 16 (2018),
  pp.~1086--1100, \url{https://doi.org/10.1029/2018SW001875}.

\bibitem{meiss2007}
{\sc J.~D. Meiss}, {\em Differential {{Dynamical Systems}}}, {SIAM}, Jan. 2007.

\bibitem{mesbahi2019}
{\sc A.~Mesbahi, J.~Bu, and M.~Mesbahi}, {\em On {{Modal Properties}} of the
  {{Koopman Operator}} for {{Nonlinear Systems}} with {{Symmetry}}}, in 2019
  {{American Control Conference}} ({{ACC}}), July 2019, pp.~1918--1923,
  \url{https://doi.org/10.23919/ACC.2019.8815342}.

\bibitem{mezic2005nd}
{\sc I.~Mezi{\'c}}, {\em Spectral properties of dynamical systems, model
  reduction and decompositions}, Nonlinear Dynamics, 41 (2005), pp.~309--325.

\bibitem{mezic2013arfm}
{\sc I.~Mezic}, {\em Analysis of fluid flows via spectral properties of the
  {K}oopman operator}, Annual Review of Fluid Mechanics, 45 (2013),
  pp.~357--378.

\bibitem{mezic2019}
{\sc I.~Mezi{\'c}}, {\em Spectrum of the {{Koopman Operator}}, {{Spectral
  Expansions}} in {{Functional Spaces}}, and {{State}}-{{Space Geometry}}},
  Journal of Nonlinear Science,  (2019),
  \url{https://doi.org/10.1007/s00332-019-09598-5}.

\bibitem{mezic2020}
{\sc I.~Mezic}, {\em On {{Numerical Approximations}} of the {{Koopman
  Operator}}}, arXiv:2009.05883 [math],  (2020),
  \url{https://arxiv.org/abs/2009.05883}.

\bibitem{mezic2004physicad}
{\sc I.~Mezi{\'c} and A.~Banaszuk}, {\em Comparison of systems with complex
  behavior}, Physica D: Nonlinear Phenomena, 197 (2004), pp.~101--133.

\bibitem{mezic2004b}
{\sc I.~Mezi{\'c} and T.~Runolfsson}, {\em Uncertainty analysis of complex
  dynamical systems}, American Control Conference,  (2004).

\bibitem{mezic2008c}
{\sc I.~Mezi{\'c} and T.~Runolfsson}, {\em Uncertainty propagation in dynamical
  systems}, Automatica, 44 (2008), pp.~3003--3013,
  \url{https://doi.org/10.1016/j.automatica.2008.04.020}.

\bibitem{mezic1999chaos}
{\sc I.~Mezi{\'c} and S.~Wiggins}, {\em A method for visualization of invariant
  sets of dynamical systems based on the ergodic partition}, Chaos: An
  Interdisciplinary Journal of Nonlinear Science, 9 (1999), pp.~213--218.

\bibitem{mizuno2011investigation}
{\sc Y.~Mizuno, D.~Duke, C.~Atkinson, and J.~Soria}, {\em Investigation of
  wall-bounded turbulent flow using dynamic mode decomposition}, in Journal of
  Physics: Conference Series, vol.~318, IOP Publishing, 2011, p.~042040.

\bibitem{moeck2013tomographic}
{\sc J.~P. Moeck, J.-F. Bourgouin, D.~Durox, T.~Schuller, and S.~Candel}, {\em
  Tomographic reconstruction of heat release rate perturbations induced by
  helical modes in turbulent swirl flames}, Experiments in Fluids, 54 (2013),
  pp.~1--17.

\bibitem{mohr2014arxiv}
{\sc R.~Mohr and I.~Mezi{\'c}}, {\em Construction of eigenfunctions for
  scalar-type operators via {L}aplace averages with connections to the
  {K}oopman operator}, arXiv preprint arXiv:1403.6559,  (2014).

\bibitem{mohr2016a}
{\sc R.~Mohr and I.~Mezi{\'c}}, {\em Koopman principle eigenfunctions and
  linearization of diffeomorphisms}, arXiv:1611.01209 [math],  (2016),
  \url{https://arxiv.org/abs/1611.01209}.

\bibitem{moore2015pnas}
{\sc C.~C. Moore}, {\em Ergodic theorem, ergodic theory, and statistical
  mechanics}, Proceedings of the National Academy of Sciences, 112 (2015),
  pp.~1907--1911.

\bibitem{moore1967}
{\sc J.~Moore and B.~Anderson}, {\em Optimal linear control systems with input
  derivative constraints}, Proceedings of the Institution of Electrical
  Engineers, 114 (1967), pp.~1987--1990.

\bibitem{morari1999model}
{\sc M.~Morari and J.~H. Lee}, {\em Model predictive control: past, present and
  future}, Computers \& Chemical Engineering, 23 (1999), pp.~667--682.

\bibitem{mukherjee2015}
{\sc A.~Mukherjee, R.~Rai, P.~Singla, T.~Singh, and A.~Patra}, {\em Laplacian
  graph based approach for uncertainty quantification of large scale dynamical
  systems}, in 2015 {{American Control Conference}} ({{ACC}}), July 2015,
  pp.~3998--4003, \url{https://doi.org/10.1109/ACC.2015.7171954}.

\bibitem{muld2012flow}
{\sc T.~W. Muld, G.~Efraimsson, and D.~S. Henningson}, {\em Flow structures
  around a high-speed train extracted using proper orthogonal decomposition and
  dynamic mode decomposition}, Computers \& Fluids, 57 (2012), pp.~87--97.

\bibitem{muld2012mode}
{\sc T.~W. Muld, G.~Efraimsson, and D.~S. Henningson}, {\em Mode decomposition
  on surface-mounted cube}, Flow, Turbulence and Combustion, 88 (2012),
  pp.~279--310.

\bibitem{nadkarni1979}
{\sc M.~G. Nadkarni}, {\em On {{Spectra}} of {{Non}}-{{Singular
  Transformations}} and {{Flows}}}, Sankhy\=a: The Indian Journal of
  Statistics, Series A (1961-2002), 41 (1979), pp.~59--66.

\bibitem{nadler2006acha}
{\sc B.~Nadler, S.~Lafon, R.~R. Coifman, and I.~G. Kevrekidis}, {\em Diffusion
  maps, spectral clustering and reaction coordinates of dynamical systems},
  Applied and Computational Harmonic Analysis, 21 (2006), pp.~113--127.

\bibitem{nakao2020}
{\sc H.~Nakao and I.~Mezi{\'c}}, {\em Spectral analysis of the {{Koopman}}
  operator for partial differential equations}, Chaos: An Interdisciplinary
  Journal of Nonlinear Science, 30 (2020), p.~113131,
  \url{https://doi.org/10.1063/5.0011470}.

\bibitem{nandanoori2020}
{\sc S.~P. Nandanoori, S.~Sinha, and E.~Yeung}, {\em Data-{{Driven Operator
  Theoretic Methods}} for {{Global Phase Space Learning}}}, in 2020 {{American
  Control Conference}} ({{ACC}}), July 2020, pp.~4551--4557,
  \url{https://doi.org/10.23919/ACC45564.2020.9147220}.

\bibitem{narasingam2019}
{\sc A.~Narasingam and J.~Kwon}, {\em {Koopman Lyapunov‐based model
  predictive control of nonlinear chemical process systems}}, AIChE Journal, 65
  (2019).

\bibitem{netto2018tps}
{\sc M.~Netto and L.~Mili}, {\em {A robust data-driven Koopman Kalman filter
  for power systems dynamic state estimation}}, IEEE Transactions on Power
  Systems, 33 (2018), pp.~7228--7237.

\bibitem{neumann1932pnas}
{\sc J.~v. Neumann}, {\em Proof of the quasi-ergodic hypothesis}, Proceedings
  of the National Academy of Sciences, 18 (1932), pp.~70--82.

\bibitem{vonneumann1932a}
{\sc J.~{\noopsort{neumann}}{von Neumann}}, {\em Zur {{Operatorenmethode}} in
  der klassischen {{Mechanik}}}, Annals of Mathematics, 33 (1932),
  pp.~587--642, \url{https://doi.org/10.2307/1968537}.

\bibitem{Noack2016jfm}
{\sc B.~R. Noack, W.~Stankiewicz, M.~Morzynski, and P.~J. Schmid}, {\em
  Recursive dynamic mode decomposition of a transient cylinder wake}, Journal
  of Fluid Mechanics, 809 (2016), pp.~843--872.

\bibitem{noe2013variational}
{\sc F.~No{\'e} and F.~Nuske}, {\em A variational approach to modeling slow
  processes in stochastic dynamical systems}, Multiscale Modeling \&
  Simulation, 11 (2013), pp.~635--655.

\bibitem{Noe2019science}
{\sc F.~No{\'e}, S.~Olsson, J.~K{\"o}hler, and H.~Wu}, {\em Boltzmann
  generators: Sampling equilibrium states of many-body systems with deep
  learning}, Science, 365 (2019), p.~eaaw1147.

\bibitem{noether1918invariante}
{\sc E.~Noether}, {\em Invariante variationsprobleme nachr. d. k{\"o}nig.
  gesellsch. d. wiss. zu g{\"o}ttingen, math-phys. klasse 1918: 235-257},
  (1918), p.~57.

\bibitem{nonomura2018dynamic}
{\sc T.~Nonomura, H.~Shibata, and R.~Takaki}, {\em Dynamic mode decomposition
  using a kalman filter for parameter estimation}, AIP Advances, 8 (2018),
  p.~105106.

\bibitem{nonomura2019extended}
{\sc T.~Nonomura, H.~Shibata, and R.~Takaki}, {\em Extended-kalman-filter-based
  dynamic mode decomposition for simultaneous system identification and
  denoising}, PloS one, 14 (2019), p.~e0209836.

\bibitem{nuske2019tensorbased}
{\sc F.~N{\"u}ske, P.~Gel{\ss}, S.~Klus, and C.~Clementi}, {\em Tensor based
  {EDMD} for the {K}oopman analysis of high-dimensional systems}, arXiv
  preprint arXiv:1908.04741,  (2019).

\bibitem{nuske2014jctc}
{\sc F.~N{\"u}ske, B.~G. Keller, G.~P{\'e}rez-Hern{\'a}ndez, A.~S. Mey, and
  F.~No{\'e}}, {\em Variational approach to molecular kinetics}, Journal of
  chemical theory and computation, 10 (2014), pp.~1739--1752.

\bibitem{nuske2016variational}
{\sc F.~N{\"u}ske, R.~Schneider, F.~Vitalini, and F.~No{\'e}}, {\em Variational
  tensor approach for approximating the rare-event kinetics of macromolecular
  systems}, J. Chem. Phys., 144 (2016), p.~054105.

\bibitem{Oppenheim2014}
{\sc A.~V. .-. Oppenheim and R.~W. .-. Schafer}, {\em Discrete-{{Time Signal
  Processing}}}, {Pearson Education Limited}, {Harlow}, 3rd ed.~ed., 2014.

\bibitem{ostoich2013interaction}
{\sc C.~M. Ostoich, D.~J. Bodony, and P.~H. Geubelle}, {\em Interaction of a
  {M}ach 2.25 turbulent boundary layer with a fluttering panel using direct
  numerical simulation}, Physics of Fluids (1994-present), 25 (2013),
  p.~110806.

\bibitem{otto2017arxiv}
{\sc S.~E. Otto and C.~W. Rowley}, {\em Linearly-recurrent autoencoder networks
  for learning dynamics}, SIAM Journal on Applied Dynamical Systems, 18 (2019),
  pp.~558--593.

\bibitem{otto2021koopman}
{\sc S.~E. Otto and C.~W. Rowley}, {\em Koopman operators for estimation and
  control of dynamical systems}, Annual Review of Control, Robotics, and
  Autonomous Systems, 4 (2021).

\bibitem{page2018a}
{\sc J.~Page and R.~R. Kerswell}, {\em Koopman analysis of {{Burgers}}
  equation}, Physical Review Fluids, 3 (2018), p.~071901,
  \url{https://doi.org/10.1103/PhysRevFluids.3.071901}.

\bibitem{page2019}
{\sc J.~Page and R.~R. Kerswell}, {\em Koopman mode expansions between simple
  invariant solutions}, Journal of Fluid Mechanics, 879 (2019), pp.~1--27,
  \url{https://doi.org/10.1017/jfm.2019.686}.

\bibitem{pan2011dynamical}
{\sc C.~Pan, D.~Yu, and J.~Wang}, {\em Dynamical mode decomposition of {G}urney
  flap wake flow}, Theoretical and Applied Mechanics Letters, 1 (2011),
  p.~012002.

\bibitem{Pan2020arxiv}
{\sc S.~Pan, N.~Arnold-Medabalimi, and K.~Duraisamy}, {\em Sparsity-promoting
  algorithms for the discovery of informative {K}oopman-invariant subspaces},
  Journal of Fluid Mechanics, 917 (2021).

\bibitem{parker2019koopman}
{\sc J.~P. Parker and J.~Page}, {\em Koopman analysis of isolated fronts and
  solitons}, SIAM Journal on Applied Dynamical Systems, 19 (2020),
  pp.~2803--2828.

\bibitem{peherstorfer2016data}
{\sc B.~Peherstorfer and K.~Willcox}, {\em Data-driven operator inference for
  nonintrusive projection-based model reduction}, Computer Methods in Applied
  Mechanics and Engineering, 306 (2016), pp.~196--215.

\bibitem{peitz2018arxiv}
{\sc S.~Peitz}, {\em Controlling nonlinear {P}{D}{E}s using low-dimensional
  bilinear approximations obtained from data}, arXiv preprint arXiv:1801.06419,
   (2018).

\bibitem{peitz2017arxiv}
{\sc S.~Peitz and S.~Klus}, {\em Koopman operator-based model reduction for
  switched-system control of {PDE}s}, Automatica, 106 (2019), pp.~184--191.

\bibitem{peitz2018feedback}
{\sc S.~Peitz and S.~Klus}, {\em Feedback control of nonlinear {PDE}s using
  data-efficient reduced order models based on the {K}oopman operator},
  (2020), pp.~257--282.

\bibitem{peitz2020arxiv}
{\sc S.~Peitz, S.~E. Otto, and C.~W. Rowley}, {\em Data-driven model predictive
  control using interpolated {K}oopman generators}, SIAM Journal on Applied
  Dynamical Systems, 19 (2020), pp.~2162--2193.

\bibitem{pendergrass2016arxiv}
{\sc S.~D. Pendergrass, J.~N. Kutz, and S.~L. Brunton}, {\em Streaming {GPU}
  singular value and dynamic mode decompositions}, arXiv preprint
  arXiv:1612.07875,  (2016).

\bibitem{eof2}
{\sc C.~Penland}, {\em \textsl{Random forcing and forecasting using Principal
  Oscillation Pattern analysis}}, Mon. Weather Rev., 117 (1989),
  pp.~2165--2185.

\bibitem{eof3}
{\sc C.~Penland and T.~Magorian}, {\em \textsl{Prediction of Ni\~{n}o 3
  sea-surface temperatures using linear inverse modeling}}, J. Climate, 6
  (1993), pp.~1067--1076.

\bibitem{perko2009}
{\sc L.~Perko}, {\em Differential Equations and Dynamical Systems},
  {Springer,}, {New York, NY ;}, 2009.

\bibitem{phan1993jota}
{\sc M.~Phan, L.~G. Horta, J.~N. Juang, and R.~W. Longman}, {\em Linear system
  identification via an asymptotically stable observer}, Journal of
  Optimization Theory and Applications, 79 (1993), pp.~59--86.

\bibitem{phan1992jas}
{\sc M.~Phan, J.~N. Juang, and R.~W. Longman}, {\em Identification of
  linear-multivariable systems by identification of observers with assigned
  real eigenvalues}, The Journal of the Astronautical Sciences, 40 (1992),
  pp.~261--279.

\bibitem{pikovsky1995}
{\sc A.~S. Pikovsky, M.~A. Zaks, U.~Feudel, and J.~Kurths}, {\em Singular
  continuous spectra in dissipative dynamics}, Physical Review E, 52 (1995),
  pp.~285--296, \url{https://doi.org/10.1103/PhysRevE.52.285}.

\bibitem{Pollicott1985}
{\sc M.~Pollicott}, {\em On the rate of mixing of {{Axiom A}} flows},
  Inventiones mathematicae, 81 (1985), pp.~413--426,
  \url{https://doi.org/10.1007/BF01388579}.

\bibitem{proctor2014epj}
{\sc J.~L. Proctor, S.~L. Brunton, B.~W. Brunton, and J.~N. Kutz}, {\em
  Exploiting sparsity and equation-free architectures in complex systems
  (invited review)}, The European Physical Journal Special Topics, 223 (2014),
  pp.~2665--2684.

\bibitem{proctor2016siads}
{\sc J.~L. Proctor, S.~L. Brunton, and J.~N. Kutz}, {\em Dynamic mode
  decomposition with control}, SIAM Journal on Applied Dynamical Systems, 15
  (2016), pp.~142--161.

\bibitem{proctor2017siads}
{\sc J.~L. Proctor, S.~L. Brunton, and J.~N. Kutz}, {\em Generalizing {K}oopman
  theory to allow for inputs and control}, SIAM Journal on Applied Dynamical
  Systems, 17 (2018), pp.~909--930.

\bibitem{proctor2015ih}
{\sc J.~L. Proctor and P.~A. Eckhoff}, {\em Discovering dynamic patterns from
  infectious disease data using dynamic mode decomposition}, International
  health, 7 (2015), pp.~139--145.

\bibitem{qian2020lift}
{\sc E.~Qian, B.~Kramer, B.~Peherstorfer, and K.~Willcox}, {\em Lift \& learn:
  Physics-informed machine learning for large-scale nonlinear dynamical
  systems}, Physica D: Nonlinear Phenomena, 406 (2020), p.~132401.

\bibitem{qin20061502}
{\sc S.~J. Qin}, {\em An overview of subspace identification}, Computers and
  Chemical Engineering, 30 (2006), pp.~1502 -- 1513,
  \url{https://doi.org/10.1016/j.compchemeng.2006.05.045}.

\bibitem{qin1997proc}
{\sc S.~J. Qin and T.~A. Badgwell}, {\em An overview of industrial model
  predictive control technology}, in AIChE Symposium Series, vol.~93, 1997,
  pp.~232--256.

\bibitem{queffelec2010}
{\sc M.~Queff{\'e}lec}, {\em Substitution {{Dynamical Systems}} - {{Spectral
  Analysis}}}, vol.~1294 of Lecture {{Notes}} in {{Mathematics}}, {Springer
  Berlin Heidelberg}, {Berlin, Heidelberg}, 2010.

\bibitem{rackauckas2020universal}
{\sc C.~Rackauckas, Y.~Ma, J.~Martensen, C.~Warner, K.~Zubov, R.~Supekar,
  D.~Skinner, and A.~Ramadhan}, {\em Universal differential equations for
  scientific machine learning}, arXiv preprint arXiv:2001.04385,  (2020).

\bibitem{raghu2017pmlr}
{\sc M.~Raghu, B.~Poole, J.~Kleinberg, S.~Ganguli, and J.~Sohl-Dickstein}, {\em
  On the expressive power of deep neural networks}, Proceedings of the 34th
  International Conference on Machine Learning, 70 (2017), pp.~2847--2854.

\bibitem{Raissi2019jcp}
{\sc M.~Raissi, P.~Perdikaris, and G.~Karniadakis}, {\em Physics-informed
  neural networks: A deep learning framework for solving forward and inverse
  problems involving nonlinear partial differential equations}, Journal of
  Computational Physics, 378 (2019), pp.~686--707.

\bibitem{raissi2020science}
{\sc M.~Raissi, A.~Yazdani, and G.~E. Karniadakis}, {\em Hidden fluid
  mechanics: Learning velocity and pressure fields from flow visualizations},
  Science, 367 (2020), pp.~1026--1030.

\bibitem{redman2021koopman}
{\sc W.~T. Redman}, {\em On koopman mode decomposition and tensor component
  analysis}, Chaos: An Interdisciplinary Journal of Nonlinear Science, 31
  (2021), p.~051101.

\bibitem{reed1978}
{\sc M.~Reed and B.~Simon}, {\em Methods of Modern Mathematical Physics.
  {{IV}}. {{Analysis}} of Operators}, {Academic Press [Harcourt Brace
  Jovanovich Publishers]}, {New York}, 1978.

\bibitem{reed1980}
{\sc M.~Reed and B.~Simon}, {\em Methods of Modern Mathematical Physics.
  {{I}}}, {Academic Press Inc. [Harcourt Brace Jovanovich Publishers]}, {New
  York}, second~ed., 1980.

\bibitem{rokhlin1966}
{\sc V.~A. Rokhlin}, {\em Selected topics from the metric theory of dynamical
  systems}, Amer. Math. Soc. Transl. Ser, 2 (1966), pp.~171--240.

\bibitem{rosenfeld2021}
{\sc J.~A. Rosenfeld, R.~Kamalapurkar, L.~F. Gruss, and T.~T. Johnson}, {\em
  Dynamic {{Mode Decomposition}} for {{Continuous Time Systems}} with the
  {{Liouville Operator}}}, arXiv:1910.03977 [math],  (2021),
  \url{https://arxiv.org/abs/1910.03977}.

\bibitem{rosenfeld2021a}
{\sc J.~A. Rosenfeld, B.~P. Russo, and R.~Kamalapurkar}, {\em Theoretical
  {{Foundations}} for the {{Dynamic Mode Decomposition}} of {{High Order
  Dynamical Systems}}}, arXiv:2101.02646 [math],  (2021),
  \url{https://arxiv.org/abs/2101.02646}.

\bibitem{rowley2009jfm}
{\sc C.~W. Rowley, I.~Mezi\'c, S.~Bagheri, P.~Schlatter, and D.~Henningson},
  {\em Spectral analysis of nonlinear flows}, J.\ Fluid Mech., 645 (2009),
  pp.~115--127.

\bibitem{roy2015pre}
{\sc S.~Roy, J.-C. Hua, W.~Barnhill, G.~H. Gunaratne, and J.~R. Gord}, {\em
  Deconvolution of reacting-flow dynamics using proper orthogonal and dynamic
  mode decompositions}, Physical Review E, 91 (2015), p.~013001.

\bibitem{Ruelle1986}
{\sc D.~Ruelle}, {\em Resonances of chaotic dynamical systems}, Physical review
  letters, 56 (1986), p.~405.

\bibitem{salova2019}
{\sc A.~Salova, J.~Emenheiser, A.~Rupe, J.~P. Crutchfield, and R.~M. D'Souza},
  {\em Koopman operator and its approximations for systems with symmetries},
  Chaos: An Interdisciplinary Journal of Nonlinear Science, 29 (2019),
  p.~093128, \url{https://doi.org/10.1063/1.5099091}.

\bibitem{sanchez2020learning}
{\sc A.~Sanchez-Gonzalez, J.~Godwin, T.~Pfaff, R.~Ying, J.~Leskovec, and
  P.~Battaglia}, {\em Learning to simulate complex physics with graph
  networks}, in International Conference on Machine Learning, PMLR, 2020,
  pp.~8459--8468.

\bibitem{sapsis2013pnas}
{\sc T.~P. Sapsis and A.~J. Majda}, {\em Statistically accurate low-order
  models for uncertainty quantification in turbulent dynamical systems},
  Proceedings of the National Academy of Sciences, 110 (2013),
  pp.~13705--13710.

\bibitem{sarkar2013mixed}
{\sc S.~Sarkar, S.~Ganguly, A.~Dalal, P.~Saha, and S.~Chakraborty}, {\em Mixed
  convective flow stability of nanofluids past a square cylinder by dynamic
  mode decomposition}, International Journal of Heat and Fluid Flow, 44 (2013),
  pp.~624--634.

\bibitem{sashidhar2021bagging}
{\sc D.~Sashidhar and J.~N. Kutz}, {\em Bagging, optimized dynamic mode
  decomposition (bop-dmd) for robust, stable forecasting with spatial and
  temporal uncertainty-quantification}, arXiv preprint arXiv:2107.10878,
  (2021).

\bibitem{sayadi2016tcfd}
{\sc T.~Sayadi and P.~J. Schmid}, {\em Parallel data-driven decomposition
  algorithm for large-scale datasets: with application to transitional boundary
  layers}, Theoretical and Computational Fluid Dynamics,  (2016), pp.~1--14.

\bibitem{sayadi2014reduced}
{\sc T.~Sayadi, P.~J. Schmid, J.~W. Nichols, and P.~Moin}, {\em Reduced-order
  representation of near-wall structures in the late transitional boundary
  layer}, Journal of Fluid Mechanics, 748 (2014), pp.~278--301.

\bibitem{Scherl2020prf}
{\sc I.~Scherl, B.~Strom, J.~K. Shang, O.~Williams, B.~L. Polagye, and S.~L.
  Brunton}, {\em Robust principal component analysis for particle image
  velocimetry}, Physical Review Fluids, 5 (2020).

\bibitem{schmid2011tcfd}
{\sc P.~Schmid, L.~Li, M.~Juniper, and O.~Pust}, {\em Applications of the
  dynamic mode decomposition}, Theoretical and Computational Fluid Dynamics, 25
  (2011), pp.~249--259.

\bibitem{schmid2009dynamic}
{\sc P.~J. Schmid}, {\em Dynamic mode decomposition of experimental data}, in
  8th International Symposium on Particle Image Velocimetry, Melbourne,
  Victoria, Australia, 2009.

\bibitem{schmid2010jfm}
{\sc P.~J. Schmid}, {\em Dynamic mode decomposition of numerical and
  experimental data}, Journal of Fluid Mechanics, 656 (2010), pp.~5--28.

\bibitem{schmid2008aps}
{\sc P.~J. Schmid and J.~Sesterhenn}, {\em Dynamic mode decomposition of
  numerical and experimental data}, in 61st Annual Meeting of the APS Division
  of Fluid Dynamics, American Physical Society, Nov. 2008.

\bibitem{schmid:2012}
{\sc P.~J. Schmid, D.~Violato, and F.~Scarano}, {\em Decomposition of
  time-resolved tomographic {PIV}}, Experiments in Fluids, 52 (2012),
  pp.~1567--1579.

\bibitem{Schmidt2009science}
{\sc M.~Schmidt and H.~Lipson}, {\em Distilling free-form natural laws from
  experimental data}, science, 324 (2009), pp.~81--85.

\bibitem{schmidt2020guide}
{\sc O.~T. Schmidt and T.~Colonius}, {\em Guide to spectral proper orthogonal
  decomposition}, Aiaa journal, 58 (2020), pp.~1023--1033.

\bibitem{sechi2021}
{\sc R.~Sechi, A.~Sikorski, and M.~Weber}, {\em Estimation of the {{Koopman
  Generator}} by {{Newton}}'s {{Extrapolation}}}, Multiscale Modeling \&
  Simulation, 19 (2021), pp.~758--774,
  \url{https://doi.org/10.1137/20M1333006}.

\bibitem{seena2011dynamic}
{\sc A.~Seena and H.~J. Sung}, {\em Dynamic mode decomposition of turbulent
  cavity flows for self-sustained oscillations}, International Journal of Heat
  and Fluid Flow, 32 (2011), pp.~1098--1110.

\bibitem{semeraro2012analysis}
{\sc O.~Semeraro, G.~Bellani, and F.~Lundell}, {\em Analysis of time-resolved
  {PIV} measurements of a confined turbulent jet using {POD} and {K}oopman
  modes}, Experiments in Fluids, 53 (2012), pp.~1203--1220.

\bibitem{Shadden2005pd}
{\sc S.~C. Shadden, F.~Lekien, and J.~E. Marsden}, {\em Definition and
  properties of {Lagrangian} coherent structures from finite-time {Lyapunov}
  exponents in two-dimensional aperiodic flows}, Physica D, 212 (2005),
  pp.~271--304.

\bibitem{sharma2016prf}
{\sc A.~S. Sharma, I.~Mezi{\'c}, and B.~J. McKeon}, {\em Correspondence between
  {K}oopman mode decomposition, resolvent mode decomposition, and invariant
  solutions of the {N}avier-{S}tokes equations}, Physical Review Fluids, 1
  (2016), p.~032402.

\bibitem{sharma2017ibpsa}
{\sc H.~Sharma, A.~D. Fontanini, U.~Vaidya, and B.~Ganapathysubramanian}, {\em
  Transfer operator theoretic framework for monitoring building indoor
  environment in uncertain operating conditions}, in 2018 Annual American
  Control Conference (ACC), IEEE, 2018, pp.~6790--6797.

\bibitem{sharma2018arxiv_a}
{\sc H.~Sharma, U.~Vaidya, and B.~Ganapathysubramanian}, {\em A transfer
  operator methodology for optimal sensor placement accounting for
  uncertainty}, Building and environment, 155 (2019), pp.~334--349.

\bibitem{sinha2020a}
{\sc S.~Sinha, B.~Huang, and U.~Vaidya}, {\em On {{Robust Computation}} of
  {{Koopman Operator}} and {{Prediction}} in {{Random Dynamical Systems}}},
  Journal of Nonlinear Science, 30 (2020), pp.~2057--2090,
  \url{https://doi.org/10.1007/s00332-019-09597-6}.

\bibitem{sinha2020arxiv}
{\sc S.~Sinha, S.~P. Nandanoori, and E.~Yeung}, {\em Data driven online
  learning of power system dynamics},  (2020), pp.~1--5.

\bibitem{sinha2020}
{\sc S.~Sinha, S.~P. Nandanoori, and E.~Yeung}, {\em Koopman operator methods
  for global phase space exploration of equivariant dynamical systems},
  IFAC-PapersOnLine, 53 (2020), pp.~1150--1155.

\bibitem{sinha2016jmaa}
{\sc S.~Sinha, U.~Vaidya, and R.~Rajaram}, {\em Operator theoretic framework
  for optimal placement of sensors and actuators for control of nonequilibrium
  dynamics}, Journal of Mathematical Analysis and Applications, 440 (2016),
  pp.~750--772.

\bibitem{slipantschuk2020}
{\sc J.~Slipantschuk, O.~F. Bandtlow, and W.~Just}, {\em Dynamic mode
  decomposition for analytic maps}, Communications in Nonlinear Science and
  Numerical Simulation, 84 (2020), p.~105179,
  \url{https://doi.org/10.1016/j.cnsns.2020.105179}.

\bibitem{son_narasingam2020}
{\sc S.~H. Son, A.~Narasingam, and J.~S.-I. Kwon}, {\em {Handling plant-model
  mismatch in Koopman Lyapunov-based model predictive control via offset-free
  control framework}}, arXiv Preprint arXiv:2010.07239,  (2020).

\bibitem{song2013global}
{\sc G.~Song, F.~Alizard, J.-C. Robinet, and X.~Gloerfelt}, {\em Global and
  {K}oopman modes analysis of sound generation in mixing layers}, Physics of
  Fluids (1994-present), 25 (2013), p.~124101.

\bibitem{sootla2016acc}
{\sc A.~Sootla and A.~Mauroy}, {\em Properties of isostables and basins of
  attraction of monotone systems}, in 2016 Annual American Control Conference
  (ACC), IEEE, 2016, pp.~7365--7370.

\bibitem{sootla2018}
{\sc A.~Sootla and A.~Mauroy}, {\em Operator-{{Theoretic Characterization}} of
  {{Eventually Monotone Systems}}}, IEEE Control Systems Letters, 2 (2018),
  pp.~429--434, \url{https://doi.org/10.1109/LCSYS.2018.2841654}.

\bibitem{steeb1980non}
{\sc W.-H. Steeb and F.~Wilhelm}, {\em Non-linear autonomous systems of
  differential equations and {C}arleman linearization procedure}, Journal of
  Mathematical Analysis and Applications, 77 (1980), pp.~601--611.

\bibitem{stewart1993early}
{\sc G.~W. Stewart}, {\em On the early history of the singular value
  decomposition}, SIAM review, 35 (1993), pp.~551--566.

\bibitem{Stoica2005}
{\sc P.~Stoica and R.~L. Moses}, {\em Spectral {{Analysis}} of {{Signals}}},
  {Pearson Prentice Hall}, 2005.

\bibitem{streif2013}
{\sc S.~Streif, P.~Rumschinski, D.~Henrion, and R.~Findeisen}, {\em Estimation
  of consistent parameter sets for continuous-time nonlinear systems using
  occupation measures and {{LMI}} relaxations}, in 52nd {{IEEE Conference}} on
  {{Decision}} and {{Control}}, Dec. 2013, pp.~6379--6384,
  \url{https://doi.org/10.1109/CDC.2013.6760898}.

\bibitem{strichartz1990}
{\sc R.~S. Strichartz}, {\em Fourier asymptotics of fractal measures}, Journal
  of Functional Analysis, 89 (1990), pp.~154--187,
  \url{https://doi.org/10.1016/0022-1236(90)90009-A}.

\bibitem{Strogatz2018book}
{\sc S.~H. Strogatz}, {\em Nonlinear dynamics and chaos with student solutions
  manual: With applications to physics, biology, chemistry, and engineering},
  CRC press, 2018.

\bibitem{suchanecki1996}
{\sc Z.~Suchanecki, I.~Antoniou, S.~Tasaki, and O.~F. Bandtlow}, {\em Rigged
  {{Hilbert}} spaces for chaotic dynamical systems}, Journal of Mathematical
  Physics, 37 (1996), pp.~5837--5847, \url{https://doi.org/10.1063/1.531703}.

\bibitem{surana2016cdc}
{\sc A.~Surana}, {\em Koopman operator based observer synthesis for
  control-affine nonlinear systems}.
\newblock 55th IEEE Conf. Decision and Control (CDC), pages 6492--6499, 2016.

\bibitem{surana2016linear}
{\sc A.~Surana and A.~Banaszuk}, {\em Linear observer synthesis for nonlinear
  systems using {K}oopman operator framework}, IFAC-PapersOnLine, 49 (2016),
  pp.~716--723.

\bibitem{surana2017cdc}
{\sc A.~Surana, M.~O. Williams, M.~Morari, and A.~Banaszuk}, {\em Koopman
  operator framework for constrained state estimation}, in Decision and Control
  (CDC), 2017 IEEE 56th Annual Conference on, IEEE, 2017, pp.~94--101.

\bibitem{susuki2011c}
{\sc Y.~Susuki and I.~Mezic}, {\em Nonlinear {{Koopman Modes}} and {{Coherency
  Identification}} of {{Coupled Swing Dynamics}}}, IEEE Transactions on Power
  Systems, 26 (2011), pp.~1894--1904,
  \url{https://doi.org/10.1109/TPWRS.2010.2103369}.

\bibitem{susuki2012nonlinear}
{\sc Y.~Susuki and I.~Mezic}, {\em Nonlinear {K}oopman modes and a precursor to
  power system swing instabilities}, IEEE Transactions on Power Systems, 27
  (2012), pp.~1182--1191.

\bibitem{susuki2015}
{\sc Y.~Susuki and I.~Mezi{\'c}}, {\em A prony approximation of {{Koopman Mode
  Decomposition}}}, in 2015 54th {{IEEE Conference}} on {{Decision}} and
  {{Control}} ({{CDC}}), Dec. 2015, pp.~7022--7027,
  \url{https://doi.org/10.1109/CDC.2015.7403326}.

\bibitem{susuki2020}
{\sc Y.~Susuki and I.~Mezi{\'c}}, {\em Invariant {{Sets}} in
  {{Quasiperiodically Forced Dynamical Systems}}}, SIAM Journal on Applied
  Dynamical Systems, 19 (2020), pp.~329--351,
  \url{https://doi.org/10.1137/18M1193529}.

\bibitem{susuki2011jns}
{\sc Y.~Susuki, I.~Mezi{\'c}, and T.~Hikihara}, {\em Coherent swing instability
  of power grids}, Journal of nonlinear science, 21 (2011), pp.~403--439.

\bibitem{stanton2016pre}
{\sc A.~Svenkeson, B.~Glaz, S.~Stanton, and B.~J. West}, {\em Spectral
  decomposition of nonlinear systems with memory}, Phys. Rev. E, 93 (2016),
  p.~022211, \url{https://doi.org/10.1103/PhysRevE.93.022211}.

\bibitem{svoronos1994discretization}
{\sc S.~Svoronos, D.~Papageorgiou, and C.~Tsiligiannis}, {\em Discretization of
  nonlinear control systems via the {C}arleman linearization}, Chemical
  engineering science, 49 (1994), pp.~3263--3267.

\bibitem{taira2017}
{\sc K.~Taira, S.~L. Brunton, S.~T.~M. Dawson, C.~W. Rowley, T.~Colonius, B.~J.
  McKeon, O.~T. Schmidt, S.~Gordeyev, V.~Theofilis, and L.~S. Ukeiley}, {\em
  Modal {{Analysis}} of {{Fluid Flows}}: {{An Overview}}}, AIAA Journal, 55
  (2017), pp.~4013--4041, \url{https://doi.org/10/gdh7zw}.

\bibitem{taira2020}
{\sc K.~Taira, M.~S. Hemati, S.~L. Brunton, Y.~Sun, K.~Duraisamy, S.~Bagheri,
  S.~T.~M. Dawson, and C.-A. Yeh}, {\em Modal {{Analysis}} of {{Fluid Flows}}:
  {{Applications}} and {{Outlook}}}, AIAA Journal, 58 (2020), pp.~998--1022,
  \url{https://doi.org/10.2514/1.J058462}.

\bibitem{takeishi2017jcai}
{\sc N.~Takeishi, Y.~Kawahara, Y.~Tabei, and T.~Yairi}, {\em Bayesian dynamic
  mode decomposition}, Twenty-Sixth International Joint Conference on
  Artificial Intelligence,  (2017).

\bibitem{Takeishi2017neurips}
{\sc N.~Takeishi, Y.~Kawahara, and T.~Yairi}, {\em Learning {K}oopman invariant
  subspaces for dynamic mode decomposition}, in Advances in Neural Information
  Processing Systems, 2017, pp.~1130--1140.

\bibitem{takeishi2017}
{\sc N.~Takeishi, Y.~Kawahara, and T.~Yairi}, {\em Subspace dynamic mode
  decomposition for stochastic {{Koopman}} analysis}, Physical Review E, 96
  (2017), p.~033310, \url{https://doi.org/10.1103/PhysRevE.96.033310}.

\bibitem{takens1981lnm}
{\sc F.~Takens}, {\em Detecting strange attractors in turbulence}, Lecture
  Notes in Mathematics, 898 (1981), pp.~366--381.

\bibitem{tallapragada2013set}
{\sc P.~Tallapragada and S.~D. Ross}, {\em A set oriented definition of
  finite-time lyapunov exponents and coherent sets}, Communications in
  Nonlinear Science and Numerical Simulation, 18 (2013), pp.~1106--1126.

\bibitem{tang2012dynamic}
{\sc Z.~Q. Tang and N.~Jiang}, {\em Dynamic mode decomposition of hairpin
  vortices generated by a hemisphere protuberance}, Science China Physics,
  Mechanics and Astronomy, 55 (2012), pp.~118--124.

\bibitem{tasaki1993}
{\sc S.~Tasaki, I.~Antoniou, and Z.~Suchanecki}, {\em Deterministic diffusion,
  {{De Rham}} equation and fractal eigenvectors}, Physics Letters A, 179
  (1993), pp.~97--102, \url{https://doi.org/10.1016/0375-9601(93)90656-K}.

\bibitem{taylor2017arxiv}
{\sc R.~Taylor, J.~N. Kutz, K.~Morgan, and B.~A. Nelson}, {\em Dynamic mode
  decomposition for plasma diagnostics and validation}, Review of Scientific
  Instruments, 89 (2018), p.~053501.

\bibitem{tibshirani1996lasso}
{\sc R.~Tibshirani}, {\em Regression shrinkage and selection via the lasso},
  Journal of the Royal Statistical Society. Series B (Methodological),  (1996),
  pp.~267--288.

\bibitem{tirunagari2017mva}
{\sc S.~Tirunagari, N.~Poh, K.~Wells, M.~Bober, I.~Gorden, and D.~Windridge},
  {\em Movement correction in {DCE-MRI} through windowed and reconstruction
  dynamic mode decomposition}, Machine Vision and Applications, 28 (2017),
  pp.~393--407.

\bibitem{Towne2018jfm}
{\sc A.~Towne, O.~T. Schmidt, and T.~Colonius}, {\em Spectral proper orthogonal
  decomposition and its relationship to dynamic mode decomposition and
  resolvent analysis}, Journal of Fluid Mechanics, 847 (2018), pp.~821--867.

\bibitem{trefethen1993science}
{\sc L.~N. Trefethen, A.~E. Trefethen, S.~C. Reddy, and T.~A. Driscoll}, {\em
  Hydrodynamic stability without eigenvalues}, Science, 261 (1993),
  pp.~578--584.

\bibitem{tu2011koopman}
{\sc J.~H. Tu, C.~W. Rowley, E.~Aram, and R.~Mittal}, {\em Koopman spectral
  analysis of separated flow over a finite-thickness flat plate with elliptical
  leading edge}, AIAA Paper 2011, 2864 (2011).

\bibitem{tu2014ef}
{\sc J.~H. Tu, C.~W. Rowley, J.~N. Kutz, and J.~K. Shang}, {\em Spectral
  analysis of fluid flows using sub-{N}yquist-rate {PIV} data}, Experiments in
  Fluids, 55 (2014), pp.~1--13.

\bibitem{tu2014jcd}
{\sc J.~H. Tu, C.~W. Rowley, D.~M. Luchtenburg, S.~L. Brunton, and J.~N. Kutz},
  {\em On dynamic mode decomposition: theory and applications}, Journal of
  Computational Dynamics, 1 (2014), pp.~391--421.

\bibitem{tukey1962future}
{\sc J.~W. Tukey}, {\em The future of data analysis}, The annals of
  mathematical statistics, 33 (1962), pp.~1--67.

\bibitem{ulam1960problems}
{\sc S.~Ulam}, {\em Problems in Modern Mathematics}, {Science Editions}, 1960.

\bibitem{vaidya2006ccc}
{\sc U.~Vaidya}, {\em Duality in stability theory: {L}yapunov function and
  {L}yapunov measure}, in 44th Allerton conference on communication, control
  and computing, 2006, pp.~185--190.

\bibitem{vaidya2007cdc}
{\sc U.~Vaidya}, {\em Observability gramian for nonlinear systems}, in Decision
  and Control, 2007 46th IEEE Conference on, IEEE, 2007, pp.~3357--3362.

\bibitem{vaidya2012jmaa}
{\sc U.~Vaidya, R.~Rajaram, and S.~Dasgupta}, {\em Actuator and sensor
  placement in linear advection {P}{D}{E} with building system application},
  Journal of Mathematical Analysis and Applications, 394 (2012), pp.~213--224.

\bibitem{vanoverschee199475}
{\sc P.~Van~Overschee and B.~De~Moor}, {\em {N4SID}: Subspace algorithms for
  the identification of combined deterministic-stochastic systems}, Automatica,
  30 (1994), pp.~75 -- 93.

\bibitem{van.1996}
{\sc P.~Van~Overschee and B.~De~Moor}, {\em Subspace Identification for Linear
  Systems: Theory - Implementation - Applications}, Springer US, 1996.

\bibitem{wang2021learning}
{\sc S.~Wang, H.~Wang, and P.~Perdikaris}, {\em Learning the solution operator
  of parametric partial differential equations with physics-informed
  deeponets}, arXiv preprint arXiv:2103.10974,  (2021).

\bibitem{wehmeyer2017arxiv}
{\sc C.~Wehmeyer and F.~No{\'e}}, {\em Time-lagged autoencoders: Deep learning
  of slow collective variables for molecular kinetics}, The Journal of Chemical
  Physics, 148 (2018), pp.~1--9.

\bibitem{wiener1941}
{\sc N.~Wiener and A.~Wintner}, {\em Harmonic analysis and ergodic theory},
  American Journal of Mathematics, 63 (1941), pp.~415--426.

\bibitem{wiggins2003}
{\sc S.~Wiggins}, {\em Introduction to Applied Nonlinear Dynamical Systems and
  Chaos}, vol.~2 of Texts in {{Applied Mathematics}}, {Springer-Verlag}, {New
  York}, second~ed., 2003.

\bibitem{williams2016ifac}
{\sc M.~O. Williams, M.~S. Hemati, S.~T.~M. Dawson, I.~G. Kevrekidis, and C.~W.
  Rowley}, {\em Extending data-driven {K}oopman analysis to actuated systems},
  IFAC-PapersOnLine, 49 (2016), pp.~704--709.

\bibitem{williams2015jnls}
{\sc M.~O. Williams, I.~G. Kevrekidis, and C.~W. Rowley}, {\em A data-driven
  approximation of the {K}oopman operator: extending dynamic mode
  decomposition}, Journal of Nonlinear Science, 6 (2015), pp.~1307--1346.

\bibitem{williams2015jcd}
{\sc M.~O. Williams, C.~W. Rowley, and I.~G. Kevrekidis}, {\em A kernel
  approach to data-driven {K}oopman spectral analysis}, Journal of
  Computational Dynamics, 2 (2015), pp.~247--265.

\bibitem{wilson2019}
{\sc D.~Wilson}, {\em Isostable reduction of oscillators with piecewise smooth
  dynamics and complex {{Floquet}} multipliers}, Physical Review E, 99 (2019),
  p.~022210, \url{https://doi.org/10.1103/PhysRevE.99.022210}.

\bibitem{wilson2020a}
{\sc D.~Wilson}, {\em A data-driven phase and isostable reduced modeling
  framework for oscillatory dynamical systems}, Chaos: An Interdisciplinary
  Journal of Nonlinear Science, 30 (2020), p.~013121,
  \url{https://doi.org/10.1063/1.5126122}.

\bibitem{wilson2020}
{\sc D.~Wilson}, {\em Phase-amplitude reduction far beyond the weakly perturbed
  paradigm}, Physical Review E, 101 (2020), p.~022220,
  \url{https://doi.org/10.1103/PhysRevE.101.022220}.

\bibitem{wilson2016}
{\sc D.~Wilson and J.~Moehlis}, {\em Isostable reduction of periodic orbits},
  Physical Review E, 94 (2016), p.~052213,
  \url{https://doi.org/10.1103/PhysRevE.94.052213}.

\bibitem{wright2009ieeetpami}
{\sc J.~Wright, A.~Yang, A.~Ganesh, S.~Sastry, and Y.~Ma}, {\em Robust face
  recognition via sparse representation}, {IEEE} Transactions on Pattern
  Analysis and Machine Intelligence {(PAMI)}, 31 (2009), pp.~210--227.

\bibitem{wu2020variational}
{\sc H.~Wu and F.~No{\'e}}, {\em Variational approach for learning markov
  processes from time series data}, Journal of Nonlinear Science, 30 (2020),
  pp.~23--66.

\bibitem{wu2017variational}
{\sc H.~Wu, F.~N{\"u}ske, F.~Paul, S.~Klus, P.~Koltai, and F.~No{\'e}}, {\em
  Variational koopman models: slow collective variables and molecular kinetics
  from short off-equilibrium simulations}, The Journal of chemical physics, 146
  (2017), p.~154104.

\bibitem{yeung2017arxiv}
{\sc E.~Yeung, S.~Kundu, and N.~Hodas}, {\em Learning deep neural network
  representations for {K}oopman operators of nonlinear dynamical systems},
  arXiv preprint arXiv:1708.06850,  (2017).

\bibitem{yeung2018acc}
{\sc E.~Yeung, Z.~Liu, and N.~O. Hodas}, {\em A {K}oopman operator approach for
  computing and balancing gramians for discrete time nonlinear systems}, in
  2018 Annual American Control Conference (ACC), IEEE, 2018, pp.~337--344.

\bibitem{you2018ifac}
{\sc P.~You, J.~Pang, , and E.~Yeung}, {\em Deep {K}oopman synthesis for
  cyber-resilient market-based frequency regulation}, IFAC-PapersOnLine, 51
  (2018), pp.~720--725.

\bibitem{liu2018acc}
{\sc L.~C. Z~Liu, S~Kundu and E.~Yeung}, {\em Decomposition of nonlinear
  dynamical systems using {K}oopman gramians}, in 2018 American Control
  Conference, IEEE, 2018.

\bibitem{zaks2017}
{\sc M.~A. Zaks and A.~Nepomnyashchy}, {\em Anomalous {{Transport}} in {{Steady
  Plane Viscous Flows}}: {{Simple Models}}}, in Advances in {{Dynamics}},
  {{Patterns}}, {{Cognition}}, Nonlinear {{Systems}} and {{Complexity}},
  {Springer, Cham}, 2017, pp.~61--76,
  \url{https://doi.org/10.1007/978-3-319-53673-6_5}.

\bibitem{zaslavsky1991}
{\sc G.~M. Zaslavsky}, {\em Lagrangian turbulence and anomalous transport},
  Fluid Dynamics Research, 8 (1991), p.~127,
  \url{https://doi.org/10.1016/0169-5983(91)90036-I}.

\bibitem{zaslavsky2002}
{\sc G.~M. Zaslavsky}, {\em Chaos, fractional kinetics, and anomalous
  transport}, Physics Reports, 371 (2002), pp.~461--580,
  \url{https://doi.org/10.1016/S0370-1573(02)00331-9}.

\bibitem{zhang2019}
{\sc H.~Zhang, C.~W. Rowley, E.~A. Deem, and L.~N. Cattafesta}, {\em Online
  {{Dynamic Mode Decomposition}} for {{Time}}-{{Varying Systems}}}, SIAM
  Journal on Applied Dynamical Systems, 18 (2019), pp.~1586--1609,
  \url{https://doi.org/10.1137/18M1192329}.

\bibitem{zhen2021}
{\sc Y.~Zhen, E.~Memin, B.~Chapron, and L.~Peng}, {\em Singular {{Values}} of
  {{Hankel Matrix}}: {{Implication}} to the {{Convergence}} of {{Hankel Dynamic
  Mode Decomposition}}}, arXiv:2107.01948 [math],  (2021),
  \url{https://arxiv.org/abs/2107.01948}.

\end{thebibliography}
\end{document}